\documentclass[a4paper,10pt]{article}
\usepackage[left=3cm,top=2cm,bottom=2cm,right=2.2cm,includehead,includefoot]{geometry}
\usepackage[ansinew]{inputenc}
\usepackage[american]{babel}
\usepackage{amsfonts} 
\usepackage{amssymb}
\usepackage{amsthm}
\usepackage{amsmath}
\usepackage{mathtools}
\usepackage{mathrsfs}
\usepackage[arrow, matrix, curve]{xy}
\usepackage{enumitem}
\usepackage{dsfont}
\usepackage{bbm}
\usepackage{url}

\title{Foliations and the cohomology \\ of moduli spaces of bounded global \textit{G}-shtukas}
\author{Stephan Neupert}

\newcommand{\id}{{\rm id}}

\newcommand{\ov}[1]{\overline{#1}}
\newcommand{\un}[1]{\underline{#1}}
\newcommand{\lmat}{\left(\begin{array}}
\newcommand{\rmat}{\end{array}\right)}
\newcommand{\M }[1]{\mathbb #1}
\newcommand{\MB}[1]{\mathbf{#1}}
\newcommand{\MC}[1]{\mathcal{#1}}
\newcommand{\MF}[1]{\mathfrak{#1}}
\newcommand{\MS}[1]{\mathscr{#1}}
\newcommand{\Fq}{\M{F}_q}
\newcommand{\F }{\M{F}}

\newcommand{\ACFq}{\ov{\M{F}}_q}
\newcommand{\unif}{\mathcal{\pi}}
\newcommand{\Sp}{{\rm Spec} \,}
\newcommand{\Spf}{{\rm Spf} \,}
\newcommand{\Spa}{{\rm Spa} \,}
\newcommand{\hooklongrightarrow}{\lhook\joinrel\relbar\joinrel\rightarrow}

\newcommand{\ENil}{E[[\zeta]]}
\newcommand{\op}[1]{\operatorname{#1}}
\newcommand{\exit}{\hfill $\square$}
\newcommand{\mmu}{\boldsymbol{\mu}}
\newcommand{\ignore}[1]{}
\newcommand{\Isom}{\op{Isom}}
\newcommand{\Aut}{\op{Aut}}
\newcommand{\End}{\op{End}}
\newcommand{\Hom}{\op{Hom}}
\newcommand{\Rep}{\op{Rep}}
\newcommand{\Res}{\op{Res}}
\newcommand{\Gal}{\op{Gal}}
\newcommand{\red}{\op{red}}

\newtheorem{lemma}{Lemma}[subsection]
\newtheorem{proposition}[lemma]{Proposition}
\newtheorem{theorem}[lemma]{Theorem}
\newtheorem{corollary}[lemma]{Corollary}
\newtheorem{definition}[lemma]{Definition}
\newtheorem{remark}[lemma]{Remark}
\newtheorem{conjecture}[lemma]{Conjecture}
\newtheorem{construction}[lemma]{Construction}
\newtheorem{notation}[lemma]{Notation}
\newtheorem{warning}[lemma]{Warning}
\newtheorem{example}[lemma]{Example}
\newtheorem{assumption}[lemma]{Assumption}
\newtheorem{sublemma}[lemma]{Sublemma}
\newtheorem{problem}[lemma]{Problem}
\newtheorem{maintheorem}{Main Theorem}

\newcommand{\lem}[3]{\begin{lemma} #1 #2 #3 \end{lemma}}
\newcommand{\prop}[3]{\begin{proposition} #1 #2 #3 \end{proposition}}
\newcommand{\thm}[3]{\begin{theorem} #1 #2 #3 \end{theorem}}
\newcommand{\cor}[3]{\begin{corollary} #1 #2 #3 \end{corollary}}
\newcommand{\defi}[3]{\begin{definition} #1 #2 #3 \end{definition}}
\newcommand{\rem}[2]{\begin{remark} #1 {\rm #2} \end{remark}}

\newcommand{\const}[3]{\begin{construction} #1 {\rm #2} $\left.\right.$ \\ {\rm #3} \end{construction}}
\newcommand{\nota}[2]{\begin{notation} #1 #2 \end{notation}}
\newcommand{\warn}[2]{\begin{warning} #1 #2 \end{warning}}
\newcommand{\ex}[2]{\begin{example} #1 {\rm #2} \end{example}}
\newcommand{\ass}[2]{\begin{assumption} #1 {\rm #2} \end{assumption}}
 
\newcommand{\prob}[3]{\begin{problem} #1 {\rm #2} \end{problem}} 
\newcommand{\mainthm}[2]{\begin{maintheorem} #1 #2 \end{maintheorem}}

\newcommand{\prooof}{\proof}

\usepackage{hyperref}

\begin{document}

\maketitle

\begin{abstract}
We decompose Newton strata in the special fiber of moduli spaces of global $G$-shtukas into a product of Rapoport-Zink spaces and Igusa varieties. This decomposition holds as well for associated adic spaces, which yields a comparison between the $\ell$-adic cohomology of these spaces together with the action by the reductive group $G$ and the Galois-group of the underlying ground field.
\end{abstract}

\section*{Introduction}
Moduli spaces of bounded global $G$-shtukas are the function field analogues of Shimura varieties over number fields. Therefore one expects that Langlands correspondences can be realized in their $\ell$-adic cohomology. In fact recently V. Lafforgue \cite{LafforgueVShtukas} made huge progress in this direction. In this article a different approach is taken: Following Harris and Taylor \cite{HarrisTaylor} and Mantovan \cite{MantoFoliation}, \cite{MantoFoliationPEL} the moduli space of bounded global $G$-shtukas is replaced by two moduli spaces parametrizing simpler objects and compare the representations found in their cohomology. \\
Let us describe the main results in slightly greater detail: Let $\Fq$ be a finite field of characteristic $p$, $G$ a reductive group over $\Fq$ and $K$ a function field over $\Fq$. By fixing a suitable ring of integers in $K$, one obtains a smooth projective geometrically connected curve $C$. Let $\sigma$ be the Frobenius on $C$. 
Based on Drinfeld's notion of elliptic modules \cite{DrinfeldElliptic}, Varshavsky \cite{Varsh} defined global $G$-shtukas as $G$-torsors $\MS{G}$ over $C$ together with an isomorphism $\varphi: \sigma^*\MS{G} \to \MS{G}$ over an open subscheme of $C$. The points $c_1, \ldots, c_n$ in the complement of this open subscheme are called characteristic places. Moreover level structures with respect to open subgroups $U \subset G(\M{A})$ and boundedness conditions exist on global $G$-shtukas, cf. \ref{subsec:GlobalBounds} and \ref{subsec:AdelicLevel} for definitions. 
Such global $G$-shtukas then admit moduli spaces $\nabla_n^{\mmu}\MC{H}^1_U(C, G)$, which are Deligne-Mumford-stacks whose (in general infinitely many) connected components are quotients of quasi-projective schemes of finite type by finite groups. They exist not only over the function field, but admit natural integral models, whose special fibers $\MB{X}^{\mmu}_U$ are characterized by fixing the locus of the characteristic places. \\
Loosely speaking, the key observation is now, that one may separate the behavior of a global $G$-shtuka locally around the characteristic places from its global structure, at least over the special fiber. 
The local part was analyzed in detail by Hartl and Viehmann \cite{HaVi}, \cite{HaVi2} and is given by the Rapoport-Zink space of local $G$-shtukas, a function field analogue of the moduli spaces of $p$-divisible groups defined by Rapoport and Zink \cite{RapoZinkSpaces}. To $G$ one may associate the loop groups $L^+G$ and $LG$ over $\Fq$ given by the sheaves
\[L^+G(S) \coloneqq G(\MC{O}_S[[z]]) \quad LG(S) \coloneqq G(\MC{O}_S((z))) \qquad \textnormal{for schemes} \; S \; \textnormal{over} \; \Fq.\]
Then the restriction of a $G$-torsor $\MS{G}$ over a curve $C$ to the formal neighborhood of a point $c_i$ is nothing else than a $L^+G$-torsor $\MC{G}$ over the point itself. Similarly the restriction of the Frobenius-isomorphism $\varphi$ is now given by a Frobenius-isomorphism $\varphi: \sigma^*(\MC{G} \times^{L^+G} LG) \to \MC{G} \times^{L^+G} LG$ between the associated $LG$-torsors. Such a pair $(\MC{G}, \varphi)$ is called a local $G$-shtuka. 
Thus fixing a fundamental alcove $b_{\nu}$ (cf. \ref{def:FundamentalAlcove}) and introducing boundedness conditions here as well, one can give the moduli problem
\[\MC{M}_{b_{\nu}}^{\preceq \mu}(S) = \begin{Bmatrix} (\MC{G}, \varphi) \; \textnormal{a local} \; G \textnormal{-shtuka bounded by} \; \mu \; \textnormal{over} \; S \; \textnormal{and} \\ \alpha: (\MC{G}, \varphi) \to (L^+G, b_{\nu}\sigma^*) \; \textnormal{a quasi-isogeny}\end{Bmatrix} \]
The formal scheme representing $\MC{M}_{b_{\nu}}^{\preceq \mu}$ is called Rapoport-Zink spaces. Denote its special fiber by $\MB{M}_{b_{\nu}}^{\preceq \mu}$. \\
Unfortunately the global counterpart, namely Igusa varieties, are not as easily constructed. Following the work of Harris and Taylor \cite{HarrisTaylor} they should parametrize global $G$-shtukas together with equivalence classes of trivializations of the associated local $G$-shtukas at the characteristic places. In other words, one should consider central leaves $\MC{C}_U^{(\nu_i)}$, i.e. the loci where the associated local $G$-shtukas lie in one specific isomorphism class, and parameterize partial trivializations over them. However whenever the chosen fundamental alcove $b_{\nu}$ is not basic, the usual definition turns out not to be representable over central leaves. 
Mantovan suggested to pull back the global $G$-shtukas along a high power of the Frobenius, which then admits a partial splitting into basic factors, which suffices to construct the Igusa varieties. We essentially follow this idea, but modify it in two ways: \\
First of all we base-change directly to the perfection of the central leaf. Then pulling pack along the Frobenius defines an isomorphism, allowing us to work directly with the given global $G$-shtuka. \\
Secondly using the description of general Igusa varieties as products of Igusa varieties for the basic case, complicates all further arguments which use the moduli description. This is mainly due to the fact, that whenever constructing partial trivializations, one has to check that they split and hence define a point in the Igusa varieties. Fortunately there is a uniform moduli description: For simplicity let us consider only the case of a single characteristic place $c_i$ in this introduction. Fix as above a fundamental alcove $b_{\nu_i}$ and consider then the central leaf, i.e. the locally closed locus in $\MB{X}^{\mmu}_U$ defined by all geometric points where the associated local $G$-shtuka is isomorphic to $(L^+G, b_{\nu_i}\sigma^*)$. Moreover define for each $d \geq 0$ the subgroup
\[I_d(b_{\nu_i}) \coloneqq \bigcap_{N \geq 0} \phi^N(K_d) \subset LG\]
where $\phi(g) = b_{\nu_i}^{-1}\sigma^{-1}(g)b_{\nu_i}$ for all $g \in LG$ and $K_d = \{g \in L^+G \,|\, g = 1 \bmod z^{d+1}\}$. Then $I_d(b_{\nu_i})$ is a closed subgroup in $L^+G$ and the universal local $G$-shtuka admits a canonical $I_0(b_{\nu_i})$-structure, i.e. a natural $I_0(b_{\nu_i})$-subtorsor inside its $L^+G$-torsor. Note that such an $I_0(b_{\nu_i})$-structure is strictly stronger than having complete slope division (once this notion is transported from $p$-divisible groups \cite[definition 1.2]{OortZinkFamilies} to local $G$-shtukas, cf. \ref{subsec:SlopeDivisible}). \\
In this situation it makes sense to talk about $I_d(b_{\nu_i})$-truncated isomorphisms, which are equivalence classes of isomorphisms that (roughly speaking) induce the same map modulo $I_d(b_{\nu_i})$. For a precise definition see \ref{subsec:DefineIgusaBasic}. Now we may define Igusa varieties as the moduli space parametrizing bounded global $G$-shtukas (with level structure) together with an $I_d(b_{\nu_i})$-truncated isomorphism between the associated local $G$-shtuka and $(L^+G, b_{\nu_i}\sigma^*)$. Such moduli spaces $\op{Ig}^{d \sharp}_{c_i, U}$ are representable by a finite \'etale cover over the perfection of the central leaf $\MC{C}_U^{(\nu_i) \sharp}$. Note that we actually construct these Igusa varieties directly, instead of using the product of Igusa varieties for the basic case and then showing the equivalence of the moduli descriptions. \\
If one uses this construction for all characteristic places, one similarly obtains an Igusa variety $\op{Ig}^{(d_i) \sharp}_U$ for tuples $(d_i)_i$. In the limit for growing $d_i$ these Igusa varieties parametrize global $G$-shtukas with actual trivializations at all characteristic places. \\
The connection between Rapoport-Zink spaces, Igusa varieties and the moduli space of global $G$-shtukas is achieved via the so-called uniformization morphism. In the setting of mixed characteristic, such morphisms were constructed in a very general setting by Rapoport and Zink \cite{RapoZinkSpaces} and Mantovan \cite{MantoFoliationPEL}, based on a large number of previous works in more special cases. For shtukas a first version was constructed by Arasteh Rad and Hartl in \cite{HarRad1}. The idea is to take the global $G$-shtuka living over the Igusa variety, cut out the torsor at the characteristic places and then fill the hole using the canonical quasi-isogeny to the universal local $G$-shtuka over the Rapoport-Zink space. This way one obtains a morphism
\[\pi_{(\infty_i)}: \prod_i \MB{M}_{b_{\nu_i}}^{\preceq \mu_i \sharp} \times \op{Ig}^{(\infty_i) \sharp}_U \to \MB{X}^{\mmu \sharp}_U.\]
Note that due to the existence of Igusa varieties only over perfect schemes, $\pi_{(\infty_i)}$ (a priori) only exist after passing to the perfection for all schemes.

\mainthm{}{
a) The morphism $\pi_{(\infty_i)}$ factors over the perfection of the Newton stratum $\MC{N}^{(\nu_i) \sharp}_U \subset \MB{X}^{\mmu \sharp}_U$, i.e. the locus where the associated local $G$-shtukas are quasi-isogenous (but not necessarily isomorphic) to $(L^+G, b_{\nu_i}\sigma^*)$. It is surjective over the Newton stratum. The fibers over geometric points of $\MC{N}^{(\nu_i) \sharp}_U$ are torsors under the group of self-quasi-isogenies $\prod J_i$ of the associated local $G$-shtukas $(L^+G, b_{\nu_i}\sigma^*)$ at all characteristic places. \\
b) $\pi_{(\infty_i)}$ can be expressed as a limit of morphisms, which are surjective on the Newton stratum themselves and are alternating finite or \'etale. 
}

Most parts of the theorem were shown by Mantovan \cite{MantoFoliation}, \cite{MantoFoliationPEL} in the case of Shimura varieties of PEL-type. Only the \'etale versions of the covering morphism are completely new and due to our change to consider everything over perfections. One more advantage of this change becomes visible, when one tries to extend this result to formal schemes: While definitions of covering morphisms like the ones used by Mantovan allow only locally and non-canonical extensions to the formal setting, the theorem above stays true (almost verbatim) when all spaces are replaced by their formal counterparts, as explained in section \ref{sec:FormalLifting}. 
In fact one can even pass to (analytic) adic spaces with only slight changes. \\
As a first application of this theorem, we prove that the dimension of leaves, i.e. the loci in $\MB{X}^{\mmu}_U$ where the associated local $G$-shtukas lie in some fixed isomorphism class, are constant inside one Newton stratum. \\
After having now a very good understanding of the geometry, we are able to deal with the cohomology: From the point of view of Langlands' correspondence we are interested in 
\[\sum_i (-1)^i \varinjlim_U H^i_c\left(\nabla^{\mmu}_n\MC{H}^1_U(C, G) \times \ov{K}, \M{Q}_\ell \right)\]
as a representation of $G(\M{A}) \times \Gamma_{\Fq}$, where $\M{A}$ is the ring of adeles and $\Gamma_{\Fq}$ is the absolute Galois group of $\Fq$ (or equivalently the Weil group of the absolute Galois group of $K$). Here the action of $G(\M{A}) \times \Gamma_{\Fq}$ is induced by the natural action of these groups on the tower of spaces $\nabla^{\mmu}_n\MC{H}^1_U(C, G) \times \ov{K}$ for varying level structure $U$.
However due to technical problems, we can only deal with it as a $G(\M{A}^{c_i}) \times \Gamma_{E'}$-representation, where $\M{A}^{c_i}$ is the ring of adeles away from the fixed characteristic places and $\Gamma_{E'}$ is the absolute Galois group of a finite extension $E'/\Fq$. Moreover for similar reasons we have to impose a projectivity condition on the moduli space of global $G$-shtukas. \\
After transporting the cohomology to the special fiber using vanishing cycles and using the decomposition into Newton strata, one is then reduced to the study of
\[\sum_i (-1)^i \varinjlim_U H^i_c\left(\MC{N}^{(\nu_i)}_U \times \ACFq, R\Psi_\eta^{an}\M{Q}_\ell\right)\]
Then a K\"unneth type formula allows us to rewrite this in terms of the cohomology of Rapoport-Zink spaces and Igusa varieties resulting in the 

\mainthm{}{
Let $\nabla^{\mmu}_n\MC{H}^1(C, G)$ be a moduli space of global $G$-shtukas, such that all connected components of $\nabla^{\mmu}_n\MC{H}^1(C, G)$ are proper over $\Sp E'[[\zeta_1, \ldots, \zeta_n]]$. \\
Then there exists a canonical isomorphism between the virtual $G(\M{A}^{c_i}) \times \Gamma_{E'}$-representations
\[\sum_i (-1)^i H^i_c\left(\nabla^{\mmu}_n\MC{H}^1(C, G) \times \ov{K}, \M{Q}_\ell\right)\]
and 
\[\sum_{(\nu_i)} \sum_{d, e, f} (-1)^{d+e+f} Tor_d^{\MC{H}(\prod J_i)} \left(H^e_c\left(\prod \MB{M}_{b_{\nu_i}}^{\preceq \mu_i} 
\times \ACFq, R\Psi_\eta^{an}\M{Q}_\ell \right), \varinjlim_U \varinjlim_{d_i} H^f_c \left(\op{Ig}^{(d_i)}_U \times \ACFq, R\Psi_\eta^{an}\M{Q}_\ell \right)\right).\]
}

This theorem implies, that instead of analyzing the cohomology of the whole moduli space of global $G$-shtukas, it essentially suffices to understand the cohomology of Rapoport-Zink spaces and Igusa varieties.
A few more remarks on this formula: While the description above suggests that the product decomposition is only needed over the special fiber, this is far from true. To apply the K\"unneth formula one has to compare vanishing cycles sheaves for Newton strata, Rapoport-Zink spaces and Igusa varieties. For this it is essential to have a good control over the situation of associated analytic adic spaces. On the other hand passing to perfections does not change the cohomology, so having the covering morphism only over perfections does not impose further difficulties. \\
When comparing this result to the corresponding one for Shimura varieties found in \cite{MantoFoliation} or \cite{MantoFoliationPEL}, the most notable difference is the appearance of vanishing cycles in our final formula. This is mainly due to problems with incompatibilities of stratifications and passage to adic spaces, which create problems when trying to argue as in \cite{MantoFoliation}. Please refer to the discussion at the end of this article, why it is even likely that removing the vanishing cycles is impossible. \vspace{2mm} \\
A few final remarks on the structure of this paper: Section \ref{sec:LocalShtuka} deals with the local theory. Much of this treatment is a straight-forward generalization of \cite{HaVi}, though the new boundedness statement \ref{Prop:BoundInverse} for inverses of quasi-isogenies simplifies many arguments. 
Moreover in \ref{subsec:TateFunctor} the description of the Tate functor is improved upon the treatment in \cite{HarRad1}. Finally we prove in \ref{subsec:TateTheorem}, that quasi-isogenies on generic fibers of normal schemes extend uniquely to the whole scheme, at least if the quasi-isogeny class of the local $G$-shtuka does not vary. The corresponding result for $p$-divisible groups is called Tate's theorem, cf. \cite[theorem 4]{TatePDivisible} and \cite[\S 4.1]{BerthelotDieu}. \\
Section \ref{sec:GlobalShtuka} is devoted to the theory of global $G$-shtukas. While most purely global statements can already be found in \cite{Varsh} or \cite{HarRad1}, \cite{HarRad2}, the focus here lies on the interplay between global and local notions. The global-local-functors $\MF{L}_{c_i}$ in \ref{subsec:GlobalLocalFunctor} were already constructed in \cite{HarRad1}, though using an alternative construction less suited for our needs. \\
Section \ref{sec:Igusa} finally introduces Igusa varieties over central leaves. The first step in this direction is to translate the notion of a complete slope division to shtukas and to see that universal local $G$-shtuka over central leaves admit them, which is done in \ref{subsec:SlopeDivisible}. In \ref{subsec:IwahoriStructure} this is strengthened to the existence of $I_0(b_{\nu_i})$-structures. This allows us to prove representability of Igusa varieties over basic strata in \ref{subsec:DefineIgusaBasic} and in general in \ref{subsec:DefineIgusaGeneral}. In these sections the basic properties of Igusa varieties are established as well. \\
Having now Igusa varieties, one can prove main theorem $1$ in section \ref{sec:FiniteCover}. \ref{subsec:Uniform} describes the construction of the uniformization morphism, removing some inaccuracies in the original proof of \cite{HarRad2}. After discussing several versions of the product composition in \ref{subsec:CoverExistence} to \ref{subsec:EtaleCoveringMorphism}, the application to dimensions of arbitrary leaves is discussed in \ref{subsec:DimensionLeaf}. \\
The next section \ref{sec:FormalLifting} focuses on extending the product decomposition to formal schemes and adic spaces, after constructing a suitable formal version of central leaves in \ref{sec:FormalIgusaVarieties}. \\
Finally we deal with the cohomology in section \ref{sec:Cohomology}. The first half deals with establishing a K\"unneth type formula over the special fiber in the flavor of theorem 5.13 of \cite{MantoFoliation}. In \ref{sec:TorsionCohomology} this is applied to the sheaves of vanishing cycles of torsion sheaves, which finally gives main theorem $2$ in \ref{sec:LimitsOfMath}.


\subsection*{Acknowledgements}
I am very grateful to my advisor E. Viehmann for suggesting to study this topic. I am grateful to A. Wei\ss{} for pointing out an error in a preliminary version and to P. Hamacher for suggesting to use perfections. I like to thank A. Ivanov and C. Liedtke for helpful discussions. \\
This author was supported by ERC starting grant 277889 'Moduli spaces of local $G$-shtukas'.

\section{Notations}\label{sec:Notations}
We fix the following data:
\begin{itemize}
\item $\Fq$ a finite field of characteristic $p$.
\item $\ACFq$ an algebraic closure of $\Fq$.
\item $E$ a finite field extension of $\Fq$ inside $\ACFq$. Any other finite field extension will always take place inside $\ACFq$.
\item $\Gamma = Gal(\ACFq/E)$ the absolute Galois group of $E$. 
\item $C$ a smooth geometrically irreducible projective curve over $\Fq$.
\item $G$ a connected reductive group over $\Fq$.
\item $B \subset G$ a Borel subgroup defined over $\Fq$.
\item $T \subset G$ a maximal (not necessarily split) torus defined over $\Fq$. 
\item $X^*(T) = \Hom_{\ACFq}(T \times \Sp \ACFq, \M{G}_m \times \Sp \ACFq)$ the character group of $T$.
\item $X_*(T) = \Hom_{\ACFq}(\M{G}_m \times \Sp \ACFq, T \times \Sp \ACFq)$ the cocharacter group of $T$.
\item $X^*(T)_{\rm dom}$ respectively $X_*(T)_{\rm dom}$ the set of dominant elements in the (co)character group of $T$ (wrt. the chosen Borel $B$).
\item $\pi_1(G) = X_*(T)/\{\textnormal{coroot lattice}\}$ the (algebraic) fundamental group of $G$. 
\item For any set $S$ with $\Gamma$-action (like $X_*(T)$ or $\pi_1(G)$), we denote the set of $\Gamma$-orbits by $S/\Gamma$.
\end{itemize}
Note that any connected reductive group over a finite field is automatically quasi-split (e.g. by \cite[4.3-4.4]{SpringerRedGroup}), hence the Borel $B$ indeed exists. \\
We will use two Frobenius morphisms depending whether we consider the local or the global situation. Both are denoted by $\sigma$ and are for $\Fq$-schemes $S$
\begin{itemize}
\item in the local context (i.e. in section \ref{sec:LocalShtuka}, etc.): $\sigma: S \to S$ the absolute $q$-Frobenius. 
\item in the global context (i.e. in most of section \ref{sec:GlobalShtuka}, etc.): $\sigma: C \times_{\Sp \Fq} S \to C \times_{\Sp \Fq} S$ the identity on $C$ and the absolute $q$-Frobenius on $S$.
\end{itemize}

Furthermore for two schemes or stacks $X, Y$ over $\Sp \Fq$ respectively $\Sp E$ we will abbreviate the fiber product 
\[X \times_{\Fq} Y \coloneqq X \times_{\Sp \Fq} Y \qquad {\rm resp.} \qquad X \times_{E} Y \coloneqq X \times_{\Sp E} Y.\] 

\rem{}{
We essentially require for $E$ to be large enough to satisfy the following two conditions: 
\begin{itemize}
 \item Every bound $\mu$ in the cocharacter group of the maximal torus of $\Res_{\M{F}/\Fq}(G)$ (for some fixed finite extension $\M{F}$ of $\Fq$) can be defined over $E$, cf. \ref{Def:BoundLocal}.
 \item The fundamental alcove $b_\nu$ can be represented by an element defined over $E$, cf. \ref{ass:FundamentalAlcoveDefined}.
\end{itemize}
In the case of split groups $G$ these two conditions are always satisfied already for $E = \Fq$.
}

\section{Local \textit{G}-shtukas}\label{sec:LocalShtuka}
We define local $G$-shtukas for connected reductive groups $G$ and study their properties. Except for the sections \ref{subsec:LocalBoundsInv}, \ref{subsec:TateFunctor} and \ref{subsec:TateTheorem}, all statements were essentially established by Hartl and Viehmann \cite{HaVi}. The main difference to their treatment is, that we no longer assume that $G$ is split, though this requires only small changes to any of their proofs. Section \ref{subsec:TateFunctor} was inspired by \cite{HarRad1} and \cite{HarRad2}. There much of the theory of $G$-shtukas is developed even for non-constant group schemes, though the statements are sometimes not strong enough for our purposes. \\
We will state all theorems needed in the other parts of this work, but mostly refer for proofs to \cite{HaVi} (at least after explaining how to reduce to the split case).

\subsection{Generalities on torsors}\label{subsec:Torsors}
Let us recall some generalities on torsors under (ind-)group schemes. These will be used frequently (without explicit reference) for various groups. In this section the base spaces $S$ will always be Deligne-Mumford stacks, from now abbreviated by DM-stacks.

\defi{}{}{
Let $H$ be an (ind-)group scheme over $\Fq$ and fix some topology $* \in \!\{fpqc, fppf, \acute{e}tale\}$. A (right) $H$-torsor on a DM-stack $S$ over $\Fq$ is a sheaf $\MC{H}$ for the $*$-topology together with a (right) action of $H$ on $\MC{H}$, such that $*$-locally on $S$ the sheaf $\MC{H}$ is isomorphic to the sheaf of points of $H$.
}

\rem{}{
i) Any torsor under an (ind-)affine (ind-)scheme is relatively representable by an (ind-)scheme. Thus whenever convenient, we will view $\MC{H}$ as an $H$-(ind-)scheme which admits after a $*$-cover $S' \to S$ an $H$-equivariant isomorphism $\MC{H} \times_S S' \cong H \times S'$. \\
ii) $H$-torsors over $S$ are classified up to isomorphism by \v{C}ech cohomology $\check{H}^1(S_*, H)$. \\
iii) If $f: S' \to S$ is any morphism of DM-stacks and $\MC{H}$ an $H$-torsor over $S$, then its pull-back $f^*\MC{H}$ is an $H$-torsor over $S'$. In particular this applies to the Frobenius morphism.
}

If $H_1 \to H_2$ is a morphism between (ind-)group schemes, then the induced morphism $\check{H}^1(S_*, H_1) \to \check{H}^1(S_*, H_2)$ associates to each $H_1$-torsor an $H_2$-torsor. This can be turned into a functor via the following explicit

\const{\label{const:GenTorsorGroup}}{}{
Let $\MC{H}_1$ be an $H_1$-torsor over a DM-stack $S$. Then consider the $*$-sheaf $\MC{H}_1 \times^{H_1} H_2$ defined as the sheafification of the presheaf
\[S' \mapsto \{H_1(S')\op{-orbits \; in \,} \MC{H}_1(S') \times H_2(S')\}\]
where $H_1(S')$ acts on $\MC{H}_1(S')$ via the given right action and on $H_2(S')$ via left multiplication by the inverse of its image under $H_1(S') \to H_2(S')$. \\
Then $\MC{H}_1 \times^{H_1} H_2$ admits a canonical (right) $H_2$-action by right multiplication on the second factor. It is easy to see that $\MC{H}_1 \times^{H_1} H_2$ admits a trivialization over $S'$ whenever this happens for $\MC{H}_1$. So $\MC{H}_1 \times^{H_1} H_2$ is indeed an $H_2$-torsor. \\
It is an easy exercise to see that this construction extends to morphisms between $H_1$-torsors and commutes with arbitrary base-changes.
}

\defi{}{\cite[definition 1]{FaltingsLoop}}{
Let $G$ be a connected reductive group over $\Fq$. \\
a) $L^+G$ is the infinite dimensional group scheme over $\Fq$ representing the fpqc-sheaf of groups
\[S \mapsto G(\Gamma(S, \MC{O}_S)[[z]]).\]
b) $LG$ is the ind-scheme (of ind-finite type) over $\Fq$ representing the sheafification of the fpqc-presheaf of groups
\[S \mapsto G(\Gamma(S, \MC{O}_S)((z))).\]
c) There is a canonical inclusion $L^+G \subset LG$. If $\MC{G}$ is a $L^+G$-torsor over a scheme $S$, then $\MC{LG} \coloneqq \MC{G} \times^{L^+G} LG$ (cf. the previous construction) is called the $LG$-torsor associated to $\MC{G}$.
}

\rem{}{
Usually $G$-torsors (and in particular global $G$-shtukas, cf. \ref{def:GlobShtuka}) are denoted by $\MS{G}$, while $\MC{G}$ is used for $L^+G$-torsors (and in particular local $G$-shtukas, cf. \ref{def:LocShtuka}).
}

As in the case of smooth affine group schemes of finite type, we can omit specifying the actual topology by the following 

\lem{\label{lem:LocShtukaTopol}}{}{
Using the notation above,
\[\check{H}^1(S_{\acute{e}t}, L^+G) = \check{H}^1(S_{fppf}, L^+G) = \check{H}^1(S_{fpqc}, L^+G)\]
i.e. $L^+G$-torsors for the \'etale, fppf- and fpqc-topology are equivalent.
}

\prooof
The proof of \cite[proposition 2.2]{HaVi} holds for any smooth affine group scheme. \exit
\vspace{3mm} \\
Similarly to the construction above, we may define functors from torsors to vector bundles using representations:

\lem{\label{const:GenTorsorVector}}{}{
Let $G$ be a connected reductive group over $\Fq$ and fix a representation $\rho: G \to GL(V)$ for a finite-dimensional $E$-vector space $V$ over a finite field extension $E/\Fq$. \\
a) The representation $\rho$ induces a canonical functor from $G$-torsors over an $E$-scheme $S$ to line bundles of rank $\dim V$ over $S$. \\
b) The representation $\rho$ induces a canonical functor from $L^+G$-torsors over an $E$-scheme $S$ to locally free $\MC{O}_S[[z]] = \MC{O}_S \otimes_E E[[z]]$-modules over $S$. \\
c) The representation $\rho$ induces a canonical functor from $LG$-torsors over an $E$-scheme $S$ to locally free $\MC{O}_S((z)) = \MC{O}_S \otimes_E E((z))$-modules over $S$. \\
All these functors are compatible for passing from $G$-torsors to $L^+G$-torsors and $LG$-torsors via the inclusions $G \subset L^+G \subset LG$.
}

\prooof
We will prove only part b) as the rest follows in the same way and the compatibilities are obvious. So let $\MC{G}$ be a $L^+G$-torsor over $S$ and consider the sheaf $\MC{G} \times^{L^+G} V$ over $S$ defined as the sheafification of the presheaf
\[S' \mapsto \{L^+G(S')\op{-orbits \; in \,} \MC{G}(S') \times (V \otimes_{E} \MC{O}_S[[z]])(S')\}\]
where we view $V$ as a constant sheaf on $S$ and take the induced $L^+G$-action on $V \otimes_{E} \MC{O}_S[[z]]$. Then $\MC{G} \times^{L^+G} V$ trivializes over the same covers as $\MC{G}$ did, hence is indeed locally free. Moreover it is an easy exercise to see that this construction extends to morphisms between $L^+G$-torsors. \exit

\subsection{Local \textit{G}-shtukas and bounded quasi-isogenies}\label{subsec:LocalDef}
We will now recall the definition of a local $G$-shtuka and morphisms between them, called quasi-isogenies. But it turns out, that most moduli spaces for local (and global) $G$-shtukas exist only as ind-schemes. To remedy this, boundedness conditions are introduced in order to specify subspaces of such moduli spaces, which exist as (formal) schemes locally of finite type. To give us more flexibility we start by defining bounded morphisms between $LG$-torsors associated to arbitrary $L^+G$-torsors, and specify then to the case of local $G$-shtukas and, in section \ref{subsec:GlobalBounds}, to global $G$-shtukas. \\
Recall that $\sigma$ denotes the absolute $q$-Frobenius and $G$ is a connected reductive group over $\Fq$, that is not necessarily split. We fix a finite field extension $E$ of $\Fq$ (as usual viewed as a subfield of $\ACFq$).

\defi{\label{def:LocShtuka}}{cf. \cite[definitions 3.1 and 3.8]{HaVi}}{
a) A local $G$-shtuka over a DM-stack $S$ over $E$ is a pair $(\MC{G}, \varphi)$ consisting of a $L^+G$-torsor $\MC{G}$ on $S$ and an isomorphism $\varphi: \sigma^*\MC{LG} \to \MC{LG}$ of associated $LG$-torsors. \\
Denote by $Sht_G$ the ind-stack representing the functor of local $G$-shtukas up to isomorphism. \\
b) A local $G$-shtuka $(\MC{G}, \varphi)$ is called \'etale if $\varphi$ is induced from a morphism $\varphi: \sigma^*\MC{G} \to \MC{G}$ between $L^+G$-torsors. The ind-stack of \'etale local $G$-shtukas is denoted by $\acute{E}tSht_G$. \\
c) A quasi-isogeny between local $G$-shtukas $\alpha: (\MC{G}, \varphi) \to (\MC{G}', \varphi')$ over $S$ is an isomorphism of the associated $LG$-torsors $\alpha: \MC{LG} \to \MC{LG}'$ satisfying $\varphi' \circ \sigma^*\alpha = \alpha \circ \varphi$. \\
The set of quasi-isogenies between two local $G$-shtukas is denoted by $QIsog((\MC{G}, \varphi), (\MC{G}', \varphi'))$.
}

\rem{\label{rem:LocShtukaFrob}}{
i) We claim above that $Sht_G$ is an ind-stack. Indeed it is not hard to see that it is given by a $LG$-bundle over a substack of the stack of $L^+G$-torsors. But for all of our applications it suffices to know that $Sht_G$ is a category fibered over the category of DM-stacks over $E$, which is obviously true. \\    
ii) The Frobenius-isomorphism $\varphi$ of a local $G$-shtuka $(\MC{G}, \varphi)$ over $S$ can be viewed as a quasi-isogeny $\varphi: (\sigma^*\MC{G}, \sigma^*\varphi) \to (\MC{G}, \varphi)$. \\
iii) Note that the definition makes sense for arbitrary connected linear algebraic groups over $\Fq$. In particular, if $P \subset G$ is any parabolic subgroup, then the notion of a local $P$-shtuka (appearing in section \ref{subsec:SlopeDivisible}) makes sense.
}

Though we will mainly be interested in local $G$-shtukas, we will need at several points some more data respectively more general notions:

\defi{\label{def:LocShtukaGeneral}}{}{
a) Let $P \subset G$ be any subgroup (which will usually be parabolic) and $H \subset L^+P$ any open subgroup. A local $G$-shtuka with $H$-structure over a DM-stack $S$ over $E$ is a triple $(\MC{G}, \varphi, \MC{H})$ consisting of a local $G$-shtuka $(\MC{G}, \varphi)$ and an $H$-torsor $\MC{H}$ over $S$ together with a fixed isomorphism $\MC{H} \times^{H} L^+G \cong \MC{G}$. \\
A quasi-isogeny is a quasi-isogeny between the local $G$-shtukas (after forgetting the $H$-structure). \\
An isomorphism is an isomorphism of the local $G$-shtukas, which is induced from an isomorphism of the $H$-torsors. \\
Denote the ind-stack of local $G$-shtukas with $H$-structure by $H\op{-}Sht_G$. \\
b) Let $H$ be any affine algebraic group over $\Fq$. Then an \'etale local $H$-shtuka over $S$ is a pair $(\MC{H}, \varphi_{\MC{H}})$ consisting of an $H$-torsor $\MC{H}$ over $S$ and a $\sigma$-linear isomorphism $\varphi_{\MC{H}}: \sigma^*\MC{H} \to \MC{H}$. \\
Denote the ind-stack of \'etale local $H$-shtukas by $H\op{-}\acute{E}tSht$.
}

\rem{}{
Local $G$-shtukas with $H$-structure will appear in the definition of Igusa-varieties, cf. section \ref{subsec:IwahoriStructure}. \'Etale local $H$-shtukas will already appear in the discussion of adelic levels structures of global $G$-shtukas, cf. section \ref{subsec:AdelicLevel}.
}

We turn now to the definition of bounds. Apart from generalizing to non-split groups, we use essentially the construction introduced in \cite{HaVi}. Nevertheless there is one slight change to the definition: Boundedness conditions for local $G$-shtukas contain certain compatibilities in the set of orbits $\pi_1(G)/\Gamma$. Nevertheless they should correspond to boundedness conditions on global $G$-shtukas (cf. section \ref{subsec:GlobalBounds}), which seem not to be able to detect differences in the torsion of $\pi_1(G)/\Gamma$. This forces us to work with $\pi_1(G)_{\M{Q}}/\Gamma = (\pi_1(G) \otimes \M{Q})/\Gamma$ rather than $\pi_1(G)/\Gamma$ itself. \vspace{2mm} \\
We will work over the category $Nilp_{\ENil}$ of DM-stacks $S$ over $\Sp \ENil$ such that $\zeta$ is locally nilpotent on $S$. Consider now two $L^+G$-torsors $\MC{G}$ and $\MC{G}'$ over a DM-stack $S \in Nilp_{\ENil}$ and $\alpha: \MC{LG} \to \MC{LG}'$ a morphism between the associated $LG$-torsors.

\const{\label{Con:LocalHodge}}{\textbf{The Hodge point of $\boldsymbol{\alpha}$} (cf. \cite[definition 3.3]{HaVi})}{
Consider a geometric point $\bar{s}: \Sp k \to S$ mapping to a point $s \in S$. Choose any trivializations $\MC{G}_{\bar{s}} \cong L^+G_{\bar{s}}$ and $\MC{G}'_{\bar{s}} \cong L^+G_{\bar{s}}$ of the restriction of the $L^+G$-torsors to $\bar{s}$. Then the restriction of $\alpha$ to $\bar{s}$ defines a morphism
\[\alpha_{\bar{s}}: LG_{\bar{s}} \cong \MC{LG}_{\bar{s}} \to \MC{LG}'_{\bar{s}} \cong LG_{\bar{s}}\]
i.e. an element $\alpha_{\bar{s}} \in LG(k)$. Note now that $G$ splits over $k$ and that $L^+G \subset LG$ is a special subgroup. Hence the Cartan decomposition yields
\[LG(k) = \bigcup_{\mu \in X_*(T)_{\rm dom}} L^+G(k) z^\mu L^+G(k)\]
cf. \cite[3.3.3]{TitsGroups}. This associates to $\alpha$ and $\bar{s}$ an element in $X_*(T)_{\rm dom} \subset X_*(T)$, which is independent of the trivialization of the two $L^+G$-torsors. Projecting down to the set $X_*(T)/\Gamma$ even eliminates the dependence on the choice of the geometric point $\bar{s}$ mapping to $s$. Hence this defines a homomorphism (of sets)
\[\mu(\alpha): S \to X_*(T)/\Gamma\]
which we call the Hodge point of $\alpha$. \\
The canonical projection $X_*(T) \to \pi_1(G) \to \pi_1(G)_{\M{Q}} = \pi_1(G) \otimes \M{Q}$ is $\Gamma$-equivariant. Hence we may compose the morphism above with $X_*(T)/\Gamma \to \pi_1(G)_{\M{Q}}/\Gamma$ to get
\[[\mu(\alpha)]: S \to \pi_1(G)_{\M{Q}}/\Gamma.\]
}

\lem{}{}{
$\pi_0(LG) \cong \pi_1(G)/\Gamma$.
}

\prooof
By \cite[theorem 5.1]{PaRa1} there is an isomorphism $\pi_0(LG \times_E \ACFq) \cong \pi_1(G)$, which is by construction $\Gamma$-equivariant. All connected components of $LG \times_E \ACFq$ in one $\Gamma$-orbit have the same image in $LG$, hence map into the same connected component of $LG$. Moreover $LG \times_E \ACFq \to LG$ is a closed morphism as a limit of finite morphism. So each connected component of $LG \times_E \ACFq$ actually maps surjectively onto a connected component of $LG$. This implies that the preimages of elements under $\pi_0(LG \times_E \ACFq) \to \pi_0(LG)$ consist of precisely one $\Gamma$-orbit, because this is true for preimages of single points under $LG \times_E \ACFq \to LG$. Thus we obtain isomorphisms
\[\pi_0(LG) \cong \pi_0(LG \times_E \ACFq)/\Gamma \cong \pi_1(G)/\Gamma.\] \exit

Using this isomorphism, the construction of the Hodge point $[\mu(\alpha)]$ can be reformulated in a slightly more conceptual way: The morphism $\alpha$ induces a map $S \to \pi_0(LG)$ by trivializing $\alpha$ locally. It is well-defined, because the image for different trivializations is contained in one connected $L^+G$-conjugacy class. Then by the previous lemma we get the map
\[S \to \pi_0(LG) \cong \pi_1(G)/\Gamma \to \pi_1(G)_{\M{Q}}/\Gamma,\]
which equals $[\mu(\alpha)]$ on geometric points by a straight-forward computation. 

\ignore{
\rem{\label{Rem:LocalHodgeRem}}{
This construction can be reformulated in a slightly more conceptual way: The morphism $\alpha$ induces a map $S \to \pi_0(LG)$ by trivializing $\alpha$ locally. It is well-defined, because the image for different trivializations is contained in one connected $L^+G$-conjugacy class. Then by the previous lemma we get the map
\[S \to \pi_0(LG) \cong \pi_1(G)/\Gamma \to \pi_1(G)_{\M{Q}}/\Gamma,\]
which equals $[\mu(\alpha)]$ on geometric points by a straight-forward computation.
\ignore{
i) Note that $[\mu(\alpha)]$ is essentially given by the Kottwitz homomorphism $\kappa$: Indeed assume for a moment that we have fixed some trivialization $\MC{G} \cong L^+G \times S$. Then for any field extension $E'/E$ (inside $\ACFq$) such that $G$ splits over $E'$ we have maps
\[S(E') \to LG(E') \xrightarrow{\kappa} \pi_1(G) \]
where the first map is given by representing $\alpha$ by an element in $LG(E')$ and the second map is given by the Kottwitz homomorphism $\kappa$. Note that due to the splitting assumption $\Gal(\ACFq/E')$ acts trivially on $\pi_1(G)$ and we may omit taking coinvariants. Then we get for each $E'$ a commutative diagram of sets
\[
 \begin{xy} \xymatrix{
   S(E') \ar^-{\kappa}[r] \ar[d] & \pi_1(G) \ar[d] \\
   S \ar^-{[\mu(\alpha)]}[r] & \pi_1(G)_{\M{Q}}/\Gamma  }
 \end{xy} 
\]
}
}
}

\lem{\label{lem:BoundLocalHodgeConstant}}{cf. \cite[proposition 3.4]{HaVi}}{
The map $[\mu(\alpha)]$ is locally constant.
} 

\prooof
Let $f: S \times_E \Sp \ACFq \to S$ be the base-change to $\ACFq$. Then by definition
\[[\mu(f^*\alpha)] = [\mu(\alpha)] \circ f: S \times_E \Sp \ACFq \to \pi_1(G)_{\M{Q}}/{\Gamma}\]
where $f^*\alpha: \MC{LG} \times_E \Sp \ACFq  \to \MC{LG}' \times_E \Sp \ACFq$ is the pullback of $\alpha$.
Thus it suffices to show that $[\mu(f^*\alpha)]$ is locally constant. But as $G$ splits over $\ACFq$, this follows from proposition 3.4 in \cite{HaVi}, which even states that the morphism to $\pi_1(G)$ itself is locally constant. \exit 

\rem{}{
Alternatively use the description in the previous remark: The map $S \to \pi_0(LG)$ factors by definition \'etale-locally over $LG$, hence it factors as $S \to \pi_0(S) \to \pi_0(LG)$, too. So the Hodge point $[\mu(\alpha)]$ factors over $\pi_0(S)$ as well.
}
%
%
Consider any dominant character $\lambda \in X^*(T)$. It is defined over some field extension $\F/\Fq$ and $(-\lambda)_{\rm dom}$, i.e. the dominant element in the orbit of $-\lambda$ under the Weyl group, defines a representation of the Borel $\ov{B}$ opposite to $B$
\[(-\lambda)_{\rm dom}: \ov{B} \times_{\Fq} \Sp \F \to \M{G}_m \times_{\Fq} \Sp \F\]
which is trivial on the unipotent radical of $\ov{B} \times_{\Fq} \Sp \F$. Thus we may consider the Weyl module of highest weight $\lambda$ 
\[V_G(\lambda) \coloneqq ({\rm Ind}^G_{\ov{B}}(-\lambda)_{\rm dom})^\vee\]
which is defined over $\F$. \\
Set now $E' = E \cdot \F$ the composite of $E$ and $\F$ in $\ACFq$. Then we may associate to any $L^+G$-torsor $\MC{G}$ over an $E$-scheme $S$ the locally free sheaves over $S \times_E \Sp E'$
\[\MC{G}_\lambda \coloneqq \MC{G} \times^{L^+G} V_G(\lambda) \qquad {\rm and} \qquad \MC{LG}_\lambda \coloneqq \MC{LG} \times^{LG} V_G(\lambda)\]
using construction \ref{const:GenTorsorVector}. Note that $\MC{G}_\lambda$ is naturally a sub-$\MC{O}_S[[z]]$-module of $\MC{LG}_\lambda$. \\
By functoriality any morphism $\alpha: \MC{LG} \to \MC{LG}'$ between two $LG$-torsors induces now an isomorphism
\[\alpha_\lambda: \MC{LG}_{\lambda} \to \MC{LG}'_{\lambda}\]
of sheaves of $\MC{O}_S((z))$-modules. 

\defi{\label{Def:BoundLocal}}{}{
Let $\mu \in X_*(T)_{\rm dom}$ be a dominant cocharacter which can be defined over $E$. \\
a) Let $\MC{G}$ and $\MC{G}'$ be two $L^+G$-torsors over DM-stack $S \in Nilp_{\ENil}$ and $\alpha: \MC{LG} \to \MC{LG}'$ a morphism between the associated $LG$-torsors. Then $\alpha$ is bounded by $\mu$ if
\begin{enumerate}[label = \arabic*.]
 \item for any dominant weight $\lambda \in X^*(T)_{\rm dom}$ defined over some finite field extension $E'/E$
\[\alpha_{\lambda}(\MC{G}_\lambda) \subseteq (z - \zeta)^{-\langle(-\lambda)_{\rm dom}, \mu\rangle}\MC{G}_\lambda' \coloneqq \MC{G}_\lambda' \otimes_{\MC{O}_{S \times_E \Sp E'}[[z]]} (z - \zeta)^{-\langle(-\lambda)_{\rm dom}, \mu\rangle}\]
inside $\MC{G}_\lambda' \otimes_{\MC{O}_{S \times_E \Sp E'}[[z]]} \MC{O}_{S \times_E \Sp E'}((z)) = \MC{LG}_\lambda'$.
 \item $[\mu(\alpha)](s) = [\mu] \in \pi_1(G)_{\M{Q}}/\Gamma$ for all points $s \in S$, where $[\mu]$ denotes the image of $\mu$ in $\pi_1(G)_{\M{Q}}/\Gamma$.
\end{enumerate}
b) A quasi-isogeny $\alpha: (\MC{G}, \varphi) \to (\MC{G}', \varphi')$ is bounded by $\mu$ if the morphism $\alpha: \MC{LG} \to \MC{LG}'$ between associated the $LG$-torsors is bounded by $\mu$. \\
c) A local $G$-shtuka $(\MC{G}, \varphi)$ is bounded by $\mu$ if $\varphi$ is bounded by $\mu$ when considered as a quasi-isogeny between $(\sigma^*\MC{G}, \sigma^*\varphi)$ and $(\MC{G}, \varphi)$. Denote by $Sht_G^{\preceq \mu}$ the stack of local $G$-shtukas bounded by $\mu$.
}

\rem{\label{Rem:BoundOldDef}}{
i) We require $\mu$ to be defined over $E$ in order to ensure that condition \textit{1.} is independent with respect to pullback along a morphism induced by some element in $\Gamma_E$, cf. remark \ref{Rem:BoundGaloisInvariant} as well. Moreover condition \textit{1.} is independent of the choice of the field extension $\F/\Fq$ over which $\lambda$ is defined. Condition \textit{2.} implies in particular that the Hodge points $[\mu(\alpha)](s)$ are $\Gamma$-orbits consisting of only one element in $\pi_1(G)_{\M{Q}}$. \\
The boundedness condition does not depend on the choice of the field $E$ (as long as the cocharacter $\mu$ can be defined over it). The only non-obvious part is condition \textit{2.}, where this follows from the fact that the $\Gamma$-orbit of $\mu$ in $\pi_1(G)_{\M{Q}}$ has always exactly one element, whatever base field we take. \\
ii) The boundedness conditions can be checked \'etale (or fpqc)-locally on $S$. For conditions \textit{1.} this is implied by part i) of this remark and for condition \textit{2.} it follows from lemma \ref{lem:BoundLocalHodgeConstant}. \\
iii) Definition $4.5$ in \cite{HarRad1} is a much more conceptual treatment of bounds via ind-subschemes $\hat{Z} \subset \widehat{Gr}$ (cf. \ref{subsec:RapoZinkSpace} for the definition of $\widehat{Gr}$), though is has two main disadvantages for our purposes: First of all it is a priori not obvious to us, that there exist sufficiently many bounds $\hat{Z}$. Nevertheless combining proposition \ref{prop:BoundClosed} and theorem \ref{thm:RapoZinkUnbounded} shows, that boundedness by $\mu$ (in our sense) can be expressed by boundedness by some $\hat{Z}$ as in \cite{HarRad1}. Secondly (and more seriously) we do not know whether an arbitrary bound $\hat{Z}$ behaves well when passing to inverse morphisms. 
}

\lem{\label{Lem:BoundWeightGenerate}}{\cite[lemma 3.7]{HaVi}}{
Condition \textit{1.} holds for every dominant weight if and only if it holds for a subset $\Lambda \subset X^*(T)_{\rm dom}$ of dominant weights that generates the monoid of dominant weights. $\Lambda$ can be chosen to be finite.
}

\prooof
The proof in \cite{HaVi} can be used verbatim for non-split $G$ as well. \exit

\rem{\label{Rem:BoundGaloisInvariant}}{
Condition \textit{1.} holds for $\lambda \in X^*(T)_{\rm dom}$ if and only if it holds for $\tau.\lambda$ for any $\tau \in \Gamma = \Gal(\ACFq/E)$: \\
Let $\F/\Fq$ be a finite field extension, such that $\lambda$ (and hence $\tau.\lambda$) can be defined over $\F$. Let $E' = E \cdot \F$ as above and $\tau: S \times_E \Sp E' \to S \times_E \Sp E'$. As $\alpha$ is defined over $E$, i.e. fixed under $\tau$, we have that
\[\alpha_{\lambda}(\MC{G}_\lambda) \subseteq (z - \zeta)^{-\langle(-\lambda)_{\rm dom}, \mu\rangle}\MC{G}_\lambda'\]
if and only if
\[\alpha_{\lambda}(\tau^*\MC{G}_\lambda) = \tau^*(\alpha_{\lambda}(\MC{G}_\lambda)) \subseteq \tau^*\left((z - \zeta)^{-\langle(-\lambda)_{\rm dom}, \mu\rangle} \MC{G}_\lambda'\right) = (z - \zeta)^{-\langle(-\lambda)_{\rm dom}, \mu\rangle} \tau^*(\MC{G}_\lambda')\]
Furthermore $\tau$ induces an isomorphism $V_G(\lambda) \cong V_G(\tau^{-1}.\lambda)$ as representations in $E$-vector spaces (but changes of course the $E'$-structure on it). This gives isomorphisms $\tau^*\MC{G}_{\lambda} \cong \MC{G}_{\tau^{-1}.\lambda}$ and $\tau^*\MC{LG}_{\lambda} \cong \MC{LG}_{\tau^{-1}.\lambda}$ as sheaves over $S \times_E \Sp E'$. In particular the inclusion above may be rewritten as
\[\alpha_{\lambda}(\MC{G}_{\tau^{-1}.\lambda}) \subseteq (z - \zeta)^{-\langle(-\lambda)_{\rm dom}, \mu\rangle} \MC{G}_{\tau^{-1}.\lambda}'.\]
Using $\langle(-\tau^{-1}.\lambda)_{\rm dom}, \mu\rangle = \langle \tau^{-1}.(-\lambda)_{\rm dom}, \mu\rangle = \langle (-\lambda)_{\rm dom}, \tau.\mu\rangle = \langle (-\lambda)_{\rm dom}, \mu\rangle$ by $\Gal(\ACFq/E)$-invariance of $\mu$, this is precisely the boundedness condition \textit{1.} for $\tau.\lambda$.
}

\subsection{Bounds for inverse morphisms}\label{subsec:LocalBoundsInv}
The next aim is to show that the inverse of a morphism bounded by $\mu$ is bounded by $(-\mu)_{\rm dom} \in X_*(T)_{\rm dom}$. If Weyl modules would be self-dual, this would be an easy exercise. But as they are not, we will first prove this fact for $G = GL_n$ and then transport the result first to arbitrary split groups and finally to any connected reductive $G$. \\
We still work with groups $G$ which are defined over $\Fq$. This restriction is hardly necessary in this section. For example (after suitable generalization of bounds) everything works as well for unramified groups over local fields. \\
The reader may replace the term '$L^+G$-torsor' by 'local $G$-shtuka' and 'morphism of associated $LG$-torsor' by 'quasi-isogeny' everywhere without changing (or simplifying) anything. In fact, except for one proof in section \ref{subsec:GlobalModuli}, these results are only used in the context of local $G$-shtukas.

\lem{\label{Lem:BoundGLEasy}}{}{
Let $G = GL_n$ with the maximal torus of diagonal matrices $T$ and the Borel of upper triangular matrices $B$, $Std = \Fq^n$ be the standard representation of $GL_n$. Let $\MC{G}$ and $\MC{G}'$ be two $GL_n$-torsors over an $E$-scheme $S$ and $\alpha: \MC{LG}\to \MC{LG}'$ a morphism between the associated $LGL_n$-torsors. Abbreviate the vector bundles constructed in \ref{const:GenTorsorVector} by $\MC{G}_{Std} = \MC{G} \times^{L^+GL_n} Std$, $\MC{LG}_{Std} = \MC{LG} \times^{LGL_n} Std$ (and similarly for $\MC{G}'$) and the induced morphism by $\alpha_{Std}: \MC{LG}_{Std} \to \MC{LG}'_{Std}$. Then $\alpha$ is bounded by a dominant element $\mu = (\mu_1, \ldots, \mu_n) \in X_*(T) \cong \M{Z}^n$ if and only if the following two conditions are satisfied: 
\begin{enumerate}[label = \arabic*.']
 \item $(\wedge^i \alpha_{Std})(\bigwedge^i \MC{G}_{Std}) \subseteq (z - \zeta)^{\mu_{n-i+1} + \ldots + \mu_n} \bigwedge^i \MC{G}_{Std}'$ for every $1 \leq i \leq n$.
 \item $[\mu(\alpha)](s) = [\mu] \in \pi_1(GL_n)$ for all points $s \in S$.
\end{enumerate}
}

\prooof
First note that $\pi_1(GL_n)$ is torsion free with trivial $\Gamma$-action, so conditions \textit{2.} and \textit{2.'} are clearly equivalent. \\
See now \cite[section 4]{HaVi} how to identify $V_{GL_n}(\delta_i) = \bigwedge^i Std$, $\MC{G}_{\delta_i} = \bigwedge^i \MC{G}_{Std}$ and $\alpha_{\delta_i} = \wedge^i \alpha_{Std}$ for the dominant weight $\delta_i = (1, \ldots, 1, 0, \ldots, 0)$ with the first $i$ entries equal to $1$. Furthermore note that $-\langle (-\delta_i)_{\rm dom}, \mu \rangle = \mu_{n-i+1} + \ldots + \mu_n$ for all $i$. Thus it is clear that the given conditions are necessary for $\alpha$ to be bounded by $\mu$. As all these $\delta_i$ for $1 \leq i \leq n$ together with $-\delta_n$ generate the monoid of dominant weights, the conditions are also sufficient by lemma \ref{Lem:BoundWeightGenerate} if they imply
\[\alpha_{-\delta_n}(\MC{G}_{-\delta_n}) \subseteq (z - \zeta)^{-N} \MC{G}_{-\delta_n}'.\]
for $N = \mu_1 + \ldots + \mu_n$. To prove this we first show that $\alpha_{\delta_n}$ induces an isomorphism on the integral level
\[\alpha_{\delta_n, 0}: \MC{G}_{\delta_n} \xrightarrow{\sim} (z - \zeta)^N \MC{G}'_{\delta_n}.\]
We only have to show surjectivity, which is a local property. Thus assume wlog. that $\MC{G}_{Std}, \MC{G}_{Std}' \cong \MC{O}_S[[z]]^n$ (for some chosen isomorphisms) and $\alpha_{Std}$ is given by an element in $LGL_n(S)$. Then $\bigwedge^n \MC{G}_{Std} \cong \MC{O}_S[[z]]$ and $\wedge^n \alpha_{Std}$, i.e. $\alpha_{\delta_n, 0}$, is simply multiplication with $\det(\alpha_{Std})$. After pulling back to any geometric point $\ov{s}$ with image $s$ we have $\det(\alpha_{Std})_{\ov{s}} = a_{\ov{s}} \cdot z^{\langle \delta_n, \mu(\alpha)(s) \rangle}$ for some unit $a_{\ov{s}}$. As the pairing of $\delta_n$ with any coroot vanishes, $\langle \delta_n, \mu(\alpha)(s) \rangle$ depends only on the image of $\mu(\alpha)(s)$ in $\pi_1(GL_n)$. Thus
\[\det(\alpha_{Std})_{\ov{s}} = a_{\ov{s}} \cdot z^{\langle \delta_n, \mu(\alpha)(s) \rangle} = a_{\ov{s}} \cdot z^{\langle \delta_n, \mu \rangle} = a_{\ov{s}} \cdot z^N.\]
Together with $\det(\alpha_{Std}) \in (z - \zeta)^N\MC{O}_S[[z]]$ this implies
\[\det(\alpha_{Std}) = a \cdot (z - \zeta)^N \qquad {\rm for \; some \;} a \in \MC{O}_S[[z]]^{\times}\]
and multiplication with $\det(\alpha_{Std})$, i.e. $\alpha_{\delta_n, 0}$, is indeed surjective. \\
Now observe that $V_{GL_n}(\delta_n)$ is one-dimensional which implies $V_{GL_n}(-\delta_n) \cong V_{GL_n}(\delta_n)^\vee$. Thus we can identify $\MC{G}_{-\delta_n} \cong \MC{G}_{\delta_n}^\vee$ and $\alpha_{-\delta_n} = (\alpha_{\delta_n}^{-1})^\vee$. Thus by dualizing $\alpha_{\delta_n}$, $\alpha_{-\delta_n}$ restricts to an isomorphisms
\[(z - \zeta)^N \MC{G}_{-\delta_n} \xrightarrow{\sim} \MC{G}_{-\delta_n}'\]
which gives after multiplying with $(z - \zeta)^{-N}$
\[\alpha_{-\delta_n}(\MC{G}_{-\delta_n}) \subseteq (z - \zeta)^{-N}\MC{G}_{-\delta_n}'\]
as desired. \exit

\rem{\label{rem:BoundGLEasyAlternative}}{
A very similar result is shown in \cite[lemma 4.3]{HaVi}. There condition \textit{2.'} is replaced by the statement that $\wedge^n \alpha_{Std}$ restricts to an isomorphism between suitable $\MC{O}_S[[z]]$-submodules. Nevertheless with view towards lemma \ref{Lem:BoundAfterRepr} the version stated above will be easier to apply. 
}

\prop{\label{Prop:BoundInverseGL}}{}{
Let $\alpha: \MC{LG} \to \MC{LG}'$ be a morphism between $LGL_n$-torsors as in the previous lemma. Assume $\alpha$ is bounded by a dominant coweight $\mu$. Then its inverse $\alpha^{-1}: \MC{LG}' \to \MC{LG}$ is bounded by $(-\mu)_{\rm dom}$.
}

\prooof
We will use the alternative description shown in the previous lemma. As both conditions are local on $S$ we may assume again that $\MC{G}, \MC{G}' \cong L^+GL_n \times S$ and $\alpha \in LGL_n(S)$. As the Kottwitz homomorphism $GL_n \to \pi_1(GL_n)$ is a group morphism we get $[\mu(\alpha^{-1})](s) = -[\mu(\alpha)](s)$ for every point $s \in S$. Thus
\[[\mu(\alpha^{-1})](s) = -[\mu(\alpha)](s) = -[\mu] = [(-\mu)_{\rm dom}] \in \pi_1(GL_n)\]
Thus the second condition is satisfied. \\
If $\mu = (\mu_1, \ldots, \mu_n)$, then $(-\mu)_{\rm dom} = (-\mu_n, \ldots, -\mu_1)$. Condition \textit{1.'} for $\alpha^{-1}$ and the coweight $(-\mu)_{\rm dom}$ can be rephrased in saying that for every $i \in \{1, \ldots, n\}$ the determinant of every $i \times i$-minor of the matrix $\alpha^{-1}$ is contained in $(z - \zeta)^{-\mu_i - \ldots - \mu_1}\MC{O}_S[[z]]$. Thus let $I, J \subset \{1, \ldots, n\}$ be two subsets with $i$ elements and let $\alpha^{-1}_{I, J}$ be the minor of $\alpha^{-1}$ which contains the rows in $I$ and the columns in $J$. Denote the $(n-i) \times (n-i)$-minor of $\alpha$ containing the rows in $I^c = \{1, \ldots, n\} \setminus I$ and the columns $J^c = \{1, \ldots, n\} \setminus J$ by $\alpha_{I^c, J^c}$. Then one has the formula
\[\det (\alpha^{-1}_{I, J}) = \pm \det (\alpha_{I^c, J^c}) \cdot \det (\alpha^{-1})\]
which is sometimes called the identity of Jacobi. But by condition \textit{1.'} for $\alpha$ we know $\det (\alpha_{I^c, J^c}) \in (z - \zeta)^{\mu_{i+1} + \ldots + \mu_n}\MC{O}_S[[z]]$. Furthermore $\det (\alpha) \in (z - \zeta)^{\mu_1 + \ldots + \mu_n}\MC{O}_S[[z]] \setminus (z - \zeta)^{\mu_1 + \ldots + \mu_n + 1}\MC{O}_S[[z]]$ using the same argument involving the Hodge point as in the previous lemma. Thus $\det (\alpha^{-1}) \in (z - \zeta)^{-\mu_1 - \ldots - \mu_n}\MC{O}_S[[z]]$, which finally gives the desired
\[\det (\alpha^{-1}_{I, J}) = \pm \det (\alpha_{I^c, J^c}) \det (\alpha^{-1}) \in (z - \zeta)^{\mu_{i+1} + \ldots + \mu_n} (z - \zeta)^{-\mu_1 - \ldots - \mu_n}\MC{O}_S[[z]] = (z - \zeta)^{-\mu_i - \ldots - \mu_1}\MC{O}_S[[z]]\] 
Thus \textit{1.'} holds for $\alpha^{-1}$ and $(-\mu)_{\rm dom}$. \exit

\rem{}{
A similar argument already appears in \cite[lemma 2.1.6]{HartlPeriod} although in the context of local shtukas (in the sense of \cite[definition 2.1.1]{HartlPeriod} or \cite[definition 4.1]{HaVi}) over formal schemes.
} 
\vspace{2mm}
We now want to extend this proposition to all groups $G$, which split over $E$. In particular all Weyl module representations are already defined over $E$. To do so we need several auxiliary lemmas. Most of the following arguments work as well in the non-split case after possibly enlarging the field in order to define the representations. Nevertheless as the non-split case follows quite easily from the split case later on, it is more convenient not to present them in the general setting. 

\lem{}{}{
Let $G$ be any reductive group over $\Fq$, which splits over $E$. Assume that $0 \to V' \to V \to V'' \to 0$ is a short exact sequence of finite-dimensional $E$-representations of $G$. Let $N$ be an integer and let $\alpha: \MC{LG} \to \MC{LG}'$ be a morphism between $LG$-torsors associated to some $L^+G$-torsors $\MC{G}$, $\MC{G}'$. Then the following are equivalent: \\
a) $\alpha_V(\MC{G} \times^{L^+G} V) \subseteq (z - \zeta)^N (\MC{G}' \times^{L^+G} V)$ \\
b) $\alpha_{V'}(\MC{G} \times^{L^+G} V') \subseteq (z - \zeta)^N (\MC{G}' \times^{L^+G} V') \quad $ and $\quad \alpha_{V''}(\MC{G} \times^{L^+G} V'') \subseteq (z - \zeta)^N (\MC{G}' \times^{L^+G} V'')$ \\
where we denote by $\alpha_V$ (and similarly for $\alpha_{V'}$ and $\alpha_{V''}$) the morphism induced by $\alpha$.
}

\prooof
First note that the functor mapping finite-dimensional $E'$-representations $V$ of $G$ to the vector bundle $\MC{G} \times^{L^+G} V$ is exact, because over an \'etale cover it simply takes $V$ to $V \otimes_{E'} \MC{O}_S[[z]]$. Thus the given short exact sequence of representations of $G$ induces a commutative diagram with exact rows
\[
 \begin{xy}
  \xymatrix{
   0 \ar[r] & \MC{G} \times^{L^+G} V' \ar[r] \ar^{\alpha_{V'}}[d] & \MC{G} \times^{L^+G} V \ar[r] \ar^{\alpha_V}[d] & \MC{G} \times^{L^+G} V'' \ar[r] \ar^{\alpha_{V''}}[d] & 0 \\
   0 \ar[r] & \MC{LG}' \times^{LG} V' \ar[r] & \MC{LG}' \times^{LG} V \ar[r] & \MC{LG}' \times^{LG} V'' \ar[r] & 0 } 
 \end{xy} 
\]
But
\[0 \to (z - \zeta)^N (\MC{G}' \times^{L^+G} V') \to (z - \zeta)^N (\MC{G}' \times^{L^+G} V) \to (z - \zeta)^N (\MC{G}' \times^{L^+G} V'') \to 0 \]
is a short exact sequence of subsheaves of the sheaves in the lower row. Thus $\alpha_V(\MC{G} \times^{L^+G} V)$ lies in $(z - \zeta)^N (\MC{G}' \times^{L^+G} V)$ if and only if the analogous statements holds for $\alpha_{V'}$ and $\alpha_{V''}$. \exit

\lem{\label{Lem:BoundAlternating}}{}{
Assume that $G$ splits over $E$. Let $\alpha: \MC{LG} \to \MC{LG}'$ be a morphism of $LG$-torsors over some $E$-scheme $S$ (associated as usual to some $L^+G$-torsors) which is bounded by $\mu$. Fix some $\lambda \in X^*(T)_{\rm dom}$ which then gives a representation $V_G(\lambda)$ over $E$. Let $1 \leq i \leq \dim V_G(\lambda)$ and let $V = \bigwedge^i V_G(\lambda)$ with $G$ acting on every factor. Finally set $N = \min_{S_i} \{\sum_{\eta \in S_i} \langle \eta, \mu \rangle\}$ where $S_i$ runs over all subsets of cardinality $i$ of the set of weights occurring in $V_G(\lambda)$ (taken with multiplicities). Then
\[\alpha_V(\MC{G} \times^{L^+G} V) \subseteq (z - \zeta)^N (\MC{G}' \times^{L^+G} V)\]
}

\rem{}{
The weights $\sum_{\eta \in S_i} \eta$ are exactly the weights occurring in the representation $V$. Thus
\[N = \min_{\eta_V} \{\langle \eta_V, \mu \rangle\}\]
taking the minimum over all weights $\eta_V$ of $V$. On the other hand, the set of weights is invariant under the action of the longest element $w_0$ of the Weyl group. Thus we may also write:
\[N = \min_{S_i} \left\{\sum\nolimits_{\eta \in S_i} \langle w_0\eta, \mu \rangle \right\} = \min_{S_i} \left\{-\sum\nolimits_{\eta \in S_i} \langle (-\eta)_{\rm dom}, \mu \rangle \right\} = \min_{\eta_V} \left\{-\langle (-\eta_V)_{\rm dom}, \mu \rangle \right\}\]
Instead of taking all weights $\eta_V$ in the rightmost expression it suffices to consider only highest weights of $V$.
}

\prooof
As $B$ is a Borel, every $E$-representation of $G$ admitting a basis of $T$-eigenvectors has a $B$-stable line. In particular it contains a quotient of a Weyl module as a subrepresentation, whose highest weight equals some weight appearing in the original representation. Thus every finite-dimensional such $E$-representation of $G$ admits a finite filtration with quotients of Weyl modules as composition factors. Let $0 = W_0 \subset W_1 \subset \ldots \subset W_n = V$ be such a filtration for $V = \bigwedge^i V_G(\lambda)$. We will prove by induction that 
\[\alpha_{W_j}(\MC{G} \times^{L^+G} W_j) \subseteq (z - \zeta)^N (\MC{G}' \times^{L^+G} W_j)\]
for all $j$. Note first that this is trivial for $j = 0$ and gives the statement of the lemma for $j = n$. Assume now it holds for $W_{j-1}$ for some $j \geq 1$ and consider the short exact sequence
\[\hspace{1.5cm} 0 \to W_{j-1} \to W_j \to W_j/W_{j-1} \to 0 \hspace{1cm} (\star)\]
By induction we have 
\[\alpha_{W_{j-1}}(\MC{G} \times^{L^+G} W_{j-1}) \subseteq (z - \zeta)^N (\MC{G}' \times^{L^+G} W_{j-1}).\]
On the other hand we have by choice of the filtration a short exact sequence
\[0 \to V'_G(\eta_V) \to V_G(\eta_V) \to W_j/W_{j-1} \to 0 \]
for some weight $\eta_V$ of $V$ (depending on $j$) and some subrepresentation $V'_G(\eta_V) \subset V_G(\eta_V)$. Therefore by the boundedness of $\alpha$ and the previous lemma
\begin{align*}
 \alpha_{W_j/W_{j-1}}(\MC{G} \times^{L^+G} W_j/W_{j-1}) & \subseteq (z - \zeta)^{-\langle (-\eta_V)_{\rm dom}, \mu \rangle} (\MC{G}' \times^{L^+G} W_j/W_{j-1}) \\
 & \subseteq (z - \zeta)^N (\MC{G}' \times^{L^+G} W_j/W_{j-1}).
\end{align*}
Thus we may apply again the previous lemma, but now for the short exact sequence $(\star)$ and using the converse implication, which gives
\[\alpha_{W_j}(\MC{G} \times^{L^+G} W_j) \subseteq (z - \zeta)^N (\MC{G}' \times^{L^+G} W_j).\] \exit

\lem{\label{Lem:BoundNiceBasis}}{}{
Assume that $G$ splits over $E$. Let $\lambda \in X^*(T)_{\rm dom}$ and $\mu \in X_*(T)_{\rm dom}$. Then there exists a basis $(v_1, \ldots, v_n)$ of the Weyl module $V_G(\lambda)$ (over $E$) such that the following properties are satisfied: \\
i) $v_i$ is a $T$-eigenvector of some weight $\lambda_i \in X^*(T)$ for all $1 \leq i \leq n$. \\
ii) $v_1$ has weight $\lambda$ and $v_n$ has weight $-(-\lambda)_{\rm dom}$. \\
iii) the vector space spanned by $v_1, v_2, \ldots, v_i$ is $B$-stable for all $1 \leq i \leq n$. \\
iv) $\langle \lambda_i, \mu \rangle \geq \langle \lambda_{i + 1}, \mu \rangle$ for all $1 \leq i \leq n - 1$.
}

\prooof
Define $v_i$ inductively for $1 \leq i \leq n$ as follows: Let $\MC{B}_i$ be the set of all $T$-eigenvectors not contained in the span of $v_1, \ldots, v_{i - 1}$. Note that $\MC{B}_i$ is non-empty because $V_G(\lambda)$ admits a basis of $T$-eigenvectors. Let $\MC{B}_i' \subset \MC{B}_i$ be the subset of eigenvectors whose weights are maximal with respect to the Bruhat order among all weights occurring as a weight of a vector in $\MC{B}_i$. Finally choose $v_i \in \MC{B}_i'$ of a weight $\lambda_i$ such that $\langle \lambda_i, \mu \rangle$ is maximal among all vectors in $\MC{B}_i'$. \\
We will check that this basis fulfills conditions i) to iv). Condition i) is trivial. Condition ii) follows from the fact that every weight $\lambda'$ of $V_G(\lambda)$ satisfies $-(-\lambda)_{\rm dom} \preceq \lambda' \preceq \lambda$ and both extremal weights occur. As the vector space in iii) is obviously $T$-stable it suffices to check that it is stable under each root subgroup $U_\alpha$ for $\alpha$ a simple root. But for any $T$-eigenvector $v_i$ of weight $\lambda_i$ and every element $u \in U_{\alpha}(\ACFq)$ the element $u \cdot v_i - v_i$ is a $T$-eigenvector of weight $\lambda_i + \alpha \succ \lambda_i$. Thus by our maximality assumption defining $\MC{B}_i$ we have $u \cdot v_i \subset \langle v_1, \ldots, v_i\rangle_{\Fq}$. Hence this subspace is $B$-stable. 
For property iv) note that by choice of $v_i$ we have $\langle \lambda_i, \mu \rangle \geq \langle \lambda', \mu \rangle$ for any weight $\lambda'$ of a vector $v' \in \MC{B}_i'$. If $v'' \in \MC{B}_i$ is any vector of weight $\lambda''$ then there is by definition of $\MC{B}_i'$ a vector $v' \in \MC{B}_i'$ of weight $\lambda'$ with $\lambda' \succeq \lambda''$. As $\mu$ as assumed to be dominant we get
\[\langle \lambda_i, \mu \rangle \geq \langle \lambda', \mu \rangle \geq \langle \lambda'', \mu \rangle\]
In particular this holds for $v'' = v_{i + 1}$ giving the desired inequality. \exit

\lem{\label{Lem:BoundAfterRepr}}{}{
Assume that $G$ splits over $E$ and let $\mu \in X_*(T)_{\rm dom}$ be a dominant coweight. Let $\lambda \in X^*(T)_{\rm dom}$, $V_G(\lambda)$ the Weyl module of $G$ (over $E$) of highest weight $\lambda$ and consider the associated representation $\rho: G \to GL \coloneqq GL(V_G(\lambda))$ over $E$. \\
Let as above $\MC{G}$ and $\MC{G}'$ be $L^+G$-torsors over some $E$-scheme $S$ and $\alpha: \MC{LG} \to \MC{LG}'$ a morphism between the associated $LG$-torsors. Let $\MC{G}_{GL} \coloneqq \MC{G} \times^{L^+G} L^+GL$, $\MC{G}'_{GL} \coloneqq \MC{G}' \times^{L^+G} L^+GL$ and let $\rho(\alpha): \MC{LG}_{GL} \to \MC{LG}'_{GL}$ be the morphism induced by $\alpha$ (cf. construction \ref{const:GenTorsorGroup}). \\
Choose now any basis $(v_1, \ldots, v_n)$ ($n = \dim V_G(\lambda)$) of $V_G(\lambda)$ satisfying the properties of lemma \ref{Lem:BoundNiceBasis}. Let $B_{GL} \subset GL$ be the Borel stabilizing the flag $E v_1 \subset E v_1 \oplus E v_2 \subset \ldots \subset V_G(\lambda)$ and let $T_{GL} \subset B_{GL}$ be the maximal torus stabilizing every line $E v_i$ for $1 \leq i \leq n$. The element $\mu \in X_*(T)_{\rm dom}$ then induces by composition with $T \to T_{GL}$ a dominant element $\mu_{GL} \in X_*(T_{GL})_{\rm dom}$.\\
a) If $\alpha$ is bounded by $\mu$ then $\rho(\alpha)$ is bounded by $\mu_{GL}$. \\
b) If $\rho(\alpha)$ is bounded by $\mu_{GL}$ then $\alpha$ fulfills condition \textit{1.} of definition \ref{Def:BoundLocal}a) for the dominant weight $\lambda$ and the dominant coweight $\mu$.
}

\prooof
a) We will check the conditions given in lemma \ref{Lem:BoundGLEasy}. \\
Consider first condition \textit{1.'}. Then for any scheme $S' \to S$ we have for the presheaves defined in construction \ref{const:GenTorsorGroup} and \ref{const:GenTorsorVector}:
\begin{align*}
 {\MC{G}_{GL}}_{Std}(S') & = \left(\MC{G}_{GL}(S') \times (Std \otimes_E \MC{O}_S[[z]])(S')\right)/L^+GL(S') \\
  & = \Bigl(\left(\MC{G}(S') \times L^+GL_n(S')\right)/L^+G(S') \times (Std \otimes_E \MC{O}_S[[z]])(S') \Bigr)/L^+GL(S') \\
  & \cong \left(\MC{G}(S') \times (Std \otimes_E \MC{O}_S[[z]])(S')\right)/L^+G(S') \\
  & = \left(\MC{G}(S') \times (V_G(\lambda) \otimes_E \MC{O}_S[[z]])(S')\right) /L^+G(S') \\
  & = \MC{G}_\lambda(S')
\end{align*}
In particular we get an isomorphism of sheaves of $\MC{O}_S[[z]]$-modules ${\MC{G}_{GL}}_{Std} \cong \MC{G}_\lambda$. Thus condition \textit{1.'} rewrites as
\[\rho(\alpha)_{\delta_i} = \wedge^i \rho(\alpha)_{Std} = \wedge^i\alpha_{\lambda}: \bigwedge\nolimits^i \MC{G}_\lambda \to \bigwedge\nolimits^i \MC{G}'_\lambda \otimes_{\MC{O}_S[[z]]} \MC{O}_S((z))\]
having image contained in $(z - \zeta)^{-\langle (-\delta_i)_{\rm dom}, \mu_{GL} \rangle} \bigwedge\nolimits^i \MC{G}'_\lambda$ for all $1 \leq i \leq n$. Here $\delta_i = (1, \ldots, 1, 0, \ldots, 0)$ with $i$ entries equal to $1$. But by lemma \ref{Lem:BoundAlternating} we already know that 
\[\wedge^i\alpha_{\lambda}\left(\bigwedge\nolimits^i \MC{G}_\lambda\right) \subseteq (z - \zeta)^N \bigwedge\nolimits^i \MC{G}'_\lambda\]
for $N = \min_{S_i} \{\langle \sum_{\eta \in S_i} \eta, \mu \rangle\}$ taking the minimum over all subsets $S_i$ of cardinality $i$ on the weights of $V_G(\lambda)$. Hence we are left to show $N = -\langle (-\delta_i)_{\rm dom}, \mu_{GL} \rangle$. \\
If $\varepsilon_j = (0, \ldots, 0, 1, 0, \ldots, 0) \in X^*(T_{GL})$ is the weight with a single entry $1$ at the $j$th coordinate, then we have for $1 \leq j \leq n$
\[z^{\langle \lambda_j, \mu \rangle} v_j = \mu(z) v_j = \mu_{GL}(z) v_j = z^{\langle \varepsilon_j, \mu_{GL} \rangle} v_j\]
where we denote the weight of $v_j$ by $\lambda_j \in X^*(T)$. Thus $\langle \lambda_j, \mu \rangle = \langle \varepsilon_j, \mu_{GL} \rangle$ for every $j$ and we get
\[-\langle (-\delta_i)_{\rm dom}, \mu_{GL} \rangle = \sum_{j = n - i + 1}^n \langle \varepsilon_j, \mu_{GL} \rangle = \sum_{j = n - i + 1}^n \langle \lambda_j, \mu \rangle\]
By property iv) of lemma \ref{Lem:BoundNiceBasis} we know that
\[\sum_{j = n - i + 1}^n \langle \lambda_j, \mu \rangle \leq \sum_{\eta \in S_i} \langle\eta, \mu \rangle\]
for any $S_i$ (and equality indeed occurs for some $S_i$). Thus 
\[-\langle (-\delta_i)_{\rm dom}, \mu_{GL} \rangle = \sum_{j = n - i + 1}^n \langle \lambda_j, \mu \rangle = \min_{S_i} \left\{\sum\nolimits_{\eta \in S_i} \langle\eta, \mu \rangle\right\} = N\]
Let us now check the condition on Hodge points, which holds for any choice of $T_{GL}$ and $B_{GL}$ containing the image of $T$ respectively $B$. Consider any geometric point $\ov{s}: \Sp k \to S$ with image $s$ where $k$ is an algebraically closed field. As the definition of the Hodge point is independent of the choice of trivializations we may assume $\ov{s}^*\MC{G} = \ov{s}^*\MC{G}' = L^+G \times \Sp k$ and $\ov{s}^*\alpha$ is given by an element $\alpha_s \in LG(k)$. Then functoriality of the Kottwitz homomorphism induces the commutative diagram:
\[
 \begin{xy}
  \xymatrix{
  LG(k) \ar^{\rho}[r] \ar_{\kappa}[d] & LGL(k) \ar_{\kappa}[d] \\
  \pi_1(G)_{\M{Q}}/\Gamma \ar[r]  & \pi_1(GL)_{\M{Q}}/\Gamma }
 \end{xy} 
\]
Hence $[\mu(\alpha)](s) = [\mu] \in \pi_1(G)_{\M{Q}}/\Gamma$ implies the desired $[\mu(\rho(\alpha))](s) = [\mu_{GL}] \in \pi_1(GL)_{\M{Q}}/\Gamma = \pi_1(GL)_{\M{Q}}$. \\
b) By choice of $\rho$ we have $\alpha_\lambda = \rho(\alpha)_{\delta_1}$. Thus 
\[\alpha_{\lambda}(\MC{G}_\lambda) \subseteq (z - \zeta)^{-\langle (-\delta_1)_{\rm dom}, \mu_{GL}\rangle} \MC{G}'_{\lambda}\] 
with
\[-\langle (-\delta_1)_{\rm dom}, \mu_{GL}\rangle = \langle \lambda_n, \mu \rangle = -\langle (-\lambda)_{\rm dom}, \mu \rangle \]
using property ii) of lemma \ref{Lem:BoundNiceBasis} and the computation done in part a). \exit

\prop{\label{Prop:BoundInverseSplit}}{}{
Let $G$ be a connected reductive group over $\Fq$ which splits over $E$. Let $\MC{G}$ and $\MC{G}'$ be two $L^+G$-torsors over a DM-stack $S$ over $E$ and $\alpha: \MC{LG} \to \MC{LG}'$ be a morphism between associated $LG$-torsors. If $\alpha$ is bounded by a dominant coweight $\mu$, then its inverse $\alpha^{-1}: \MC{LG}' \to \MC{LG}$ is bounded by $(-\mu)_{\rm dom}$.
}

\prooof
Let $\lambda$ be any dominant weight and consider the representation $\rho: G \times_{\Fq} \Sp E \to GL \coloneqq GL(V_G(\lambda))$ acting on the Weyl module of $G$ (over $E'$) of highest weight $\lambda$. By the previous proposition $\rho(\alpha)$ is then bounded by $\mu_{GL}$. Applying proposition \ref{Prop:BoundInverseGL} it follows that $\rho(\alpha)^{-1} = \rho(\alpha^{-1})$ is bounded by $(-\mu_{GL})_{\rm dom}$. 
Unfortunately lemma \ref{Lem:BoundAfterRepr}b) cannot be applied directly, because in general (unless further conditions on the $v_i$ are imposed) $(-\mu_{GL})_{\rm dom}$ does not equal $((-\mu)_{\rm dom})_{GL}$. However the only property used in its proof is $\langle (-\lambda)_{\rm dom}, (-\mu)_{\rm dom} \rangle = \langle (-\delta_1)_{\rm dom}, (-\mu_{GL})_{\rm dom}\rangle$. But this simply follows from
\[\langle (-\lambda)_{\rm dom}, (-\mu)_{\rm dom} \rangle = \langle \lambda, \mu \rangle = \langle \lambda_1, \mu \rangle = \langle \varepsilon_1, \mu_{GL} \rangle = \langle \delta_1, \mu_{GL} \rangle = \langle (-\delta_1)_{\rm dom}, (-\mu_{GL})_{\rm dom} \rangle.\]
Hence indeed
\[\alpha^{-1}_{\lambda}(\MC{G}'_\lambda) \subseteq (z - \zeta)^{-\langle(-\lambda)_{\rm dom}, (-\mu)_{\rm dom}\rangle} \MC{G}_\lambda\]
and condition \textit{1.} of definition \ref{Def:BoundLocal} is satisfied. \\
The argument for condition \textit{2.} is exactly the same as in the case of $GL_n$ (c.f. proposition \ref{Prop:BoundInverseGL}). \exit

\prop{\label{Prop:BoundInverse}}{}{
Let $G$ be any connected reductive group over $\Fq$. Let again $\MC{G}$ and $\MC{G}'$ be two $L^+G$-torsors over a DM-stack $S$ over $E$ and $\alpha: \MC{LG} \to \MC{LG}'$ be a morphism between associated $LG$-torsors. If $\alpha$ is bounded by an element $\mu \in X_*(T)_{\rm dom}$ defined over $E$, then its inverse $\alpha^{-1}: \MC{LG}' \to \MC{LG}$ is bounded by $(-\mu)_{\rm dom}$. \\
In particular any quasi-isogeny $\alpha: (\MC{G}, \varphi) \to (\MC{G}', \varphi')$ of local $G$-shtukas which is bounded by $\mu$ has an inverse which is bounded by $(-\mu)_{\rm dom}$.
}

\prooof
Choose any finite field extension $E'/E$ such that $G$ splits over $E'$. As the assertion can be checked \'etale locally by \ref{Rem:BoundOldDef}ii), we may replace $S$ by $S \times_E \Sp E'$. Moreover by \ref{Rem:BoundOldDef}i) we may replace $\Gamma_E$ by the absolute Galois group of $E'$ in condition \textit{1.} But then we can apply proposition \ref{Prop:BoundInverseSplit} over the ground field $E'$. \exit

\subsection{Further properties}\label{subsec:LocalFurther}
This section establishes several properties of bounded quasi-isogenies needed in the proofs of the representability results in the next sections. Contrary to the previous section, we start encountering arguments requiring local $G$-shtukas and not only $L^+G$-torsors. Throughout this section $G$ is any reductive group without any splitting hypothesis.

\lem{\label{Lem:BoundTrivia}}{}{
Let $\alpha: (\MC{G}, \varphi) \to (\MC{G}', \varphi')$ and $\alpha': (\MC{G}', \varphi') \to (\MC{G}'', \varphi'')$ be quasi-isogenies between local $G$-shtukas over a DM-stack $S$ over $E$ which are bounded by $\mu$ respectively $\mu'$. \\
a) The Frobenius-pullback $\sigma^*\alpha: (\sigma^*\MC{G}, \sigma^*\varphi) \to (\sigma^*\MC{G}', \sigma^*\varphi')$ is bounded by $\sigma^*\mu \in X_*(T)_{\rm dom}$. \\
b) The composite $\alpha' \circ \alpha$ is bounded by $\mu + \mu' \in X_*(T)_{\rm dom}$. 
}

\prooof
a) By definition of the Hodge point we have $[\mu(\sigma^*(\alpha))] = \sigma^*[\mu(\alpha)]$. To check \textit{1.} we proceed as in remark \ref{Rem:BoundGaloisInvariant}: We have $\sigma^*(\alpha_\lambda(\MC{G}_\lambda)) = \sigma^*\alpha_\lambda(\sigma^* \MC{G}_\lambda) = (\sigma^*\alpha)_{\sigma^{*-1}\lambda}(\sigma^*\MC{G}_{\sigma^{*-1}\lambda})$ and $\langle \sigma^{*-1}\lambda, \mu \rangle = \langle \lambda, \sigma^*\mu \rangle$. Thus pulling back the inclusions in \textit{1.} for $\alpha$, $\mu$ and some $\lambda$ along $\sigma^*$, gives precisely the corresponding condition \textit{1.} for $\sigma^*\alpha$, $\sigma^*\mu$ and $\sigma^{-1}\lambda$. Thus we are done by varying $\lambda$. \\
b) Condition \textit{1.} is immediate from the definition and \textit{2.} can be checked on geometric points and after trivializations and thus follows from the fact, that the Kottwitz homomorphism $L^+G(\ACFq) \to \pi_1(G)$ is a group morphism.
\exit

\defi{}{}{
Let $S = \varinjlim S_m$ be an admissible formal DM-stack over $\ENil$ such that $\zeta$ is locally nilpotent. \\
a) A formal local $G$-shtuka over $S$ is a tuple $(\MC{G}, \varphi) = \varinjlim (\MC{G}_m, \varphi_m)$ of local $G$-shtukas $(\MC{G}_m, \varphi_m)$ over the DM-stack $S_m$, together with isomorphisms between $(\MC{G}_m, \varphi_m)$ and the pullback of $(\MC{G}_{m+1}, \varphi_{m+1})$ to $S_m$ for every $m$. \\
b) A quasi-isogeny between formal local $G$-shtukas over $S$ is a tuple $\alpha = \varinjlim \alpha_m: (\MC{G}, \varphi) = \varinjlim (\MC{G}_m, \varphi_m) \to (\MC{G}', \varphi') = \varinjlim (\MC{G}'_m, \varphi'_m)$ of quasi-isogenies $\alpha_m: (\MC{G}_m, \varphi_m) \to (\MC{G}'_m, \varphi'_m)$, such that the pullback of $\alpha_{m + 1}$ to $S_n$ equals $\alpha_m$ for every $m$. \\
c) A quasi-isogeny $\alpha$ between formal local $G$-shtukas over $S$ is bounded by $\mu$ if every $\alpha_m$ is bounded by $\mu$.
A formal local $G$-shtuka $(\MC{G}, \varphi)$ over $S$ is bounded by $\mu$ if every local $G$-shtuka $(\MC{G}_m, \varphi_m)$ is bounded by $\mu$. 
}

\rem{}{
We will hardly make use of formal local $G$-shtukas, but they are needed for one argument in the proof of theorem \ref{Thm:RapoZinkBounded}.
}

\prop{\label{Prop:BoundRigid}}{\textbf{(Rigidity of quasi-isogenies)} cf. \cite[proposition 3.9]{HaVi}}{
Let $(\MC{G}, \varphi)$ and $(\MC{G}', \varphi')$ be two local $G$-shtukas over $S \in Nilp_{\ENil}$. Let $\iota: \ov{S} \hookrightarrow S$ be a closed immersion defined by a sheaf of ideals $\MC{I}$ which is locally nilpotent. Then there is a bijection (of sets)
\[\iota^*: QIsog_S((\MC{G}, \varphi), (\MC{G}', \varphi')) \to QIsog_{\ov{S}}(\iota^*(\MC{G}, \varphi), \iota^*(\MC{G}', \varphi')) \qquad \alpha \mapsto \iota^*\alpha\]
Let $S$ now be quasi-compact and $(\MC{G}, \varphi)$ and $(\MC{G}', \varphi')$ bounded local $G$-shtukas. If a quasi-isogeny of the left-hand side is bounded by $\ov{\mu}$, then the corresponding quasi-isogeny $\alpha$ on the left-hand side is bounded by some $\mu$, too. This $\mu$ can be chosen to depend only on $\ov{\mu}$, $\iota$ and the bounds for the two local $G$-shtukas, but not on the actual local $G$-shtukas or the quasi-isogeny $\alpha$. \\
The same holds for formal local $G$-shtukas over admissible formal DM-stacks $S$.
}

\prooof
The construction of the map between the sets of quasi-isogenies can be found in \cite[proof of proposition 3.9]{HaVi}. There a slightly weaker statement about the behavior of bounds is shown as well. Their arguments may be strengthened as follows: \\
Assume that $(\MC{G}, \varphi)$ is bounded by $\mu_0$ and $(\MC{G}', \varphi')$ is bounded by $\mu_0'$. 
Let now $\alpha \in QIsog_S((\MC{G}, \varphi), (\MC{G}', \varphi'))$ be a quasi-isogeny such that $\iota^*\alpha$ is bounded by $\ov{\mu}$. As $S$ is quasi-compact there is some $n > 0$ with $\MC{I}^{q^n} = 0$. Now consider the commutative diagram of isomorphisms
\[
 \begin{xy}
  \xymatrix{
  \MC{LG} \ar^{\alpha}[r] \ar_{\sigma^{n-1 \, *}\varphi \circ \ldots \circ \sigma^*\varphi \circ \varphi}[d] & \MC{LG}' \ar^{\sigma^{n-1 \, *}\varphi' \circ \ldots \circ \sigma^*\varphi' \circ \varphi'}[d]\\
  \sigma^{n \, *}\MC{LG} \ar^{\sigma^{n \, *}\alpha}[r]  & \sigma^{n \, *}\MC{LG}' }
 \end{xy} 
\]
As explained in \cite{HaVi} the $q^n$-Frobenius on $S$ factors over $\ov{S}$. Thus ${\sigma^n}^*\alpha$ is the pullback of $\iota^*\alpha$ via some morphism $j: S \to \ov{S}$ and is thus bounded by $\ov{\mu}$. Then \ref{Prop:BoundInverse} and \ref{Lem:BoundTrivia} imply that 
\[\alpha = (\sigma^{n-1 \, *}\varphi \circ \ldots \circ \sigma^*\varphi \circ \varphi)^{-1} \circ \sigma^{n \, *}\alpha \circ (\sigma^{n-1 \, *}\varphi' \circ \ldots \circ \sigma^*\varphi' \circ \varphi')\]
is bounded by $\mu = \sum_{i = 0}^{n-1} \sigma^{i \, *}(-\mu_0)_{\rm dom} + \ov{\mu} + \sum_{i = 0}^{n-1} \sigma^{i \, *}\mu_0'$. \\
Consider now the case of formal DM-stacks. Let $S = \varinjlim S_m$. Then we have correspondingly $\qquad$ $\ov{S} = \varinjlim \ov{S}_m$, where $\ov{S}_m \subset S_m$ is defined by the nilpotent sheaf of ideals $\MC{I}$ on $S$. Let $\ov{\alpha} = \varinjlim \ov{\alpha}_m \in QIsog_{\ov{S}}(\iota^*(\MC{G}, \varphi), \iota^*(\MC{G}', \varphi'))$. Then we may lift each $\ov{\alpha}_m$ to a quasi-isogeny $\alpha_m$ over $S_m$. By uniqueness of this lift, the $\alpha_m$ are compatible and define indeed a quasi-isogeny between formal local $G$-shtukas over $S$. \\
To get the boundedness assertion, assume again that the formal local $G$-shtukas $(\MC{G}, \varphi)$ resp. $(\MC{G}', \varphi')$ are bounded by $\mu_0$ resp. $\mu_0'$. If $\alpha = \varinjlim \alpha_m \in QIsog_S((\MC{G}, \varphi), (\MC{G}', \varphi'))$ is a quasi-isogeny such that $\iota^*\alpha$ is bounded by $\ov{\mu}$, we have just shown that every $\alpha_m$ is bounded by $\mu = \sum_{i = 0}^{n-1} \sigma^{i \, *}(-\mu_0)_{\rm dom} + \ov{\mu} + \sum_{i = 0}^{n-1} \sigma^{i \, *}\mu_0'$ (where $n$ is a fixed integer with $\MC{I}^{q^n} = 0$), which is independent of $m$. Hence $\alpha$ is indeed bounded. \exit

\rem{}{
a) It is in general not true that any quasi-isogeny over $\ov{S}$ which is bounded by $\mu$ lifts to a quasi-isogeny that is again bounded by $\mu$. Take for example $G = GL_2$, the closed immersion $i: \ov{S} = \Sp \Fq \to S = \Sp \Fq[\zeta]/(\zeta^2)$, $b_n = \lmat{cc} 1 & \zeta z^{-n-1} \\ 0 & 1 \rmat$ for some $n > 0$, $(\MC{G}, \varphi) = (L^+GL_2, \sigma^*)$, $(\MC{G}', \varphi') = (L^+GL_2, b_n \sigma^*(b_n)^{-1} \sigma^*)$ and $\alpha = b_n \in LGL_2$. Then $i^*\alpha = \id$ is bounded by $\mu = (0, 0)$, but $\alpha$ itself is only bounded by $\mu' = (n, -n)$. \\
b) The assumption that $(\MC{G}, \varphi)$ and $(\MC{G}', \varphi')$ are bounded is necessary. To see this take for $(\MC{G}, \varphi)$ any local $G$-shtuka over a scheme $S$, which is not bounded by any dominant coweight but restricts to a bounded local $G$-shtuka over the reduced subscheme $S_{red}$. An example for this can be found in \cite[example 3.12]{HaVi}. Let $(\MC{G}', \varphi')$ be a bounded local $G$-shtuka over $S$ and $\alpha$ any quasi-isogeny between them. Then $\alpha$ cannot be bounded by any dominant coweight, because otherwise exactly three of the quasi-isogenies in the commutative square $\varphi' \circ \sigma^*\alpha = \alpha \circ \varphi$ are bounded, contradicting \ref{Prop:BoundInverse} and \ref{Lem:BoundTrivia}. But the pullback of $\alpha$ to $S_{red}$ is automatically bounded by \ref{Lem:BoundFields}b).
}

\prop{\label{prop:BoundClosed}}{cf. \cite[lemma 3.10]{HaVi}}{
Let $\alpha: (\MC{G}, \varphi) \to (\MC{G}', \varphi')$ be a quasi-isogeny between local $G$-shtukas over a DM-stack $S$ over $E$ and let $\mu \subset X_*(T)_{\rm dom}$ be a cocharacter that can be defined over $E$. Then the condition that $\alpha$ is bounded by $\mu$ is representable by a closed immersion into $S$, which is locally finitely presented. \\
In particular under the assumption that $S$ is reduced, $\alpha$ is bounded by $\mu$ if and only if this holds for the pullback to every generic point of $S$.
}

\prooof 
Let $E'/E$ be some finite field extension, such that $G$ splits over $E'$ or equivalently such that each weight $\lambda$ can be defined over $E'$. We will first show that the $\mu$-bounded locus is representable after base-change to $E'$. The argument for this is essentially the same as in \cite{HaVi}: \\
As the Hodge point of $\alpha$ is locally constant, condition \textit{2.} is trivially representable by a finitely presented closed immersion into $S \times_E \Sp E'$. So assume wlog. condition \textit{2.} holds on all of $S \times_E \Sp E'$. Fix a finite set of weights $\Lambda$ which generate the monoid of dominant weights and note that it suffices by \ref{Lem:BoundWeightGenerate} to check condition \textit{1.} only for the weights in $\Lambda$. 
As condition \textit{1.} is local we may assume $S \times_E \Sp E'$ is an affine scheme isomorphic to $\Sp R$ for some $E'$-algebra $R$ and that the $R[[z]]$-modules $\MC{G}_\lambda$ and $\MC{G}'_\lambda$ are free for every dominant weight $\lambda \in \Lambda$ (and indeed exist over $R$ by our choice of $E'$). As $\Sp R$ is quasi-compact we may choose for every $\lambda \in \Lambda$ an $N_\lambda \gg 0$ such that $\alpha_\lambda(\MC{G}_\lambda) \subset (z - \zeta)^{-N_\lambda} \MC{G}_\lambda'$. So $\alpha$ is bounded by $\mu$ if and only if $\alpha_\lambda$ maps all generators of $\MC{G}_\lambda$ to zero in $M_\lambda \coloneqq (z - \zeta)^{-N_\lambda} \MC{G}_\lambda'/(z - \zeta)^{-\langle(-\lambda)_{\rm dom},\mu\rangle} \MC{G}_\lambda'$. Since $M \coloneqq \bigoplus_{\lambda \in \Lambda} M_\lambda$ is a free $R$-module of finite rank, this condition (for all $\lambda \in \Lambda$) is represented by a finitely presented closed immersion. \\
Thus we get a closed immersion $Z' \subset S \times_E \Sp E'$ representing the locus of $\mu$-boundedness over $E'$. We have to see, that it descends to $E$, i.e. that $Z'$ is invariant under $\Gal(E'/E)$. Fix some $\tau \in \Gal(E'/E)$ and consider the induced isomorphism $\tau: S \times_E \Sp E' \to S \times_E \Sp E'$ over $E$. Then as in remark \ref{Rem:BoundGaloisInvariant} we see that $\alpha|_{Z'}$ fulfilling condition \textit{1.} for some $\lambda \in X^*(T)_{\rm dom}$ implies that $\alpha|_{\tau^{-1}(Z')} = \tau^*\alpha|_{Z'}$ fulfills condition \textit{1.} for $\tau^{-1}\lambda \in X^*(T)_{\rm dom}$. As condition \textit{2.} is obviously $\Gal(E'/E)$-invariant, we get that $\alpha|_{\tau^{-1}(Z')}$ is again bounded by $\mu$, i.e. $\tau^{-1}(Z') = Z'$. Thus $Z'$ is $\Gal(E'/E)$-invariant as desired and descends to a closed immersion into $S$. \exit

\rem{}{
The same result is obviously true if one replaces the local $G$-shtukas by $L^+G$-torsors and $\alpha$ by any morphism between associated $LG$-torsors. 
}

\defi{}{}{
Recall the Bruhat order $\preceq$ on $X_*(T)_{\M{Q}, dom}$: For two dominant cocharacters $\mu, \mu' \in X_*(T)_{\M{Q}}$ the relation $\mu' \preceq \mu$ holds if $\mu - \mu'$ is a non-negative rational linear combination of simple coroots.
}

\lem{\label{Lem:BoundFields}}{cf. \cite[lemma 3.11 and 3.13]{HaVi}}{
a) Let $k$ be an algebraically closed field and $\alpha$ a quasi-isogeny between local $G$-shtukas over $\Sp k \in Nilp_{\ENil}$. Then $\alpha$ is bounded by $\mu$ if and only if its Hodge point $\mu(\alpha)(\Sp k) \in X_*(T)_{\rm dom}$ satisfies $\mu(\alpha)(\Sp k) \preceq \mu$. \\
b) Let $S$ be quasi-compact and connected and let $\alpha$ be a morphism between local $G$-shtukas over $S$. If either
\begin{itemize}
 \item $G$ is split semi-simple or
 \item $S$ reduced and the (constant) image of $[\mu(\alpha)]$ in $\pi_1(G)_{\M{Q}}/\Gamma$ is an $1$-element orbit,
\end{itemize}
then $\alpha$ is bounded by a $\Gamma$-invariant element $\mu \in X_*(T)_{\rm dom}$.
}

\prooof
The proofs are very similar to the ones in \cite{HaVi}. As in \cite[lemma 3.11]{HaVi} one can see that over a geometric point boundedness by $\mu$ is equivalent to $\langle \lambda, \mu - \mu(\alpha)(\Sp k) \rangle \geq 0$ for all $\lambda \in X^*(T)_{\rm dom}$ and $[\mu(\alpha)(\Sp k)] = [\mu] \in \pi_1(G)_{\M{Q}}$. But this is equivalent to $\mu(\alpha)(\Sp k) \preceq \mu$. The proof of part b) given in \cite{HaVi} uses only the existence of some $\Gamma$-invariant $\mu$ with $[\mu(\alpha)] = [\mu] \in \pi_1(G)_{\M{Q}}/\Gamma$ and the fact that for any $\alpha$ and any $\lambda$ there is some $N \gg 0$ such that the image of the $L^+G$-torsor under $\alpha_\lambda$ is contained in the subsheaf defined by multiplication with $(z - \zeta)^{-N}$. \exit 

\rem{}{
i) In \cite[lemma 3.11]{HaVi} a slightly different partial order is used: There $\mu' \preceq \mu$ if $\mu - \mu'$ is a non-negative \textit{integral} linear combination of simple coroots. The change made in the partial order precisely reflects the weakening of condition \textit{2.} by allowing differences in the torsion part of $\pi_1(G)$. \\
ii) As we require $\mu$ to be $\Gamma$-invariant, there might be no such $\mu$ satisfying the assertion in part a). Similarly if in part b) $[\mu(\alpha)]$ is not an $1$-element orbit in $\pi_1(G)_{\M{Q}}$ under the $\Gamma$-action, then we cannot hope to bound $\alpha$ by any $\mu$ (at least for this ground field $E$).
}

\subsection{Rapoport-Zink spaces for local \textit{G}-shtukas}\label{subsec:RapoZinkSpace}
We construct a formal scheme classifying pairs consisting of a bounded local $G$-shtuka and a quasi-isogeny between it and some fixed local $G$-shtuka. These spaces are the analog of Rapoport-Zink spaces for $p$-divisible groups as introduced in \cite{RapoZinkSpaces}. As most of the results below follow rather formally from the results of the previous section, we will mostly refer for proofs to \cite{HaVi}. \vspace{4mm} \\
We continue in the use of the notations in the previous sections. Let $Gr = LG / L^+G \times_{\Fq} \Sp E$ be the affine Grassmannian (base-changed to the ground field $E$). It exists as an ind-scheme over $E$ of ind-finite type, as shown e.g. in \cite[after corollary $3$]{FaltingsLoop}. Let $\widehat{Gr}$ be the completion of $Gr \times_E \Sp \ENil$ along the special fiber defined by $\zeta$. Equivalently $\widehat{Gr}$ is the fiber product $Gr \times_E \Spf \ENil$ in the category of ind-schemes. \\
Fix now an element $b_0 \in LG(E)$, which gives a local $G$-shtuka $(L^+G, b_0\sigma^*)$ over $E$. 

\thm{\label{thm:RapoZinkUnbounded}}{\cite[theorem 6.2]{HaVi}}{
The ind-scheme $\widehat{Gr}$ pro-represents the functor
\[\begin{array}{ccc}
  Nilp_{\ENil} & \to  & Sets \hspace{7cm}  \\
  S   & \mapsto & \left\{ (\MC{G}, \varphi), \alpha) \left| \begin{array}{c} (\MC{G}, \varphi) \; \textnormal{a local} \; G \textnormal{-shtuka over} \; S, \\ \alpha: (\MC{G}, \varphi) \to (L^+G, b_0\sigma^*) \; \textnormal{a quasi-isogeny over} \; S \end{array} \right. \right\}
\end{array}\]
Here two triples $((\MC{G}, \varphi), \alpha)$ and $((\MC{G}', \varphi'), \alpha')$ are isomorphic, if $\alpha^{-1} \circ \alpha'$ comes from an isomorphism of $L^+G$-torsors $\MC{G}' \to \MC{G}$ over $S$.
}

\prooof
The proof is exactly the same as in \cite{HaVi}. \exit 

\rem{}{
Note that by proposition \ref{Prop:BoundRigid}, the quasi-isogeny $\alpha$ is uniquely determined by its restriction to $S \times_{\Spf \ENil} \Sp \ENil/(\zeta)$. Hence the functor above indeed coincides with the one considered in \cite{HaVi}.
}

\defi{\label{def:RapoZinkDecent}}{}{
A local $G$-shtuka $(L^+G, \varphi)$ with trivial $L^+G$-torsor over an $E$-scheme $S$ is called decent if there exists an integer $s > 1$ and a cocharacter $\tilde{\mu} \in X_*(T)$ such that
\[\varphi \circ \sigma^*\varphi \circ \ldots \circ \sigma^{* (s-1)}\varphi = z^{\tilde{\mu}} \cdot \sigma^{*s}\]
}

\thm{\label{Thm:RapoZinkBounded}}{cf. \cite[theorem 6.3]{HaVi}}{
Let $(L^+G, b_0\sigma^*)$ be a decent local $G$-shtuka over $\Sp E$. Then the functor
\[\begin{array}{cccc}
 \MC{M}_{b_0}^{\preceq \mu}: & Nilp_{\ENil} & \to     & Sets  \hspace{7cm} \\
  & S      & \mapsto & \left\{((\MC{G}, \varphi), \alpha) \left| 
	\begin{array}{c} 
	  (\MC{G}, \varphi) \; \textnormal{a local} \; G \textnormal{-shtuka over} \; S \; \textnormal{bounded by} \; \mu, \\
          \alpha: (\MC{G}, \varphi)_{\ov{S}} \to (L^+G_E, b_0\sigma^*)_{\ov{S}} \; \textnormal{a quasi-isogeny} 
	\end{array} \right.\right\}
\end{array}\]
is representable by a closed immersion into $\widehat{Gr}$. It is a (non-$\zeta$-adic) formal scheme over $\Spf \ENil$, whose underlying reduced subscheme is locally of finite type.
}

\prooof
Consider first the case of split groups: The proof given in \cite{HaVi} applies to all definitions of bounds, for which the statements \ref{Prop:BoundRigid} - \ref{Lem:BoundFields} hold. \\
Nevertheless its proof uses at one step \cite[proposition 3.16]{HaVi} (respectively \cite[lemma 4.18]{HarRad1}). Unfortunately we were not able to understand the descent argument in the proof of this proposition (cf. the more detailed discussion in remark \ref{Rem:UnifProblem}i)). Therefore let us present a slightly modified argument here, that the functor $\MC{M}_n$ is representable. Please refer to \cite[proof of theorem 6.3, first claim]{HaVi} for definitions of all the objects used here.  \\
Let $\MC{M}^m_n$ be the formal completion of $\MC{M}^m$ along $(\MC{M}_n)_{\rm red}$. Any open affine subset $U \subset (\MC{M}_n)_{\rm red}$ defines formal affine subschemes $\Spf R_m \subset \MC{M}^m_n$. Let $R = \varprojlim R_m$ and $J$ the inverse image in $R$ of the largest ideal of definition of $R_n$. Then $R$ is a $J$-adic ring and we fix some integer $c > 0$. To show that $\MC{M}_n$ is a formal scheme, it suffices to show that there is an integer $m_0$ such that for all $m \geq m_0$ the natural map $R_m/J^cR_m \to R_{m_0}/J^cR_{m_0}$ is an isomorphism. 
To do so let $(\MC{G}_m, \varphi_m)$ be the universal local $G$-shtuka over $\Sp R_m/J^cR_m$ and $(\MC{G}, \varphi) = \varinjlim (\MC{G}_m, \varphi_m)$ the formal local $G$-shtuka over $\Spf R/J^c = \varinjlim \Sp R_m/J^cR_m$. By rigidity for quasi-isogenies between formal local $G$-shtukas we may lift the universal quasi-isogeny $\alpha$ from $\Spf R/J = \Spf R_n/JR_n$ to $\Spf R/J^c$. By assumption $\rho_*\alpha$ and $\rho_*\alpha^{-1}$ are bounded by $2n\varrho^\vee$ (as quasi-isogenies between formal local $G$-shtukas) over $\Spf R/J$. 
As lifting commutes with applying the functor $\rho_*$ the boundedness assertion in \ref{Prop:BoundRigid} for formal local $G$-shtukas implies, that $\rho_*\alpha$ and $\rho_*\alpha^{-1}$ are bounded by some dominant coweight over $\Spf R/J^c$, hence also by $2m_0\varrho^\vee$ for some $m_0 \gg n$. Hence by the universal property of $\MC{M}^{m_0}_n$ there is a unique map $\Spf R/J^c \to \Spf R_{m_0}$ inducing the given point $((\MC{G}, \varphi), \alpha)$. This gives the desired isomorphism. \\
Consider now the case of arbitrary connected reductive groups. As in the split case the only non-trivial part is to show that $\MC{M}_{b_0}^{\preceq \mu}$ exists as a formal scheme (and not only as a formal ind-scheme). This may be checked after passing to a finite field extension. But then we are again in the split case and may use the arguments above. \exit 

\rem{}{
i) The decency condition can essentially be removed: Let $(L^+G, b\sigma^*)$ be any local $G$-shtuka over $\Sp E$. Again it is clear that $\MC{M}_b^{\preceq \mu}$ exists as a closed immersion into $\widehat{Gr}$. All other assertions can be checked after passing to some sufficiently large field extension such that $G$ is split and there exists some quasi-isogeny $\alpha_0: (L^+G, b\sigma^*) \to (L^+G, b_0\sigma^*)$ to some decent local $G$-shtuka defined over $\Sp E$. As $G$ splits, this $\alpha_0$ is bounded by some $\mu_0 \in X_*(T)_{\rm dom}$. Then composition with $\alpha_0$ gives an immersion
\[\MC{M}_b^{\preceq \mu} \to \MC{M}_{b_0}^{\preceq \mu + \mu_0}\]
Its image is exactly the subspace where the universal quasi-isogeny over $\MC{M}_{b_0}^{\preceq \mu + \mu_0}$ is bounded by $\mu$ after composition with $\alpha_0^{-1}$. In particular it is given by a closed subscheme locally of finite presentation. Hence $\widehat{X}^{\preceq \mu + \mu_0}_{(L^+G_E, b_0\sigma^*)}$ being locally formally of finite type implies the same result for $\widehat{X}^{\preceq \mu}_{(\MC{G}, \varphi)}$. \\
ii) In accordance with the usual notation for Rapoport-Zink spaces of $p$-divisible groups (i.e. in the case of mixed characteristic), we usually refer to the formal scheme representing the functor above as $\MC{M}_{b_0}^{\preceq \mu}$ instead of $\widehat{X}^{\preceq \mu}_{(L^+G_E, b_0\sigma^*)}$.
Its reduced fiber will be denoted by $\MB{M}_{b_0}^{\preceq \mu}$. It equals closed affine Deligne-Lusztig varieties using the same arguments as in \cite[theorem 6.3]{HaVi}.
}

\cor{\label{cor:RapoZinkProper}}{}{
Each irreducible component of $\MB{M}_{b_0}^{\preceq \mu}$ is proper and of finite type. In particular $\MB{M}_{b_0}^{\preceq \mu}$ satisfies the valuation criterion for properness.
}

\prooof
In the case of split groups, this was shown in \cite[corollary 6.5]{HaVi} that they are projective and of finite type. Thus even in the non-split case, this holds after passing to a finite field extension, giving the desired result. \exit

Let us finally remark, that the construction of the formal deformation space of a local $G$-shtuka given in \cite[theorem 5.6]{HaVi} works as well in the case of non-split reductive groups $G$ in the very same way.

\subsection{The Tate functor for \'etale local \textit{G}-shtukas}\label{subsec:TateFunctor}
We construct the (dual) Tate functor for \'etale local $G$-shtukas, i.e. an equivalence of categories between \'etale local $G$-shtukas over $S$ (cf. definition \ref{def:LocShtuka}b)) and morphisms $\pi_1(S, \ov{s}) \to \Aut(\un{L^+G})$, where $\un{L^+G}$ denotes a trivial $L^+G(\Fq)$-torsor and $\pi_1(S, \ov{s})$ the \'etale fundamental group. This will play an important role in defining level structures of global $G$-shtukas (cf. section \ref{subsec:AdelicLevel}). \\
Let us briefly sketch the idea behind this correspondence: 
An \'etale local $G$-shtuka $(\MC{G}, \varphi)$ admits over the universal cover $\widetilde{S}$ an isomorphism to $(\MC{G}_{\widetilde{S}}, \sigma^*)$, where $\MC{G}_{\widetilde{S}}$ is a trivial $L^+G$-torsor. The descent data is given by a $\pi_1(S, \ov{s})$-action on the torsor $\MC{G}_{\widetilde{S}}$, which has to be compatible with the Frobenius-isomorphism. Hence it preserves the set of $\sigma$-invariant elements, which is a trivial torsor under the finite group $L^+G(\Fq)$. 
This allows us to encode the action of $\gamma \in \pi_1(S, \ov{s})$ by its induced isomorphism $\Isom(\MC{G}_{\widetilde{S}}^{\sigma^*}|_{\tilde{\ov{s}}}, \MC{G}_{\widetilde{S}}^{\sigma^*}|_{\gamma(\tilde{\ov{s}})}) \cong L^+G(\Fq)$ (where $\tilde{\ov{s}} \in \widetilde{S}$ is a fixed point over $\ov{s}$). Therefore we get the desired functor as
\[(\MC{G}, \varphi) \mapsto \left(\pi_1(S, \ov{s}) \to L^+G(\Fq) = \Aut(\MC{G}^{\sigma^*}|_{\tilde{\ov{s}}})\right).\]
The inverse is just given by descending a trivial $L^+G$-torsor together with the trivial Frobenius-isomorphism via the given $\pi_1(S, \ov{s})$-action. \\
We stress the point, that almost every result in this section is already contained in \cite{HarRad1} after rephrasing the morphism $\pi_1(S, \ov{s}) \to \Aut(\un{L^+G})$ as a tensor functor from the category of $G$-representations to the category of $\pi_1(S, \ov{s})$-representations on finite free $\Fq[[z]]$-modules. Nevertheless the slightly different point of view presented here, hopefully helps in gaining a deeper understanding of these constructions. \vspace{2mm} \\
Fix a DM-stack $S$ over $E$ which is (for simplicity) connected. Furthermore fix a closed point $s \in S$, a geometric point $\ov{s}$ with image $s$ (and residue field isomorphic an algebraic closure of the residue field of $s$) and a (geometric) point $\tilde{\ov{s}}$ in the universal cover $\widetilde{S}$ of $S$ lying over $\ov{s}$. We denote the category of \'etale local $G$-shtukas by $\acute{E}tSht_G(S)$. 
Note at this point that the definition of the \'etale fundamental group given in \cite{GroSGAI} can be extended to DM-stacks. For more details (and the occurring subtleties) see \cite{ZoonekyndPi1Stack}.
Let us start by defining some more categories:

\defi{\label{def:TateEtaleShtuka}}{}{
For every $n \geq 0$ let $\acute{E}tSht_G(S)[n]$ be the category of pairs consisting of a $L^+G/K_n$-torsor $\MC{G}[n]$ over $S$ and a $\sigma$-linear automorphism $\varphi[n]: \sigma^*\MC{G}[n] \to \MC{G}[n]$ of this torsor. Then there is a canonical functor 
\[\begin{array}{cccc}
 -[n]: & \acute{E}tSht_G(S) & \to & \acute{E}tSht_G(S)[n] \\
       & (\MC{G}, \varphi) & \mapsto & (\MC{G}[n], \varphi[n]) \coloneqq (\MC{G} \times^{L^+G} L^+G/K_n, \varphi \times \id)  \\
       &                   &         & \hspace{1.7cm} \cong (\MC{G} \otimes_{\MC{O}_S[[z]]} \MC{O}_S[[z]]/(z^{n+1}), \varphi \bmod z^{n+1}) \hspace{-1.7cm}
\end{array}\]
}

\defi{}{}{
a) A set-theoretic torsor for some abstract group $G_0$ is a set $\un{G_0}$ together with a simply transitive action by $G_0$. If $T$ is any scheme, then $\un{G_0}_T$ denotes the trivial $G_0$-torsor on the \'etale site of $T$ given by $\un{G_0}_T(T') = \un{G_0}$ for any $T' \to T$ \'etale. \\
b) For every $n \geq 0$ denote by $\Rep(\pi_1(S), L^+G/K_n(\Fq))$ the category of representations $\rho[n]: \pi_1(S, \ov{s}) \to \Aut(\un{L^+G/K_n})$ where $\un{L^+G/K_n}$ is a set-theoretic torsor under the finite group $L^+G/K_n(\Fq)$. A morphism $\alpha: \rho[n]_1 \to \rho[n]_2$ in $\Rep(\pi_1(S), L^+G/K_n(\Fq))$ is given by an isomorphism between the trivial set-theoretic $L^+G/K_n(\Fq)$-torsors conjugating $\rho[n]_1$ into $\rho[n]_2$. \\ 
Then define the limit category 
\[\Rep(\pi_1(S), L^+G(\Fq)) = \varprojlim_n \Rep(\pi_1(S), L^+G/K_n(\Fq)),\]
i.e. the category of continuous representations $\rho: \pi_1(S, \ov{s}) \to \Aut(\un{L^+G})$ where $\un{L^+G}$ is a set-theoretic torsor under the pro-finite group $L^+G(\Fq)$.
}

\rem{}{
We will use $L^+G/K_n(\Fq)$ both for the abstract group of $\Fq$-valued points of $L^+G/K_n$ and for the finite group scheme (over base scheme $S$) representing the sheaf of groups given by the constant sheaf with value the abstract group $L^+G/K_n(\Fq)$.
}

Before we are able to define functors between these categories we need two preparatory lemmas. These statements will be generalized later on by the construction of Igusa varieties.

\lem{\label{lem:TatePartTrivial}}{}{
Let $\varphi \in L^+G/K_n(S)$ be any $S$-valued point (for some $n \geq 0$). Then there is a finite \'etale cover $S' \to S$ such that there is an element $g \in L^+G/K_n(S')$ satisfying $g^{-1} \varphi \sigma(g) = 1$ inside $L^+G/K_n(S')$. 
}

\prooof
The quotient $L^+G/K_n$ is a linear algebraic group. Hence Lang's morphism $L: L^+G/K_n \to L^+G/K_n$, $L(g) = g\sigma(g)^{-1}$ is surjective and a finite \'etale cover. Thus consider the cartesian diagram
\[
 \begin{xy} \xymatrix{
   S' \ar^-{g}[r] \ar[d] & L^+G/K_n \ar^{L}[d] \\
   S \ar^-{\varphi}[r]  & L^+G/K_n 
 } \end{xy} 
\]
which defines a finite \'etale cover $S' \to S$ and an element $g \in L^+G/K_n(S')$ with $\varphi = g\sigma(g)^{-1}$ as desired. \exit

\lem{\label{lem:TateTrivialEtale}}{}{
Let $(\MC{G}, \varphi)$ be an \'etale local $G$-shtuka over $S$. Then for every $n \geq 0$ the $\varphi$-invariants $\MC{G}[n]^{\varphi}$ define an \'etale $L^+G/K_n(\Fq)$-torsor over $S$. Furthermore there is a canonical isomorphism of $L^+G/K_n$-torsors
\[\MC{G}[n]^{\varphi} \times^{L^+G/K_n(\Fq)} L^+G/K_n \cong \MC{G}[n].\]
Under this isomorphism $\varphi$ acts on the left-hand side as $\sigma^* \times \id_{L^+G/K_n}$.
}

\prooof
Obviously the $L^+G$-action defines a $L^+G/K_n(\Fq)$-action on the \'etale sheaf $\MC{G}[n]^{\varphi}$. We have to see that it admits a trivialization over an \'etale cover. For this consider any connected $U \subset S$ such that $\MC{G}$ trivializes over a finite \'etale cover $U' \to U$. Then $\varphi$ may be represented by an element in $L^+G(U')$ and by the previous lemma there is a finite \'etale cover $U'_n \to U'$ admitting an isomorphism $(\MC{G}[n], \varphi[n])_{U'_n} \cong (L^+G/K_n, \sigma^*)_{U'_n}$. But the $\sigma^*$-invariants of the trivial torsor give exactly the trivial $L^+G/K_n(\Fq)$-torsor. \\
In particular this shows that the canonical morphism $\MC{G}[n]^{\varphi} \times^{L^+G/K_n(\Fq)} L^+G/K_n \to \MC{G}[n]$ is an isomorphism after pulling back to $U'_n$. Hence it is an isomorphism already over $S$. The restriction of $\varphi$ to $\MC{G}[n]^{\varphi}$ acts by definition simply as $\sigma^*$. But by $L^+G/K_n$-equivariance, any isomorphism $\sigma^*\MC{G}[n]^{\varphi} \times^{L^+G/K_n(\Fq)} L^+G/K_n \to \MC{G}[n]^{\varphi} \times^{L^+G/K_n(\Fq)} L^+G/K_n$ is uniquely determined by its action on $\MC{G}[n]^{\varphi}$ and the last assertion follows. \exit

\prop{\label{prop:TateTrivialTorsor}}{}{
Let $(\MC{G}, \varphi)$ be an \'etale local $G$-shtuka over $S$. Then for every $n \geq 0$, the functor on the \'etale (or fppf, fpqc) site of $S$
\[S' \mapsto \Isom(\MC{G}[n]^{\varphi} \times_S S', L^+G/K_n(\Fq) \times_{\Fq} S')\]
is representable by a finite \'etale cover $S_n \to S$. Over $S_n$ the universal isomorphism induces a trivialization $(\MC{G}[n], \varphi)_{S_n} \cong (L^+G/K_n, \sigma^*)_{S_n}$. Furthermore $(\MC{G}, \varphi)$ trivializes canonically over the pro-finite \'etale cover $\varprojlim_n S_n$.
}

\prooof
As $L^+G/K_n(\Fq)$ is a finite group scheme, the functor defined in the corollary is representable by a scheme $S_n \to S$. By the previous lemma and its proof, $S_n \to S$ is surjective and locally on $S$, it is a $L^+G/K_n(\Fq)$-torsor. Thus $S_n \to S$ is a finite \'etale cover with \'etale group $L^+G/K_n(\Fq)$. Then the trivialization of $\MC{G}[n]^{\varphi}$ by the universal isomorphism together with the last assertions of the previous lemma, yields
\begin{align*}
 (\MC{G}[n], \varphi[n]) & \cong (\MC{G}[n]^{\varphi} \times^{L^+G/K_n(\Fq)} L^+G/K_n, \sigma^* \times \id_{L^+G/K_n}) \\
   & \cong (L^+G/K_n(\Fq) \times^{L^+G/K_n(\Fq)} L^+G/K_n, \sigma^* \times \id_{L^+G/K_n}) = (L^+G/K_n, \sigma^*)
\end{align*}
over $S_n$. As all constructions are compatible for different $n$, we may take the projective limit to get the corresponding result for $(\MC{G}, \varphi)$ itself. \exit

\rem{}{
In the case of $G = GL_d$ (or rather in the case of vector bundles) this was already stated in \cite[theorem 2.5]{BoeHarUnif}, although over rigid analytic varieties as base spaces $S$. Let us explain this relationship in the case of schemes $S$ (for simplicity over $\Fq$): 
Consider an \'etale local shtuka $(\MC{F}, \varphi)$ consisting of a locally free $\MC{O}_S \otimes_{\Fq} \Fq[[z]]$-module of rank $d$ together with an isomorphism $\varphi: \sigma^*\MC{F} \to \MC{F}$. Let $(\MC{F}[n], \varphi)$ be the associated locally free $\MC{O}_S \otimes_{\Fq} \Fq[[z]]/ z^{n+1}$-module with the induces $\sigma$-linear isomorphism $\varphi$ (in \cite{BoeHarUnif} $(\MC{F}[n], \varphi)$ is called $(\MC{F}/\MF{a}\MC{F}, \tau)$). Consider $\MC{F}[n]$ as a geometric vector bundle $F[n]$ over $S$ of rank $d(n+1)$ and $\MC{F}$ as a vector bundle $F$. These come together with a group structure making it locally isomorphic to $(\M{G}_a^{n+1})^d$ resp. $(\M{G}_a^{\M{Z}_{\geq 0}})^d$ and with a compatible action of $\Fq[[z]]/z^{n+1}$ resp. $\Fq[[z]]$. 
One has a subgroup-scheme $F[n]^{\varphi} \subset F[n]$ consisting of the fix-points of $\varphi$ (called $_{\MF{a}}\un{\MC{F}}$ in \cite{BoeHarUnif}). Then \cite[theorem 2.5(a)]{BoeHarUnif} states that $F[n]^{\varphi}$ trivializes over a finite \'etale cover and \cite[theorem 2.5(b)]{BoeHarUnif} gives an isomorphism $\MC{F}[n]^{\varphi} \otimes_{\Fq} \MC{O}_S \cong \MC{F}[n]$ (when viewed as the corresponding sheaves of modules on $S$). \\
Consider now the equivalence between locally free sheaves $\MC{F}$ and $L^+GL_d$-torsors (where we give torsors by the sheaf they represent). $\MC{F}$ corresponds to the torsor
\[\MC{G} \coloneqq \MC{I}som_{\MC{O}_S \otimes_{\Fq} \Fq[[z]]}(\MC{F}, (\MC{O}_S \otimes_{\Fq} \Fq[[z]])^d) \cong \MC{I}som_{z, grp}(L^+\M{G}_a^d, F)\]
where $\MC{I}som_{z, grp}$ denotes the sheaf of isomorphisms of group schemes, which are compatible with the $z$-action. Similarly we have
\[\MC{G}[n] \coloneqq \MC{I}som_{\MC{O}_S \otimes_{\Fq} \Fq[[z]]/z^{n+1}}(\MC{F}[n], (\MC{O}_S \otimes_{\Fq} \Fq[[z]]/z^{n+1})^d) \cong \MC{I}som_{z, grp}(L^+\M{G}_a[n]^d, F[n])\]
(where $L^+\M{G}_a[n]^d$ is the group scheme associated to the locally free sheaf $(\MC{O}_S \otimes_{\Fq} \Fq[[z]]/z^{n+1})^d$ or equivalently the group scheme representing the functor $S' \mapsto \M{G}_a^d(\MC{O}_{S'} \otimes_{\Fq} \Fq[[z]]/z^{n+1})$) and 
\[\MC{G}[n]^{\varphi} \coloneqq \MC{I}som_{\Fq[[z]]/z^{n+1}}(\MC{F}[n]^{\varphi}, (\Fq[[z]]/z^{n+1})^d) \cong \MC{I}som_{z, grp}(L^+\M{G}_a[n]^d(\Fq), F[n]^{\varphi}).\]
Hence \cite[theorem 2.5(a)]{BoeHarUnif} is equivalent to the trivialization of the $L^+GL_d/K_n(\Fq)$-torsor $\MC{G}[n]^{\varphi}$ after a finite \'etale cover and \cite[theorem 2.5(b)]{BoeHarUnif} translates into the isomorphism given in lemma \ref{lem:TateTrivialEtale}. 
}

\const{}{Tate-functor for local $G$-shtukas}{
Fix any $n \geq 0$ and let $(\MC{G}, \varphi) \in \acute{E}tSht_G(S)$. Then we have just seen that the torsor $\MC{G}[n]^{\varphi}$ trivializes over the universal cover $\widetilde{S}$ of $S$ (and even over a finite \'etale cover). To fix a canonical trivialization, let $\un{\MC{G}[n]^{\varphi}}$ be the set-theoretic $L^+G/K_n(\Fq)$-torsor consisting of all $\kappa(\ov{s})$-valued points of the torsor $\MC{G}[n]^{\varphi}|_{\tilde{\ov{s}}}$ over $\Sp \kappa(\ov{s})$ (the residue field of the geometric point $\tilde{\ov{s}}$). 
Now two isomorphisms $\MC{G}[n]_{\widetilde{S}}^{\varphi} \cong \un{\MC{G}[n]^{\varphi}}_{\widetilde{S}}$ differ by an automorphism of the right-hand side, i.e. by an element of $\Aut(\un{\MC{G}[n]^{\varphi}}) \cong L^+G/K_n(\Fq)$. Hence there is a uniquely determined isomorphism $\xi[n]$, which restricts to the identity over $\tilde{\ov{s}} \in \widetilde{S}$. Then for any $\gamma \in \pi_1(S, \ov{s})$, viewed as an element $\gamma \in \Aut_S(\widetilde{S})$, we have an isomorphism $\MC{G}[n]_{\widetilde{S}} \to \gamma^*(\MC{G}[n]_{\widetilde{S}})$ coming from the $\gamma$-invariance of the torsor $\MC{G}$ over $S$. 
This induces the following diagram
\[\begin{aligned}
 \begin{xy} \xymatrix{
   \MC{G}[n]^{\varphi}_{\widetilde{S}} \ar^-{\xi[n]}[rrr] \ar[d] &&& \un{\MC{G}[n]^{\varphi}}_{\widetilde{S}} \ar@{.>}[d] \\
   \gamma^*(\MC{G}[n]^{\varphi}_{\widetilde{S}}) \ar^-{\gamma^*(\xi[n])}[r] & \gamma^*(\un{\MC{G}[n]^{\varphi}}_{\widetilde{S}}) \ar@{=}[r] & \un{\MC{G}[n]^{\varphi}}_{  \gamma^*\widetilde{S}} \ar@{=}[r] & \un{\MC{G}[n]^{\varphi}}_{\widetilde{S}}
 } \end{xy}
\end{aligned} \qquad (\star) \]
Hence using the induced arrow on the right, we get a morphism
\[\rho_{\MC{G}}[n]: \pi_1(S, \ov{s}) \to \Aut(\un{\MC{G}[n]^{\varphi}})\]
and therefore the functor
\[\begin{array}{cccc}
 T[n]: & \acute{E}tSht_G(S) & \to & \Rep(\pi_1(S), L^+G/K_n(\Fq)) \\
         & (\MC{G}, \varphi) & \mapsto & \rho_{\MC{G}}[n].
\end{array}\]
Now define the (dual) Tate-functor as
\[\begin{array}{cccccc}
 T: & \acute{E}tSht_G(S) & \to & \Rep(\pi_1(S), L^+G(\Fq)) \\
     & (\MC{G}, \varphi) & \mapsto & \rho_{\MC{G}} \coloneqq \varprojlim_{n} \rho_{\MC{G}}[n].
\end{array}\]
We call $T(\MC{G})$ the (dual) Tate module associated to the \'etale local $G$-shtuka $(\MC{G}, \varphi)$.
}

\rem{\label{rem:TateLimit}}{
First note that by definition $\MC{G} \cong \varprojlim_n \MC{G}[n]$ (and the same for the pullback to $\widetilde{S}$). Furthermore writing $\un{\MC{G}^\varphi}$ for the set-theoretic $L^+G(\Fq)$-torsor consisting of all $\kappa(\ov{s})$-valued points of the torsor $\MC{G}[n]^{\varphi}|_{\tilde{\ov{s}}}$ over $\Sp \kappa(\ov{s})$, we have 
\[\xi: \MC{G}_{\widetilde{S}}^{\varphi} = (\varprojlim_n \MC{G}[n])_{\widetilde{S}}^{\varphi} = \varprojlim_n (\MC{G}[n]_{\widetilde{S}}^{\varphi}) \xrightarrow{\varprojlim \xi[n]} \varprojlim_n (\un{\MC{G}[n]^{\varphi}}_{\widetilde{S}}) = \un{\MC{G}^{\varphi}}_{\widetilde{S}}\] 
which is compatible with all the actions of $\pi_1(S, \ov{s})$ on the torsors $\un{\MC{G}[n]^{\varphi}}$. In fact this action $\rho_\MC{G}$ on the limit describes precisely how the trivial $L^+G(\Fq)$-torsor over $\widetilde{S}$ descends to the $L^+G(\Fq)$-torsor $\MC{G}^{\varphi}$ over $S$. Hence an alternative way to describe the dual Tate-functor is
\[\begin{array}{cccccc}
 T: & \acute{E}tSht_G(S) & \to & \Rep(\pi_1(S), L^+G(\Fq))\\
     & (\MC{G}, \varphi) & \mapsto & \rho_{\MC{G}}: \pi_1(S, \ov{s}) \to \Aut(\un{\MC{G}^{\varphi}})
\end{array}\]
where for any $\gamma \in \pi_1(S, \ov{s})$ the element $\rho_{\MC{G}}(\gamma)$ is defined via the commutative diagram $(\star)$ but for $\MC{G}$ instead of $\MC{G}[n]$ and $\xi$ instead of $\xi[n]$.
}

The next construction follows the ideas presented in \cite[proposition 1.3.7]{HartlPeriod} or \cite[proposition 3.6]{HarRad1}:

\const{}{The inverse of the Tate-functor}{
Let $\rho \in \Rep(\pi_1(S), L^+G(\Fq))$ acting on the $L^+G(\Fq)$-torsor $\un{L^+G}_{\widetilde{S}}$. For any $n \geq 0$ consider the trivial $L^+G/K_n$-torsor $\widetilde{\MC{G}}[n] \coloneqq (\un{L^+G}_{\widetilde{S}} \times^{L^+G(\Fq)} L^+G/K_n(\Fq)) \times^{L^+G/K_n(\Fq)} L^+G/K_n$ over $\widetilde{S}$ together with the Frobenius-linear morphism $\sigma^*[n]: \sigma^*\widetilde{\MC{G}}[n] \to \widetilde{\MC{G}}[n]$ acting only on the last factor. Then for every $\gamma \in \pi_1(S, \ov{s})$ the automorphism $\rho(\gamma)$ defines a unique isomorphism of $L^+G/K_n$-torsors $\rho[n](\gamma): \widetilde{\MC{G}}[n] \to \widetilde{\MC{G}}[n] = \gamma^*\widetilde{\MC{G}}[n]$. As $\widetilde{\MC{G}}[n]^{\sigma^*} = \un{L^+G}_{\widetilde{S}} \times^{L^+G(\Fq)} L^+G/K_n(\Fq)$ by construction, the isomorphism $\rho[n](\gamma)$ commutes with the $\sigma^*$-action. \\
Let us now descend $(\widetilde{\MC{G}}[n], \sigma^*[n])$ to $S$: First note that the kernel of $\pi_1(S, \ov{s}) \to \Aut(\un{L^+G}_{\widetilde{S}}) \to \Aut(\un{L^+G}_{\widetilde{S}} \times^{L^+G(\Fq)} L^+G/K_n(\Fq))$ has finite index in $\pi_1(S, \ov{s})$. Thus there is a finite \'etale cover $S_n \to S$, such that $(\widetilde{\MC{G}}[n], \sigma^*[n])$ descends to the trivial $L^+G/K_n$-torsor over $S_n$ (together with the Frobenius-linear isomorphism $\sigma^*[n]$) and the $\pi_1(S, \ov{s})$-action factors over $\Aut_S(S_n)$. Thus we obtain an \'etale descent datum for the finite cover $S_n$ to $S$. In this way we obtain a unique pair $(\MC{G}[n], \varphi[n])$ over $S$. \\
Finally note that we have canonical isomorphisms $\MC{G}[n] \times^{L^+G/K_n} L^+G/K_{n-1} \cong \MC{G}[n-1]$ compatible with the morphisms $\varphi[n]$. Hence we may define the quasi-inverse of the Tate-functor as
\[T^{-1}(\rho) \coloneqq \varprojlim_n(\MC{G}[n], \varphi[n])\]
}

\rem{}{
After defining $\widetilde{\MC{G}} \coloneqq \un{L^+G}_{\widetilde{S}} \times^{L^+G(\Fq)} L^+G$, we have a canonical isomorphism
\[\varprojlim_n (\MC{G}[n], \varphi[n])_{\widetilde{S}} \cong (\widetilde{\MC{G}}, \sigma^*)\]
of \'etale local $G$-shtukas over $\widetilde{S}$. Furthermore the canonical descent datum on the left-hand side coincides with the action of $\pi_1(S, \ov{s})$ defined by $\rho$. Thus $T^{-1}(\rho)$ is nothing else than $(\widetilde{\MC{G}}, \sigma^*)$ after descending it along the pro-finite \'etale cover $\widetilde{S} \to S$. 
}

\thm{}{}{
The functor $T$ is an equivalence of categories.
}

\prooof
Denote the category of trivial $L^+G$-torsors $\widetilde{\MC{G}}$ over $\widetilde{S}$ together with a descent datum of $(\widetilde{\MC{G}}, \sigma^*)$ from $\widetilde{S}$ to $S$ by $\MC{H}^1(\widetilde{S}, L^+G)^{\op{descent}}$. The triviality of the $L^+G(\Fq)$-torsors on $\widetilde{S}$ implies that we have an equivalence of categories
\[\Rep(\pi_1(S), L^+G(\Fq)) \cong \MC{H}^1(\widetilde{S}, L^+G)^{\op{descent}}\]
But by proposition \ref{prop:TateTrivialTorsor} every \'etale local $G$-shtuka trivializes over $\widetilde{S}$ and the canonical (and obviously fully faithful) functor
\[\MC{H}^1(\widetilde{S}, L^+G)^{\op{descent}} \to \acute{E}tSht_G(S)\]
is essentially surjective. Hence $T$ is indeed an equivalence. \exit

\subsection{Tate's theorem on extending quasi-isogenies}\label{subsec:TateTheorem}
The main result of the section is theorem \ref{thm:ThmTateGShtuka}: On connected normal schemes every quasi-isogeny between local $G$-shtukas over the generic fiber can be extended uniquely to a quasi-isogeny over the whole scheme. The analogue for $p$-divisible groups is called the Theorem of Tate and was shown by Tate \cite[theorem 4]{TatePDivisible} in characteristic $0$ and by Berthelot \cite[section 4.1]{BerthelotDieu} in positive characteristic. 
The following argumentation follows mainly the ideas of Berthelot. \\
The results of this section will not be needed until section \ref{sec:FormalLifting}.

\defi{}{}{
Let $A$ a $\Fq$-algebra with its Frobenius $\sigma$ (over $\Fq$). A shtuka over $A$ is a pair $(\Lambda, \varphi)$ consisting of a locally free $A[[z]] = A \widehat{\otimes}_{\Fq} \Fq[[z]]$-module $\Lambda$ (of finite rank $n$) and a $\sigma$-linear morphism $\varphi: \sigma^*\Lambda \to \Lambda$ inducing an isomorphism after inverting $z$. \\
A homomorphism $\alpha: (\Lambda_1, \varphi_1) \to (\Lambda_2, \varphi_2)$ is a morphism between the $A[[z]]$-modules $\Lambda_i$ satisfying $\alpha \circ \varphi_1 = \varphi_2 \circ \sigma^*\alpha$. It is called an isogeny if it induces an isomorphism after inverting $z$.
}

\rem{}{
i) Although we defined shtukas in this generality, we will only use it for complete valuation rings $A$ with algebraically closed residue field $k$, or quotient fields thereof. Note that in these cases $\Lambda$ is necessarily free. \\
ii) If one replaces $A[[z]]$ with the ring of Witt vectors $W(A)$, we get almost the definition of crystals in mixed characteristic: There we have to give in addition to $\Lambda$ and $\varphi$ a connection on $\Lambda$, which satisfies certain compatibility conditions. This connection is mainly needed to have a canonical comparison between the crystals over the perfection of $A$ for different choices of liftings of the Frobenius $\sigma$ (cf. \cite{BerthelotDieu}, explanation after theorem 4.6.1). As this problem disappears in equal characteristic, we may dispense with the connection as we did in the definition. Note that whenever we refer for proofs to the case of mixed characteristic (notably the next two statements), these do not use the connection except for passing to the perfection of $A$. \\
iii) There is a canonical functor from shtukas of rank $n$ to local $GL_n$-shtukas over $\Sp A$: To $(\Lambda, \varphi)$ associate the trivial $L^+GL_n$-torsor $Isom_{A[[z]]}(\Lambda, A[[z]]^n)$ and the Frobenius isomorphism induced by $\varphi$. For further information cf. \cite[\S 4]{HaVi}. Thus we may translate all definitions made above for local $G$-shtukas to the case of shtukas. In particular we have Newton points of shtukas, which here determine the quasi-isogeny class of the shtuka, i.e. the isomorphism class of $(\Lambda[\frac 1z], \varphi)$.
}

\prop{\label{prop:IsogenyTheorem}}{}{
Let $A$ be a complete valuation ring of characteristic $p$ with perfect residue field $k$, and $(\Lambda, \varphi)$ a shtuka over $A$ with constant Newton point. Let $\lambda_1, \ldots, \lambda_s$ be the distinct slopes which appear with multiplicities $d_1, \ldots, d_s$. Then $(\Lambda, \varphi)$ is isogenous to a shtuka $(\Lambda', \varphi')$ which is completely slope divisible, i.e. there exists a filtration
\[0 \subset (\Lambda_1', \varphi_1') \subset \ldots \subset (\Lambda_s', \varphi_s') = (\Lambda', \varphi')\]
via subshtukas such that each composition factor $(\Lambda_i'/\Lambda_{i-1}', \ov{\varphi}_i')$ is a shtuka of rank $d_i$ and constant Newton point isoclinic of slope $\lambda_i$. Furthermore on each $\Lambda_i'/\Lambda_{i-1}'$ the morphism $z^{-\lambda_i N}\ov{\varphi}_i'^N$ defines an isogeny (and hence an isomorphism) for $N \gg 0$ sufficiently divisible. \\
This filtration admits a unique splitting over the perfection $A^{\sharp} = k[[t^{1/p^{\infty}}]]$ of $A$.
}

\prooof
This proposition is the analog of a result of Katz \cite[corollary 2.6.3]{KatzSlope}. His proof, which includes the proofs of theorem 2.4.2, theorem 2.5.1, theorem 2.6.1 and corollary 2.6.2 in loc.cit., translates word for word to our situation. \exit

\rem{}{
i) In the case of $p$-divisible groups the same result holds by Oort and Zink \cite[corollary 2.2]{OortZinkFamilies} over any normal base scheme. It seems likely that the same holds in equal characteristic for shtukas and also for arbitrary local $G$-shtukas. Nevertheless the result of Katz is sufficient for our purposes and his proofs are written in the language of crystals rather than $p$-divisible groups which allows a much easier translation to our situation. \\
ii) We will define completely slope divisible local $G$-shtukas in the next section. One easily checks that it matches the definition above in the case $G = GL_n$. 
}

\prop{\label{prop:ThmTateGShtukaDVR}}{}{
Let $A$ be a complete valuation ring of characteristic $p$ with uniformizer $\unif$ and algebraically closed residue field $k$. Let $(\MC{G}_i, \varphi_i)$ (for $i = 1, 2$) be two local $G$-shtukas over $S = \Sp A$ with constant Newton points. Let $\eta = \Sp A[\unif^{-1}]$ be the generic point and denote the generic fibers by $(\MC{G}_{i, \eta}, \varphi_{i, \eta})$. Then the canonical map
\[\op{QIsog}((\MC{G}_1, \varphi_1), (\MC{G}_2, \varphi_2)) \to \op{QIsog}((\MC{G}_{1, \eta}, \varphi_{1, \eta}), (\MC{G}_{2, \eta}, \varphi_{2, \eta}))\]
is an isomorphism.
}

\prooof
Note first that every $L^+G$-torsor over $S = \Sp A$ is trivial. Indeed by isotriviality (cf. \cite[lemma XIV 1.4]{Ray119}) any $G$-torsor over $\Sp A[z]/(z^n)$ (for some $n > 0$) trivializes over a finite \'etale cover. But there are no non-trivial ones, as by \cite[IX, proposition 1.7]{GroSGAI} they correspond to finite unramified extensions of $k$ (see lemma \ref{lem:UnifLem2}a) for more details). Hence any $G$-torsor over $\Spf A[[z]]$ is trivial, or equivalently any $L^+G$-torsor over $S$ is trivial. Thus we may identify $\varphi_i = b_i\sigma^*$ with $b_i \in LG(A) = G(A((z)))$ (for both $i = 1, 2$). Thus we have to show:
\[\{g \in LG(A) \;|\; g^{-1}b_1\sigma^*(g) = b_2\} = \{g_{\eta} \in LG(A[\unif^{-1}]) \;|\; g_{\eta}^{-1}b_{1, \eta}\sigma^*(g_{\eta}) = b_{2, \eta}\}\]
where $b_{i, \eta}$ is the image of $b_i$ in $LG(A[\unif^{-1}])$, i.e. defines the local $G$-shtuka at the generic point. Note that the inclusion $\subseteq$ is obvious and we will concentrate on the other one. \\
Pick now any faithful representation $G \to GL_n$ (defined over $k$, thus ignoring all subtleties coming from non-splitness). Then $LG(A) = LGL_n(A) \cap LG(A[\unif^{-1}])$. Hence we may wlog. assume that $G = GL_n$. In particular we may view elements in $LGL_n(A)$ respectively $LGL_n(A[\unif^{-1}])$ as isomorphisms of the free module $\Lambda' \coloneqq A((z))^n$ respectively $\Lambda'_\eta \coloneqq A[\unif^{-1}]((z))^n$ of rank $n$. Thus the assertion rewrites as
\begin{align*}
 & \{g \in \Aut(\Lambda') \;|\; g^{-1}b_1\sigma^*(g) = b_2: \sigma^*\Lambda' \to \Lambda'\} =  \\
 & \hspace{1,5cm} = \{g_{\eta} \in \Aut(\Lambda'_\eta) \;|\; g_{\eta}^{-1}{b_1}_\eta\sigma^*(g_{\eta}) = {b_2}_\eta: \sigma^*\Lambda'_{\eta} \to \Lambda'_{\eta}\}.
\end{align*}
Let $\Lambda = A[[z]]^n$ and $\Lambda_\eta \coloneqq A[\unif^{-1}][[z]]^n$. Then after replacing $b_1$ and $b_2$ by $z^Nb_1$ and $z^Nb_2$ for some $N \gg 0$ we may assume $b_i\sigma^*: \sigma^*\Lambda \to \Lambda$ are homomorphisms already defined on $\Lambda$. Furthermore by the same replacement for $g$ and $g_\eta$ (although we may not find such an $N$ for all $g$ or $g_\eta$ simultaneously) we are left to show
\begin{align*}
 & \{g \in \End^*(\Lambda) \;|\; g^{-1}b_1\sigma^*(g) = b_2: \sigma^*\Lambda \to \Lambda\} =  \\
 & \hspace{1,5cm} = \{g_{\eta} \in \End^*(\Lambda_\eta) \;|\; g_{\eta}^{-1}{b_1}_\eta\sigma^*(g_{\eta}) = {b_2}_\eta: \sigma^*\Lambda_{\eta} \to \Lambda_{\eta}\}
\end{align*}
where $End^*$ denotes all endomorphisms which induce isomorphisms after inverting $z$. In other words the left-hand side denotes the isogenies between the shtukas $(\Lambda, b_1\sigma^*)$ and $(\Lambda, b_2\sigma^*)$ and the right-hand side denotes those over the generic fiber. \\
In this situation we can copy Berthelot's proof of theorem 4.7.1 in \cite{BerthelotDieu} (after the usual translations from mixed to equal characteristic) which uses proposition \ref{prop:IsogenyTheorem} (or rather its analog in mixed characteristic) to reduce to the isoclinic case. \exit

\thm{\label{thm:ThmTateGShtuka}}{{\rm (Theorem of Tate for local $G$-shtukas on base schemes in equal characteristic)}}{
Let $S$ be a normal connected scheme with function field $K$ of characteristic $p$. Let $(\MC{G}_i, \varphi_i)$ (for $i = 1, 2$) be two local $G$-shtukas over $S$ with constant Newton points. Denote their generic fibers by $(\MC{G}_{i, \eta}, \varphi_{i, \eta})$. Then the canonical map
\[\op{QIsog}((\MC{G}_1, \varphi_1), (\MC{G}_2, \varphi_2)) \to \op{QIsog}((\MC{G}_{1, \eta}, \varphi_{1, \eta}), (\MC{G}_{2, \eta}, \varphi_{2, \eta}))\]
is an isomorphism.
}

\prooof
Having a canonical isomorphism allows us to use \'etale decent. Hence it suffices to prove the corollary after passing to an \'etale cover and we may assume that $S = \Sp A$ is affine with function field $K$ and both $L^+G$-torsors are trivializable. In this situation we may identify quasi-isogenies over all of $S$ with elements in $LG(A) = G(A((z)))$ and quasi-isogenies over the generic fiber with elements in $LG(K) = G(K((z)))$. Hence the map above is obviously injective. Surjectivity immediately follows from the \\
\textit{Claim:} Any $g \in G(K((z)))$ defining a quasi-isogeny over the generic fiber is defined over $A((z))$. \\ 
As $A$ was assumed to be normal, $A$ is the intersection of all localizations $A_{\MF{p}}$ for prime ideals $\MF{p}$ of height $1$. As the same is true for $A((z))$ we may replace $A$ by $A_{\MF{p}}$ and have to show the claim for locals rings of dimension $1$, i.e. valuation rings. 
If $\widehat{A}$ is the completion of $A$ at the unique maximal ideal, we have $A = \widehat{A} \cap K$. Thus we are reduced to complete valuation rings. If $A^{nr}$ is the maximal unramified extension of $A$, then $A = A^{nr} \cap K$. Thus we may actually assume that $A$ has an algebraically closed residue field. Now the claim is actually equivalent to the previous proposition. \exit

\rem{\label{rem:IsogenyTheoremRemark}}{
i) Note that in the previous two statements we require only the constancy of the Newton point, as the Kottwitz point is constant on the connected schemes anyway. Hence it is automatic, that the quasi-isogeny class of the local $G$-shtukas is constant, too. \\
ii) The assumption on the constancy of the Newton point is essential. de Jong showed similar statements in \cite{JongHomosBTGroup} for $p$-divisible groups with varying Newton point, though he allows degeneration phenomena, i.e. extensions as simple homomorphisms that are no longer necessarily quasi-isogenies. Such degenerated homomorphisms have no direct analogue as morphisms of $G$-shtukas, so it would be hard to even formulate a conjecture mirroring de Jong's result in the world of $G$-shtukas.
}

\section{Global \textit{G}-shtukas and their moduli space}\label{sec:GlobalShtuka}
We turn now to the global situation and study global $G$-shtukas: $G$-torsors $\MS{G}$ over the fixed curve $C$ together with a $\sigma$-linear morphism
\[\varphi: \sigma^*\MS{G}|_{(C \times_{\Fq} S) \setminus \bigcup_i \Gamma_{c_i}} \to \MS{G}|_{(C \times_{\Fq} S) \setminus \bigcup_i \Gamma_{c_i}}\]
of $G$-torsors, which is defined outside a finite number of $S$-valued points $c_i \in C(S)$, called the characteristic places of the global $G$-shtuka. Here we abbreviate $\Gamma_{c_i}$ for the graph of the point $c_i$ in $C \times_{\Fq} S$. \\
We first prove that the local behavior of a global $G$-shtuka can be described by local $G$-shtukas, which will be most useful when the $c_i$ are constant along $S$. This allows us to define two boundedness conditions for global $G$-shtukas:
\begin{itemize}
 \item For general characteristic places there are boundedness conditions in the flavor of Varshavsky \cite[definition 2.4b)]{Varsh}, cf. definition \ref{def:GlobalBoundVarsh}.
 \item If the characteristic places are constant along $S$ coming from $\F$-valued points in $C$, then one may call a global $G$-shtuka bounded if the associated local $\Res_{\F/\Fq}(G)$-shtukas are bounded in the sense of \ref{Def:BoundLocal}.
\end{itemize}
Although this second definition turns out to be a special case of the first boundedness condition, its local nature makes it more useful for analyzing the corresponding moduli space. Furthermore level structures away from the characteristic places will be introduced. 
Finally we consider the moduli space of bounded global $G$-shtukas with level structures, which exists in general only as a DM-stack. However any open compact substack admits a finite \'etale cover defined by some increased level structure, which is representable as a scheme. But it is not known to us, whether there exists one level structure such that the corresponding moduli space is globally representable by a scheme.

\subsection{Global \textit{G}-shtukas without level structure}\label{subsec:GlobalShtukaDef}
In the following we will define global $G$-shtukas, ignoring for now the more subtle notions of boundedness conditions and level structures. \\
As usual $C$ is the fixed smooth projective curve over $\Fq$ and as already mentioned above, we denote by $\Gamma_{c_i} \subset C \times_{\Fq} S$ the graph of a $S$-valued point $c_i \in C(S)$. As in the previous sections $S$ is any DM-stack over $E$.

\defi{\label{def:Torsors}}{}{
Fix some topology $* \in \{fpqc, fppf, \acute{e}tale\}$ (though the category will not depend on this choice, cf. \cite[section 6]{GroFGAI}). Then let $\MC{H}^1(C, G)$ be the category fibered in groupoids over the category of DM-stacks over $E$, whose objects $\MC{H}^1(C, G)(S)$ over some DM-stack $S$ (over $E$) are the $G$-torsors $\MS{G}$ (for the chosen topology) over $C \times_{\Fq} S$. The morphisms are given by isomorphisms of $G$-torsors.
}

\rem{}{
i) $\MC{H}^1(C, G)$ is in fact a smooth Artin stack locally of finite type. This is shown in \cite[proposition 1]{Heinloth2010} in far greater generality. \\
ii) The Frobenius $\sigma: C \times_{\Fq} S \to C \times_{\Fq} S$ (which acts as the identity on $C$ and as the absolute $q$-Frobenius on $S$) acts on $\MC{H}^1(C, G)$ by pullback along this morphism. 
}

\defi{\label{def:GlobShtuka}}{}{
a) A global $G$-shtuka with $n$ characteristic places (but without level structure) over a DM-stack $S$ over $E$ is a triple consisting of
\begin{itemize}
 \item a $G$-torsor $\MS{G} \in \MC{H}^1(C, G)(S)$ over $C \times_{\Fq} S$,
 \item an $n$-tuple $c_1, \ldots, c_n \in C(S)$ of distinct $S$-valued points, 
 \item and an isomorphism 
\[\varphi: \sigma^*\MS{G}|_{C \times_{\Fq} S \setminus \bigcup_i \Gamma_{c_i}} \xrightarrow{\;\sim\;} \MS{G}|_{C \times_{\Fq} S \setminus \bigcup_i \Gamma_{c_i}}\]
\end{itemize}
The points $c_i \in C(S)$ are called the characteristic places and $\varphi$ is called the Frobenius-isomorphism on $\MS{G}$. \\
b) The moduli stack of global $G$-shtukas (without level structure) is the stack $\nabla_n\MC{H}^1(C, G)$ over $\Sp E$ classifying (up to isomorphism) global $G$-shtukas with $n$ characteristic places. Two global $G$-shtukas $(\MS{G}, \varphi)$ and $(\MS{G}', \varphi')$ are isomorphic if their characteristic places coincide (as ordered tuples) and there is an isomorphism $f: \MS{G} \to \MS{G}'$ of $G$-torsors such that
\[f \circ \varphi = \varphi' \circ \sigma^*f.\]
Let $\Delta \subset C^n$ be the complement of the locus of distinct points in $C$. \ignore{be the locus where in at least two components the points in $C$ coincide.} Then there is a canonical morphism $\nabla_n\MC{H}^1(C, G) \to C^n \setminus \Delta$ mapping a global $G$-shtuka to its characteristic points. \\
c) Let $c_1, \ldots, c_n \in C(\F)$ be any distinct points on $C$ (for some finite field $\F$). Then denote by $\nabla_{(c_i)}\MC{H}^1(C, G)$ the moduli stack of global $G$-shtukas (without level structure) with characteristic places exactly the $c_i$. \\
d) Similarly denote by $\nabla_{(\hat{c_i})}\MC{H}^1(C, G)$ the moduli stack of global $G$-shtukas, such that the characteristic places vary only in the formal neighborhood of the points $c_i$. 
}

\rem{\label{rem:DefiGlobalShtuka}}{
i) When no confusion is possible, we will omit mentioning the characteristic places and denote a global $G$-shtuka (without level structure) simply by $(\MS{G}, \varphi)$. \\
ii) An alternative description of $\nabla_n\MC{H}^1(C, G)$ (even with some kind of level structure) is given in \cite[definition 3.4]{HarRad2}. \\
iii) Fixing local coordinates $\zeta_i$ at $c_i \in C$, the canonical morphism $\nabla_{(\hat{c}_i)}\MC{H}^1(C, G) \to C^n \setminus \Delta$ factors over the formal completion at the $c_i$, which by choice of the local coordinates equals $\Spf \F[[\zeta_1, \ldots, \zeta_n]]$. This way $\nabla_{(\hat{c}_i)}\MC{H}^1(C, G)$ can be viewed as a moduli stack over $\Sp E \times_{\Fq} \Spf \F[[\zeta_1, \ldots, \zeta_n]]$ or more conveniently as a moduli stack over $\Spf E[[\zeta_1, \ldots, \zeta_n]]$ (by forgetting the $\F$-action). \\
iv) The symmetric group $S_n$ acts on $C^n \setminus \Delta$ and $\nabla_n\MC{H}^1(C, G)$ by permuting the characteristic places. As this action gives canonical isomorphisms between the fibers of $\nabla_n\MC{H}^1(C, G) \to C^n \setminus \Delta$ over the points in one $S_n$-orbit, we might consider the moduli stack $\nabla_n\MC{H}^1(C, G)/S_n$ over $(C^n \setminus \Delta)/S_n$ parametrizing now $G$-torsors, a set (and not a tuple) of $n$ distinct characteristic places and a Frobenius-isomorphism away from them. All further statements are valid (after appropriate modifications) in this unordered setup as well.
} 

\defi{\label{def:QIsog}}{}{
Let $(\MS{G}, (c_i), \varphi), (\MS{G}', (c_i'), \varphi') \in \nabla_n\MC{H}^1(C, G)(S)$ be two global $G$-shtukas over $S$. \\
a) A quasi-isogeny is an equivalence class of pairs $(\alpha, U)$ consisting of a dense open substack $U \subset C \times_{\Fq} S$ and an isomorphism $\alpha: \MC{G}|_U \to \MC{G}'|_U$ of $G$-torsors over $U$ satisfying
\[\alpha \circ \varphi = \varphi' \circ \sigma^*(\alpha)\]
after restricting each of these morphisms to $U \cap \setminus (\bigcup_i \Gamma_{c_i} \cup \Gamma_{c_i'})$. \\
Two such pairs $(\alpha_1, U_1)$ and $(\alpha_2, U_2)$ are equivalent if there is an open substack $U \subset U_1 \cap U_2$ such that $\alpha_1|_U = \alpha_2|_U$. \\
b) Let $D \subset C \times_{\Fq} S$ be some closed substack. A $D$-quasi-isogeny is an isogeny such that there exists a representative $(\alpha, U)$ with $U = (C \times_{\Fq} S) \setminus D$, i.e. the isomorphism $\alpha$ of $G$-torsors can be extended outside of $D$.
}

\rem{\label{rem:QIsogFrob}}{
Consider a global $G$-shtuka $(\MS{G}, (c_i), \varphi) \in \nabla_n\MC{H}^1(C, G)(S)$. Then it is easy to check that 
\[(\varphi, C \times_{\Fq} S \setminus \bigcup\nolimits_i \Gamma_{c_i}): (\sigma^*\MS{G}, (\sigma^*c_i), \sigma^*\varphi) \to (\MS{G}, (c_i), \varphi)\]
is a quasi-isogeny.
}

The next lemma will be needed in section \ref{subsec:GlobalBounds}.

\lem{\label{lem:QIsogExtend}}{}{
Let $D \subset C \times_{\Fq} S$ be an effective divisor and $(\alpha, U): (\MS{G}, (c_i), \varphi) \to (\MS{G}', (c_i'), \varphi')$ be (a representative of) a $D$-quasi-isogeny between global $G$-shtukas over a quasi-compact DM-stack $S$ over $E$. 
Let $E'$ be a finite extension of $E$ and $\rho$ be any representation of $G$ on a finite dimensional $E'$-vector space $V$. Consider the vector bundles $\MS{G}_V$ and $\MS{G}'_V$ as constructed in \ref{const:GenTorsorVector}. The quasi-isogeny $\alpha$ induces now a canonical morphism 
\[\alpha_V|_U: \MS{G}_V|_U \to \MS{G}'_V|_U\]
over the open substack $U$. Then for each sufficiently big integer $N \gg 0$ the morphism $\alpha_V|_U$ extends to
\[\alpha_V: \MS{G}_V \to \MS{G}'_V \otimes_{\MC{O}_{C \times_{\Fq} S \times_E \Sp E'}} \MC{O}_{C \times_{\Fq} S \times_E \Sp E'}(-N \cdot (D \times_E \Sp E'))\]
of vector bundles over all of $C \times_{\Fq} S \times_E \Sp E'$. This $\alpha_V$ is equivariant for the $G$-action defined via $\rho$.
}

\prooof
As $\alpha_V$ is obviously unique, we may check the extension property \'etale locally. Therefore assume for simplicity that $S$ is a scheme. As $\MS{G}_V$ is a vector bundle of finite rank, $\alpha_V$ is determined locally around a point $x \in D \times_E \Sp E'$ by its value on finitely many sections of $\MS{G}_V$. Hence for $N_x \gg 0$ sufficiently large $\alpha_V$ extends to a morphism $\alpha_V: \MS{G}_V \to \MS{G}'_V \otimes_{\MC{O}_{C \times_{\Fq} S \times_E \Sp E'}} \MC{O}(-N_x \cdot (D \times_E \Sp E'))$ on an open neighborhood of $x$. 
As $S$ is quasi-compact, there is a finite set of points on $D$ such that their open neighborhoods cover $D \times_E \Sp E'$. Thus we may choose one $N \gg 0$ sufficiently large for all these points $x$ and we get an effective divisor $N \cdot (D \times_E \Sp E')$ such that $\alpha_V|_U$ extends to a morphism of sheaves
\[\alpha_V: \MS{G}_V \to \MS{G}'_V \otimes_{\MC{O}_{C \times_{\Fq} S \times_E \Sp E'}} \MC{O}_{C \times_{\Fq} S}(-N \cdot (D \times_E \Sp E'))\]
on the whole scheme $C \times_{\Fq} S \times_E \Sp E'$. \\
It remains to show that $\alpha_V$ is $G$-equivariant, i.e. the following diagram commutes
\[
\begin{xy}
 \xymatrix{
   G \times_{\Fq} \MS{G}_V \ar^{\id \times \alpha_V}[r] \ar[d] & G \times_{\Fq} \MS{G}'_V \ar[d] \\
   \MS{G}_V \ar^{\alpha_V}[r] & \MS{G}'_V}
\end{xy} 
\]
where the vertical maps are given by the $G$-action. By definition this holds over $U$. But any morphism $(G \times_{\Fq} \MS{G}_V)|_U \to \MS{G}'_V|_U$ extends uniquely to $C \times_{\Fq} S$ (if it extends at all). Hence the diagram commutes indeed over all of $C \times_{\Fq} S \times_E \Sp E'$. \exit

\rem{}{
This proof nowhere needs the existence of the Frobenius-morphisms $\varphi$ and $\varphi'$. Thus the lemma (and its proof) remains valid, if one replaces the global $G$-shtukas by $G$-torsors and the quasi-isogeny $\alpha$ by any morphism between $G$-torsors defined over an open dense subset $U$.
}

\subsection{The local description of global \textit{G}-shtukas}\label{subsec:GlobalLocalFunctor}
Fix a finite extension $\F$ of $\Fq$ and any point $v \in C(\F)$. Our aim is to construct a morphism of stacks
\[\MF{L}_v: \nabla_n\MC{H}^1(C, G) \to Sht_{\Res_{\F/\Fq}(G)}\]
where $\Res_{\F/\Fq}(G)$ is the restriction of scalars. This functor should remember the local behavior of both the $G$-torsor and the Frobenius-isomorphism at the point $v$. \\
In the case of a non-characteristic place $v$, the image is contained in the substack of \'etale local $\Res_{\F/\Fq}(G)$-shtukas and we will use this construction in section \ref{subsec:AdelicLevel} to define level structures. If $v$ equals one of the characteristic places (which is then necessarily constant along $S$), the image of the functor $\MF{L}_{c_i}$ is very close to be surjective, though it will take us until section \ref{subsec:Uniform} to see this. This already indicates that much of the information about the global $G$-shtuka is preserved under the functor 
\[\MF{L} = \prod_{c_i} \MF{L}_{c_i}: \nabla_{(c_i)}\MC{H}^1(C, G) \to \prod_{c_i} Sht_{\Res_{\F/\Fq}(G)}.\]
Note that the functor $\MF{L}_v$ was first constructed by Hartl and Rad in \cite[section 5.2]{HarRad1}. Nevertheless we hope that the reader appreciates the slightly more conceptual approach on this construction presented here, including an explanation how to get the local Frobenius-isomorphism. \\
Throughout this section we will never (except in lemma \ref{lem:LoGloRes}a)) specify the underlying Grothendieck topology. So (except for definition \ref{def:LoGloTor}, where this does not make much sense) in all statements it is understood that we fix some topology $* \in \{\acute{e}t, fppf, fpqc\}$. We start with some more general nonsense about $G$-torsors and their behavior under restriction of scalars, working at first purely over the category of $E$-schemes $S$:

\lem{\label{lem:LoGloRes}}{}{
Let $R$ be a finite $\Fq$-algebra, $G$ a smooth group scheme over $R$ and $S$ a scheme over $E$. \\
a) For every \'etale morphism $S' \to S \times_{\Fq} \Sp R$ (over $E$), there exists an \'etale morphism $S'' \to S$ (again over $E$), such that $S'' \times_{\Fq} R \to S \times_{\Fq} \Sp R$ factors over $S'$. If $S' \to S \times_{\Fq} \Sp R$ is a cover, then one can choose $S''$ in such a way that $S'' \to S'$ is surjective and $S'' \times_{\Fq} R$ is a refinement of the cover $S'$. \\
b) There are canonical equivalences of categories
\[\MC{R}^+: \MC{H}^1((S \times_{\Fq} \Sp R)_*, L^+G) \xrightarrow{\sim} \MC{H}^1(S_*, L^+\Res_{R/\Fq}(G))\]
\[\MC{R}: \MC{H}^1((S \times_{\Fq} \Sp R)_*, LG) \xrightarrow{\sim} \MC{H}^1(S_*, L\Res_{R/\Fq}(G))\]
which are compatible with associating $LG$-torsors to $L^+G$-torsors.
}

\prooof
a) By \cite[\href{http://stacks.math.columbia.edu/tag/05YD}{Tag 05YD}]{stacks-project}(currently lemma 69.11.3) $\Res_{R/\Fq}(f): \Res_{R/\Fq}(S') \to \Res_{R/\Fq}(S \times_{\Fq} \Sp R)$ is again \'etale and it is surjective, if $f: S' \to S \times_{\Fq} \Sp R$ was. After base change via the canonical morphism $S \to \Res_{R/\Fq}(S \times_{\Fq} \Sp R)$ the same is true for the morphism over $E$
\[S'' \coloneqq \Res_{R/\Fq}(S') \times_{\Res_{R/\Fq}(S \times_{\Fq} \Sp R)} S \xrightarrow{pr_2} S.\]
This gives the commutative diagram
\[
\begin{xy} \xymatrix{
    S'' \times_{\Fq} \Sp R \ar[r] \ar[d] & \Res_{R/\Fq}(S') \times_{\Fq} \Sp R \ar^-{p}[r] \ar^{\Res_{R/\Fq}(f) \times \id}[d] & S' \ar^{f}[d] \\
    S \times_{\Fq} \Sp R \ar^-{\iota}[r] & \Res_{R/\Fq}(S \times_{\Fq} \Sp R) \times_{\Fq} \Sp R \ar[r] & S \times_{\Fq} \Sp R
} \end{xy}
\]
where $p$ comes from the identity on $\Res_{R/\Fq}(S)$ by adjunction and the composite of the two lower horizontal morphisms is the identity. Hence we get the desired factorization. \\ 
Finally we have to show that $S'' \times_{\Fq} \Sp R$ is a refinement of $S'$ assuming that $f: S' \to S \times_{\Fq} \Sp R$ is surjective. It is easy to see that on geometric points (or for points over local artinian rings with algebraically closed residue field): Let $x: \Sp k \to S'$ be any geometric point (over $\Sp E \times_{\Fq} \Sp R$). It has the image $f(x): \Sp k \to S' \to S \times_{\Fq} \Sp R$ which gives an $E \otimes_{\Fq} R$-linear morphism $\Sp k \times_{\Fq} \Sp R \to S \times_{\Fq} \Sp R$. This defines a point in $\Res_{R/\Fq}(S \times_{\Fq} \Sp R)(\Sp k) = (S \times_{\Fq} \Sp R)(\Sp k \times_{\Fq} \Sp R)$, which is $\iota(f(x))$ after base change to $R$. Then we have by definition a commutative diagram over $E \otimes_{\Fq} R$
\[
\begin{xy} \xymatrix{
    \Sp k  \ar^{x}[r] \ar[d] & S' \ar^{f}[d]\\
    \Sp k \times \Sp R \ar^-{\iota(f(x))}[r] & S \times \Sp R
} \end{xy}
\]
where the left vertical map is given by the identity on the first factor and the canonical morphism to $\Sp R$ (as $x$ was a geometric point over $\Sp R$). Hence we may lift the point $\iota(f(x)) \in \Res_{R/\Fq}(S \times_{\Fq} \Sp R)(\Sp k) = (S \times_{\Fq} \Sp R)(\Sp k \times_{\Fq} \Sp R)$ (over $E$) by surjectivity and smoothness of $f$ to a point $y \in \Res_{R/\Fq}(S')(\Sp k) = S'(\Sp k \times_{\Fq} \Sp R)$ (over $E$) such that the corresponding point in $\Res_{R/\Fq}(S') \times \Sp R$ (now over $E \otimes_{\Fq} R$) satisfies $p(y) = x$. Note here that $y$ is far from unique. But because we defined $S''$ as the fiber product, the triple $(y, f(x), \Res_{R/\Fq}(f)(y) = \iota(f(x)))$ defines a $k$-valued point in $S'' \times \Sp R$ (over $E$) which is by $p(y) = x$ a preimage of $x$. \\
b) As we have equivalences between torsors for the respective topologies, we may restrict ourselves to $* = \acute{e}t$. Furthermore we will only consider $\MC{R}^+$ as the arguments for $\MC{R}$ are analogous and the compatibility assertion will be obvious. \\
Given any $L^+G$-torsor $\MC{G}$ over $S \times_{\Fq} \Sp R$, the element $\MC{R}^+\MC{G}$ in $\MC{H}^1(S_{\acute{e}t}, L^+\Res_{R/\Fq}(G))$ is simply given by $\MC{R}^+\MC{G}(S'') \coloneqq \MC{G}(S'' \times_{\Fq} \Sp R)$ for every \'etale $S'' \to S$. If $\MC{G}$ trivializes over some $S'$, take any $S'' \to S$ as in a) and $\MC{G}$ trivializes again over $S'' \times_{\Fq} \Sp R$. Hence $\MC{R}^+\MC{G}$ trivializes over $S''$, i.e. it is indeed a torsor. \\
For the inverse construction start with $\MC{R}^+\MC{G} \in \MC{H}^1(S_{\acute{e}t}, L^+\Res_{R/\Fq}(G))$. Note that it suffices by the refinement statement on covers in a) to define the $L^+G$-torsor $\MC{G}$ over schemes of the form $S'' \times_{\Fq} \Sp R$ for $S'' \to S$ \'etale. But there we may just set $\MC{G}(S'' \times_{\Fq} \Sp R) \coloneqq \MC{R}^+\MC{G}(S'')$. It is obvious that this defines a $L^+G$-torsor and that the constructions are mutually inverse. \exit

\rem{}{
i) In all applications the group scheme $G$ will be a constant group scheme coming by base-change from a reductive group over $\Fq$. \\
ii) It is not known to us, whether $\Res_{R/\Fq}(S')$ (or $\Res_{R/\Fq}(S \times_{\Fq} \Sp R)$ for that matter) exists as a scheme if $S'$ is not quasi-projective. Nevertheless it is still an algebraic space by \cite[\href{http://stacks.math.columbia.edu/tag/05YF}{Tag 05YF}]{stacks-project}(currently proposition 69.11.5). Thus $S''$ is an algebraic space \'etale over the scheme $S$, hence by \cite[corollary 6.17]{KnutAlgSpace} indeed a scheme. \\
If $S$ is quasi-projective, then proofs of all claims concerning restriction of scalars (including existence results) can be found in \cite[section 7.6]{BoschNeron}. \\
iii) Note one subtlety concerning restriction of scalars and the surjectivity result in a): If $\{S_i\}_i$ form an \'etale cover of $S \times \Sp R$ (i.e. are a finite set of jointly surjective \'etale morphisms), then $\{\Res_{R/\Fq}(S_i)\}$ needs not to be such a cover for $S$. Nevertheless $\Res_{R/\Fq}(\coprod_i S_i) \to S$ is still surjective and \'etale. For a more detailed discussion refer to \cite[appendix A.5]{CGPRedGroup}. There one can also find a proof of part b) for linear algebraic groups in the case where $R$ is a field extension of $\Fq$. \\
iv) If one tries to prove part a) for fppf- or fpqc-morphisms, one runs into two problems: First of all we can no longer apply Knutson's algebraization result to ensure to have a scheme $S''$. Secondly (and more seriously), we do not know whether $\Res_{R/\Fq}$ preserves flatness or surjectivity of morphisms, although this seems to hold at least if $\Fq \to R$ is \'etale. 
}

\defi{\label{def:LoGloTor}}{}{
Let $X$ be any locally ringed space with a morphism $X \to \Sp \Fq$ and a fixed Grothendieck topology. Let $G$ be a reductive group over $\Fq$ and denote its pullback to $X$ by $G_X$. This is a group object in the category of locally ringed spaces over $X$. Then a $G$-torsor over $X$ is a locally ringed space $\MS{G}$ over $X$ together with a morphism $G_X \times_X \MS{G} \to \MS{G}$ such that locally (for the fixed Grothendieck topology) $\MS{G}$ admits a trivialization, i.e. a $G_X$-equivariant isomorphism $\MS{G} \cong G_X$ to the trivial torsor. \\
Analogously to definition \ref{def:Torsors}, let $\MC{H}^1(\Fq, G)(X)$ be the category of $G$-torsors over $X$.
}

The next proposition is (at least partially) well-known: Similar statements can be found in \cite[proposition 2.2a)]{HaVi} or \cite[proposition 2.4]{HarRad1}. 
\prop{\label{prop:LoGloCompare}}{}{
Let $S$ be any scheme over $\Fq$ and write $S[[z]]/z^{n+1} \coloneqq S \times_{\Fq} \Sp \Fq[[z]]/z^{n+1}$ for $n \geq 0$. Consider furthermore the formal scheme $S[[z]] \coloneqq S \widehat{\times}_{\Fq} \Spf \Fq[[z]]$ and the locally ringed space $S((z))$ with underlying topological space $|S[[z]]| = |S|$ and structure sheaf $\MC{O}_{S((z))}$ the sheafification of $\MC{O}_{\widehat{S}}(U)[\frac 1z] = \MC{O}_{S}(U) \widehat{\otimes}_{\Fq} \Fq((z))$ (for $U \subset |S|$ open). We view $S[[z]]$ and $S((z))$ with the topology induced from the one on $S$, i.e. if $S' \to S$ is a cover in $*$, then $S'[[z]] \to S[[z]]$ respectively $S'((z)) \to S((z))$ are again covers. \\
a) For every $n \geq 0$ there is a canonical equivalence of categories
\[\MC{H}^1(\Fq, G)(S[[z]]/z^{n+1}) \xrightarrow{\; \cong \;} \MC{H}^1(\Fq, L^+G/K_n)(S)\]
b) There is a canonical equivalence of categories
\[\MC{H}^1(\Fq, G)(S[[z]]) \xrightarrow{\; \cong \;} \MC{H}^1(\Fq, L^+G)(S)\]
c) There is a canonical functor
\[\MC{H}^1(\Fq, G)(S((z))) \longrightarrow \MC{H}^1(\Fq, LG)(S)\]
All these functors are compatible in the sense that we have commutative diagrams
\ignore{
\[
 \begin{xy} \xymatrix{
    G\op{-}Tor_*(S[[z]]) \ar[r] \ar[d] & \MC{H}^1(S_*, L^+G) \ar[d] \\
    G\op{-}Tor_*(S[[z]]/z^{n+1}) \ar[r] & \MC{H}^1(S_*, L^+G/K_n)
 } \end{xy}  \qquad
 \begin{xy} \xymatrix{
    G\op{-}Tor_*(S[[z]]) \ar[r] \ar[d] & \MC{H}^1(S_*, L^+G) \ar[d] \\
    G\op{-}Tor_*(S((z))) \ar[r] & \MC{H}^1(S_*, LG)
 } \end{xy} 
\]
}
\[
 \begin{xy} \xymatrix{
    \MC{H}^1(\Fq, G)(S[[z]]) \ar[r] \ar[d] & \MC{H}^1(\Fq, L^+G)(S) \ar[d] \\
    \MC{H}^1(\Fq, G)(S[[z]]/z^{n+1}) \ar[r] & \MC{H}^1(\Fq, L^+G/K_n)(S)
 } \end{xy}  \qquad
 \begin{xy} \xymatrix{
    \MC{H}^1(\Fq, G)(S[[z]]) \ar[r] \ar[d] & \MC{H}^1(\Fq, L^+G)(S) \ar[d] \\
    \MC{H}^1(\Fq, G)(S((z))) \ar[r] &  \MC{H}^1(\Fq, LG)(S)
 } \end{xy} 
\]
where the vertical arrows on the left-hand sides are given by base-change along the canonical morphism between the base spaces.
}

\prooof
a) The same arguments as in lemma \ref{lem:LoGloRes}b) but for $G$ instead of $L^+G$ give for $R = \Fq[[z]]/z^{n+1}$ an equivalence
\[\MC{H}^1(\Fq, G)(S[[z]]/z^{n+1}) = \MC{H}^1(\Fq, G)(S \times_{\Fq} \Sp R) \cong \MC{H}^1(\Fq, \Res_{R/\Fq}(G))(S).\]
But $\Res_{R/\Fq}(G) = L^+G/K_n$ as they represent the same functor of groups. \\
b) This follows directly from a):
\begin{align*}
 \MC{H}^1(\Fq, G)(S[[z]]) & \cong \varprojlim_n \MC{H}^1(\Fq, G)(S \times_{\Fq} \Sp \Fq[[z]]/(z^{n+1})) \\
 & \cong \varprojlim_n \MC{H}^1(\Fq, L^+G/K_n)(S) \\
 & \cong \MC{H}^1(\Fq, L^+G)(S).
\end{align*}
c) We use essentially the same construction as Hartl and Rad: If $\MS{G}$ is a $G$-torsor over $S((z))$, then consider the sheaf of sets
\[\MC{LG}(S') \coloneqq Hom_{S((z))}(S'((z)), \MS{G})\]
for every cover $S' \to S$. 
It is easy to see that applied to the trivial $G$-torsor over $S((z))$, this defines the trivial $LG$-torsor over $S$, hence this construction defines indeed an element in $\MC{H}^1(\Fq, LG)(S)$. \\
The diagram on the left-hand side is obviously commutative from our construction in b). To see the commutativity on the right-hand side just observe that (as explained in \cite[proposition 2.4]{HarRad1}) we may explicitly describe the functor in b) as mapping $\MS{G} \in \MC{H}^1(\Fq, G)(S((z)))$ to the torsor representing 
\[\MC{G}(S') \coloneqq Hom_{S[[z]]}(S'[[z]], \MS{G}).\]
\exit

\rem{\label{rem:LoGloTorRemark}}{
i) Note that we do not claim to have an equivalence in c), although it seems reasonable to assume it is. The main problem in proving this is to establish representability results or a descent theory for locally ringed spaces of the form $S((z))$. However we will never need this equivalence. \\
ii) An alternative way to describe the functors in a) and b) can be found in \cite[section 5.1]{HarRad1}: Let $R = \Fq[[z]]$ or $R = \Fq[[z]]/(z^{n+1})$. Then we may apply the functor $\Res_{R/\Fq}$ to go from $G$-torsors over $S \widehat{\times} \Spf R$ to $\Res_{R/\Fq}(G \widehat{\times} \Spf R)$-torsors over $\Res_{R/\Fq}(S \widehat{\times} \Spf R)$ (working in the category of formal algebraic spaces). Then the pull-back along $S \to \Res_{R/\Fq}(S \widehat{\times}_{\Fq} \Spf R)$ gives the desired torsor over $S$ (at least as algebraic space). \\
A similar construction can be done in the situation of lemma \ref{lem:LoGloRes}, though it would be harder to show that the resulting object is indeed a torsor under the desired group. \\
iii) A third way to construct the equivalences in a) and b) is to construct a correspondence of \'etale descend data on both sides. For more details see \cite[proof of proposition 2.2(a)]{HaVi}.
}

After these preparations let us tackle the construction of the global-local functor $\MF{L}_v$: \\
In all constructions let $v \in C(\F)$ be some point over a finite field extension $\F$ of $\Fq$ and fix a local coordinate $z$ of $C$ at the point $v$. Abbreviate $G_v = \Res_{\F/\Fq}(G)$ for the restriction of scalars of $G$. \\
Fix now a global $G$-shtuka $(\MS{G}, \varphi) \in \nabla_{(c_i)}\MC{H}^1(C, G)(S)$ over an $E$-scheme $S$ with some characteristic places $c_i$.

\const{\label{Con:LoGloTorsor}}{\textbf{The local} $\mathbf{L^+G_v}$ \textbf{\!-torsor}}{
Let $Z = \{v\} \times_{\Fq} S \subset C \times_{\Fq} S$ be the graph of $v$. Then the completion $\widehat{Z}$ of $C \times_{\Fq} S$ along $Z$ is isomorphic to $\Spf \F[[z]] \widehat{\times}_{\Fq} S = (\Sp \F \times_{\Fq} S)[[z]]$ by using the chosen local coordinate $z$. Restricting the $G$-torsor $\MS{G}$ to $\widehat{Z}$, i.e. considering the base-change $\MS{G} \times_{C \times_{\Fq} S} \widehat{Z}$, defines an element $\MS{G}|_{\widehat{Z}}$ in $\MC{H}^1(\Fq, G)(\widehat{Z}) = \MC{H}^1(\Fq, G)((\Sp \F \times_{\Fq} S)[[z]])$. 
By proposition \ref{prop:LoGloCompare}b), this corresponds to an element $\widehat{\MS{G}} \in \MC{H}^1(\Fq, L^+G)(\Sp \F \times_{\Fq} S)$ and by lemma \ref{lem:LoGloRes} this is equivalent to a $L^+G_v$-torsor $\MC{R}^+\widehat{\MS{G}} \in \MC{H}^1(\Fq, L^+G_v)(S)$.
}

\rem{}{
i) Note that we could have applied directly the functor $\Res_{\F/\Fq}$ to the $G$-torsor $\MS{G}|_{\widehat{Z}}$ in order to produce an element in $\MC{H}^1(\Fq, G_v)(\Res_{\F/\Fq} \widehat{Z})$ which can be pulled back to an element in $\MC{H}^1(\Fq, G_v)(S[[z]])$. Then proposition \ref{prop:LoGloCompare}b) produces the desired element in $\MC{H}^1(\Fq, L^+G_v)(S)$. \\
ii) By construction we may reconstruct the torsor $\MS{G}|_{\widehat{Z}}$ out of the $L^+G_v$-torsor over $S$.
}

\const{\label{Con:LoGloFrobEtale}}{\textbf{The local Frobenius-isomorphism at non-characteristic places}}{
Assume for now that $v$ is disjoint from all characteristic places $c_i$. Then $\varphi$ is defined on an open neighborhood of $v$ and we may restrict $\varphi$ to $\widehat{Z}$. Thus we obtain a Frobenius-linear isomorphism $\varphi: \sigma^*\MS{G}|_{\widehat{Z}} \to \MS{G}|_{\widehat{Z}}$. As both proposition \ref{prop:LoGloCompare}b) and lemma \ref{lem:LoGloRes}b) were equivalences of categories, we first obtain a morphism $\widehat{\varphi}: \sigma^*\widehat{\MS{G}} \to \widehat{\MS{G}}$ and then $\MC{R}^+\widehat{\varphi}: \sigma^* \MC{R}^+\widehat{\MS{G}} = \MC{R}^+(\sigma^*\widehat{\MS{G}}) \to \MC{R}^+\widehat{\MS{G}}$.
}

This allows us to define the global-local functor $\MF{L}_v$ (over the category of $E$-schemes $S$) for points $v$ not meeting any characteristic place as
\begin{align*}
 \MF{L}_v: \nabla_{(c_i)}\MC{H}^1(C, G) & \to \acute{E}tSht_{G_v} \\
 (\MS{G}, \varphi) & \mapsto(\MC{R}^+\widehat{\MS{G}}, \MC{R}^+\widehat{\varphi})
\end{align*}

\rem{}{
i) Actually the definition of $\MF{L}_v$ obviously extends to all global $G$-shtukas, whose characteristic places do not meet $v$. \\
ii) Assume that $\F'/\F$ is another finite field extension and $v' \in C(\F')$ denotes the point $v$ viewed as a $\F'$-valued point. Then we can factor $\MF{L}_{v'}$ as
\[\MF{L}_{v'} = (- \times^{L^+G_v} L^+G_{v'}) \circ \MF{L}_{v}\] 
where we use the canonical morphism $G_v = \Res_{\F/\Fq}(G) \to \Res_{\F'/\F}(\Res_{\F/\Fq}(G)) = \Res_{\F'/\Fq}(G) = G_{v'}$ to get
\begin{align*}
 (- \times^{L^+G_v} L^+G_{v'}): \acute{E}tSht_{G_v} & \to \acute{E}tSht_{G_{v'}} \\
 (\MC{G}, \varphi) & \mapsto (\MC{G} \times^{L^+G_v} L^+G_{v'}, \varphi \times id_{LG_{v'}})
\end{align*}
The proof of this is an easy exercise by going through all the constructions. 
The same factorization will hold as well if $v$ coincides with a characteristic place.
}

The definition of the Frobenius-isomorphism is more complicated at points $v = c_i$ coinciding with some characteristic place. 
We first construct the $LG_v$-torsor in a convenient way:

\const{\label{Con:LoGloTorsorRat}}{\textbf{The local} $\mathbf{LG_v}$ \textbf{\!-torsor}}{
Recall that we have a morphism $\iota: (\Sp \F \times_{\Fq} S)[[z]] \cong \widehat{Z} \hookrightarrow C \times_{\Fq} S$. Fix an arbitrary open immersion $j: U \to C$. We claim that we can extend $\iota$ to a commutative diagram of morphisms between locally ringed spaces
\[
 \begin{xy} \xymatrix{
    (\Sp \F \times_{\Fq} S)((z)) \ar^-{\iota'}[r] \ar[d] & (C \times_{\Fq} S, (j \times \id)_*\MC{O}_{U \times_{\Fq} S}) \ar[d] \\
    (\Sp \F \times_{\Fq} S)[[z]] \ar^-{\iota}[r] & C \times_{\Fq} S.
 } \end{xy} 
\]
Indeed if $v \in U$, then even $\iota$ factors over $(U, \MC{O}_U) \hookrightarrow (C \times_{\Fq} S, (j \times \id)_*\MC{O}_{U \times_{\Fq} S}) \hookrightarrow (C \times_{\Fq} S, \MC{O}_{C \times_{\Fq} S})$. If $v \notin U$ this argument shows, that it suffices to see that a local coordinate at the point $v$ gets mapped to an invertible element in $\MC{O}_{(\Sp \F \times_{\Fq} S)((z))}$, which is obviously true as $z$ represents such a local coordinate. Hence we have indeed such a morphism $\iota'$. \\
Thus we may pull back the $G$-torsor $\MS{G}$ over $C \times_{\Fq} S$ to a $G$-torsor $\MS{G}|_{(\Sp \F \times_{\Fq} S)((z))}$ over $(\Sp \F \times_{\Fq} S)((z))$ via either of the two ways. As $\MS{G}|_{(\Sp \F \times_{\Fq} S)((z))}$ comes by pull-back from a torsor over $(\Sp \F \times_{\Fq} S)[[z]]$, it trivializes indeed over a locally ringed space of the form $S'((z))$ for $S' \to \Sp \F \times_{\Fq} S$ \'etale. Proposition \ref{prop:LoGloCompare}c) gives now an element $\widehat{\MC{L}\MS{G}} \in \MC{H}^1(\Fq, LG)(\Sp \F \times_{\Fq} S)$, which by lemma \ref{lem:LoGloRes} is equivalent to a $LG_v$-torsor $\MC{R}\widehat{\MC{L}\MS{G}} \in \MC{H}^1(\Fq, LG_v)(S)$. \\
As $\MS{G}|_{(\Sp \F \times_{\Fq} S)((z))}$ can be constructed as the pull-back of $\MS{G}|_{(\Sp \F \times_{\Fq} S)[[z]]}$, the compatibility assertions in proposition \ref{prop:LoGloCompare} and lemma \ref{lem:LoGloRes} yield a canonical identification
\[\MC{R}\widehat{\MC{L}\MS{G}} \cong \MC{L}(\MC{R}^+\widehat{\MS{G}}).\]
}

\const{\label{Con:LoGloFrobChar}}{\textbf{The local Frobenius-isomorphism at characteristic places}}{
By definition $\varphi$ can be seen as an isomorphism of $G$-torsors
\[\varphi: \sigma^*\MS{G}|_{(C \times_{\Fq} S, (j \times \id)_*\MC{O}_{U \times_{\Fq} S})} \to \MS{G}|_{(C \times_{\Fq} S, (j \times \id)_*\MC{O}_{U \times_{\Fq} S})}\]
for $U = C \setminus \bigcup_i c_i$ the open complement of all characteristic places. Hence its pull-back to $(\Sp \F \times_{\Fq} S)((z))$ via $\iota'$ (cf. the previous construction \ref{Con:LoGloTorsorRat}) defines an isomorphism
\[\varphi|_{(\Sp \F \times_{\Fq} S)((z))}: \sigma^*\MS{G}|_{(\Sp \F \times_{\Fq} S)((z))} \to \MS{G}|_{(\Sp \F \times_{\Fq} S)((z))}.\]
As in construction \ref{Con:LoGloFrobEtale}, we get morphisms $\widehat{\varphi}: \sigma^*\widehat{\MC{L}\MS{G}} \to \widehat{\MC{L}\MS{G}}$ and then $\MC{R}\widehat{\varphi}: \sigma^* \MC{R}\widehat{\MC{L}\MS{G}} = \MC{R}(\sigma^*\widehat{\MC{L}\MS{G}}) \to \MC{R}\widehat{\MC{L}\MS{G}}$. Then the last identification obtained in construction \ref{Con:LoGloTorsorRat} allows to rewrite this as the desired
\[\MC{R}\widehat{\varphi}: \sigma^*\MC{L}(\MC{R}^+\widehat{\MS{G}}) \to \MC{L}(\MC{R}^+\widehat{\MS{G}}).\]
} $\left. \right.$ \\
With these constructions we can define even for a characteristic place $v = c_i$ (again over the category of $E$-schemes $S$):
\begin{align*}
 \MF{L}_{c_i}: \nabla_{(c_i)}\MC{H}^1(C, G) & \to Sht_{G_v} \\
 (\MS{G}, \varphi) & \mapsto(\MC{R}^+\widehat{\MS{G}}, \MC{R}\widehat{\varphi})
\end{align*}

\rem{\label{rem:GlobLocGeneral}}{
i) Note that for any non-characteristic place, the local Frobenius-isomorphism $\MC{R}\widehat{\varphi}$ constructed in \ref{Con:LoGloFrobChar} coincides with the Frobenius-isomorphism induced from $\MC{R}^+\widehat{\varphi}$ as constructed in \ref{Con:LoGloFrobEtale}. \\
ii) This construction also works in the case of non-constant group schemes $G$ like the ones considered in \cite{HarRad1} and \cite{HarRad2}. One only has to be slightly more careful what group objects to consider in order to obtain torsors everywhere. \\
iii) Probably the proper way to define the associated Frobenius-isomorphism at characteristic places, would have been to consider the generic fiber of $\widehat{Z}$ as an adic space (or at least a rigid analytic space), see that $\varphi$ induces an isomorphism of the $G$-torsor over this adic space and finally prove a correspondence between $G$-torsors over such adic spaces and $LG$-torsors. \\
iv) In the very same way one can associate to any two $G$-torsors $\MS{G}$, $\MS{G}'$ over $C \times S$, any $U \subset C$ open and any morphism $\alpha: \MS{G}|_{U \times_{\Fq} S} \to \MS{G}|_{U \times_{\Fq} S}$ the local datum at a point $v$ as $\MC{R}\widehat{\alpha}: \MC{L}\MC{R}^+\widehat{\MS{G}} \to \MC{L}\MC{R}^+\widehat{\MS{G}'}$.
}

View now $\nabla_{(\hat{c}_i)}\MC{H}^1(C, G)$ as a stack over $\Spf E[[\zeta_1, \ldots, \zeta_n]]$ by remark \ref{rem:DefiGlobalShtuka}iii). Assume wlog. that the $\zeta_i$ coincide with the local coordinate $z$ whenever $v = c_i$ equals a characteristic place. Note that this convention will only become important in the next section when comparing boundedness conditions.

\thm{}{}{
Both functors $\MF{L}_v$ and $\MF{L}_{c_i}$ extend to functors
\[\MF{L}_v: \nabla_{(\hat{c}_i)}\MC{H}^1(C, G) \to \acute{E}tSht_{G_v} \times_E \Spf E[[\zeta_1, \ldots, \zeta_n]]\]
and
\[\MF{L}_{c_i}: \nabla_{(\hat{c}_i)}\MC{H}^1(C, G) \to Sht_{G_{c_i}} \times_E \Spf E[[\zeta_1, \ldots, \zeta_n]]\]
of stacks over $\Spf E[[\zeta_1, \ldots, \zeta_n]]$. \\
$\MF{L}_v$ again remembers the local behavior at the constant point $v$, while $\MF{L}_{c_i}$ remembers the local behavior at the (slightly varying) characteristic place $c_i$. 
}

\prooof
To get $\MF{L}_v$ as a morphism of stacks, the main part is to check compatibility with fpqc-descent data. This however is immediate from the construction. \\
For $\MF{L}_{c_i}$ we have to see in addition, that the constructions above work as well for places varying in formal neighborhoods. But it suffices to ensure, that the local coordinate $z = \zeta_i$ defines an isomorphism between the formal completion $\widehat{Z}$ of $C \times_{\Fq} S$ along the graph $Z$ of the characteristic place $c_i$ and $(\Sp \F \times_{\Fq} S)[[z]]$. But this certainly holds as long the characteristic place $c_i$ varies only in a formal neighborhood of a point. \exit

\subsection{Bounded global \textit{G}-shtukas}\label{subsec:GlobalBounds}
We discuss now ways to describe boundedness conditions on quasi-isogenies $\alpha: \MS{G} \to \MS{G}'$ between global $G$-shtukas (and hence on global $G$-shtukas themselves as well), all given by an $n$-tuple $\mmu = (\mu_i)_i$ of dominant coweights, invariant under $\Gamma = Gal(\ACFq/E)$. \\
The most natural way (motivated by \cite[definition 2.4b)]{Varsh}) to do so is definition \ref{def:GlobalBoundVarsh} in the case of $\F = \Fq$: For any dominant coweight $\lambda$ of $G$ we write $\MS{G}_\lambda = \MS{G} \times^G V_G(\lambda)$ (for the Weyl module $V_G(\lambda)$) and require
\[\MS{G}'_\lambda \otimes_{\MC{O}_{C \times_{\Fq} S}} \MC{O}({\textstyle \sum_i} \langle \lambda, \mu_i \rangle \cdot \Gamma_{c_i}) \subseteq \alpha_\lambda(\MS{G}_\lambda) \subseteq \MS{G}'_\lambda \otimes_{\MC{O}_{C \times_{\Fq} S}} \MC{O}({\textstyle \sum_i} -\langle (-\lambda)_{\rm dom}, \mu_i \rangle \cdot \Gamma_{c_i})\]
where we consider $\sum_i \langle \lambda, \mu_i \rangle \cdot \Gamma_{c_i}$ and $\sum_i -\langle (-\lambda)_{\rm dom}, \mu_i \rangle \cdot \Gamma_{c_i}$ as divisors on $C \times_{\Fq} S$ and take the corresponding line bundles. \\
If the characteristic places are constant $\Fq$-valued points on $C$, then this boundedness condition can be expressed using associated local $G$-shtukas. Unfortunately we do not know any such description if the characteristic places are only defined over larger fields. \\
To remedy this problem, we give a slightly more complicated version of this boundedness condition, which now depends on a fixed finite extension $\F$ of $\Fq$. This modification allows us in theorem \ref{thm:CompareLocalGlobalBound} to describe the boundedness condition locally in terms of associated local $\Res_{\F/\Fq}(G)$-shtukas if (and only if) the characteristic places are constant coming from $\F$-valued points on $C$. This reformulation provides the key to study special fibers of the moduli space of bounded global $G$-shtukas in a rather local manner.

\defi{\label{def:GlobalBoundVarsh}}{}{
Fix a finite field extension $\F$ over $\Fq$ and an $n$-tuple $\mmu = (\mu_i)_i$ where each $\mu_i$ is a $\Gamma$-invariant dominant coweight of $\Res_{\F/\Fq}(G)$. \\
a) Let $(\alpha, U): (\MS{G}, (c_i), \varphi) \to (\MS{G}', (c_i), \varphi')$ be a quasi-isogeny between two global $G$-shtukas over a DM-stack $S$ over $E$, which have the same characteristic places. 
Consider for every dominant character $\lambda \in X^*(\Res_{\F/\Fq}(G))$, which we view as a morphism defined over some finite field $E'$ containing $E$, the representation 
\[G \times_{\Fq} \Sp E' \hookrightarrow \Res_{\F/\Fq}(G) \times_{\Fq} \Sp E' \to GL(V_{\Res_{\F/\Fq}(G)}(\lambda))\]
where $V_{\Res_{\F/\Fq}(G)}(\lambda)$ denotes the Weyl module associated to $\Res_{\F/\Fq}(G)$ and $\lambda$. We abbreviate $\MS{G}_\lambda \coloneqq \MS{G} \times^{G} V_{\Res_{\F/\Fq}(G)}(\lambda)$ and similarly for $\MS{G'}$ (both of them are defined over $C \times_{\Fq} S \times_E \Sp E'$). 
Moreover consider the divisors over $C \times_{\Fq} S \times_E \Sp E'$
\[D_{\mmu}(\lambda) = \sum_i -\langle (-\lambda)_{\rm dom}, \mu_i \rangle \cdot \Gamma_{c_i}\]
and
\[D_{\mmu}'(\lambda) = \sum_i \langle \lambda, \mu_i \rangle \cdot \Gamma_{c_i}\]
and the corresponding line bundles $\MC{O}(D_{\mmu})$ and $\MC{O}(D_{\mmu}')$ over $C \times_{\Fq} S \times_E \Sp E'$.
Then $\alpha$ is globally $\F$-bounded by $\mmu$ if for every dominant character $\lambda \in X^*(\Res_{\F/\Fq}(G))$, the associated morphism 
\[\alpha_\lambda|_U: \MS{G}_\lambda|_U \to \MS{G}'_\lambda|_U\]
between associated vector bundles satisfies
\[\MS{G}'_\lambda \otimes_{\MC{O}_{C \times_{\Fq} S \times_E \Sp E'}} \MC{O}(D_{\mmu}') \subseteq \alpha_\lambda(\MS{G}_\lambda) \subseteq \MS{G}'_\lambda \otimes_{\MC{O}_{C \times_{\Fq} S \times_E \Sp E'}} \MC{O}(D_{\mmu})\]
over all of $C \times_{\Fq} S \times_E \Sp E'$. \\
To make sense of this condition, note that by lemma \ref{lem:QIsogExtend}, $\alpha_\lambda$ extends (at least locally on $S$) to all of $C \times_{\Fq} S \times_E \Sp E'$ after tensoring the target bundle by some line bundle given by some sufficiently large divisor. Then all three sheaves make sense as subsheaves of this twisted target bundle. \\
b) A global $G$-shtuka $(\MS{G}, (c_i), \varphi)$ is $\F$-globally bounded by $\mmu$ if $\varphi$ considered as a quasi-isogeny (cf. remark \ref{rem:QIsogFrob}) is $\F$-globally bounded by $\mmu$. \\
c) The stack of global $G$-shtukas which are $\F$-globally bounded by $\mmu$ is denoted by $\nabla_n^{\mmu}\MC{H}^1(C, G)$. It admits again a canonical morphism to $C^n \setminus \Delta$. If $c_1, \ldots, c_n \in C(\F)$ are fixed $\F$-valued points of $C$ then denote $\nabla_{(c_i)}^{\mmu} \MC{H}^1(C, G) \coloneqq \nabla_n^{\mmu}\MC{H}^1(C, G) \times_{C^n \setminus \Delta} (c_i)_i$, where $(c_i)_i$ is seen as a $\F$-valued point in $C^n \setminus \Delta$. 
}

\rem{}{
Note that the same approach works to define boundedness conditions for morphisms between arbitrary $G$-torsors on $C \times_{\Fq} S$, not necessarily having the Frobenius-morphism turning them into a global $G$-shtuka.
}

Before we start the promised comparison with local boundedness conditions, let us state the global counterpart of proposition \ref{prop:BoundClosed}:

\prop{\label{prop:GlobalBoundClosed}}{}{
Let $(\alpha, U): (\MS{G}, (c_i), \varphi) \to (\MS{G}', (c_i), \varphi')$ be a quasi-isogeny between two global $G$-shtukas over a DM-stack $S$ over $E$, which have the same characteristic places. Moreover fix an $n$-tuple $\mmu = (\mu_i)_i$ of $\Gamma$-invariant dominant coweights of $\Res_{\F/\Fq}(G)$. Then the locus where $\alpha$ is globally $\F$-bounded by $\mmu$ is closed in $S$. 
}

\prooof
As for local $G$-shtukas it suffices to check the condition
\[\MS{G}'_\lambda \otimes_{\MC{O}_{C \times_{\Fq} S \times_E \Sp E'}} \MC{O}(D_{\mmu}') \subseteq \alpha_\lambda(\MS{G}_\lambda) \subseteq \MS{G}'_\lambda \otimes_{\MC{O}_{C \times_{\Fq} S \times_E \Sp E'}} \MC{O}(D_{\mmu})\]
only on a generating set of the monoid of dominant characters. Then the very same arguments as in the proof of proposition \ref{prop:BoundClosed} show that the locus where all these inclusions hold is a closed subscheme inside $C \times_{\Fq} S$. By applying upper semi-continuity of fiber dimensions of the morphism $C \times_{\Fq} S \to S$ and the fact that this fiber dimension is $1$ if and only if its fiber is all of $C$, we see that there is even a closed subscheme $Z_0 \subset S$ parametrizing all points $x \in S$ such that the inclusions hold over all of $C \times_{\Fq} x$. This is the desired locus where $\alpha$ is globally $\F$-bounded by $\mmu$. \exit
\vspace{2mm} \\
Assume for the rest of this section that the $c_i$ are $\F$-valued points in $C$. Our aim is now to describe the condition to be globally $\F$-bounded locally in terms of associated local $\Res_{\F/\Fq}(G)$-shtukas, i.e. we want to compare it to 

\defi{}{}{
Let $\F$ be some finite field extension of $\Fq$ and $\mmu = (\mu_i)_i$ an $n$-tuple of $\Gamma$-invariant dominant coweights of $\Res_{\F/\Fq}(G)$ as above. We fix an $n$-tuple of characteristic points $c_1, \ldots, c_n \in C(\F)$. Consider a quasi-isogeny $(\alpha, U): (\MS{G}, \varphi) \to (\MS{G}', \varphi')$ in the category $\nabla_{(\hat{c}_i)}\MC{H}^1(C, G)$. Then $\alpha$ is locally bounded by $\mmu$ if $\alpha$ extends to an isomorphism outside the characteristic places and for each characteristic place $c_i$ the associated quasi-isogeny
\[\MF{L}_{c_i}\alpha: \MC{L}\MF{L}_{c_i}\MS{G} \to \MC{L}\MF{L}_{c_i}\MS{G}'\]
between local $\Res_{\F/\Fq}(G)$-shtukas is bounded by $\mu_i$ (in the sense of \ref{Def:BoundLocal}).
}

Let us now start to reformulate the condition of being locally bounded:

\prop{\label{prop:LocalBoundReformulate}}{}{
Let $\alpha: (\MC{G}, \varphi) \to (\MC{G}', \varphi')$ be a quasi-isogeny between two local $\Res_{\F/\Fq}(G)$-shtukas over $S \in Nilp_{\ENil}$. Then $\alpha$ is bounded by a $\Gamma$-invariant dominant cocharacter $\mu$ if and only if $\alpha$ satisfies for every dominant character $\lambda \in X^*(\Res_{\F/\Fq}(G))$ defined over $E'$ (containing $E$)
\begin{align}
 (z - \zeta)^{\langle \lambda, \mu\rangle} \MC{G}_\lambda' \subseteq \alpha_\lambda(\MC{G}_\lambda) \subseteq (z - \zeta)^{-\langle(-\lambda)_{\rm dom}, \mu\rangle} \MC{G}_\lambda'.
\end{align}
}

\prooof
Assume $\alpha$ is bounded by $\mu$. Then $\alpha^{-1}$ is bounded by $(-\mu)_{\rm dom}$. Hence we get for each $\lambda$:
\[\alpha_\lambda^{-1}(\MC{G}'_\lambda) \subseteq (z - \zeta)^{-\langle(-\lambda)_{\rm dom}, (-\mu)_{\rm dom}\rangle} \MC{G}_\lambda = (z - \zeta)^{-\langle \lambda, \mu \rangle} \MC{G}_\lambda\]
which is nothing else than the inclusion on the left-hand side. \\
Conversely assume now that both inclusions in $(1)$ are valid. We treat first the case that $S = \{s\}$ is the spectrum of an algebraically closed field. Then $\alpha$ is bounded by $\mu$ if and only if $\mu(\alpha)(s) \preceq \mu$ for the partial order on $X_*(T)_{\M{Q}}$ (cf. lemma \ref{Lem:BoundFields}a)). Consider now a Weyl module representation $\rho_{\lambda_0}: \Res_{\F/\Fq}(G) \to GL(V_{\lambda_0})$. 
Lemma \ref{Lem:BoundAfterRepr}a) shows that the right-hand side of $(1)$ is preserved under passage to a Weyl module representation (note that the proof of loc. cit. does not use the equality of Hodge points to show this). The other inclusion follows via the same arguments applied to $\alpha^{-1}$. But for $GL_n$-shtukas we have seen in remark \ref{rem:BoundGLEasyAlternative} that condition $(1)$ implies (and is actually equivalent) to boundedness. Thus we get $\rho_{\lambda_0}(\mu(\alpha)(s)) \preceq \rho_{\lambda_0}(\mu)$ for each $\lambda_0$. Now by \cite[lemma 2.2]{RapoRich} (and after noting that the given proof of the equivalence of part iv) in loc. cit. uses only Weyl module representations) this suffices to conclude $\mu(\alpha)(s) \preceq \mu$. \\
Over arbitrary bases $S$, it suffices to see that for every geometric point $s \in S$ we have $[\mu(\alpha)](s) = [\mu] \in \pi_1(\Res_{\F/\Fq}(G))_{\M{Q}}/\Gamma$, because $[\mu(\alpha)]$ is locally constant anyway. But to check this condition we may restrict everything to $s$. In that situation we just proved that $\alpha|_s$ is bounded by $\mu$ implying in particular the desired $[\mu(\alpha)](s) = [\mu]$.  \exit \vspace{2mm} \\
The next proposition shows how to translate the restriction of scalars appearing in the definition of an associated local $\Res_{\F/\Fq}(G)$-shtuka into the global context:

\prop{\label{prop:LocalBoundReformulate2}}{}{
Let $\MC{G}$ be a $L^+G$-torsor over $S \times_{E} E' \times_{\Fq} \F$ and $\MC{LG}$ the corresponding $LG$-torsor. By lemma \ref{lem:LoGloRes}b) these torsors correspond to a $L^+\Res_{\F/\Fq}(G)$-torsor $\MC{R}^+\MC{G}$ and a $L\Res_{\F/\Fq}(G)$-torsor $\MC{R}\MC{LG} = \MC{L}\MC{R}^+\MC{G}$ over $S \times_{E} E'$. Fix a representation $\rho: \Res_{\F/\Fq}(G) \to V$ defined over $E'$ and consider it via the canonical morphism $\rho_0: G \to \Res_{\F/\Fq}(G) \to V$ as a representation of $G$. 
Then there is a canonical isomorphism
\[\MC{R}\MC{LG} \times^{L\Res_{\F/\Fq}(G)} (V \otimes_{E'} \MC{O}_{S \times_{E} E'}((z))) \times_{\Fq} \Sp \F \to \MC{LG} \times^{LG} (V \otimes_{E'} \MC{O}_{S \times_{E} E' \times_{\Fq} \Sp \F}((z))).\]
Furthermore consider for every $N \in \M{Z}$ the canonical inclusions
\begin{align*}
 \MC{R}^+\MC{LG} \times^{L^+\Res_{\F/\Fq}(G)} (V \otimes_{E'} (z - \zeta)^N\MC{O}_{S \times_{E} E'}[[z]]) & \subset \MC{R}\MC{LG} \times^{L\Res_{\F/\Fq}(G)} (V \otimes_{E'} \MC{O}_{S \times_{E} E'}((z))) \\
\MC{G} \times^{L^+G} (V \otimes_{E'} (z - \zeta)^N\MC{O}_{S \times_{E} E' \times_{\Fq} \Sp \F}[[z]]) & \subset \MC{LG} \times^{LG} (V \otimes_{E'} \MC{O}_{S \times_{E} E' \times_{\Fq} \Sp \F}((z))).
\end{align*}
Then the isomorphism above restricts to an isomorphism
\[\MC{R}^+\!\MC{LG} \times^{L^+\Res_{\F/\Fq}(G)} (V \otimes_{E'} (z - \zeta)^N\MC{O}_{S \times_{E} E'}[[z]]) \times_{\Fq} \Sp \F \to \MC{G} \times^{L^+G} (V \otimes_{E'} (z - \zeta)^N\MC{O}_{S \times_{E} E' \times_{\Fq} \Sp \F}[[z]]).\]
}

\prooof
By the very definition of the functor $\MC{R}$ there is a canonical morphism
\[\MC{R}\MC{LG} \times_{\Fq} \Sp \F \to \MC{LG}\]
over $S \times_{E} E' \times_{\Fq} \F$. It is equivariant with respect to the canonical morphism $L\Res_{\F/\Fq}(G) \times_{\Fq} \Sp \F \to LG$ (coming by adjunction from the identity on $\Res_{\F/\Fq}(G)$). Therefore by the universal property of the equivariant fiber products it induces a canonical morphism 
\begin{align*}
 & \MC{R}\MC{LG} \times^{L\Res_{\F/\Fq}(G)} (V \otimes_{E'} \MC{O}_{S \times_{E} E'}((z))) \times_{\Fq} \Sp \F \\
 & \qquad = (\MC{R}\MC{LG} \times_{\Fq} \Sp \F)\times^{L\Res_{\F/\Fq}(G) \times_{\Fq} \Sp \F} (V \otimes_{E'} \MC{O}_{S \times_{E} E' \times_{\Fq} \Sp \F}((z))) \\
 & \qquad \to \MC{LG} \times^{LG} (V \otimes_{E'} \MC{O}_{S \times_{E} E' \times_{\Fq} \Sp \F}((z))).
\end{align*}
But over a trivializing \'etale cover both sides are isomorphic to $V \otimes_{E'} \MC{O}_{S \times_{E} E' \times_{\Fq} \Sp \F}((z))$. Thus this morphism is in fact an isomorphism. \\
For the second assertion note that all previous arguments work as well if one replaces $\MC{LG}$ by $\MC{G}$, $\MC{R}$ by $\MC{R}^+$ and $\MC{O}_{S \times_{E} E'}((z))$ by $(z - \zeta)^N \MC{O}_{S \times_{E} E'}[[z]]$. Moreover it is clear that the isomorphism defined in this way is just the restriction of the isomorphism defined above. \exit

\rem{}{
The morphisms 
\[\MC{R}\MC{LG} \times_{\Fq} \Sp \F \to \MC{LG}\]
considered in the proof coincides with the composition $\MC{R}\MC{LG} \times_{\Fq} \Sp \F \to \Res_{\F/\Fq}(\MC{LG}) \times_{\Fq} \Sp \F \to \MC{LG}$. Here we identify $\MC{R}\MC{LG}$ with the pullback of $\Res_{\F/\Fq}(\MC{LG})$ along $S \times_{E} E' \to \Res_{\F/\Fq}(S \times_{E} E' \times_{\Fq} \Sp \F)$ (cf. remark \ref{rem:LoGloTorRemark}ii) giving the morphism on the left-hand side. The other morphism comes by adjunction from the identity on $\Res_{\F/\Fq}(\MC{LG})$. 
Note again that due to difficulties to define $\Res_{\F/\Fq}(\MC{LG})$, all morphisms in this remark have to be considered in the category of (ind-)algebraic spaces. This makes the direct definition coming from the definition of $\MC{R}$ the easier and more useful choice.
}

\lem{}{}{
Let $\MC{G}$ be a $L^+G$-torsor over $S \times_{E} E' \times_{\Fq} \Sp \F$ and $\MS{G}$ the corresponding $G$-torsor over $S \times_{E} E' \times_{\Fq} \Spf \F[[z]]$ (cf. proposition \ref{prop:LoGloCompare}b). Let furthermore $\rho: G \to H$ be any morphism between reductive groups. Then 
\[\MC{G} \times^{L^+G} L^+H \qquad \op{and} \qquad \MS{G} \times^G H\]
correspond under the equivalence given in proposition \ref{prop:LoGloCompare}b) (of course applied to $H$ instead of $G$ now).
}

\prooof
This follows immediately from the definition of the equivalence \ref{prop:LoGloCompare}b). \exit

\thm{\label{thm:CompareLocalGlobalBound}}{}{
Let $\alpha$ be a quasi-isogeny in the category $\nabla_{(\hat{c}_i)}\MC{H}^1(C, G)$, where the $c_i \in C(\F)$. Then $\alpha$ is globally $\F$-bounded by $\mmu$ if and only if it locally bounded by the same tuple $\mmu$.
}

\prooof
Both boundedness definitions require $\alpha$ to be an isomorphism outside the characteristic places. So it suffices to see that boundedness of the associated local $\Res_{\F/\Fq}(G)$-shtuka by some $\mu_i$ is equivalent to have for all $\lambda$
\[\MS{G}'_\lambda \otimes_{\MC{O}_{C \times_{\Fq} S \times_E \Sp E'}} \MC{O}(\langle \lambda, \mu_i \rangle \cdot \Gamma_{c_i}) \subseteq \alpha_\lambda(\MS{G}_\lambda) \subseteq \MS{G}'_\lambda \otimes_{\MC{O}_{C \times_{\Fq} S \times_E \Sp E'}} \MC{O}(-\langle (-\lambda)_{\rm dom}, \mu_i \rangle \cdot \Gamma_{c_i})\]
fpqc-locally around $c_i$. To see this let us start with a quasi-isogeny $(\alpha, U): (\MS{G}, (c_i), \varphi) \to (\MS{G}', (c_i), \varphi')$ which is locally bounded by $\mmu = (\mu_i)$. Then by proposition \ref{prop:LocalBoundReformulate} this is (at the characteristic place $c_i$) equivalent to 
\[(z - \zeta)^{\langle \lambda, \mu_i\rangle} (\MF{L}_{c_i}\MS{G}')_\lambda \subseteq \alpha_\lambda((\MF{L}_{c_i}\MS{G})_\lambda) \subseteq (z - \zeta)^{-\langle(-\lambda)_{\rm dom}, \mu_i\rangle} (\MF{L}_{c_i}\MS{G}')_\lambda\]
for each $\lambda$. By proposition \ref{prop:LocalBoundReformulate2} this in turn is equivalent to
\[(z - \zeta)^{\langle \lambda, \mu_i\rangle} {\widehat{\MS{G}}'}_\lambda \subseteq \alpha_\lambda({\widehat{\MS{G}}}_\lambda) \subseteq (z - \zeta)^{-\langle(-\lambda)_{\rm dom}, \mu_i\rangle} {\widehat{\MS{G}}'}_\lambda\]
over $\Sp \F \times_{\Fq} S \times_E \Sp E'$. Here $\widehat{\MS{G}}$ denotes (as in construction \ref{Con:LoGloTorsor}) the $L^+G$-torsor associated to the $G$-torsor over the formal completion $\widehat{\Gamma}_{c_i} \subset C \times_{\Fq} S$ of the graph $\Gamma_{c_i}$ of the characteristic place $c_i$ and ${\widehat{\MS{G}}}_\lambda$ abbreviates $\widehat{\MS{G}} \times^{L^+G} V_{\Res_{\F/\Fq}(G)}(\lambda)$. 
Using now the previous lemma for $H = GL(V_{\Res_{\F/\Fq}(G)}(\lambda))$ and the canonical equivalence between $GL_n$-torsors and vector bundles, we see that under the equivalence of proposition \ref{prop:LoGloCompare}b) ${\widehat{\MS{G}}}_\lambda$ corresponds to $\MS{G}|_{\widehat{\Gamma}_{c_i}} \times^{G} V_{\Res_{\F/\Fq}(G)}(\lambda)$ (using the representation $G \hookrightarrow \Res_{\F/\Fq}(G) \to GL(V_{\Res_{\F/\Fq}(G)}(\lambda))$). Moreover as we chose the same local coordinate for $z$ and $\zeta$, the ideal $(z - \zeta)^N$ corresponds to the line bundle $\MC{O}_{\widehat{\Gamma}_{c_i}}(N \cdot \Gamma_{c_i})$ for any integer $N$. Thus the inclusions of $\MC{O}_{\Sp \F \times_{\Fq} S \times_E \Sp E'}[[z]]$-vector bundles above are equivalent to 
\begin{align*}
 & (\MS{G}'|_{\widehat{\Gamma}_{c_i}} \times^{G} V_{\Res_{\F/\Fq}(G)}(\lambda)) \otimes_{\MC{O}_{\Spf \F[[z]] \times_{\Fq} S \times_E \Sp E'}} \MC{O}_{\widehat{\Gamma}_{c_i}}(\langle \lambda, \mu_i\rangle \cdot \Gamma_{c_i}) \\
 & \qquad \subseteq \alpha_\lambda(\MS{G}|_{\widehat{\Gamma}_{c_i}} \times^{G} V_{\Res_{\F/\Fq}(G)}(\lambda)) \\
 & \qquad \subseteq (\MS{G}'|_{\widehat{\Gamma}_{c_i}} \times^{G} V_{\Res_{\F/\Fq}(G)}(\lambda)) \otimes_{\MC{O}_{\Spf \F[[z]] \times_{\Fq} S \times_E \Sp E'}} \MC{O}_{\widehat{\Gamma}_{c_i}}(-\langle(-\lambda)_{\rm dom}, \mu_i\rangle \cdot \Gamma_{c_i})
\end{align*}
But this is nothing else than the condition coming from being globally $\F$-bounded considered fpqc-locally around $c_i$. \exit

\subsection{Adelic level structures}\label{subsec:AdelicLevel}
We define now level structures for global $G$-shtukas. For this we fix some finite reduced subscheme $D_0 \subset C$ and let $\M{A}_{D_0}^{int}$ be the ring of integral adeles of $C$ at the places in $D_0$. We will always assume that the characteristic places of the global $G$-shtukas are contained in $(C \setminus D_0) \times_{\Fq} S$.
Then fix an open subgroup $U = \prod_{v \in D_0} U_v \subset G(\M{A}_{D_0}^{int})$. For any such $v$, denote its residue field by $\kappa(v)$ and abbreviate $G_v = \Res_{\kappa(v)/\Fq}(G)$. Each $U_v$ acts then on the $L^+G(\Fq)$-torsor appearing in the Tate module of the associated local $G_v$-shtuka at the place $v$ (cf. section \ref{subsec:TateFunctor}). We then define $U$-level structures by fixing for each $v$ some suitable $U_v$-orbit under this action. \\
Furthermore we give (for suitable $U$) a ``naive'' characterization of a $U$-level structure (cf. proposition \ref{prop:AdelicCompareLocal}b) and \ref{prop:AdelicCompareNaive}) and identify it with the one Hartl and Rad defined in \cite[section 6]{HarRad2} (cf. proposition \ref{prop:AdelicCompareRad}b)). \\
If the characteristic places are concentrated in (formal neighborhoods of) some fixed $\F$-valued points, then we may replace $\M{A}_{D_0}^{int}$ by $\M{A}^{int (c_i)_i}$ the ring of all integral adeles of $C$ away from the characteristic places $c_i$. All arguments in this section work as well for $\M{A}^{int (c_i)_i}$ if one replaces $D_0$ by $C \setminus \bigcup_i c_i$ whenever that makes sense (even if this is not explicitly stated).
\vspace{1mm} \\
We start by defining adelic level structures on \'etale local $G_v$-shtukas, where $v$ is a fixed point of $C$.

\defi{\label{def:AdelicLocal}}{}{
Let $U_v \subset L^+G_v(\Fq)$ be an open subgroup and $(\MC{G}, \varphi) \in \acute{E}tSht_{G_v}(S)$ an \'etale local $G_v$-shtuka over a DM-stack $S$. Let $T(\MC{G}) = (\rho_{\MC{G}}: \pi_1(S, \ov{s}) \to \Aut(\un{\MC{G}^{\varphi}}))$ be its image under the Tate functor. Then an adelic $U$-level structure on $(\MC{G}, \varphi)$ consists of a $\pi_1(S, \ov{s})$-invariant $U_v$-orbit $\psi$ in $\un{\MC{G}^{\varphi}}$. \\
An isomorphism between two \'etale local $G$-shtukas with adelic $U_v$-level structure is an isomorphism of the \'etale local $G$-shtukas, such that the induced morphism on their Tate modules gives a bijection between the respective level structures.
}

\rem{}{
i) If $U_v = L^+G_v(\Fq)$, then there is only one choice of an adelic $U$-level structure, namely taking $\psi = \un{\MC{G}^{\varphi}}$. This way we get an equivalence between the category of \'etale local $G_v$-shtukas without level structure and the category of \'etale local $G_v$-shtukas with adelic $L^+G_v(\Fq)$-level structure. \\
ii) Let $(\MC{G}, \varphi) \in \acute{E}tSht_{G_v}(S)$ be an \'etale local $G_v$-shtuka and consider the quasi-isogeny $\varphi: (\sigma^*\MC{G}, \sigma^*\varphi) \to (\MC{G}, \varphi)$. Let now $\psi$ be a $U_v$-level structure on $(\MC{G}, \varphi)$. Then $\sigma^*\psi \subset \sigma^*(\un{\MC{G}^{\varphi}}) = \un{(\sigma^*\MC{G})^{\sigma^*\varphi}}$ is a $U_v$-level structure on $(\sigma^*\MC{G}, \sigma^*\varphi)$. Moreover $\varphi: (\sigma^*\MC{G}, \sigma^*\varphi) \to (\MC{G}, \varphi)$ becomes a quasi-isogeny between local $G_v$-shtukas with $U_v$-level structure, because it induces the morphism $\sigma^*: \un{(\sigma^*\MC{G})^{\sigma^*\varphi}} = \sigma^*(\un{\MC{G}^{\varphi}}) \to \un{\MC{G}^{\varphi}}$ on the Tate module.
}

\lem{\label{lem:AdelicLocalNormalQuotient}}{}{
Let $U_v \subset U_v' \subset L^+G_v(\Fq)$ be two open subgroups such that $U_v$ is normal in $L^+G_v(\Fq)$. Let $(\MC{G}, \varphi) \in \acute{E}tSht_{G_v}(S)$. Then we have a bijection between $U_v'$-level structures on $(\MC{G}, \varphi)$ and $\pi_1(S, \ov{s})$-invariant $U_v'/U_v$-orbits $\ov{\psi}$ in the $L^+G_v(\Fq)/U_v$-torsor $\un{\MC{G}^{\varphi}} \times^{L^+G_v(\Fq)} L^+G_v(\Fq)/U_v$.
}

\prooof
Note first that as $\pi_1(S, \ov{s})$ acts through elements of $L^+G_v(\Fq)$ on $\un{\MC{G}^{\varphi}}$ and $U_v$ is normal in $L^+G_v(\Fq)$, we get indeed a $\pi_1(S, \ov{s})$-action on $\un{\MC{G}^{\varphi}} \times^{L^+G_v(\Fq)} L^+G_v(\Fq)/U_v$. Thus the canonical bijection between $U_v'$-orbits in $\un{\MC{G}^{\varphi}}$ and $U_v'/U_v$-orbits in $\un{\MC{G}^{\varphi}} \times^{L^+G_v(\Fq)} L^+G_v(\Fq)/U_v$ preserves the property of being $\pi_1(S, \ov{s})$-invariant. \exit

\prop{\label{prop:AdelicCompareLocal}}{}{
Let $U_v \subset L^+G_v(\Fq)$ be an open subgroup and $(\MC{G}, \varphi) \in \acute{E}tSht_{G_v}(S)$ an \'etale local $G_v$-shtuka. Assume that there exists a subgroup-scheme $\MB{U}_v \subset L^+G_v$ (defined over $\Fq$) whose $\Fq$-valued points are exactly $U_v \subset L^+G_v(\Fq)$.  \\
a) Then there is a bijection (using the notations of definition \ref{def:LocShtukaGeneral}b))
\[
 \begin{Bmatrix} \textnormal{adelic} \; U_v\textnormal{-level} \\ \textnormal{structures} \; \psi \end{Bmatrix} \xleftrightarrow{\;1:1\;} \begin{Bmatrix} (\MC{U}, \varphi_{\MC{U}}) \in \MB{U}\op{-}\acute{E}tSht(S) \; \textnormal{and} \\ \psi^\sharp: (\MC{U}, \varphi_{\MC{U}}) \times^{\MB{U}_v} L^+G \cong (\MC{G}, \varphi) \end{Bmatrix}
\]
where we consider elements the right-hand side only up to the equivalence class defined by isomorphisms in $\MB{U}_v\op{-}\acute{E}tSht(S)$ respecting $\psi^\sharp$, i.e. $(\MC{U}_1, \varphi_{\MC{U} 1}, \psi^\sharp_1) \cong (\MC{U}_2, \varphi_{\MC{U} 2}, \psi^\sharp_2)$, if there exists an isomorphism $\xi: (\MC{U}_1, \varphi_{\MC{U} 1}) \cong (\MC{U}_2, \varphi_{\MC{U} 2})$ satisfying $\psi^\sharp_2 = \psi^\sharp_1 \circ (\xi \times \id)$. \\
b) Assume furthermore that $\MB{U}_v \subset L^+G_v$ is normal. Let $\ov{\MC{G}} = \MC{G} \times^{L^+G_v} L^+G_v/\MB{U}_v$ and $\ov{\varphi}: \sigma^*\ov{\MC{G}} \to \ov{\MC{G}}$ be the morphism induced by $\varphi$. Then there is a bijection
\[
 \begin{Bmatrix} \textnormal{adelic} \; U_v\textnormal{-level} \\ \textnormal{structures} \; \psi \end{Bmatrix} \xleftrightarrow{\;1:1\;}  \begin{Bmatrix} \textnormal{isomorphisms in} \; L^+G_v/\MB{U}_v\op{-}\acute{E}tSht(S) \\ \ov{\psi^\sharp}: (\ov{\MC{G}}, \ov{\varphi}) \cong (L^+G/\MB{U}, \sigma^*) \end{Bmatrix}
\]
where $(L^+G_v/\MB{U}_v, \sigma^*)$ is the trivial element in $L^+G_v/\MB{U}_v\op{-}\acute{E}tSht(S)$.
}

\prooof
a) Assume we are given an adelic $U_v$-level structure $\psi$ on $(\MC{G}, \varphi)$, i.e. a $\pi_1(S, \ov{s})$-invariant $U_v$-subtorsor $\un{U} \subset \un{\MC{G}^\varphi}$. This defines a trivial $U_v = \MB{U}_v(\Fq)$-torsor $\un{U}_{\tilde{S}} \subset \un{\MC{G}^{\varphi}}_{\tilde{S}} \subset \MC{G}_{\tilde{S}}$. Then via the morphism $\un{\MB{U}_v(\Fq)}_{\tilde{S}} \subset {\MB{U}_v}_{\tilde{S}}$ we get a $\MB{U}_v$-subtorsor $\tilde{\MC{U}} \subset \MC{G}_{\tilde{S}}$. By definition $\tilde{\MC{U}}$ is invariant under the $\pi_1(S, \ov{s})$-action on $\MC{G}_{\tilde{S}}$ and $\varphi_{\tilde{S}}$ restricts to a Frobenius-isomorphism on $\tilde{\MC{U}}$. Hence $\tilde{\MC{U}}$ descends to a $\MB{U}_v$-subtorsor $\MC{U} \subset \MC{G}$ over $S$ which admits the Frobenius-isomorphism $\varphi|_{\MC{U}}$. \\
Conversely assume we have $\psi^\sharp: (\MC{U}, \varphi_{\MC{U}}) \times^{\MB{U}_v} L^+G_v \cong (\MC{G}, \varphi)$ or equivalently a $\MB{U}_v$-subtorsor $\MC{U} \subset \MC{G}$ with $\varphi_{\MC{U}} = \varphi|_{\MC{U}}$. This gives $\MC{U}^{\varphi_{\MC{U}}} \subset \MC{G}^{\varphi}$ over $S$ and a similar inclusion over the universal cover. But this is nothing else than a $U_v$-orbit $\un{\MC{U}^{\varphi_{\MC{U}}}} \subset \un{\MC{G}^{\varphi}}$ (when restricted to $\tilde{\ov{s}}$). As this subset was already defined over $S$, it is indeed $\pi_1(S, \ov{s})$-invariant. Hence it gives the desired adelic $U_v$-level structure. \\
It is clear that these constructions are mutually inverse. \\
b) Assume now that $\MB{U}_v$ is a normal subgroup. Start with an adelic $U_v$-level structure $\psi$ on $(\MC{G}, \varphi)$. By a) this corresponds to an inclusion $\psi^\sharp: (\MC{U}, \varphi_{\MC{U}}) \hookrightarrow (\MC{G}, \varphi)$ of an element $(\MC{U}, \varphi_{\MC{U}}) \in \MB{U}_v\op{-}\acute{E}tSht(S)$. This defines an isomorphism
\begin{align*}
 \ov{\psi^\sharp}: \ov{\MC{G}} & = \MC{G} \times^{L^+G_v} L^+G_v/\MB{U}_v \stackrel{\psi^\sharp}{\cong} \MC{U} \times^{\MB{U}_v} L^+G_v \times^{L^+G_v} L^+G_v/\MB{U}_v \\
 & = \MC{U} \times^{\MB{U}_v} L^+G_v/\MB{U}_v = (\MC{U} \times^{\MB{U}_v} \MB{1}) \times^{\MB{1}} L^+G_v/\MB{U}_v.
\end{align*}
where $\MB{1}$ denotes the trivial group (as a group scheme over $S$). But any $\MB{1}$-torsor is trivial, i.e. isomorphic to the structure sheaf $\MC{O}_S$ on $S$. Hence we have a canonical isomorphism of $\MB{1}$-torsors $\MC{U} \times^{\MB{U}_v} \MB{1} \cong \MB{1}$. This defines the desired morphism
\[\ov{\psi^\sharp}: \ov{\MC{G}} \to L^+G_v/\MB{U}_v.\]
Using that there are no non-trivial automorphisms of $\MB{1}$-torsors, we get for the Frobenius-isomorphism
\[
 \begin{xy} \xymatrix {
  \sigma^*\ov{\MC{G}} \ar@{=}[r] \ar^{\ov{\varphi}}[d]  & \sigma^*\MC{U} \times^{\MB{U}_v} \! L^+G_v \times^{L^+G_v} \! L^+G_v/\MB{U}_v \ar@{=}[r] \ar^{\varphi_U \times \id \times \id}[d] & (\sigma^*\MC{U} \times^{\MB{U}_v} \! \MB{1}) \times^{\MB{1}} \! L^+G_v/\MB{U}_v \ar@{=}[r] \ar^{(\varphi_U \times \id) \times \id}[d] & \sigma^*\MB{1} \times^{\MB{1}} \! L^+G_v/\MB{U}_v \ar^{\sigma^* \times \id}[d] \\
  \ov{\MC{G}} \ar@{=}[r]  & \MC{U} \times^{\MB{U}_v} \! L^+G_v \times^{L^+G_v} \! L^+G_v/\MB{U}_v \ar@{=}[r] & (\MC{U} \times^{\MB{U}_v} \! \MB{1}) \times^{\MB{1}} \! L^+G_v/\MB{U}_v \ar@{=}[r] & \MB{1} \times^{\MB{1}} \! L^+G_v/\MB{U}_v
 } \end{xy} 
\]
Hence we get indeed an isomorphism
\[\ov{\psi^\sharp}: (\ov{\MC{G}}, \ov{\varphi}) \xrightarrow{\sim} (L^+G_v/\MB{U}_v, \sigma^*)\]
Conversely assume we have a local $G$-shtuka $(\MC{G}, \varphi)$ and an isomorphism $\ov{\psi^\sharp}: (\ov{\MC{G}}, \ov{\varphi}) \xrightarrow{\sim} (L^+G_v/\MB{U}_v, \sigma^*)$. Then consider the sequence of local systems on $S$
\[\MC{G}^\varphi \to \MC{G}^\varphi \times^{L^+G_v(\Fq)} L^+G_v(\Fq)/U_v \cong (\MC{G} \times^{L^+G_v} L^+G_v/\MB{U}_v)^\varphi = \ov{\MC{G}}^{\ov{\varphi}} \xrightarrow{\ov{\psi^\sharp}} (L^+G_v/\MB{U}_v)^{\sigma^*} = L^+G_v(\Fq)/U_v\]
Now the preimage of $1 \in L^+G_v(\Fq)/U_v$ defines a $U_v$-subtorsor of $\MC{G}^\varphi$. Hence we get a $\pi_1(S, \ov{s})$-invariant $U_v$-subtorsor of $\MC{G}^\varphi_{\tilde{S}}$ over the universal cover, which defines the desired adelic $U_v$-level structure inside $\un{\MC{G}^\varphi}$. \\
It is clear from the constructions that they are inverse to each other (just note that the trivialization of $\ov{\MC{G}}$ is uniquely defined by the preimage of $1 \in L^+G_v/\MB{U}_v$).
\exit
$\left. \right.$ \vspace{3mm} \\
We extend the definition of adelic level structures now to the global setting: \\
Let $D_0 \subset C$ be any finite reduced subscheme and $\M{A}_{D_0}^{int}$ be the ring of integral adeles of $C$ at all places in ${D_0}$. Let $U = \prod_{v \in {D_0}} U_v \subset G(\M{A}_{D_0}^{int})$ be an open subgroup. Consider now any global $G$-shtuka $(\MS{G}, (c_i)_i, \varphi) \in \nabla_n\MC{H}^1(C, G)(S)$ over $S$ with characteristic places $c_i$ contained in $(C \setminus D_0) \times_{\Fq} S$. Then recall for every non-characteristic place $v \in D_0$ with residue field $\kappa(v)$ the functor
\[\MF{L}_v: \nabla_n\MC{H}^1(C, G)(S) \to \acute{E}tSht_{G_v}(S)\]
associating to the global $G$-shtuka $(\MS{G}, (c_i)_i, \varphi)$ its local $G_v = \op{Res}_{\kappa(v)/\Fq}(G)$-shtuka. Furthermore recall that we fixed a local coordinate $z$ at $v$ in the definition of $\MF{L}_v$, i.e. an isomorphism $\M{A}_v^{int} \cong \kappa(v)[[z]]$. Thus we may view
\[U_v \subset G(\M{A}_v^{int}) \cong G(\kappa(v)[[z]]) = G_v(\Fq[[z]]) = L^+G_v(\Fq).\]
This allows us to make the

\defi{\label{def:AdelicGlobal}}{}{
a) Let $(\MS{G}, (c_i)_i, \varphi) \in \nabla_n\MC{H}^1(C, G)(S)$ be a global $G$-shtuka with characteristic places contained in $(C \setminus D_0) \times_{\Fq} S$ and $U = \prod_{v \in D_0} U_v \subset G(\M{A}_{D_0}^{int})$ an open subgroup. Then an adelic $U$-level structure consists of an adelic $U_v$-level structure of the local $G_v$-shtuka $\MF{L}_v(\MS{G}, \varphi)$ for every place $v \in D_0$. \\
b) The stack (over the category DM-stacks over $S$) of global $G$-shtukas with $U$-level structure is denoted by $\nabla_n\MC{H}^1_U(C, G)$. It has a canonical morphism to $(C \setminus D_0)^n \setminus \Delta$ by forgetting everything except the characteristic places. \\
c) A quasi-isogeny between global $G$-shtukas with $U$-level structure $(\MS{G}, (c_i)_i, \varphi, \psi)$ and $(\MS{G}', (c_i')_i, \varphi', \psi')$ is a quasi-isogeny $(\alpha, D'): (\MS{G}, (c_i)_i, \varphi) \to (\MS{G}', (c_i')_i, \varphi')$ such that (for some representative of $\alpha$) $D'$ and $D_0$ are disjoint and for every place $v \in D_0$ the associated local morphism $\MF{L}_v(\alpha): \MF{L}_v(\MS{G}, \varphi) \to \MF{L}_v(\MS{G}', \varphi')$ is an isomorphism of \'etale local $G_v$-shtukas with $U_v$-level structure (as in definition \ref{def:AdelicLocal}).
}

\rem{}{
i) As usual, we leave it to the interested reader to check that a change of the local coordinates induces a canonical isomorphism between the sets of level structures. \\
ii) As in the local case the Frobenius-isomorphism $\varphi$ in a global $G$-shtuka with $U$-level structure $(\MS{G}, (c_i)_i, \varphi, \psi)$ defines a quasi-isogeny $\varphi: (\sigma^*\MS{G}, \sigma^*\varphi, \sigma^*\psi) \to (\MS{G}, \varphi, \psi)$ respecting the $U$-level structure. \\
iii) For global $G$-shtukas in $\nabla_{\hat{c}_i}\MC{H}^1(C, G)(S)$, the same definition gives $U$-level structures for any open subgroup $U \subset G(\M{A}^{int (c_i)_i})$
}

\prop{\label{prop:AdelicCompareNaive}}{}{
Let $D \subset C$ be a finite (not necessarily reduced) subscheme contained in the formal completion of $C$ along $D_0$. Consider the open subgroup $U = \ker\left(G(\M{A}_{D_0}^{int}) \to G(\MC{O}_D)\right)$. Let $(\MS{G}, (c_i)_i, \varphi) \in \nabla_n\MC{H}^1(C, G)(S)$ be a global $G$-shtuka satisfying the usual condition on characteristic places. Then there is a bijection
\[
 \begin{Bmatrix} \textnormal{adelic} \; U\textnormal{-level} \\ \textnormal{structures} \; \psi \\ \textnormal{on} \; (\MS{G}, \varphi) \end{Bmatrix} \xleftrightarrow{\;1:1\;} \begin{Bmatrix} G\textnormal{-equivariant isomorphisms} \\ \psi^\sharp: \MS{G}|_{D \times S} \to G_{D \times S} \, \textnormal{with} \\  \psi^\sharp  \circ \varphi|_{D \times S} = \sigma^* \circ \sigma^*\psi^\sharp: \sigma^*\MS{G}|_{D \times S} \to G_{D \times S} \end{Bmatrix}
\]
where on the right-hand side $\sigma: D \times_{\Fq} S \to D \times_{\Fq} S$ is the restriction of the Frobenius on $C \times_{\Fq} S$, i.e. it is the identity on $D$ and the absolute $q$-Frobenius on $S$.
}

\prooof
Wlog. assume that $D$ is concentrated in one place $v$, i.e. after fixing a local coordinate $z$ at $v$, we have $D = \Sp \kappa(v)[[z]]/(z^{n+1}) \hookrightarrow \Spf \M{A}_v^{int}$ in the completion of $C$ at $v$. Then $U = \prod_{v' \in D} U_{v'}$ with $U_{v'} = G(\M{A}_{v'}^{int})$ for $v' \neq v$ and $U_v = K_n = \{g \in L^+G_v(\Fq) \,|\, g \equiv 1 \bmod z^n\} \subset L^+G_v(\Fq) = G(\M{A}_v^{int})$. Hence both sides of the bijection are trivial outside the place $v$. \\
But in the place $v$, we have due to proposition \ref{prop:AdelicCompareLocal}b) a bijection
\[
 \begin{Bmatrix} \textnormal{adelic} \; U_v\textnormal{-level structures} \; \psi \\ \textnormal{on} \; \MF{L}_v(\MS{G}, \varphi) \end{Bmatrix} \xleftrightarrow{\;1:1\;}  \begin{Bmatrix} \textnormal{isomorphisms in} \, L^+G_v/K_n\op{-}\acute{E}tSht(S) \\ \ov{\psi^\sharp}: (\ov{\MF{L}_v\MS{G}}, \ov{\MF{L}_v\varphi}) \cong (L^+G_v/K_n, \sigma^*) \end{Bmatrix}
\]
By proposition \ref{prop:LoGloCompare}a) we have a canonical equivalence
\[
 \begin{Bmatrix} G\textnormal{-torsors over} \, D \times S \end{Bmatrix} \xleftrightarrow{\;1:1\;}  \begin{Bmatrix} L^+G_v/K_n\textnormal{-torsors over} \, S \end{Bmatrix}
\]
compatible with Frobenius-linear isomorphisms and the functor $\MF{L}_v$. Hence there is a canonical bijection between the set of isomorphisms $(\ov{\MF{L}_v\MS{G}}, \ov{\MF{L}_v\varphi}) \cong (L^+G_v/K_n, \sigma^*)$ and the set of isomorphisms $\psi^\sharp: \MS{G}|_{D \times S} \to G_{D \times S}$ with $\psi^\sharp  \circ \varphi|_{D \times S} = \sigma^* \circ \sigma^*\psi^\sharp$. This proves the proposition. \exit

\rem{\label{rem:AdelicNaive}}{
i) Such trivializations of the global $G$-shtuka $(\MS{G}, (c_i)_i, \varphi)$ over $D \times_{\Fq} S$ were used in \cite{Varsh} as level structures. In the following we call them ``naive level structures'' to distinguish them from the adelic ones defined above. \\
ii) If $U = \ker\left(G(\M{A}_{D_0}^{int}) \to G(\MC{O}_D)\right)$ as in the proposition, we will write $\nabla_n\MC{H}^1_D(C, G) \coloneqq \nabla_n\MC{H}^1_U(C, G)$. 
}
%
The next aim is to compare adelic $U$-level structures to integral $U$-level structures as defined by Hartl and Rad. We refer to \cite[section 3]{HarRad1} and \cite[section 6]{HarRad2} for a more detailed description of integral $U$-level structures. 

\lem{\label{lem:AdelicRad}}{}{
Let $v \in C$ be a place with residue field $\kappa(v)$ and let $G_v = Res_{\kappa(v)/\Fq}(G)$ as above. \\
a) There is a canonical equivalence of categories
\[Funct^\otimes_0(\Rep_{\kappa(v)[[z]]}G, FMod_{\kappa(v)[[z]][\pi_1(S, \ov{s})]}) \cong \Rep(\pi_1(S), L^+G_v(\Fq)) \]
where 
\begin{itemize}
 \item $\Rep_{\kappa(v)[[z]]}G$ denotes the category of representations of the group scheme $G$ in free $\kappa(v)[[z]]$-modules of finite rank, i.e. the category of morphisms $G \to GL_d$ over $\Spf \kappa(v)[[z]]$ (for variable $d$);
 \item $FMod_{\kappa(v)[[z]][\pi_1(S, \ov{s})]}$ is the category of $\pi_1(S, \ov{s})$-representations on free $\kappa(v)[[z]]$-modules of finite rank;
 \item $Funct^\otimes_0(\Rep_{\kappa(v)[[z]]}G, FMod_{\kappa(v)[[z]][\pi_1(S, \ov{s})]})$ denotes the category of tensor functors $F$ between them, which admit an isomorphism (of tensor functors) 
\[Forget \circ F \cong \omega^0: Rep_{\kappa(v)[[z]]}G \to FMod_{\kappa(v)[[z]]},\]
where $Forget$ denotes the functor forgetting the $\pi_1(S, \ov{s})$-action and $\omega^0$ denotes the usual fiber functor forgetting the $G$-action.
\end{itemize}
b) Let $F \in Funct^\otimes_0(Rep_{\kappa(v)[[z]]}G, FMod_{\kappa(v)[[z]][\pi_1(S, \ov{s})]})$ and $\rho: \pi_1(S, \ov{s}) \to \Aut(\un{L^+G_v})$ the associated element in $\Rep(\pi_1(S), L^+G_v(\Fq))$. There is a canonical bijection
\[\Isom^\otimes(Forget \circ F, \omega^0) \cong \un{L^+G_v}\]
Furthermore this bijection is compatible with the $\pi_1(S, \ov{s})$-action. 
}

\rem{}{
Note that contrary to \cite{HarRad1} we use group schemes over $\Spf \kappa(v)[[z]]$ instead of group schemes over $\Sp \kappa(v)[[z]]$, as local $G$-shtuka give  naturally tensor functors defined on representations over this formal scheme.
}

\prooof
a) Let $F: \Rep_{\kappa(v)[[z]]}G \to FMod_{\kappa(v)[[z]][\pi_1(S, \ov{s})]}$ be a tensor functor in the category on the left-hand side. Then $Forget \circ F: \Rep_{\kappa(v)[[z]]}G \to FMod_{\kappa(v)[[z]]}$ defines a $G$-torsor $\MC{G}$ over $\Spf \kappa(v)[[z]]$. Indeed for every $n \geq 0$ we have an induced tensor functor $(Forget \circ F)[n]: \Rep_{\kappa(v)[[z]]/(z^{n+1})}G \to FMod_{\kappa(v)[[z]]/(z^{n+1})}$. This is equivalent to a tensor functor from the category of $\kappa(v)$-representations of the linear algebraic group $\Res_{\left. ^{\kappa(v)[[z]]}\!/\!_{(z^{n+1})} \right/ \kappa(v)}(G)$ into $\kappa(v)$-vector spaces.
By classical Tannakian formalism, $(Forget \circ F)[n]$ defines a $\Res_{\left. ^{\kappa(v)[[z]]}\!/\!_{(z^{n+1})} \right/ \kappa(v)}(G)$-torsor over $\kappa(v)$ or equivalently a $G$-torsor over $\Sp \kappa(v)[[z]]/(z^{n+1})$. Thus passing to the limit gives the desired $G$-torsor over $\Spf \kappa(v)[[z]]$. The condition $Forget \circ F \cong \omega^0$ implies now, that this $G$-torsor $\MC{G}$ is trivial. Furthermore the $\pi_1(S, \ov{s})$-action on the image of $F$ now gives an action on $\MC{G}$. Hence after the canonical equivalences
\begin{align*}
 \begin{Bmatrix} \pi_1(S, \ov{s})\textnormal{-actions on trivial} \\ G\textnormal{-torsors over} \; \Spf \kappa(v)[[z]] \end{Bmatrix} & \xleftrightarrow{\;1:1\;}  \begin{Bmatrix} \pi_1(S, \ov{s})\textnormal{-actions on trivial} \\ L^+G_v\textnormal{-torsors over} \; \Sp \Fq \end{Bmatrix}   \\
 & \xleftrightarrow{\;1:1\;}  \begin{Bmatrix} \pi_1(S, \ov{s})\textnormal{-actions on} \\ \textnormal{set-theoretic} \; L^+G_v(\Fq)\textnormal{-torsors} \end{Bmatrix}
\end{align*}
this defines the desired element in $\Rep(\pi_1(S), L^+G_v(\Fq))$. \\
Conversely start with an element $\rho: \pi_1(S, \ov{s}) \to \Aut(\un{L^+G_v})$ in $\Rep(\pi_1(S), L^+G_v(\Fq))$. This is equivalent to a $\pi_1(S, \ov{s})$-action on the $G$-torsor $\MS{G}$ over $\Spf \kappa(v)[[z]]$ defined by $\un{L^+G_v}$. If $\rho: G \to GL(V)$ is now any representation in $\Rep_{\kappa(v)[[z]]}G$, then $\MS{G} \times^G GL(V)$ together with the induced $\pi_1(S, \ov{s})$-action is an element in $FMod_{\kappa(v)[[z]][\pi_1(S, \ov{s})]}$. This defines the desired tensor functor (and one easily checks that it indeed lies in $Funct^\otimes_0(\Rep_{\kappa(v)[[z]]}G, FMod_{\kappa(v)[[z]][\pi_1(S, \ov{s})]})$). \\
b) In a) we constructed an equivalence between such functors $F$ and $L^+G_v$-torsors $\un{L^+G_v}$. Furthermore an isomorphism between tensor functors from $\Rep_{\kappa(v)[[z]]}G$ to $FMod_{\kappa(v)[[z]]}$ induces an isomorphism of the respective torsors. Hence we have a canonical bijection
\[\Isom^\otimes(Forget \circ F, \omega^0) \cong \Isom(\un{L^+G_v}, L^+G_v(\Fq))\]
where the $z$-adic group $L^+G_v(\Fq)$ is seen as the trivial torsor. But the isomorphisms on the right-hand side are uniquely determined by the preimage of $1 \in L^+G_v(\Fq)$, which is the desired element of the set-theoretic torsor $\un{L^+G_v}$. \\
By definition each element $\gamma \in \pi_1(S, \ov{s})$ defines an automorphism of the functor $Forget \circ F$, hence an automorphism $\un{\gamma}$ of $\un{L^+G_v}$. As $\gamma$ acts on $\Isom^\otimes(Forget \circ F, \omega^0)$ by precomposition of the automorphism defined by $\gamma^{-1}$, it also acts on $\Isom(\un{L^+G_v}, L^+G_v)$ by precomposition with $\un{\gamma}^{-1}$. Thus $\gamma$ acts on the preimages of $1$ simply by applying $\un{\gamma}$. \exit
$\left. \right.$ \vspace{3mm} \\
Recall now that Hartl and Rad construct in \cite[definition 3.5]{HarRad1} a functor
\[\check{\MC{T}}: \acute{E}tSht_{G_v}(S) \to Funct^\otimes_0(\Rep_{\kappa(v)[[z]]}G, FMod_{\kappa(v)[[z]][\pi_1(S, \ov{s})]}).\]
This allows one to define an (a priori) different adelic version of level structures, as done in \cite[definition 6.3a)]{HarRad2}, which we recall for the reader's convenience.

\defi{}{}{
a) Let $(\MC{G}_v, \varphi_v) \in \acute{E}tSht_{G_v}(S)$ be an \'etale local $G_v$-shtuka and $U_v \subset G(\M{A}_v^{int}) \cong L^+G_v(\Fq)$ be an open subgroup. Then an integral $U_v$-level structure for $(\MC{G}_v, \varphi_v)$ consists of a $\pi_1(S, \ov{s})$-invariant $U_v$-orbit in $\Isom^\otimes(Forget \circ \check{\MC{T}}(\MC{G}_v, \varphi_v), \omega^0)$. \\
b) Let $(\MS{G}, (c_i)_i, \varphi) \in \nabla_n\MC{H}^1(C, G)(S)$ be a global $G$-shtuka whose characteristic places lie in $(C \setminus D_0) \times_{\Fq} S$ and let $U = \prod_{v \in D_0} U_v \in G(\M{A}_{D_0}^{int})$ an open subgroup. Then an integral $U$-level structure consists of an integral $U_v$-level structure of the local $G_v$-shtuka $\MF{L}_v(\MS{G}, \varphi)$ for every place $v \in D_0$.
}

\rem{}{
i) Note that in \cite[definition 6.3]{HarRad2} Hartl and Rad define only rational $U$-level structures. The only difference to the definition above is, that they consider tensor functors into the category $FMod_{\kappa(v)((z))[\pi_1(S, \ov{s})]}$ by taking the tensor functors above and inverting $z$. Then the set $\Isom^\otimes(Forget \circ F, \omega^0)$ comes naturally with a $G(\M{A}_v) = LG_v(\Fq)$-action (where $\M{A}_v$ denotes the ring of adeles at $v$) and one may allow any compact open $U \in G(\M{A}_{D_0})$.
We sketch the slight differences between integral and rational $U$-level structures at the end of this section. \\
ii) Furthermore note that in \cite[section 6]{HarRad2} a global version of the functor $\check{\MC{T}}$ is used to define the level structures. Nevertheless Hartl and Rad remark that the global version of $\check{\MC{T}}$ is just the product of all composites $\check{\MC{T}} \circ \MF{L}_v$. This way one may easily see that the definition above is just a reformulation of the actual definition in \cite{HarRad2}. 
}

\prop{\label{prop:AdelicCompareRad}}{}{
a) Let $(\MC{G}_v, \varphi_v) \in \acute{E}tSht_{G_v}(S)$ be an \'etale local $G_v$-shtuka. Then the elements 
\[\check{\MC{T}}(\MC{G}_v, \varphi_v) \in Funct^\otimes_0(\Rep_{\kappa(v)[[z]]}G, FMod_{\kappa(v)[[z]][\pi_1(S, \ov{s})]})\] 
and 
\[T(\MC{G}_v, \varphi_v) \in \Rep(\pi_1(S), L^+G_v(\Fq))\]
coincide via the equivalence defined in lemma \ref{lem:AdelicRad}a). In other words, there is a commutative diagram of categories
\[
 \begin{xy} \xymatrix @C=4pc {
  \acute{E}tSht_{G_v}(S) \ar^-{\check{\MC{T}}}[r] \ar@{=}[d] & Funct^\otimes_0(Rep_{\kappa(v)[[z]]}G, FMod_{\kappa(v)[[z]][\pi_1(S, \ov{s})]}) \ar^{\cong}[d] \\
  \acute{E}tSht_{G_v}(S) \ar^-{T}[r]  & \Rep(\pi_1(S), L^+G_v(\Fq))
 } \end{xy} 
\]
b) Let $(\MS{G}, (c_i)_i, \varphi) \in \nabla_n\MC{H}^1(C, G)(S)$ be a global $G$-shtuka with characteristic places in $(C \setminus D_0) \times_{\Fq} S$ and $U \subset G(\M{A}_{D_0}^{int})$ an open subgroup. Then there is a canonical bijection
\[
  \begin{Bmatrix} \textnormal{integral} \; U\textnormal{-level} \\ \textnormal{structures on} \; (\MS{G}, \varphi) \end{Bmatrix} \xleftrightarrow{\;1:1\;} \begin{Bmatrix} \textnormal{adelic} \; U\textnormal{-level} \\ \textnormal{ structures on} \; (\MS{G}, \varphi) \end{Bmatrix}
\]
}

\prooof
a) This is clear, because both functors $\check{\MC{T}}$ and $T$ encode the $\pi_1(S, \ov{s})$-action on $\MC{G}_v^{\varphi_v}|_{\tilde{\ov{s}}}$ and the equivalence constructed in lemma \ref{lem:AdelicRad}a) gives precisely the translation between the two ways to encode this information. \\
b) It suffices to show that we have such a bijection for level structures on \'etale local $G_v$-shtukas $(\MC{G}_v, \varphi_v)$. Then part a) of this proposition and lemma \ref{lem:AdelicRad}b) imply that there is a $\pi_1(S, \ov{s})$-equivariant isomorphism
\[\Isom^\otimes(Forget \circ \check{\MC{T}}(\MC{G}_v, \varphi_v), \omega^0) \cong \un{\MC{G}_v^{\varphi_v}}.\]
Hence the set of $\pi_1(S, \ov{s})$-invariant $U_v$-orbits are the same, which on the left-hand side encode integral $U_v$-level structures and on the right-hand side adelic $U_v$-level structures. \exit \vspace{4mm} \\
Finally we sketch the comparison between integral and rational $U$-level structures. First note that any compact open subgroup $U \in G(\M{A}_{D_0})$ can be conjugated into $G(\M{A}_{D_0}^{int})$, so we will always assume that $U$ lies already in this smaller group. For simplicity we restrict ourselves to the case of only one place with non-trivial level structure, i.e. we assume $D_0 = \{v\}$ and $U = U_v$. The general case then boils down to do the following observations at each place separately. \\
Now we claim to have a diagram, where $\Rep_{\kappa(v)((z))}G$ denotes representations over $\Sp \kappa(v)((z))$,
\[
 \begin{xy} \xymatrix  {
  Funct^\otimes_0(\Rep_{\kappa(v)[[z]]}G, FMod_{\kappa(v)[[z]][\pi_1(S, \ov{s})]}) \ar^-{\sim}[r] \ar[d] & \Rep(\pi_1(S), L^+G_v(\Fq)) \ar[dd] \\
  Funct^\otimes_0(\Rep_{\kappa(v)[[z]]}G, FMod_{\kappa(v)((z))[\pi_1(S, \ov{s})]}) \ar^-{\sim}[d] & \\
  Funct^\otimes_0(\Rep_{\kappa(v)((z))}G, FMod_{\kappa(v)((z))[\pi_1(S, \ov{s})]}) \ar^-{\sim}[r] & \Rep(\pi_1(S), LG_v(\Fq))
 } \end{xy} 
\]
Let us explain the various maps: The upper horizontal map comes from lemma \ref{lem:AdelicRad}a) and the lower horizontal map is shown in a similar way (though the construction of the associated $G$-torsor can now be done directly by applying the Tannakian formalism over the field $\kappa(v)((z))$). The first map on the left-hand side is just postcomposition with the canonical $FMod_{\kappa(v)[[z]][\pi_1(S, \ov{s})]} \to FMod_{\kappa(v)((z))[\pi_1(S, \ov{s})]}$ obtained by inverting $z$. The construction of the second map is similar to \cite[II.3.10-II.3.11]{DMOS} (and having a constant group scheme $G$ is essential for this to be an equivalence). Finally the morphism on the right-hand side is given by taking the induced action on $\un{L^+G_v} \times^{L^+G_v} LG_v \cong \un{LG_v}$. \\
Consider now an element $F_{\M{Q}} \in Funct^\otimes_0(Rep_{\kappa(v)((z))}G, FMod_{\kappa(v)((z))[\pi_1(S, \ov{s})]})$ and let $\rho_{\M{Q}}: \pi_1(S, \ov{s}) \to \Aut(\un{LG_v})$ the associated element in $\Rep(\pi_1(S), LG_v(\Fq))$ (similarly to lemma \ref{lem:AdelicRad}). Then there is again a canonical bijection
\[\Isom^\otimes(Forget \circ F_{\M{Q}}, \omega_Q^0) \cong \un{LG_v}\]
compatible with the $\pi_1(S, \ov{s})$-action. Here $\omega_{\M{Q}}^0$ denotes the usual fiber functor over $\kappa(v)((z))$. In particular a rational $U_v$-level structure, which is by \cite[definition 6.3]{HarRad2} a $\pi_1(S, \ov{s})$-invariant $U_v$-orbit in $\Isom^\otimes(Forget \circ F_Q, \omega_Q^0)$, is the same as a $\pi_1(S, \ov{s})$-invariant $U_v$-orbit in $\un{LG_v}$. Actually \cite{HarRad2} uses tensor functors $F_{\M{Q}} \in Funct^\otimes_0(\Rep_{\kappa(v)[[z]]}G, FMod_{\kappa(v)((z))[\pi_1(S, \ov{s})]})$ but this makes no difference by the equivalence displayed in the diagram above.\\
To compare integral and rational $U_v$-level structures, start with a global $G$-shtuka giving a tensor functor $F \in Funct^\otimes_0(\Rep_{\kappa(v)[[z]]}G, FMod_{\kappa(v)[[z]][\pi_1(S, \ov{s})]})$ corresponding to some $\rho: \pi_1(S, \ov{s}) \to \Aut(\un{LG_v})$. We denote their rational analogues again by $F_{\M{Q}}$ and $\rho_{\M{Q}}$. Then we have
\[
 \begin{xy} \xymatrix  {
  \Isom^\otimes(Forget \circ F, \omega^0) \ar^-{\sim}[r] \ar[d] & \un{L^+G_v} \ar[d] \\
  \Isom^\otimes(Forget \circ F_{\M{Q}}, \omega_{\M{Q}}^0) \ar^-{\sim}[r]  & \un{LG_v}
 } \end{xy} 
\]
where the vertical morphism on the left-hand side is given by base-changing everything to $\kappa(v)((z))$ and the one on the right-hand side is again given by $\un{L^+G_v} \subset \un{L^+G_v} \times^{L^+G_v} LG_v \cong \un{LG_v}$. Using the translation of integral and rational $U_v$-level structures to the choice of subtorsors of $\un{L^+G_v}$ respectively $\un{LG_v}$, this shows that the space of rational $U_v$-level structures is strictly bigger than the space of integral $U_v$-level structures. Note that together with proposition \ref{prop:AdelicCompareNaive}, this reproves the strict inclusion shown in \cite[theorem 6.4]{HarRad2}. There even equality is asserted, which however contradicts the statements proven in \ref{prop:AdelicCompareRad}. \\
Nevertheless a global $G$-shtuka admits an integral $U_v$-level structure if it admits a rational one. Indeed any $\pi_1(S, \ov{s})$-invariant $U_v$-orbit in $\un{LG_v}$ can be translated by an element in $LG_v(\Fq)$ to an invariant $U_v$-orbit in $\un{L^+G_v}$ (recall for this our assumption $U_v \subset L^+G_v(\Fq)$).

\subsection{The moduli space of bounded global \textit{G}-shtukas with level structure}\label{subsec:GlobalModuli}
We will finally define the notion of a globally $\F$-bounded global $G$-shtuka with $U$-level structure. It is shown that the moduli space of such bounded global $G$-shtukas is a DM-stack which can be covered by quotients of quasi-projective schemes by finite groups. In particular we will show that any quasi-compact open substack is representable by a quasi-projective scheme after enlarging the level structure (cf. proposition \ref{prop:RepresentGlobalShtukaLocally} for the precise statement). For this proof we follow the ideas of Varshavsky \cite{Varsh}. 

\defi{}{}{
Fix a finite reduced subscheme $D_0 \subset C$ and an open subgroup $U \subset G(\M{A}_{D_0}^{int})$. Assume wlog. that $U_v \neq G(\M{A}_v^{int})$ for any $v \in D_0$. A globally $\F$-bounded global $G$-shtuka with $U$-level structure is a tuple $(\MS{G}, (c_i), \varphi, \psi)$ consisting of a globally $\F$-bounded global $G$-shtuka $(\MS{G}, (c_i), \varphi)$ (as in definition \ref{def:GlobalBoundVarsh}) and a $U$-level structure on it (as in definition \ref{def:AdelicGlobal}) such that $D_0$ is disjoint from the characteristic places. \\
We denote by $\nabla_n^{\mmu}\MC{H}^1_U(C, G)$ the stack of global $G$-shtuka with $U$-level structure which are globally $\F$-bounded by some tuple $\mmu$. It has a canonical morphism to $(C \setminus D_0)^n \setminus \Delta$, mapping as usual a global $G$-shtuka to its characteristic places. 
}

\rem{}{
In short we have
\[\nabla_n^{\mmu}\MC{H}^1_U(C, G) = \nabla_n^{\mmu}\MC{H}^1(C, G) \times_{\nabla_n\MC{H}^1(C, G)} \nabla_n\MC{H}^1_U(C, G)\]
}

\prop{\label{Prop:BoundGlobHecke}}{c.f. \cite[lemma 3.1b)]{Varsh}}{
Fix a finite field extension $\F/\Fq$ and an $n$-tuple  $\mmu$ of dominant $\Gamma$-invariant cocharacters of $\Res_{\F/\Fq}(G)$. Consider the stack over the category of $E$-schemes
\[Hecke^{\mmu}(S) = \left\{(\MS{G}', \MS{G}, (c_i)_i, \varphi)  \left| \begin{array}{c} 
  \MS{G}', \MS{G} \in \MC{H}^1(C, G)(S) \;,\; (c_i)_i \in (C^n \setminus \Delta)(S) \\
 \varphi: \MS{G}'|_{(C \times_{\Fq} S) \setminus \bigcup\nolimits_i \Gamma_{c_i}} \to \MS{G}|_{(C \times_{\Fq} S) \setminus \bigcup\nolimits_i \Gamma_{c_i}} \, \textnormal{isomorphism} \\
  \textnormal{such that} \; \varphi \; \textnormal{is globally} \; \M{F}\textnormal{-bounded by} \, \mmu
\end{array} \right.\right\}
\]
Then the forgetful morphism of stacks
\[Hecke^{\mmu} \to \MC{H}^1(C, G) \times_{\Fq} (C^n \setminus \Delta) \qquad,\qquad (\MS{G}', \MS{G}, (c_i)_i, \varphi) \mapsto (\MS{G},  (c_i)_i)\]
is representable by a projective scheme.
}

\prooof
We follow the proof in \cite{Varsh}, which in turn was inspired by \cite{Gaits}. 
Fix a faithful representation $\Res_{\F/\Fq}(G) \to GL(V)$ where $V = \bigoplus_{\lambda \in \Lambda} V(\lambda)$ is a finite direct sum of Weyl modules of $\Res_{\F/\Fq}(G)$. We may choose $\Lambda$ to be invariant under $Gal(\ACFq/\Fq)$ implying that $V$ is invariant under the absolute Galois group. In particular $V$ is already defined over $\Fq$. Moreover we view $V$ as a faithful representation of $G$ via the canonical inclusion $G \hookrightarrow \Res_{\F/\Fq}(G)$.
Set
\[N_i = \min_{\lambda \in \Lambda} -\langle (-\lambda)_{\rm dom}, \mu_i \rangle \quad {\rm and} \quad N_i' = -\min_{\lambda \in \Lambda} -\langle (-\lambda)_{\rm dom}, (-\mu_i)_{\rm dom} \rangle = \max_{\lambda \in \Lambda} \langle \lambda, \mu_i \rangle\]
for all $1 \leq i \leq n$. 
Define now two more stacks, namely
\[
Hecke'(S) = \left\{(\MS{G}', \MS{G}, (c_i)_i, \varphi) \left| \begin{array}{c} 
  \MS{G}', \MS{G} \in \MC{H}^1(C, G)(S) \;,\; (c_i)_i \in (C^n \setminus \Delta)(S) \\    
  \varphi: \MS{G}'|_{(C \times_{\Fq} S) \setminus \bigcup_i \Gamma_{c_i}} \to \MS{G}|_{(C \times_{\Fq} S) \setminus \bigcup_i \Gamma_{c_i}} \textnormal{isomorphism such that} \\ 
  \MS{G}_V \otimes_{\MC{O}_{C \times S}} \MC{O}(\sum_i N_i' \cdot \Gamma_{c_i}) \subseteq \varphi_V(\MS{G}'_V) \subseteq \MS{G}_V \otimes_{\MC{O}_{C \times S}} \MC{O}(\sum_i N_i \cdot \Gamma_{c_i}) \end{array} \right.\right\}
\]
and
\[Hecke''(S) = \left\{(\MC{E}, \MS{G}, (c_i)_i, \varphi_V) \left|\begin{array}{c} 
  \MC{E} \; \textnormal{a vector-bundle over} \; C \times_{\Fq} S \; \textnormal{of rank} \; \dim V \\
  \MS{G} \in \MC{H}^1(C, G)(S) \;,\; (c_i)_i \in (C^n \setminus \Delta)(S) \\
  \varphi_V: \MC{E}|_{(C \times_{\Fq} S) \setminus \bigcup_i \Gamma_{c_i}} \to \MS{G}_V|_{(C \times_{\Fq} S) \setminus \bigcup_i \Gamma_{c_i}} \textnormal{isomorphism such that}  \\ 
  \MS{G}_V \otimes_{\MC{O}_{C \times S}} \MC{O}(\sum_i N_i' \cdot \Gamma_{c_i}) \subseteq \varphi_V(\MC{E}) \subseteq \MS{G}_V \otimes_{\MC{O}_{C \times S}} \MC{O}(\sum_i N_i \cdot \Gamma_{c_i}) \end{array} \right.\right\}
\]
Here we make sense of the last assumptions in the definition of each stack as follows: Locally on $S$, we may extend $\varphi$ respectively $\varphi'$ to $C \times_{\Fq} S$ after tensoring the target bundle with $\MC{O}(-N_0 \sum_i Z_i)$ for some $N_0 \gg 0$ (cf. lemma \ref{lem:QIsogExtend}). Then we may ask that the image of this morphism has the desired properties. \\
Our aim is to prove the properties indicated in the diagram
\[Hecke^{\mmu} \xhookrightarrow{closed} Hecke' \xhookrightarrow{closed} Hecke'' \xrightarrow{proj.} \MC{H}^1(C, G) \times_{\Fq} (C^n \setminus \Delta)\]
\textbf{Claim 1:} $Hecke''$ is representable as a projective scheme over the stack $\MC{H}^1(C, G) \times_{\Fq} (C^n \setminus \Delta)$ \\
Consider any element in $\MC{H}^1(C, G)(S) \times (C^n \setminus \Delta)(S)$, i.e. let $S$ be an $E$-scheme together with a $G$-torsor $\MS{G}$ on $C \times_{\Fq} S$ and fix an $n$-tuple of distinct $S$-valued points $c_i$ in $C$. Define the coherent sheaf $\MC{F} \coloneqq \MS{G}_V(\sum_i N_i \cdot \Gamma_{c_i})/\MS{G}_V(\sum_i N_i' \cdot \Gamma_{c_i})$ on $C \times S$ supported in $\bigcup_i \Gamma_{c_i}$. Then the Quot-scheme $\MF{Quot}_{\MC{F}/C \times_{\Fq} S/S}$ is representable by a projective scheme over $S$ (as $C/\Fq$ is projective) by \cite[theorem 3.1]{GroFGAIV}. Identify $S \times_{\MC{H}^1(C, G) \times_{\Fq} (C^n \setminus \Delta)} Hecke''$ with $\MF{Quot}_{\MC{F}/C \times S/S}$ as follows: If $(\MC{E}, \MS{G}, \varphi_V) \in S \times_{\MC{H}^1(C, G) \times_{\Fq} (C^n \setminus \Delta)} Hecke''$, then the condition on $\varphi_V(\MC{E})$ in the definition of $Hecke''$ gives a quotient map
\[\MC{F} = \MS{G}_V({\textstyle\sum_i} N_i \cdot \Gamma_{c_i})/\MS{G}_V({\textstyle\sum_i} N_i' \cdot \Gamma_{c_i}) \twoheadrightarrow \MS{G}_V({\textstyle\sum_i} N_i \cdot \Gamma_{c_i}) / \varphi_V(\MC{E})  \]
(abbreviating $\MS{G}_V(D') \coloneqq \MS{G}_V \otimes_{\MC{O}_{C \times_{\Fq} S}} \MC{O}(D')$ for any divisor $D'$). As the sheaf on the right-hand side is flat over $S$, it defines an element in $\MF{Quot}_{\MC{F}/C \times S/S}$. As both the vector bundle $\MC{E}$ and the morphism $\varphi_V$ can be reconstructed from this quotient morphism, it follows immediately that this defines the desired identification. \\
\textbf{Claim 2:} $Hecke'$ is a substack of $Hecke''$ \\
Consider the morphism of stacks $Hecke' \to Hecke''$, $(\MS{G}', \MS{G}, \varphi) \mapsto (\MS{G}'_V, \MS{G}, \varphi_V)$. To show that this map is injective, consider the natural morphism $\MS{G}' \to \MS{G}' \times^G GL(V) \cong GL(\MS{G}'_V)$ of sheaves over $C \times_{\Fq} S$. By choice of $V$ this map is injective. Dividing out the $G$-action gives a section
\[l: G \backslash \MS{G}' \to G \backslash GL(\MS{G}'_V)\]
Note that the $G$-equivariant preimage in $GL(\MS{G}'_V)$ of the image of $l$ is canonically isomorphic to $\MS{G}'$. But $l$ is uniquely determined by its restriction to $(C \times_{\Fq} S) \setminus \bigcup_i \Gamma_{c_i}$, where it can be alternatively described as
\[G \backslash \MS{G}' \cong GL(V) \backslash GL(\MS{G}'_V) \stackrel{\varphi_V}{\longrightarrow} GL(V) \backslash GL(\MS{G}_V) \cong G \backslash \MS{G} \hooklongrightarrow G \backslash GL(\MS{G}_V) \stackrel{\varphi_V^{-1}}{\longrightarrow} G \backslash GL(\MS{G}'_V) \qquad (\star)\]
Thus $l$ and hence $\MS{G}'$ are determined by $(\MS{G}'_V, \MS{G}, \varphi_V)$. $\varphi$ itself is simply the restriction of $\varphi_V: GL(\MS{G}'_V) \to GL(\MS{G}_V)$ to $\MS{G}'$ and can be reconstructed, too. Note that as $\varphi_V$ is only defined on $(C \times_{\Fq} S) \setminus \bigcup_i \Gamma_{c_i}$ its restriction automatically maps $\MS{G}'$ into $\MS{G}$ as required. \\
\textbf{Claim 3:} $Hecke'$ is closed in $Hecke''$ \\
Consider again an element in $\MC{H}^1(C, G)(S) \times (C^n \setminus \Delta)(S)$. By claim 1, $S'' = S \times_{\MC{H}^1(C, G) \times_{\Fq} (C^n \setminus \Delta)} Hecke''$ is a projective scheme over $S$. Let $\MC{E}^{univ}$ be the universal vector bundle corresponding to $\MC{E}$ over $C \times S''$ and $\varphi_V^{univ}$ the universal isomorphism of vector bundles over $(C \setminus \{c_i\}_i) \times_{\Fq} S''$. Then $\varphi_V^{univ}$ defines as in $(\star)$ a section
\[l^{univ}: GL(V) \backslash GL(\MC{E}^{univ}) \to G \backslash GL(\MC{E}^{univ})\] 
over $(C \times_{\Fq} S'') \setminus \bigcup_i \Gamma_{c_i}$. As explained during claim 2, a point in $S''$ lies in the image of $S' = S'' \times_{Hecke''} Hecke'$ if and only if the section $l$ extends to all of $C$ over it. In other words $S'$ is the largest substack of $S''$ such that $l^{univ}$ extends to a section on $C \times_{\Fq} S' \subset C \times_{\Fq} S''$. But this locus is closed, cf. e.g. \cite[appendix A.4]{Varsh}. \\
\textbf{Claim 4:} $Hecke^{\mmu}$ is a closed substack of $Hecke'$ \\
Let $\varphi$ be any morphism globally $\F$-bounded by $\mmu$. Then by choice of $N_i$ we have
\[\MS{G}'_V = \bigoplus_{\lambda \in \Lambda} \MS{G}'_\lambda \xhookrightarrow{\oplus_i \varphi_\lambda} \bigoplus_{\lambda \in \Lambda} \MS{G}_\lambda \otimes_{\MC{O}_{C \times_{\Fq} S}} \MC{O}({\textstyle\sum_i} -\langle (-\lambda)_{\rm dom}, \mu_i\rangle \cdot \Gamma_{c_i}) \subseteq \MS{G}_\lambda \otimes_{\MC{O}_{C \times_{\Fq} S}} \MC{O}({\textstyle\sum_i} N_i \cdot \Gamma_{c_i})\]
and similarly for the other inclusion. Thus $Hecke^{\mmu}$ is naturally a substack of $Hecke'$. As the locus where a quasi-isogeny is globally $\F$-bounded by $\mmu$ is closed in the base scheme by proposition \ref{prop:GlobalBoundClosed}, $Hecke^{\mmu}$ is closed in $Hecke'$. 
\exit

\rem{}{
A very similar proof can be found in \cite[propositions 3.8 and 3.10]{HarRad2}. The main difference in their setup is, that in \cite{HarRad2} one considers only the case of fixed $\F$-valued points as characteristic places. Moreover their definition of boundedness is only similar, but does not coincide, with our ``globally $\M{F}$-bounded by $\mmu$'' or the corresponding ``locally bounded by $\mmu$''.
}

We include now the level structures. As boundedness conditions do no longer appear in the argumentation, we will abbreviate ``global $G$-shtukas'' for ``global $G$-shtukas globally $\F$-bounded by $\mmu$''. Before treating the case of general bounds though, we discuss the situation of naive level structures as defined in remark \ref{rem:AdelicNaive}. We start with the following well-known

\lem{\label{Lem:RepresentTorsorLevel}}{}{
For any finite closed subscheme $D' \subset C$ denote the category (fibered over DM-stacks $S$ over $E$) of $G$-torsors $\MS{G}$ over $C \times_{\Fq} S$ together with a trivialization $\psi: \MS{G}|_{D' \times_{\Fq} S} \to G \times_{\Fq} (D' \times_{\Fq} S)$ by $\MC{H}^1_{D'}(C, G)$. \\
Let $\MC{X} \subset \MC{H}^1(C, G)$ be a quasi-compact open substack and fix a (non-empty) finite closed subscheme $D \subset C$. Then there is an integer $N \gg 0$ such that fiber product $\MC{X} \times_{\MC{H}^1(C, G)} \MC{H}^1_{N \cdot D}(C, G)$ is representable by a smooth quasi-projective scheme over $E$. 
}

\prooof
A detailed proof is given e.g. in \cite{Wang}. Proposition 5.0.9 in loc. cit. shows representability, while proposition 6.0.18 proves smoothness and theorem 5.0.14 quasi-projectiveness. The same methods can be found in \cite[appendix A.4 and proof 3.2]{Varsh}. \exit

\prop{\label{prop:RepresentGlobalShtukaLocally}}{}{
Let $\MC{X} \subset \MC{H}^1(C, G)$ be a quasi-compact open substack stable under the action of the Frobenius $\sigma$ and $D \subset C$ be any finite closed subscheme. Let $N \gg 0$ be an integer such that $\MC{X} \times_{\MC{H}^1(C, G)} \MC{H}^1_{N \cdot D}(C, G)$ is representable by a scheme. Then $\MC{X} \times_{\MC{H}^1(C, G)} \nabla_n^{\mmu}\MC{H}^1_{N \cdot D}(C, G)$ is representable by a quasi-projective scheme locally of finite type over $\Sp E \times_{\Fq} ((C \setminus D)^n \setminus \Delta)$.
}

\prooof
We follow the proof in \cite[proof 3.2]{Varsh} and use proposition \ref{prop:AdelicCompareNaive} to identify the level structures appearing in $\nabla_n^{\mmu}\MC{H}^1_{N \cdot D}(C, G)$ with trivializations of the restriction of the global $G$-shtuka to $(N \cdot D) \times_{\Fq} S$. For notational reasons we will abbreviate $\MC{X} \times \MC{Z} \coloneqq \MC{X} \times_{\MC{H}^1(C, G)} \MC{Z}$ for any stack $\MC{Z}$ over $\MC{H}^1(C, G)$. Furthermore assume wlog. $N = 1$. Consider now the stack (over the category of $E$-schemes)
\[Hecke^{\mmu}_D(S) = \left\{(\MS{G}', \MS{G}, (c_i)_i, \varphi, \psi) \left| \begin{array}{c} (\MS{G}', \MS{G}, (c_i)_i, \varphi) \in Hecke^{\mmu}(S) \; \op{with} \; c_i \in (C \setminus D)(S) \; \forall i \\
  \psi: \MS{G}|_{D \times_{\Fq} S} \to G \times_{\Fq} (D \times_{\Fq} S) \; \textnormal{trivialization} \end{array} \right.\right\} \]
which maps to $\MC{H}^1(C, G)$ by forgetting everything except $\MS{G}$. 
Then $\MC{X} \times Hecke^{\mmu}_D$ can be described as the fiber product
\[
\begin{xy}
 \xymatrix{
   \MC{X} \times Hecke^{\mmu}_D \ar[d] \ar[r] & Hecke^{\mmu} \ar[d]  \\
   \MC{X} \times \MC{H}^1_D(C, G) \times_{\Fq} ((C \setminus D)^n \setminus \Delta) \ar[r] & \MC{H}^1(C, G) \times_{\Fq} ((C \setminus D)^n \setminus \Delta)}
\end{xy} 
\]
Hence the previous proposition and our assumption on $\MC{X} \times \MC{H}^1_D(C, G)$ imply that $\MC{X} \times Hecke^{\mmu}_D$ is representable as a quasi-projective scheme over $\Sp E \times_{\Fq} ((C \setminus D)^n \setminus \Delta)$. Now consider the fiber product diagram of schemes
\[
\begin{xy}
 \xymatrix@C=0.4pc @R=0.4pc {
  &  (\MS{G}, (c_i)_i, \varphi, \psi) \ar@{|->}[rr] && (\sigma^*\MS{G}, \MS{G}, (c_i)_i, \varphi, \psi) & \\
  (\MS{G}, (c_i)_i, \varphi, \psi) \ar@{|->}[dd] & \MC{X} \times \nabla_n^{\mmu}\MC{H}^1_D(C, G) \ar[dd] \ar[rr] && \MC{X} \times Hecke^{\mmu}_D \ar[dd] & (\MS{G}', \MS{G}, (c_i)_i, \varphi, \psi)  \ar@{|->}[dd] \\
  &  & & & \\
 (\MS{G}, \psi)  &  \MC{X} \times \MC{H}^1_D(C, G) \ar[rr]^-{\Gamma_{\sigma}} && (\MC{X} \times \MC{H}^1_D(C, G)) \times (\MC{X} \times \MC{H}^1_D(C, G)) & ((\MS{G}', \psi \circ \varphi|_D), (\MS{G}, \psi)) \\
  & (\MS{G}, \psi) \ar@{|->}[rr] && ((\sigma^*\MS{G}, \sigma^*\psi), (\MS{G}, \psi)) &
 }
\end{xy} 
\]
where the $\Gamma_\sigma$ denotes the graph morphism of $\sigma$, which is closed by separatedness of source and target. Thus $\MC{X} \times \nabla_n^{\mmu}\MC{H}^1(C, G)$ is as a closed subscheme of $\MC{X} \times Hecke^{\mmu}_D$ indeed quasi-projective over $\Sp E \times_{\Fq} ((C \setminus D)^n \setminus \Delta)$. \exit

\rem{}{
A explicit cover of $\MC{H}^1(C, G)$ by quasi-compact open substacks stable under the action of the Frobenius $\sigma$, was described in \cite[2.1]{Varsh} using $B$-structures (as discussed in \cite{DrinfeldSimpson}): For any element $\eta \in X_*(T)_{dom}$ let
\[\MC{H}^1(C, G)^{\leq \eta}(S) = \left\{\MC{G} \in \MC{H}^1(C, G)(S) \left|\begin{array}{c} 
  \deg \MC{B}_\lambda \leq \langle \lambda, \eta \rangle \; \textnormal{for all} \; \lambda \in X^*(T)_{dom} \; \textnormal{and all}\\
   B\textnormal{\!-structures} \; \MC{B} \subset \ov{s}^{*}\MC{G} \; \textnormal{for geometric points} \; \ov{s}
   \end{array} \right.\right\}\]
Then $\MC{H}^1(C, G)^{\leq \eta}$ is a quasi-compact open substack of $\MC{H}^1(C, G)$ by \cite[lemma A.3]{Varsh}.
} 

\prop{\label{prop:RepresentLevelStructuresGlobal}}{}{
Let $U \subset U' \subset G(\M{A}_D^{int})$ be two open subgroups with $U$ normal in $U'$ (and as usual $D \subset C$ a finite subscheme). Consider the canonical morphism 
\[\nabla_n^{\mmu}\MC{H}^1_U(C, G) \to \nabla_n^{\mmu}\MC{H}^1_{U'}(C, G)\]
obtained by replacing the $U$-orbit (defined by the level structure) by the corresponding $U'$-orbit. Then this morphism is relatively representable in the category of schemes by a torsor under the finite group $U'/U$.
}

\prooof
First of all note, that if $(\MS{G}, \varphi, \psi) \in \nabla_n^{\mmu}\MC{H}^1_U(C, G)$, $v \in D$ a non-characteristic place, $\psi_v \subset (\MF{L}_v\MS{G})^{\MF{L}_v\varphi}$ the $\pi_1(S, \ov{s})$-invariant subset defined by the level structure at $v$ and finally $g_v \in U'_v \subset G(\M{A}_v^{int})$ any element, then the subset $g_v \cdot \psi_v \subset (\MF{L}_v\MS{G})^{\MF{L}_v\varphi}$ is again a $\pi_1(S, \ov{s})$-stable $U$-orbit. This follows by normality of $U$ in $U'$ and the fact that $\pi_1(S, \ov{s})$ acts through elements in $U$. In particular we have an action of $U'$ on the set of $U$-level structures on $(\MS{G}, \varphi)$: An element $g = (g_v) \in U \subset G(\M{A}_D^{int}) = \prod_{v \in D} G(\M{A}_v^{int})$ maps a $U$-level structure $\psi = (\psi_v)_{v \in D}$ to $g \cdot \psi \coloneqq (g_v \cdot \psi_v)_{v \in D}$. \\
Let us check now the representability assertion: For this we may assume wlog. that $U$ and $U'$ differ only at one place $v \in D$. Denote their components at $v$ by $U_v \subset U_v' \subset G(\M{A}_v^{int}) = L^+G_v(\Fq)$ (abbreviating again $G_v = \Res_{\kappa(v)/\Fq}(G)$). Let $(\MS{G}, \varphi, \psi') \in \nabla_n^{\mmu}\MC{H}^1_{U'}(C, G)(S)$ be a global $G$-shtuka with $U'$-level structure over a DM-stack $S$. Then we have to consider the stack $S \times_{\nabla_n^{\mmu}\MC{H}^1_{U'}(C, G)} \nabla_n^{\mmu}\MC{H}^1_{U}(C, G)$ given on a DM-stack $S'$ over $S$ by $(\MS{G} \times_S S', \varphi \times \id_{S'}, \psi)$ with a $U$-level structures $\psi$ such that $U' \cdot \psi = \psi'$. This is obviously a fppf-sheaf, hence we may check representability after an \'etale cover. 
Choose now a positive integer $N \gg 0$ such that $U_{N \cdot D} = \ker(G(\M{A}_D^{int}) \to G(\MC{O}_{N \cdot D}))$ is contained in $U$ and thus automatically normal in it. Similarly to proposition \ref{prop:TateTrivialTorsor} in the local setting, there is an \'etale cover $\widetilde{S} \to S$ such that $(\MS{G}, \varphi)|_{N \cdot D \times_{\Fq} \widetilde{S}} \cong (G \times_{\Fq} (N \cdot D \times_{\Fq} \widetilde{S}), \sigma^*)$ after passing to this \'etale cover $\widetilde{S}$. 
Thus assume wlog. that the global $G$-shtuka together with the Frobenius-isomorphism is trivial over $N \cdot D$. But for the place $v \in D$ a trivialization of $(\MS{G}, \varphi)$ over $N \cdot D$ induces a trivialization of $\MF{L}_v\MS{G} \times^{L^+G_v} L^+G_v/U_{ND \, v}$, where $U_{ND \, v}$ is the $v$-component of $U_{N \cdot D}$ (viewed here as the corresponding subgroup scheme of $L^+G_v$). As we have a trivialization of the Frobenius-isomorphism as well, we also get a trivialization of the $L^+G_v(\Fq)/U_{Dv}$-torsor 
\[(\MF{L}_v\MS{G} \times^{L^+G_v} L^+G_v/U_{Dv})^{\MF{L}_v\varphi} = (\MF{L}_v\MS{G})^{\MF{L}_v\varphi} \times^{L^+G_v(\Fq)} L^+G_v(\Fq)/U_{Dv}\]
over $S$. In particular $\pi_1(S, \ov{s})$ acts trivially on the set-theoretic torsor 
\[\un{(\MF{L}_v\MS{G})^{\MF{L}_v\varphi}} \times^{L^+G_v(\Fq)} L^+G_v(\Fq)/U_{Dv}.\]
By lemma \ref{lem:AdelicLocalNormalQuotient} the $v$-component of the $U'$-level structure $\psi'$ corresponds to some $U'_v/U_{Dv}$-orbit in $\un{(\MF{L}_v\MS{G})^{\MF{L}_v\varphi}} \times^{L^+G_v(\Fq)} L^+G_v(\Fq)/U_{Dv}$. By triviality of the $\pi_1(S, \ov{s})$-action and again lemma \ref{lem:AdelicLocalNormalQuotient}, the $U$-level structures $\psi$ with $U' \cdot \psi = \psi'$ correspond bijectively to $U_v/U_{Dv}$-orbits in the $U'_v/U_{Dv}$-orbit defined by $\psi'$. But this is obviously representable by the trivial $U_v'/U_v \cong U'/U$-torsor over $\widetilde{S}$. \exit

\rem{}{
It is easy to see directly that we have an isomorphism of stacks
\begin{align*}
 \nabla_n^{\mmu}\MC{H}^1_U(C, G) \times_E U'/U & \cong \nabla_n^{\mmu}\MC{H}^1_U(C, G) \times_{\nabla_n^{\mmu}\MC{H}^1_{U'}(C, G)} \nabla_n^{\mmu}\MC{H}^1_U(C, G) \\
 ((\MS{G}, \varphi, \psi) , g \cdot U) & \mapsto ((\MS{G}, \varphi, \psi), (\MS{G}, \varphi, g \cdot \psi))
\end{align*}
cf. \cite[proposition 6.5]{HarRad2}. However this does not help at all to prove representability. 
}

\thm{\label{thm:RepresentGlobalShtuka}}{cf. \cite[proposition 2.16a)]{Varsh}}{
a) Let $U \subset G(\M{A}_D^{int})$ be an open subgroup. Then $\nabla_n^{\mmu}\MC{H}^1_U(C, G)$ is a DM-stack locally of finite type over $E$. \\
b) Let $U \subset U' \subset G(\M{A}_D^{int})$ be any two open subgroups. Then the canonical morphism
\[\nabla_n^{\mmu}\MC{H}^1_U(C, G) \to \nabla_n^{\mmu}\MC{H}^1_{U'}(C, G)\]
is finite \'etale.
}

\prooof
Part a) follows from the previous proposition: As this is a local property it suffices to check that $\MC{X} \times_{\MC{H}^1(C, G)} \nabla_n^{\mmu}\MC{H}^1_U(C, G)$ is a DM-stack locally of finite type for each quasi-compact open substack $\MC{X} \subset \MC{H}^1(C, G)$. Consider now an integer $N \gg 0$ such that the corresponding subgroup $U_{N \cdot D} = \ker(G(\M{A}_D^{int}) \to G(\MC{O}_{N \cdot D}))$ is contained in $U$. Then after enlarging $N$ we may even assume that $\MC{X} \times_{\MC{H}^1(C, G)} \nabla_n^{\mmu}\MC{H}^1_{N \cdot D}(C, G)$ is representable by a scheme with the desired properties. But $U_{N \cdot D} = \ker(G(\M{A}_D^{int}) \to G(\MC{O}_{N \cdot D})) \subset U_D \subset U$ is a subgroup, which is normal in $G(\M{A}_D^{int})$. Hence the morphism 
\[\MC{X} \times_{\MC{H}^1(C, G)} \nabla_n^{\mmu}\MC{H}^1_{N \cdot D}(C, G) \to \MC{X} \times_{\MC{H}^1(C, G)} \nabla_n^{\mmu}\MC{H}^1_U(C, G)\]
is representable by a torsor under a finite group. This shows that $\MC{X} \times_{\MC{H}^1(C, G)} \nabla_n^{\mmu}\MC{H}^1_U(C, G)$ is indeed a DM-stack. \\
Part b) follows immediately from the previous proposition. \exit

\nota{\label{nota:ModuliSpaces}}{
For fixed places $c_i$ in $C \setminus D$, abbreviate
\[\MB{X}^{\mmu}_U \coloneqq \nabla_n^{\mmu}\MC{H}^1_U(C, G) \times_{(C \setminus D)^n} (c_i)_i \]
and
\[\MF{X}^{\mmu}_U \coloneqq \nabla_n^{\mmu}\MC{H}^1_U(C, G) \times_{(C \setminus D)^n} \widehat{(c_i)}_i \]
where $(c_i)_i \in (C \setminus D)^n$ is the point given by the $c_i$ and $\widehat{(c_i)}_i$ is its completion. \\
$\MB{X}^{\mmu}_U$ parameterizes global $G$-shtukas in $\nabla_{(c_i)}\MC{H}^1(C, G)$ together with $U$-level structures and satisfying the boundedness conditions. In the same way $\MF{X}^{\mmu}_U$ parameterizes bounded global $G$-shtukas in $\nabla_{(\hat{c}_i)}\MC{H}^1(C, G)$ together with levels structures.
}

\section{Igusa varieties on central leaves}\label{sec:Igusa}
We now analyze the special fiber $\MB{X}^{\mmu}_U$ in the moduli space of bounded global $G$-shtukas in greater detail. Central leaves, defined in section \ref{subsec:CentralLeaf} as the locus where the associated local $G$-shtukas are (point-wise) isomorphic to a suitably chosen fixed local $G$-shtuka, are of particular interest, because their universal local $G$-shtuka (associated to some fixed characteristic place) has a very easy form: The Frobenius-isomorphism lies (after a suitable trivialization) in a certain parabolic subgroup as explained in section \ref{subsec:SlopeDivisible}. Moreover the $L^+G$-torsor actually comes from a torsor for a certain subgroup $I_0(b_\nu)$ of $L^+G$, cf. section \ref{subsec:IwahoriStructure}. 
In sections \ref{subsec:DefineIgusaBasic} and \ref{subsec:DefineIgusaGeneral} we finally turn to the definition of Igusa varieties: These are moduli spaces parametrizing equivalence classes of isomorphisms, which are representable at least over basic strata. The situation over central leaves in arbitrary Newton strata is more complicated: First of all one needs a very detailed structure theory for local $G$-shtukas over these central leaves before being able to write down a useful moduli problem. Secondly one has to pass to the perfection of the central leaf, before the moduli spaces become representable. But once existence of Igusa varieties is shown, it is almost immediate that they are finite \'etale covers of the central leaf. \\
Similar theorems in the world of abelian varieties and $p$-divisible groups (with additional structures) can e.g. be found in \cite{HarrisTaylor}, \cite{MantoFoliation}, \cite{MantoFoliationPEL} and \cite{OortFoliation}, though their moduli descriptions of Igusa varieties are usually less natural than ours in the non-basic case. Many of the following definitions and proofs were inspired by ones found in these references.

\subsection{Newton strata}\label{subsec:NewtonPoint}
In this section we define a stratification of the moduli space of global $G$-shtukas via the quasi-isogeny classes of the associated local $G$-shtukas. Essentially everything here is already contained in \cite[section 7]{HaVi}, though the description of the set of $\sigma$-conjugacy classes goes back to Kottwitz \cite{KottIso1}, \cite{KottIso2}. \\
Fix for the moment an algebraically closed field $k$ over $\Fq$. Recall that any local $G$-shtuka over $k$ is of the form $(L^+G, b\sigma^*)$ for some $b \in LG(k)$ and that quasi-isogenies correspond to $\sigma$-conjugation of the element $b$. This gives a canonical bijection
\[
  \begin{Bmatrix} \textnormal{quasi-isogeny classes of local} \; G\textnormal{-shtukas over} \; k \end{Bmatrix} \xleftrightarrow{\;1:1\;} \begin{Bmatrix} \sigma\textnormal{-conjugacy classes in} \; LG(k) \end{Bmatrix}.
\]
By \cite[lemma 1.3]{RapoRich} the right-hand side does not depend on the field $k$ (as long as it is algebraically closed). Thus we may denote these sets by $\MC{B}(G)$, independently of $k$. 

\const{\label{const:LocalB(G)}}{}{
Let $S$ be some DM-stack over $E$ (or even over $\Fq$) and $(\MC{G}, \varphi) \in Sht_G(S)$ be a local $G$-shtuka over $S$. Let $\ov{s} \in S(k)$ be any geometric point (for some algebraically closed field $k$) and consider the restriction $(\MC{G}, \varphi)_{\ov{s}}$ of the local $G$-shtuka to it. This defines by the discussion above an element $[\MC{G}]_{\ov{s}} \in \MC{B}(G)$. We claim that it depends only on the image $s \in S$ of the geometric point. Indeed $\sigma$-conjugacy classes are trivially stable under the action of $\sigma$ and hence they are in particular stable under $Gal(k/E)$. 
Thus we write simply $[\MC{G}]_s$ for this element and obtain a map
\begin{align*}
 Sht_G(S) & \to Hom_{sets}(S, \MC{B}(G)) \\
 (\MC{G}, \varphi) & \mapsto (s \mapsto [\MC{G}]_s)
\end{align*}
}

\rem{}{
A more explicit description of $\MC{B}(G)$ is given in \cite[4.13]{KottIso2} or alternatively \cite[theorem 1.15iii)]{RapoRich} by the canonical embedding
\[\MC{B}(G) \hookrightarrow (X_*(T)_{\M{Q}}/W_G)^{Gal(k/\Fq)} \times \pi_1(G)_{Gal(k/\Fq)}\]
(where $W_G$ denotes the Weyl group of $G$) such that any representative $b \in LG(k)$ of a $\sigma$-conjugacy class maps to 
\begin{itemize}
 \item its Newton point in the first factor: See \cite[4.3]{KottIso1} for its definition in the general setting and \cite[section 1.3]{KatzSlope} or \cite[section 6.4]{ZinkCartier} for the case of $G = GL_n$.
 \item its Kottwitz point in the second factor as defined in \cite[lemma 6.1]{KottSVRep}, although written down using the notation of \cite[theorem 1.15i)]{RapoRich}. 
\end{itemize}
In \cite{KottIso2} and \cite{RapoRich} even the image of this morphism is described. 
}

\defi{}{}{
Let $(\MS{G}, \varphi, \psi) \in \MB{X}^{\mmu}_U(S)$ (cf. \ref{nota:ModuliSpaces}) be a bounded global $G$-shtuka over some DM-stack $S$ with the fixed characteristic places $c_i$. Then its image under the composite
\begin{align*}
 \MB{X}^{\mmu}_U(S) & \xrightarrow{(\MF{L}_{c_i})_i} \prod_i Sht_{G_{c_i}}(S) \longrightarrow \prod_i  Hom_{sets}(S, \MC{B}(G_{c_i})) = Hom_{sets}\left(S, \textstyle{\prod_i} \MC{B}(G_{c_i})\right) \\
  (\MS{G}, \varphi, \psi) \qquad & \xmapsto{\hspace*{75mm}} (s \mapsto ([\MF{L}_i\MS{G}]_s)_i)
\end{align*}
gives a canonical map $[\MS{G}]_{-}: S \to \prod_i \MC{B}(G_{c_i})$, where we assign to each point of $S$ and each associated local $G_{c_i}$-shtuka its quasi-isogeny class.
}

\prop{\label{prop:NewtStrataLocClosed}}{}{
Fix a DM-stack $S$ over $E$ again. \\
a) Let $(\MC{G}, \varphi) \in Sht_G(S)$ be a local $G$-shtuka over $S$ and $\nu \in \MC{B}(G)$. Then the subset 
\[\MC{N}^\nu_{(\MC{G}, \varphi)} = \{s \in S \;|\; [\MC{G}]_s = \nu\} \subset S\]
is locally closed. \\
b) Let $(\MS{G}, \varphi, \psi) \in \MB{X}^{\mmu}_U(S)$ be a global $G$-shtuka over $S$ and $(\nu_i)_i \in \prod_i \MC{B}(G_{c_i})$. Then the subset \[\MC{N}^{(\nu_i)}_{(\MS{G}, \varphi, \psi)} = \{s \in S \;|\; [\MS{G}]_s = (\nu_i)_i\} \subset S\]
is locally closed. \\
In both cases these subsets are called Newton strata of $S$ for the respective $G$-shtuka.
}

\prooof
The first part is an immediate consequence of the analogue of \cite[theorem 3.6ii)]{RapoRich} in equal characteristic and can be shown (after the obvious translations) in the very same way. Note at this point, that it suffices to check the assertion after an \'etale cover, allowing us to reduce to the case of schemes $S$. The second part is a direct consequence of the first and the previous definition. \exit

\rem{}{
For more properties of Newton strata see \cite[section 7]{HaVi}, where \ref{prop:NewtStrataLocClosed}a) can be found as well.
}

\subsection{Central leaves in the moduli space of global \textit{G}-shtukas}\label{subsec:CentralLeaf}
Given a global $G$-shtuka over a scheme $S$, we define the central leaf of $S$ to be the locus inside a Newton stratum where the associated local $G$-shtukas have fixed isomorphism class. This isomorphism class will be given by $P$-fundamental alcoves.

\nota{\label{nota:Iwahori}}{
For any integer $n \geq 0$ define the subgroup scheme $I_n \subset L^+G$ by
\[I_n(S) = \{g \in K_n(S) \;|\: g \bmod z^{n+1} \in B(\MC{O}_S(S)[[z]]/(z^{n+1}))\} \subset K_n(S)\]
Furthermore for any subgroup $H \subset G$ we denote
\[I_H = I_0 \cap L^+H\]
}

\defi{\label{def:FundamentalAlcove}}{}{
Let $P \subset G$ be a semi-standard parabolic subgroup, i.e. a parabolic subgroup containing $T$ but not necessarily $B$. Let $M$ be its Levi subgroup and $N$ its unipotent radical. Let $\ov{N}$ be the unipotent radical of the parabolic subgroup opposite to $P$. An element $b \in \widetilde{W}$ in the extended affine Weyl group of $G$ is called a $P$-fundamental alcove if
\[\phi_b(I_M) = I_M \; , \; \phi_b(I_N) \subseteq I_N \quad and \quad \phi_b(I_{\ov{N}}) \supseteq I_{\ov{N}}\]
where 
\[\phi_b: LG \to LG \quad , \quad g \mapsto b \cdot \sigma(g) \cdot b^{-1}\]
Here we implicitly identify elements of $\widetilde{W}$ with representatives in $LG$. 
}

\warn{}{
This definition does not coincide with the definition \cite[definition 6.1]{VieTruncLevel1} (if $G$ does not split). There $\phi_b$ is replaced by the map
\[g \mapsto \sigma(b \cdot g \cdot b^{-1})\]
Therefore an element $b$ is a $P$-fundamental alcove in our sense if and only if $\sigma^{-1}(b)$ is a $P$-fundamental alcove as in \cite{VieTruncLevel1}. Nevertheless all properties of $P$-fundamental alcoves given in \cite{VieTruncLevel1} are also valid in our situation, either using the correspondence above or the very same arguments given in \cite{VieTruncLevel1}. \\
Of course one could formulate all following results using the definition \cite[definition 6.1]{VieTruncLevel1}, but then one would need a twist by $\sigma$ in the definition of the local $G$-shtuka in \ref{nota:CentralLeaf}d), i.e. considering now $(L^+G_{c_i}, \sigma(b_{\nu_i})\sigma^*)$. This would result in (even more) cumbersome notations later on as most of sections \ref{subsec:DefineIgusaBasic} and \ref{subsec:DefineIgusaGeneral} uses this local $G$-shtuka, but none of its specific properties. Moreover with view towards the construction in the proof of proposition \ref{prop:I_0Structure}, the given definition of a $P$-fundamental alcove seems to be the more natural.
}

\rem{}{
i) Neither $P$ nor $b$ can in general be defined over $\Fq$. We call an element a fundamental alcove, if we do not want to specify the parabolic subgroup. \\
ii) By \cite[theorem 6.5]{VieTruncLevel1} (or for split groups by \cite[corollary 13.2.4]{GHKR2}) there is a $P$-fundamental alcove (for some parabolic $P$) in each $\sigma$-conjugacy class in $LG(k)$ for some algebraically closed field $k$. \\
iii) For further properties of $P$-fundamental alcoves see \cite[section 6]{VieTruncLevel1}. In the case of split groups even more can be said, cf. \cite[section 13]{GHKR2} and \cite[remark 4.2]{HaVi2}.
}

\ass{\label{ass:FundamentalAlcoveDefined}}{
From now on assume always that for each $\nu_i \in \MC{B}(G_{c_i})$ appearing in the definition of the Newton stratum in consideration, there exists a fundamental alcove $b_{\nu_i} \in \widetilde{W}$ in the $\sigma$-conjugacy class of $\nu_i$ which admits a representative in $LG_{c_i}(E)$ defined over $E$. We will usually fix such a representative and call it by abuse of notation $b_{\nu_i}$, too.
}

\nota{\label{nota:CentralLeaf}}{
Throughout this section we fix the following data: \\
a) an $n$-tuple of characteristic places $(c_i)_i$ in $C$. \\
b) an open subgroup $U \subset G(\M{A}^{int (c_i)})$ with support in $D$.  \\
c) an $n$-tuple $(\nu_i)_i \in \prod_i \MC{B}(G_{c_i})$ of $\sigma$-conjugacy classes. \\
d) for each $\nu_i$ a local $G_{c_i}$-shtuka $(L^+G_{c_i}, b_{\nu_i}\sigma^*)$ over $E$ with $b_{\nu_i}$ a fundamental alcove lying in the $\sigma$-conjugacy class $\nu_i$. \\ 
Furthermore we make the following conventions: \\
e) We abbreviate $\mathbf{X}^{\mmu}_U \coloneqq \nabla_n^{\mmu}\MC{H}^1_U(C, G)$. Note that it is by theorem \ref{thm:RepresentGlobalShtuka}b) a DM-stack locally of finite type over $\Sp E$. \\
f) $(\MS{G}^{univ}_U, \varphi^{univ}_U, \psi^{univ}_U)$ is the universal global $G$-shtuka over the stack $\mathbf{X}^{\mmu}_U$. It is given by a sheaf of global $G$-shtukas over the category $\acute{E}tSch/\mathbf{X}^{\mmu}_U$ of schemes which are \'etale over $\mathbf{X}^{\mmu}_U$ (which form a basis for the Grothendieck topology). We will denote the evaluation of this sheaf on some $S \in \acute{E}tSch/\mathbf{X}^{\mmu}_U$ by $(\MS{G}^{univ}_U, \varphi^{univ}_U, \psi^{univ}_U)|_S$. \\
g) $\MC{N}^{(\nu_i)}_U \subset \mathbf{X}^{\mmu}_U$ denote the Newton strata as defined in proposition \ref{prop:NewtStrataLocClosed}.
}

\defi{\label{def:CentralLeaf}}{}{
Let $S$ be a DM-stack over $E$ and $(\MS{G}, \varphi, \psi)$ a global $G$-shtuka over $S$. Then consider the subset
\[\MC{C}^{(L^+G_{c_i}, b_{\nu_i}\sigma^*)_i} = \begin{Bmatrix} \ov{s}: \Sp k \to S \,\textnormal{geometric point, such that for all}\, i \\ \MF{L}_{c_i}(\MS{G}, \varphi) \times_S \ov{s} \cong (L^+G_{c_i}, b_{\nu_i}\sigma^*) \times_{\Sp E} \Sp k \end{Bmatrix}\]
Then the reduced stack with geometric points $\MC{C}^{(L^+G_{c_i}, b_{\nu_i}\sigma^*)_i}$ (whose existence is shown in the next proposition) is called the central leaf in $S$ corresponding to the tuple $(b_{\nu_i})_i$. 
}

In the case of $S = \mathbf{X}^{\mmu}_U$ together with the universal global $G$-shtuka over it, we use the abbreviattion $\MC{C}_U^{(\nu_i)}$ for the corresponding central leaf. 
Note that central leaves depend on the chosen $b_{\nu_i}$, which are not unique, not even (in general) in the fixed conjugacy class.

\rem{}{
i) The isomorphism condition depends only on the image of $\ov{s}$, but not on the actual choice of the geometric point. \\
ii) As Newton strata fix the quasi-isogeny class of the local $G_{c_i}$-shtukas at the characteristic places, we have
\[\MC{C}^{(L^+G_{c_i}, b_{\nu_i}\sigma^*)_i}\subset \MC{N}^{(\nu_i)}\]
where $\MC{N}^{(\nu_i)} \subset S$ is the Newton stratum to $(\nu_i)_i \in \prod_i \MC{B}(G_{c_i})$. \\
iii) Although we stated the definition of a central leaf in great generality, we will mostly use it in the universal case. \\
iv) The condition for a global $G$-shtuka to lie in the central leaf depends only on the local behavior at the characteristic places. In particular it does not depend on the choice of the level structure. 
v) A similar definition of a central leaf for local $G_{c_i}$-shtukas was studied by Hartl and Viehmann in \cite[\S 6]{HaVi2}. As the definition for global $G$-shtukas is based on the one for local $G_{c_i}$-shtukas, we will frequently invoke their results.
}

\prop{\label{prop:CentralLeafProperties}}{}{
a) $\MC{C}_U^{(\nu_i)}$ is closed in the Newton stratum $\MC{N}^{(\nu_i)}_U \subset \mathbf{X}^{\mmu}_U$. In particular it is locally closed in $\mathbf{X}^{\mmu}_U$. 
In the following we endow $\MC{C}_U^{(\nu_i)}$ with the reduced substack structure over $E$. \\
b) $\MC{C}_U^{(\nu_i)}$ is a smooth DM-stack locally of finite type over $E$. 
}

\prooof
a) There is at most one reduced substack whose geometric points equal to the set given in \ref{def:CentralLeaf}. Hence it suffices to show this over an \'etale cover. But over schemes, the arguments given by Hartl and Viehmann in \cite[corollary 6.7]{HaVi} (together with their proofs of previous statements) translate verbatim from split groups $G_{c_i}$ to the general case. \\
%
b) As already noted above $\mathbf{X}^{\mmu}_U$ is representable by a DM-stack locally of finite type over $E$. By a) $\MC{C}_U^{(\nu_i)}$ is locally closed in $\mathbf{X}^{\mmu}_U$, hence (endowed with the reduced substack structure) a reduced DM-stack locally of finite type over $E$. Thus by generic smoothness (over perfect base schemes as in \cite[\href{http://stacks.math.columbia.edu/tag/04QM}{Tag 04QM}]{stacks-project}, currently lemma 29.15.7) $\MC{C}_U^{(\nu_i)}$ contains a smooth point (and even an open dense subset of smooth points). \\
Now take any two closed points $x_j = (\MS{G}_j, \varphi_j, \psi_j) \in (\MC{C}_U^{(\nu_i)} \times_E \Sp \ACFq)(\ACFq)$ (for $j = 1, 2$) and consider for each the completion $(\MC{C}_U^{(\nu_i)} \times_E \Sp \ACFq)^{\wedge x_j}$ of $(\MC{C}_U^{(\nu_i)} \times_E \Sp \ACFq)$ at $x_j$ as a formal DM-stack. Taking an \'etale cover by a scheme and comparing the completions, shows that $(\MC{C}_U^{(\nu_i)} \times_E \Sp \ACFq)^{\wedge x_j}$ actually exists as a honest formal scheme. By \cite[theorem 5.10]{HarRad1} we have isomorphisms of local deformation spaces (considered as formal schemes over $\ACFq$)
\[Def(\MS{G}_j, \varphi_j, \psi_j) \cong \prod_i Def(\MF{L}_{c_i}(\MS{G}_j, \varphi_j)).\]
By definition of the central leaf this gives isomorphisms
\[Def(\MC{G}_1, \varphi_1, \psi_1) \cong \prod_i Def(L^+G_{c_i}, b_{\nu_i}\sigma^*) \cong Def(\MC{G}_2, \varphi_2, \psi_2)\]
which restricts to an isomorphism of formal schemes
\[(\MC{C}_U^{(\nu_i)} \times_E \Sp \ACFq)^{\wedge x_1} \cong (\MC{C}_U^{(\nu_i)} \times_E \Sp \ACFq)^{\wedge x_2}\]
Thus $\MC{C}_U^{(\nu_i)} \times_E \Sp \ACFq$ is regular at every point if it is regular at one point, which was just shown above. In particular $\MC{C}_U^{(\nu_i)}$ is smooth. 
\exit

\rem{}{
i) For split groups $G$ and $G_{c_i}$ Newton-strata and central leaves are already defined over $\Fq$. Although the Newton strata can be defined over $\Fq$ even for non-split groups $G$, this is no longer true for central leaves. \\
ii) Similar arguments as for a) show that central leaves exist as locally closed subschemes for any family of global $G$-shtukas over an arbitrary DM-stack. 
}

\subsection{Completely slope divisible local \textit{G}-shtukas}\label{subsec:SlopeDivisible}
We recall the definition of a completely slope divisible local $G$-shtuka. Our notion (almost) coincides with the notion `completely slope divisible in the sense of Zink' as found in \cite[section 5]{HaVi2}. Furthermore we show that over a connected normal scheme a local $G$-shtuka is completely slope divisible if its generic fiber is completely slope divisible. The results of this section imply that the universal family over a central leaf is completely slope divisible, cf. theorem \ref{thm:UnivComplSlopeDiv}. \\
For a parabolic subgroup $P \subset G$, we will use the notion of local $P$-shtukas. These are essentially given by replacing the reductive group $G$ by $P$ in the definition of a local $G$-shtuka. The main use of local $P$-shtukas will be to pose natural conditions on isomorphisms between $LG$-torsors and thus on quasi-isogenies between local $G$-shtukas. This allows us to isolate the ones that behave nicely.

\defi{}{}{
Let $G$ be a connected reductive group over $\Fq$ and $P \subset G$ a parabolic subgroup defined over $E$ (but not necessarily over $\Fq$). A local $P$-shtuka over an $E$-scheme $S$ is a pair $(\MC{P}, \varphi)$ consisting of a $L^+P$-torsor $\MC{P}$ over $S$ and a $LP$-equivariant $\sigma$-linear isomorphism $\varphi: \sigma^*\MC{LP} \to \MC{LP}$. Here $LP$-equivariance means the following: $\sigma^*\MC{P}$ has an induced action of the group $\sigma^*(L^+P) = L^+(\sigma^*P)$ (which may differ from $L^+P$). The same holds for induced $LP$-torsors. Hence we may as that the diagram of ind-schemes
\[
\begin{xy}
 \xymatrix {
    L(\sigma^*P) \times_E \MC{L}\sigma^*\MC{P} \ar[r] \ar^{\sigma^* \times \varphi}[d] & \MC{L}\sigma^*\MC{P} \ar^{\varphi}[d]  \\
    LP \times_E \MC{L}\MC{P} \ar[r] & \MC{L}\MC{P}  }
\end{xy} 
\]
commutes.
}

\defi{\label{def:ComplSlopeDivi}}{\cite[definition 5.3]{HaVi2}}{
A local $G$-shtuka $(\MC{G}, \varphi)$ over an $E$-scheme $S$ is called completely slope divisible, if there exists a standard parabolic subgroup $P \subset G$ (defined over $E$) with Levi subgroup $M \subset P$ and opposite parabolic $\ov{P}$ together with a local $\ov{P}$-shtuka $(\ov{\MC{P}}, \varphi_{\ov{\MC{P}}})$ over $S$, an integer $s > 0$ and a $G$-dominant $M$-central cocharacter $\nu': \M{G}_m \to T$ (defined over $\Fq$) such that: \\
a) There is an isomorphism of local $G$-shtukas $\alpha: (\MC{G}, \varphi) \to (\ov{\MC{P}} \times^{L^+\ov{P}} L^+G, \varphi_{\ov{\MC{P}}} \times^{L\ov{P}} \id_{LG})$ over $S$. \\
b) $z^{-\nu'}\varphi_{\ov{\MC{P}}}^s$ restricts to an isomorphism $\sigma^{* s}\ov{\MC{P}} \to \ov{\MC{P}}$ of $L^+\ov{P}$-torsors over $S$. Here we view $z^{-\nu'}$ via $LM(E) \subset LM(S) \subset L\ov{P}(S) = Hom_{L\ov{P}}(\MC{L\ov{P}}, \MC{L\ov{P}})$ as an isomorphism of the $L\ov{P}$-torsor associated to $\ov{\MC{P}}$. \\
c) $M$ is the centralizer of the cocharacter $\nu'$ in $G$. \\
The data $(P, s, \nu', (\ov{\MC{P}}, \varphi_{\ov{\MC{P}}}), \alpha)$ is called a complete slope division of $(\MC{G}, \varphi)$.
}

\rem{\label{rem:ComplSlopeDiviGeneral}}{
i) Via $\ov{\MC{P}} \hookrightarrow \ov{\MC{P}} \times^{L^+\ov{P}} L^+G \xrightarrow{\alpha^{-1}} \MC{G}$ we will view $\ov{\MC{P}}$ as a subtorsor of $\MC{G}$. Furthermore note that under this inclusion $\varphi_{\ov{\MC{P}}}$ is just the restriction of $\alpha \circ \varphi \circ \sigma^*(\alpha)^{-1}$ to $\ov{\MC{P}}$. \\
ii) Consider a complete slope division $(P, s, \nu', (\ov{\MC{P}}, \varphi_{\ov{\MC{P}}}), \alpha)$ of $(\MC{G}, \varphi)$ and some positive integer $N$. We claim that $(P, Ns, N\nu', (\ov{\MC{P}}, \varphi_{\ov{\MC{P}}}), \alpha)$ is again a complete slope division. Indeed we only have to check that $z^{-N\nu'} \varphi_{\ov{\MC{P}}}^{Ns}$ restricts to an isomorphism $\sigma^{* Ns}\ov{\MC{P}} \to \ov{\MC{P}}$. To do so, write
\[z^{-N\nu'}\varphi_{\ov{\MC{P}}}^{Ns} = \prod_{j = 1}^N z^{-(N-j)\nu'} (z^{-\nu'}\varphi_{\ov{\MC{P}}}^s) z^{(N-j)\nu'}\]
Because $\nu'$ is $G$-dominant, conjugation with $z^{-\nu'}$ fixes $\ov{P}$. Thus we get isomorphisms of local $\ov{P}$-shtukas
\[z^{-(N-j)\nu'} (z^{-\nu'}\varphi_{\ov{\MC{P}}}^s) z^{(N-j)\nu'}: \sigma^{* (N-j+1)s}\ov{\MC{P}} \to \sigma^{* (N-j)s}\ov{\MC{P}}\]
and their composite $z^{-N\nu'}\varphi_{\ov{\MC{P}}}^{Ns}$ has the desired properties. \\
iii) Condition c) can always be satisfied by enlarging $M$ and $P$ and replacing the subtorsor $\ov{\MC{P}}$ by its image under the action of the new opposite parabolic subgroup. Nevertheless it seems convenient to require c) as it is necessary for the following uniqueness result and will imply that the isoclinic constituents (cf. \ref{const:ComplSlopeDiviFiltration}) are in some sense maximal. \\
iv) Note that the definition makes sense for every Borel defined over $E$ (and not necessarily over $\Fq$). Consider therefore a Borel $B' = w.B$ obtained from $B$ by the action of the Weyl group element $w \in W$, which we assume to be defined over $E$. Then a complete slope division $(P, s, \nu', (\ov{\MC{P}}, \varphi_{\ov{\MC{P}}}), \alpha)$ w.r.t. $B$ yields one w.r.t $B'$ by conjugation of all constituents by $w$. \\
v) In the case of $G = GL_n$ the analogy to complete slope divisibility of $p$-divisible groups is explained in \cite[proposition 5.4]{HaVi2}. \\
vi) We try to distinguish between $\sigma$-conjugacy classes, usually denoted by $\nu$, and (multiples of) the corresponding Newton point (which was defined in \ref{const:LocalB(G)} as the image of $\nu$ under the map $\MC{B}(G) \to (X_*(T)_{\M{Q}}/W_G)^{Gal(k/\Fq)}$), usually denoted by $\nu'$. \\
vii) As $\nu'$ is a constant multiple of the Newton point, which is defined over $\Fq$, there is no loss of generality to assume that $\nu'$ is defined over $\Fq$ and not only over $E$. Moreover as $M$ is the centralizer of $\nu'$, $M$ is automatically defined over $\Fq$. Nevertheless we do not require this for $P$ or for $N$.
}

\lem{\label{lem:ComplSlopeDiviUnique}}{}{
Let $(\MC{G}, \varphi)$ be a completely slope divisible local $G$-shtuka over a scheme $S$. Then for fixed $s$ the complete slope division is unique up to unique isomorphism.
}

\prooof
By condition b), $\frac 1s\nu'$ equals the constant Newton point of $(\MC{G}, \varphi)$. As we require $\nu'$ to be $G$-dominant, it is uniquely determined. By c) the Levi $M$ is defined by $\nu'$ and as $P$ is standard we have $P = M \cdot B$. Thus it remains to check the uniqueness of the $L^+\ov{P}$-torsor $\ov{\MC{P}}$ and the isomorphism $\alpha$: \\ 
It suffices to do this over each trivializing \'etale cover. Therefore we may assume $(\MC{G}, \varphi) = (L^+G, b \sigma^*)$. Consider now two complete slope divisions $(P, s, \nu', (\ov{\MC{P}}_i, \varphi_{\ov{\MC{P}} i}), \alpha_i)$ for $i = 1, 2$. By viewing $\ov{\MC{P}}_i$ (for each $i$) via $\alpha_i$ as a subtorsor of $\MC{G}$, we may identify it with $L^+\ov{P} \xhookrightarrow{g_i} L^+G$ for some element $g_i \in L^+G(S)$. Then $\varphi_{\ov{\MC{P}} i}$ is identified with $p_i\sigma^* \coloneqq g_i^{-1}b\sigma^*(g_i)\sigma^*: \sigma^*L\ov{P} \to L\ov{P}$, where $p_i$ lies in $L\ov{P}(S)$. Define $g_0 = g_2^{-1}g_1 \in L^+G(S)$. Then $p_2 = g_0 p_1 \sigma^*(g_0)^{-1}$ and our main goal is to prove $g_0 \in L^+\ov{P}(S)$. For any integer $N > 0$ remark \ref{rem:ComplSlopeDiviGeneral}ii) enables us to write $(p_i\sigma^*)^{Ns} = z^{N\nu'} \cdot \tilde{p}_{i, N} \cdot \sigma^{* Ns}$ for some $\tilde{p}_{i, N} \in L^+\ov{P}(S)$. Thus we get
\[z^{N\nu'} \tilde{p}_{2, N} \sigma^{* Ns} = (p_2\sigma^*)^{Ns} = (g_0 p_1 \sigma^*(g_0)^{-1}\sigma^*)^{Ns} = g_0 (p_1\sigma^*)^{Ns} g_0^{-1} = g_0 z^{N\nu'} \tilde{p}_{1, N} \sigma^{* Ns}(g_0)^{-1} \sigma^{* Ns}\]
which rewrites as
\[z^{-N\nu'} g_0 z^{N\nu'} = \tilde{p}_{2, N} \sigma^{*Ns}(g_0) \tilde{p}_{1, N}^{-1} \in L^+G(S).\]
Because $\nu'$ is $G$-dominant, $\ov{P}$ is opposite to a standard parabolic and its Levi $M$ is the centralizer of $\nu'$, we get
\[\bigcap_{N \geq 0} z^{N\nu'} L^+G(S) z^{-N\nu'} = L^+\ov{P}(S)\]
Hence in particular $g_0 \in L^+\ov{P}(S)$ and we may define the isomorphism $g_0: (L^+\ov{P}, p_1\sigma^*) \to (L^+\ov{P}, p_2\sigma^*)$ of local $\ov{P}$-shtukas. This gives the desired isomorphism between the complete slope divisions. Uniqueness of the isomorphism follows from the fact, that it has to fit into a commutative diagram
\[
\begin{xy}
 \xymatrix @C=+0mm {
    \qquad (\ov{\MC{P}}_1, \varphi_{\ov{\MC{P}} 1}) \ar[rr] \ar@{^{(}->}[dr] & & (\ov{\MC{P}}_2, \varphi_{\ov{\MC{P}} 2}) \qquad \ar@{_{(}->}[dl]  \\
    & (\MC{G}, \varphi) &  }
\end{xy} 
\]
\exit

\lem{}{}{
Let $(\MC{G}, \varphi)$ be a completely slope divisible local $G$-shtuka over $S$. Then the torsor $\MC{G}$ trivializes over a pro-finite \'etale cover. 
}

\prooof
By definition $(\MC{G}, z^{-\nu'}(\varphi_{\ov{\MC{P}}} \times^{L\ov{P}} \id_{LG})^s)$ is an \'etale local $G$-shtuka (where the Frobenius is taken relative to $\M{F}_{q^s}$). Thus this is a direct consequence of proposition \ref{prop:TateTrivialTorsor}.

\rem{}{
In fact over central leaves, this pro-finite \'etale cover may be taken to be the tower of Igusa varieties, cf. section \ref{subsec:DefineIgusaGeneral} and proposition \ref{prop:IgusaTowerModuli}. 
}

\prop{\label{prop:ComplSlopeDiviField}}{}{
Let $S' \to S$ be any fpqc-morphism of schemes over $E$. Then a local $G$-shtuka $(\MC{G}, \varphi)$ over $S$ is completely slope divisible if and only if its base-change $(\MC{G}_{S'}, \varphi_{S'})$ to $S'$ is completely slope divisible. \\
In particular this applies to $S' =  \Sp \ov{K} \to S = \Sp K$ for any inclusion of a field $K$ (containing $E$) into an algebraic closure $\ov{K}$.
}

\prooof
Any complete slope division of $(\MC{G}, \varphi)$ induces by base-change of the local $\ov{P}$-shtuka to $S'$ one for $(\MC{G}_{S'}, \varphi_{S'})$. \\
Thus assume now that $(\MC{G}_{S'}, \varphi_{S'})$ has a complete slope division $\un{\ov{\MC{P}}} \coloneqq (P, s, \nu', (\ov{\MC{P}}, \varphi_{\ov{\MC{P}}}), \alpha)$. We will apply fpqc-decent to get a complete slope division over $S$: By the previous lemma there is an isomorphism $\xi: pr_1^*\un{\ov{\MC{P}}} \to pr_2^*\un{\ov{\MC{P}}}$ over $S' \times_S S'$ between the pullbacks along the two projections $pr_i: S' \times_S S' \to S'$ (onto the $i$th factor). Let $pr_{ij}: S' \times_S S' \times_S S' \to S' \times_S S'$ be the projection on the $i$th and $j$th factor. Then 
\[pr_{23}^*(\xi) \circ pr_{12}^*(\xi): pr_{12}^* pr_1^*\un{\ov{\MC{P}}} \to pr_{12}^* pr_2^*\un{\ov{\MC{P}}} \cong pr_{23}^* pr_1^*\un{\ov{\MC{P}}} \to pr_{23}^* pr_2^*\un{\ov{\MC{P}}}\] 
and 
\[\hspace{13mm} pr_{13}^*(\xi): pr_{12}^* pr_1^*\un{\ov{\MC{P}}} \cong pr_{13}^* pr_1^*\un{\ov{\MC{P}}} \to pr_{13}^* pr_2^*\un{\ov{\MC{P}}} \cong pr_{23}^* pr_2^*\un{\ov{\MC{P}}}\]
are two isomorphisms between complete slope divisions over $S' \times_S S' \times_S S'$. Hence by the uniqueness assertion of the previous lemma we get indeed $pr_{23}^*(\xi) \circ pr_{12}^*(\xi) = pr_{13}^*(\xi)$. Hence the decent datum $\xi$ is effective and defines a complete slope division over $S$. \exit

\rem{}{
For the analogous result by Oort and Zink for $p$-divisible groups (at least in the case of field extensions) see the remark in front of proposition $1.3$ in \cite{OortZinkFamilies}. There it is noted that the statement follows directly from the existence of slope filtrations for $p$-divisible groups over arbitrary fields. We may define slope divisions for local $G$-shtukas by replacing conditions b) and c) in definition \ref{def:ComplSlopeDivi} by requiring that $M$ centralizes the $G$-dominant Newton point. Then the same implication would be true for local $G$-shtukas, but unfortunately the existence of slope filtrations is in general unknown to us.
}

For our application to Igusa varieties we need to know whether complete slope divisions defined over the generic fiber of a scheme extend over all of it. We first prove that at least the $L^+\ov{P}$-torsor extends:

\lem{\label{lem:ExtendPTorsor}}{}{
Let $S$ be a connected normal noetherian scheme over $E$ and $\MC{G}$ be a $L^+G$-torsor over $S$ with generic fiber $\MC{G}_\eta$. Let $P \subset G$ be a standard parabolic subgroup, $\ov{P}$ its opposite parabolic subgroup and $\ov{\MC{P}}_\eta$ a $L^+\ov{P}$-torsor over the generic fiber of $S$ together with a $\ov{P}$-equivariant inclusion $\ov{\MC{P}}_\eta \hookrightarrow \MC{G}_\eta$. Then the closure of $\ov{\MC{P}}_\eta$ in $\MC{G}$ defines a $\ov{P}$-torsor over $S$.
}

\prooof
Obviously there is at most one $\ov{P}$-torsor over $S$ extending $\ov{\MC{P}}_\eta$. Hence it suffices to prove the assertion over an \'etale cover. \\
Assume first that $S = \Sp A$ is affine with $A$ a DVR. As $A$ is local and $G$ isotrivial we may find a connected finite \'etale cover $S' = \Sp A'$ of $A$ on which $(\MC{G}, \varphi)$ trivializes. Then $A'$ is again a DVR and we may assume wlog. that $\MC{G} = L^+G$ over $S$ and $\ov{\MC{P}}_\eta = L^+\ov{P}_\eta \hookrightarrow L^+G_\eta$ over the generic point. 
Let $n > 0$ and consider now $L^+\ov{P}[n] = L^+\ov{P}/(K_n \cap L^+\ov{P}) \subset L^+G[n] = L^+G/K_n$. Then the data above induce a $L^+\ov{P}[n]$-equivariant inclusion $L^+\ov{P}[n]_\eta \hookrightarrow L^+G[n]_\eta$ into the generic fiber of the trivial $L^+G[n]$-torsor over $S$. This is equivalent to a morphism $\eta \to L^+\ov{P}[n] \backslash L^+G[n]$ over $\Fq$, where the target is a quotient of group schemes over $\Fq$. This quotient is well-defined as $L^+G[n] = Res_{\Fq[[z]]/(z^n)/\Fq}(G \times \Fq[[z]]/(z^n))$ is a linear algebraic group and the same holds for $L^+\ov{P}[n]$. 
Furthermore it is a projective variety as being a parabolic subgroup is stable under restriction of scalars. Hence by the valuative criterion for properness, this morphism extends to $S \to L^+\ov{P}[n] \backslash L^+G[n]$, i.e. gives a canonical $L^+\ov{P}[n]$-subtorsor of $L^+G[n]$ over all of $S$ which extends $L^+\ov{P}[n]_\eta$. As these subtorsors are by construction compatible for different $n$, they give an inclusion $L^+\ov{P} \hookrightarrow L^+G$ extending $L^+\ov{P}_\eta \hookrightarrow L^+G_\eta$. It is now easy to see, that the defined $L^+\ov{P}$-torsor is indeed the closure of $\ov{\MC{P}}_\eta$. \\
Let us now consider the case of general connected normal noetherian base schemes $S$: By the previous considerations over DVRs, $\ov{\MC{P}}_\eta$ extends to each point of codimension $1$. Then by normality of $S$, the sections defining $\ov{\MC{P}}$ in codimension at most $1$ extend to sections over all of $S$. Hence $\ov{\MC{P}}$ is indeed defined over all of $S$. \exit

\rem{}{
The idea of the first part is to consider $\ov{\MC{P}}_\eta \hookrightarrow \MC{G}_\eta$ as a point in the ``proper scheme'' $L^+\ov{P} \backslash L^+G$ to apply the valuative criterion. Nevertheless we do not know, whether the quotient $L^+\ov{P} \backslash L^+G = \varprojlim_n L^+\ov{P}[n] \backslash L^+G[n]$ exists as an infinite-dimensional quasi-compact proper scheme.
}

\thm{\label{thm:ComplSlopeDiviGeneric}}{}{
Let $(\MC{G}, \varphi)$ be a local $G$-shtuka with constant Newton point (and hence constant quasi-isogeny class) over a connected normal locally noetherian scheme $S$. Then $(\MC{G}, \varphi)$ is completely slope divisible if and only if it is completely slope divisible over the generic fiber.
}

\prooof
As the 'only if' part is trivial, let us consider a complete slope division $(P, s, \nu', (\ov{\MC{P}}_\eta, \varphi_{\ov{\MC{P}} \eta}), \alpha_\eta)$ over the generic fiber. Cover $S$ by connected normal noetherian open subschemes. Then the uniqueness statement lemma \ref{lem:ComplSlopeDiviUnique} implies that we may glue complete slope divisions over these open subschemes. Hence we may assume wlog. that $S$ is noetherian. \\
View again $\ov{\MC{P}}_\eta$ as a subtorsor of $\MC{G}_\eta$ (the generic fiber of $\MC{G}$). Then by the previous lemma, the closure $\ov{\MC{P}}$ of $\ov{\MC{P}}_\eta$ defines a $\ov{P}$-torsor over $S$. The remaining data are now easy to construct: We keep $P$, $s$ and $\nu'$ unchanged. The action of $L^+G$ on $\MC{G}$ gives (when restricted to $\ov{\MC{P}}$) a morphism $\ov{\MC{P}} \times L^+G \to \MC{G}$, which factors over a morphism $\alpha^{-1}: \ov{\MC{P}} \times^{L^+\ov{P}} L^+G \to \MC{G}$. This $\alpha^{-1}$ extends by definition $\alpha_\eta^{-1}$. Consider now the morphism 
\[\alpha \circ \varphi \circ \sigma^*{\alpha^{-1}}: \sigma^*\MC{L}\ov{\MC{P}} \times^{L\ov{P}} LG \to \MC{L}\ov{\MC{P}} \times^{L\ov{P}} LG.\]
Because its restriction to the generic fiber gives a morphism between $L\ov{P}$-torsors and because $\ov{P} \subset G$ is closed, this defines an isomorphism $\varphi_{\ov{P}}: \sigma^*\MC{L}\ov{\MC{P}} \to \MC{L}\ov{\MC{P}}$. By definition of $\varphi_{\ov{P}}$ the morphism $\alpha$ gives indeed an isomorphism of local $G$-shtukas and not only of $L^+G$-torsors. \\ 
We are left to show that $z^{-\nu'}\varphi_{\ov{\MC{P}}}^s$ restricts to an isomorphism between $L^+\ov{P}$-torsors: It suffices to do this on some trivializing \'etale cover. Hence assume that $z^{-\nu'}\varphi_{\ov{\MC{P}}}^s$ is represented by $z^{-\nu'}p\sigma^{*s}$ for some element $p \in L\ov{P}(S)$ (with $S = \Sp A$ affine) whose generic fiber $p_\eta$ lies in $L^+\ov{P}(K)$ for the fraction field $K$ of $A$. 
Furthermore choosing a faithful representation $\ov{P} \to GL_n$ lets us consider everything in $GL_n$ and we have to prove $z^{-\nu'}p \in L^+GL_n(S)$ (for the precise argument see the proof of proposition \ref{prop:ThmTateGShtukaDVR}). Thus writing $\Lambda = A[[z]]^n$, $\Lambda' = A((z))^n$ and $\Lambda_\eta = K[[z]]^n$ we may view $z^{-\nu'}p\sigma^{*s}: \sigma^{*s}\Lambda' \to \Lambda'$ and know that it defines on the generic fiber an isomorphism $(z^{-\nu'}p)_\eta\sigma^{*s}: \sigma^{*s}\Lambda_\eta \stackrel{\sim}{\to} \Lambda_\eta$. \\
We first show that $z^{-\nu'}p\sigma^{*s}$ restricts to a morphism $\sigma^{*s}\Lambda \to \Lambda$: Let $b_1 = z^{-\nu'} p \sigma^*(z^{-\nu'})^{-1}$ and $b_2 = \sigma^*(p)$. As the Newton point of $p$ was assumed to be constant, the same holds for $b_1$ and $b_2$. Then we saw in the proof of proposition \ref{prop:ThmTateGShtukaDVR} that any element in $End^*(\Lambda_\eta)$ which $\sigma$-conjugates ${b_1}_\eta$ into ${b_2}_\eta$ extends uniquely to an element in $End^*(\Lambda)$ (actually we showed this over a complete DVR, but the reduction steps done in the proof of corollary \ref{thm:ThmTateGShtuka} work here as well). Applied to $(z^{-\nu'}p)_\eta$ we see that it extends to a morphism of $\Lambda$ which coincides with $z^{-\nu'}p$ by uniqueness of the extension. \\
The bijectivity of $z^{-\nu'}p\sigma^{*s}: \sigma^{*s}\Lambda \to \Lambda$ follows now quickly: We know that $\det(z^{-\nu'}p) \in A[[z]] \cap A((z))^\times$ as $z^{-\nu'}p$ lies in $End^*(\Lambda)$. Furthermore $(z^{-\nu'}p)_\eta$ is by definition an automorphism of $\Lambda_\eta$. Therefore $\det(z^{-\nu'}p) \in K[[z]]^\times$. Together this implies $\det(z^{-\nu'}p) \in A[[z]]^\times$, i.e. $z^{-\nu'}p\sigma^*$ defines indeed an isomorphism. \exit

\rem{}{
i) The analogous result for $p$-divisible groups (even over more general base schemes) can be found in \cite[proposition 2.3]{OortZinkFamilies}. \\
ii) Note that the constancy of the Newton point is essential, because otherwise it would follow that central leaves are not only closed in Newton strata, but closed in the whole moduli space. 
}

\thm{\label{thm:UnivComplSlopeDiv}}{}{
The local $G$-shtuka $(\MC{G}^{univ}_{c_i, U}, \varphi^{univ}_{c_i, U})$ associated (at the characteristic place $c_i$) to the universal global $G$-shtuka over the central leaf $\MC{C}_U^{(\nu_i)}$ is completely slope divisible.
}

\prooof
As complete slope divisions are unique, we may check this over an \'etale cover. So let $S$ be any irreducible scheme \'etale over $\MC{C}_U^{(\nu_i)}$, let $\eta \in S$ be its generic point and $\ov{\eta}$ a geometric point with image $\eta$. The the pullback of $(\MC{G}^{univ}_{c_i, U}, \varphi^{univ}_{c_i, U})|_S$ to $\ov{\eta}$ is completely slope divisible by definition of the central leaf. Hence by proposition \ref{prop:ComplSlopeDiviField} it is already completely slope divisible over the generic point. Now theorem \ref{thm:ComplSlopeDiviGeneric} implies that $(\MC{G}^{univ}_{c_i, U}, \varphi^{univ}_{c_i, U})|_S$ is completely slope divisible over the smooth scheme $S$. \exit
$\quad$ \vspace{3mm} \\
In the theory of completely slope divisible $p$-divisible groups, one can define associated basic $p$-divisible groups as the composition factors in the filtration given by the complete slope division. In the following we explain the analogous constructions for completely slope divisible local $G$-shtukas. They will only be needed as an auxiliary tool to deduce the representability of arbitrary Igusa varieties from the basic case, cf. section \ref{subsec:DefineIgusaGeneral}.

\defi{}{}{
a) A local $G$-shtuka $(\MC{G}, \varphi)$ over a scheme $S$ is called basic if there is an integer $s > 0$ and a central cocharacter $\nu': \M{G}_m \to G$ such that its Newton point is constant and equal to $\frac 1s \nu'$. \\
b) A local $G$-shtuka $(\MC{G}, \varphi)$ over a scheme $S$ is called completely basic if there is an integer $s > 0$ and a central cocharacter $\nu': \M{G}_m \to G$ such that $z^{-\nu'}\varphi^s$ restricts to an isomorphism $\sigma^{* s}\MC{G} \to \MC{G}$ over $S$.
}

\rem{}{
a) In the case of $G = GL_n$ and $S = \Sp k$ for an algebraically closed field $k$, a local $G$-shtuka is basic if and only if all its Newton slopes are equal. \\ 
b) Any completely basic local $G$-shtuka is basic: Indeed consider any completely basic $(\MC{G}, \varphi)$ and choose a geometric point $\ov{x} \in S(k)$ ($k$ some algebraically closed field). Then there is a trivialization $(\MC{G}, \varphi) \cong (L^+G, b_{\ov{x}}\sigma^*)$ such that $z^{-\nu'}(b_{\ov{x}}\sigma^*)^s = g_{\ov{x}}\sigma^{*s}$ for some $g_{\ov{x}} \in L^+G(k)$. By the Theorem of Lang-Steinberg, which applies to the group scheme $L^+G$ by \cite[lemma 2.1]{VieTruncLevel1} (for $H = K = L^+G$ and $g = 1$), there is some $g'_{\ov{x}} \in L^+G(k)$ with $g_{\ov{x}} = g'_{\ov{x}}\sigma^{*s}(g'_{\ov{x}})^{-1}$. Using that $z^{-\nu'}$ is central in $G$ we get 
\[{g'_{\ov{x}}}^{-1}(b_{\ov{x}}\sigma^*)^s g'_{\ov{x}} = z^{\nu'}.\]
In particular the Newton point of $b_{\ov{x}}$ equals $\frac 1s \nu'$.
}

\const{\label{const:ComplSlopeDiviFiltration}}{}{
Let $(\MC{G}, \varphi)$ be a local $G$-shtuka over a scheme $S$ and $(P, s, \nu', (\ov{\MC{P}}, \varphi_{\ov{\MC{P}}}), \alpha)$ a complete slope division for it. Let $M$ be the Levi subgroup of $P$ or equivalently the Levi subgroup of $\ov{P}$. Then the canonical projection $\ov{P} \to M$ induces the local $M$-shtuka
\[(\MC{M}, \varphi_{\MC{M}}) \coloneqq (\ov{\MC{P}} \times^{L^+\ov{P}} L^+M, \varphi_{\ov{\MC{P}}} \times^{L\ov{P}} \id_{LM})\]
This $(\MC{M}, \varphi_{\MC{M}})$ is called the basic constituent of $(\MC{G}, \varphi)$.
}

\prop{}{}{
Let $(\MC{G}, \varphi)$ be a completely slope divisible local $G$-shtuka over $S$. Then the basic constituent $(\MC{M}, \varphi_{\MC{M}})$ is a completely basic local $M$-shtuka.
}

\prooof
Fix a complete slope division $(P, s, \nu', (\ov{\MC{P}}, \varphi_{\ov{\MC{P}}}), \alpha)$. Then $\nu'$ defines a $M$-central cocharacter and by assumption $z^{-\nu'}\varphi_{\MC{M}}^s = z^{-\nu'}\varphi_{\ov{\MC{P}}}^s \times^{L\ov{P}} \id_{LM}: \sigma^{*s}\MC{LM} \to \MC{LM}$ induces an isomorphism between $L^+M$-torsors. \exit

\subsection{Iwahori-type structures on completely slope divisible local \textit{G}-shtukas}\label{subsec:IwahoriStructure}
Recall that a fundamental $P$-alcove $b$ satisfies $\phi_b(I_M) = I_M$, $\phi_b(I_N) \subseteq I_N$ and $\phi_b(I_{\ov{N}}) \supseteq I_{\ov{N}}$ for $M \subset P$ the Levi, $N \subset P$ the unipotent radical and 
\[\phi_b: LG \to LG \quad , \quad g \mapsto \sigma(b \cdot g \cdot b^{-1}).\]
Assume that $b_\nu$ is a fundamental alcove in the $\sigma$-conjugacy class $\nu \in \MC{B}(G)$ and that $(\MC{G}, \varphi)$ is a completely slope divisible local $G$-shtuka over a scheme $S$ with the property that its restriction to any geometric point is isomorphic to the local $G$-shtuka $(L^+G, b_\nu \sigma^*)$. Our aim is to construct a subgroup $I_0(b_\nu) \subset L^+G$ with the following two properties:
\begin{itemize}
 \item There is a canonical $I_0(b_\nu)$-torsor $\MC{I}_0$ over $S$ together with an isomorphism $\MC{I}_0 \times^{I_0(b_\nu)} L^+G \cong \MC{G}$.
 \item There is a sequence $\ldots \subset I_n(b_\nu) \subset I_{n-1}(b_\nu) \subset \ldots \subset I_0(b_\nu)$ of normal subgroups of $I_0(b_\nu)$ with $\bigcap_n I_n(b_\nu) = \{1\}$ and $b_\nu^{-1} \cdot I_n(b_\nu) \cdot b_\nu \subset I_n(b_\nu)$ for all $n$. 
\end{itemize}
We do not require that the $I_n(b_\nu)$ are open in $L^+G$; in fact the $I_n(b_\nu)$ will only be open in $L^+\ov{P}_{b_\nu}$ for some parabolic subgroup $\ov{P}_{b_\nu}$. For a comparison of $\MC{I}_0$ and its induced $\ov{P}_{b_\nu}$-torsor and complete slope divisions, see proposition \ref{prop:I_0Shtuka}. \\
We will show later on in proposition \ref{prop:IgusaRepComplSlope}, that the universal local $G$-shtuka over the perfection of a central leaf has a canonical $I_0(b_\nu)$-structure.

\defi{\label{def:IwahoriDefinition}}{}{
Let $b_\nu \in LG(E)$ be a $P$-fundamental alcove (for some parabolic $P$). Then define for every $n \geq 0$ the subgroup of $L^+G$ defined by
\[I_n(b_\nu) = \bigcap_{N \geq 0} \phi_{b_\nu}^N(K_n) \subset LG.\]
}

\lem{\label{lem:IwahoriDefinition2}}{}{
Let $I_n(b_\nu)$ be as in the definition. Then there is a parabolic subgroup $\ov{P}_{b_\nu} \subset G$ and a subgroup $Q_{b_\nu} \subset \ov{P}_{b_\nu}$ containing $T$ such that (on each $E$-scheme $S$)
\[I_n(b_\nu)(S) = \{g \in (K_n \cap L^+\ov{P}_{b_\nu})(S) \;|\; g \bmod z^{n+1} \in Q_{b_\nu}(\MC{O}_S(S)[[z]]/(z^{n+1}))\}.\]
In particular $I_n(b_\nu)$ is representable as a group scheme over $E$. Furthermore $\ov{P}_{b_\nu}$ and $Q_{b_\nu}$ can both be chosen independently of $n$.
}

\prooof
Let $r > 0$ be an integer such that $G$ splits over the extension $\M{F}_{q^r}$ of $\Fq$ of degree $r$. Then $\sigma$ acts on $\widetilde{W}$ (of order at most $r$) and we may view $b_\nu \cdot \sigma \in \widetilde{W} \rtimes \langle \sigma \rangle$. Then $\phi_{b_\nu}$ is conjugation with this element, taking the non-trivial $\sigma$-action into account. 
Then by \cite[remark 6.2]{VieTruncLevel1} $b_\nu' \coloneqq (b_\nu \cdot \sigma)^r = b_\nu \cdot \sigma(b_\nu) \cdot \sigma^2(b_\nu) \cdot \ldots \cdot \sigma^{r-1}(b_\nu) \cdot \sigma^r \in \widetilde{W} \rtimes \langle \sigma \rangle/\sigma^r$ lies in the extended Weyl group $\widetilde{W}_M$ of $M$ and its $M$-dominant Newton point (wrt. the Frobenius $\sigma^r$) is central in $M$. Let $M_{b_\nu}$ be the centralizer of this Newton point and let $P_{b_\nu} = P \cdot M_{b_\nu}$ with Levi $M_{b_{\nu}}$ and unipotent radical $N_{b_\nu}$. Then $b_\nu$ is a $P_{b_\nu}$-fundamental alcove (cf. \cite[lemma 6.3]{VieTruncLevel1} and its proof). \\
Choose now a further integer $s > 0$ such that ${b_\nu'}^s \subset X_*(T) \subset \widetilde{W}$, i.e. we may write $(b_\nu \cdot \sigma)^{rs} = z^{\nu'} \sigma^{rs}$ for some $\sigma^r$-invariant element $\nu' \in X_*(T)$, which is then by definition a multiple of the Newton point of $b_\nu$. In particular $\phi_{b_\nu}^{rs}(g) = z^{\nu'} \cdot \sigma^{rs}(g) \cdot z^{-\nu'}$. Thus denoting the parabolic opposite to $P_{b_\nu}$ by $\ov{P}_{b_\nu}$ and recalling that $M_{b_\nu}$ is the centralizer of $\nu'$, we get
\begin{align*}
 \bigcap_{N \geq 0} \phi_{b_\nu}^{Nrs}(K_n) & = \bigcap_{N \geq 0} z^{N\nu'} \cdot \sigma^{Nrs}(K_n) \cdot z^{-N\nu'} = \bigcap_{N \geq 0} z^{N\nu'} \cdot K_n \cdot z^{-N\nu'} \\
 & = K_n \cap L^+\ov{P}_{b_\nu} = K_n \cap L\ov{P}_{b_\nu}
\end{align*}
We claim
\[I_n(b_\nu) = \bigcap_{N = 0}^{rs-1} \phi_{b_\nu}^{N}(K_n \cap L\ov{P}_{b_\nu}).\]
Indeed if $g \in I_n(b_\nu)$, then $g \in L\ov{P}_{b_\nu}$ by the computation above. But as $b_\nu$ is a $P_{b_\nu}$-fundamental alcove, we have $\phi_{b_\nu^{-1}}(L\ov{P}_{b_\nu}) = L\ov{P}_{b_\nu}$. Thus for every $N \geq 0$ we have $\phi_{b_\nu}^{-N}(g) \subset K_n \cap L\ov{P}_{b_\nu} = K_n \cap L^+\ov{P}_{b_\nu}$. Conversely if $g$ is an element of the right hand side and $N = \delta rs + \varepsilon$ for $\delta \geq 0$ and $\varepsilon \in \{0, 1, \ldots rs-1\}$, then $\phi_{b_\nu}^{-\varepsilon}(g) \subset K_n \cap L\ov{P}_{b_\nu}$ by assumption and $\phi_{b_\nu}^{-\delta rs}(\phi_{b_\nu}^{-\varepsilon}(g)) \subset K_n$ by the previous computation. \\
Therefore it suffices to show that for every $N \in \{0, 1, \ldots, rs - 1\}$ there is a subgroup $Q_{b_\nu}^N \subset \ov{P}_{b_\nu}$ satisfying
\begin{align*}
 I_n^N(b_\nu)(S) & \coloneqq K_n(S) \cap \phi_{b_\nu}^{N}(K_n \cap L\ov{P}_{b_\nu})(S) \\
 & \, = \{g \in (K_n \cap L\ov{P}_{b_\nu})(S) \;|\; g \bmod z^{n+1} \in Q_{b_\nu}^N(\MC{O}_S(S)[[z]]/(z^{n+1}))\}
\end{align*}
on $E$-schemes $S$. 
First note that it follows from $I_{\ov{P}_{b_\nu}} \subseteq \phi_{b_\nu}(I_{\ov{P}_{b_\nu}})$ that also $I_n \cap L\ov{P}_{b_\nu} \subseteq \phi_{b_\nu}(I_n \cap L\ov{P}_{b_\nu})$ (for the subgroup $I_n$ as defined in \ref{nota:Iwahori}). In particular
\[K_{n+1} \cap L\ov{P}_{b_\nu} \subset I_n \cap L\ov{P}_{b_\nu} \subset K_n \cap \phi_{b_\nu}^{N}(K_n \cap L\ov{P}_{b_\nu}) = I_n^N(b_\nu)\]
Hence it suffices to show
\[I_n^N(b_\nu)/K_{n+1} = (K_n \cap L^+Q_{b_\nu}^N)/K_{n+1}\]
for some $Q_{b_\nu}^N \subset \ov{P}_{b_\nu}$. With $\widetilde{W} \rtimes \langle \sigma \rangle = X_*(T) \rtimes W \rtimes \langle \sigma \rangle$ write now $(b_\nu \cdot \sigma)^N = z^{\nu_N'} \cdot w_N \cdot \sigma^N$ for some $\nu_N' \in X_*(T)$ and $w_N \in W$. Then
\begin{align*}
 K_n \cap \phi_{b_\nu}^{N}(K_n \cap L\ov{P}_{b_\nu}) & = K_n \cap z^{\nu_N'} \cdot w_N \cdot \sigma^N(K_n \cap L\ov{P}_{b_\nu}) \cdot w_N^{-1} \cdot z^{-\nu_N'} \\
 & = K_n \cap z^{\nu_N'} \cdot (K_n \cap L\ov{P}_{b_\nu}^N) \cdot z^{-\nu_N'} \\
 & = K_n \cap (z^{\nu_N'} \cdot K_n \cdot z^{-\nu_N'}) \cap L\ov{P}_{b_\nu}^N
\end{align*}
with $\ov{P}_{b_\nu}^N \coloneqq w_N \cdot \sigma^N\ov{P}_{b_\nu} \cdot w_N^{-1}$ being another parabolic subgroup defined over $E$. But there is obviously a parabolic subgroup $Q^N$ such that
\[(K_n \cap z^{\nu_N'} \cdot K_n \cdot z^{-\nu_N'})/K_{n+1} = (K_n \cap L^+Q^N)/K_{n+1}.\]
Then $Q_{b_\nu}^N = Q^N \cap \ov{P}_{b_\nu}^N$ has all desired properties wrt. $I_n^N(b_\nu)$. Thus setting $Q_{b_\nu} = \bigcap_{N = 0}^{rs - 1} Q_{b_\nu}^N$ indeed has all the properties asserted in the lemma. \exit

\rem{\label{rem:IwahoriSlopeDivision}}{
Fix any Borel $B_{b_\nu}$ inside $P_{b_\nu}$ and consider the completely slope divisible local $G$-shtuka $(L^+G, b_\nu \sigma^*)$. By remark \ref{rem:ComplSlopeDiviGeneral}iv) it admits a complete slope division $(P, s, \nu, (P, b_\nu \sigma^*), \alpha)$ with respect to $B_{b_\nu}$. Then by construction of $P_{b_\nu}$ and $P$ (or rather their opposite parabolic subgroups), they coincide as subgroups of $G$. We will show a much stronger comparison result in proposition \ref{prop:I_0Shtuka}b).
}

\ex{}{
Assume that $G = GL_n$ with the Borel $B$ of upper triangular matrices and that the fundamental alcove $b_\nu \in LG(\Fq)$ is superbasic in the sense of \cite[section 5.9]{GHKR1}, i.e. its $\sigma$-conjugacy class does not intersect any proper Levi subgroup. In particular $b_\nu$ is isoclinic and a fundamental alcove for the parabolic $P_{b_\nu} = G$. Then $b_\nu$ has the form $w \cdot z^{\nu'}$ for $\nu' = (a, \ldots, a, a-1, \ldots, a-1) \in \M{Z}^n \cong X_*(T)$ with some integer $a$ appearing $d$ times, and $w \in W = S_n$ the permutation defined by shifting the indices of the canonical basis of $\M{Z}^n$ by $d$. Furthermore the condition superbasic implies $\gcd(d, n) = 1$ and that $w$ is a cyclic permutation of order $n$. Now an easy computation shows that $Q_{b_\nu} = B$ and $I_n(b_\nu) = I_n$ are the subgroups defined in \ref{nota:Iwahori} with respect to the Borel $B$. 
}

\ex{\label{ex:PseudoIwahori}}{
Assume now $G = GL_n$ with $B$ the upper triangular matrices and $b_\nu = (13)(254) \cdot z^{(1, 1, 0, 0, 0)}$ (identifying $W \cong S_5$ via permutation matrices and $X_*(T) \cong \M{Z}^5$). Then an explicit computation shows:
\[Q_{b_\nu} = \lmat{ccccc} * &  & * &  &  \\ * & * & * & * & * \\  &  & * &  &  \\ * &  & * & * & * \\  &  & * &  & * \rmat \quad \subset \quad 
\ov{P}_{b_\nu} = \lmat{ccccc} * &  & * &  &  \\ * & * & * & * & * \\ * &  & * &  &  \\ * & * & * & * & * \\ * & * & * & * & * \rmat\]
In particular $P_{b_\nu}$ may be non-standard and the $Q_{b_\nu}$ defined above may be not parabolic.
}

\lem{\label{lem:PseudoIwahoriInclusion}}{}{
Let $b_\nu$ be as before. Then $I_n(b_\nu)$ is a normal open subgroup of $I_0(b_\nu)$ for every $n$. Furthermore
\[\phi_{b_\nu}^{-1} (I_n(b_\nu)) \subseteq I_n(b_\nu).\]
}

\prooof
Consider any elements $g \in I_n(b_\nu)$ and $h \in I_0(b_\nu)$. Then for any $N \geq 0$
\[\phi_{b_\nu}^{-N} (h^{-1}gh) = \phi_{b_\nu}^{-N}(h)^{-1} \cdot \phi_{b_\nu}^{-N}(g) \cdot \phi_{b_\nu}^{-N}(h) \in \phi_{b_\nu}^{-N}(h)^{-1} \cdot K_n \cdot \phi_{b_\nu}^{-N}(h) \subset K_n\] 
because $\phi_{b_\nu}^{-N}(h) \in K_0 = L^+G$ and $K_n$ is normal in $K_0$. Hence $h^{-1}gh \in I_n(b_\nu)$. \\
The second assertion is obvious from the definition of $I_n(b_\nu)$.  \exit

\defi{\label{defi:StronglyComplSlopeDivi}}{}{
Let $(\MC{G}, \varphi)$ be a local $G$-shtuka over a scheme $S$. Then $(\MC{G}, \varphi)$ is called strongly completely slope divisible, if it has a complete slope division $(P, s, \nu', (\MC{\ov{P}}, \varphi_{\MC{\ov{P}}}), \alpha)$ such that after passing to a fpqc-cover of $S$, there exists an isomorphism $(\MC{\ov{P}}, \varphi_{\MC{\ov{P}}}) \cong (L^+P, b_\nu\sigma^*)$ for some $P$-fundamental alcove $b_\nu$ (whose $\sigma$-conjugacy class $\nu \in \MC{B}(G)$ has automatically a $\M{Q}$-multiple of $\nu'$ as Newton point).
}

\rem{}{
We will frequently use pro-finite \'etale covers to construct such isomorphisms. In fact any pro-finite \'etale cover is fpqc, because it is affine hence quasi-compact and faithfully flat as a limit of faithfully flat morphisms. 
}

\prop{\label{prop:I_0Structure}}{}{
Every strongly completely slope divisible local $G$-shtuka $(\MC{G}, \varphi)$ has a canonical $I_0(b_\nu)$-structure, for the element $b_\nu$ appearing in definition \ref{defi:StronglyComplSlopeDivi}. 
}

\prooof
We define the underlying $I_0(b_\nu)$-torsor as follows: Identify $\MC{G} \cong Isom_{L^+G}(L^+G, \MC{G})$. Then each $g \in \MC{G}$ defines for every $N \geq 0$ a morphism
\[g^\sharp_N: \sigma^{*N}LG \xrightarrow{(b_\nu\sigma^*)^N} LG \xrightarrow{\;g\;} \MC{LG} \xrightarrow{\varphi^{-N}} \sigma^{*N}\MC{LG}.\]
Define $\MC{I}_0 \subset \MC{G}$ as the fpqc-sheaf consisting of all elements $g$, such that $g^\sharp_N$ restricts to an isomorphism of $L^+G$-torsors for every $N$. We have to see that $\MC{I}_0$ is indeed an (\'etale) $I_0(b_\nu)$-torsor. As $I_0(b_\nu)$ is smooth, it suffices to check that $\MC{I}_0$ is a torsor for the fpqc-topology. But fpqc-locally we may fix an isomorphism $(\MC{\ov{P}}, \varphi_{\MC{\ov{P}}}) \cong (L^+P, b_\nu\sigma^*)$ and under this isomorphism
\[g^\sharp_N = (b_\nu \cdot \sigma^*)^{-N} \cdot g \cdot (b_\nu \cdot \sigma^*) = \sigma^{-1}(b_\nu) \cdot \ldots \cdot \sigma^{-N}(b_\nu) \cdot \sigma^{-N}(g)  \cdot \sigma^{-N}(b_\nu^{-1}) \cdot \ldots \cdot \sigma(b_\nu^{-1}) = \phi_{b_{\nu}}^{-N}(g)\]
Thus $g^\sharp_N \in \MC{I}_0$ if and only if $g \in \phi_{b_{\nu}}^N(L^+G)$ for every $N \geq 0$. Thus under this trivialization, $\MC{I}_0$ gets identified with $I_0(b_\nu) \subset L^+G$ (viewed as trivial torsors for the respective groups).  \exit

\rem{\label{rem:I_0Canonical}}{
The word canonical in the statement means that for every two such local $G$-shtukas $(\MC{G}_i, \varphi_i)$ (for $i = 1, 2$) with constant quasi-isogeny class $\nu \in \MC{B}(G)$ and each isomorphism $\xi: (\MC{G}_1, \varphi_1) \to (\MC{G}_2, \varphi_2)$ the $I_0(b_\nu)$-structures $\MC{I}_{0 i} \subset \MC{G}_i$ satisfy $\xi(\MC{I}_{0 1}) = \MC{I}_{0 2}$.  
}

\prop{\label{prop:I_0Shtuka}}{}{
Let $(\MC{G}, \varphi)$ be a strongly completely slope divisible local $G$-shtuka over a scheme $S$. \\
a) The $I_0(b_\nu)$-subtorsor $\MC{I}_0$ of $\MC{G}$ from proposition \ref{prop:I_0Structure} induces a canonical local $\ov{P}_{b_\nu}$-shtuka $(\ov{\MC{P}}_{b_\nu}, \varphi_{\ov{\MC{P}}_{b_\nu}})$ together with an isomorphism of local $G$-shtukas
\[\alpha: (\MC{G}, \varphi) \to (\ov{\MC{P}}_{b_\nu} \times^{L^+\ov{P}_{b_\nu}} L^+G, \varphi_{\ov{\MC{P}}_{b_\nu}} \times^{L\ov{P}_{b_\nu}} \id_{LG}).\]
b) If $B_{b_\nu}$ is any Borel in $P_{b_\nu}$, then $(\ov{\MC{P}}_{b_\nu}, \varphi_{\ov{\MC{P}}_{b_\nu}})$ coincides with the local $\ov{P}_{b_\nu}$-shtuka defined by the complete slope division with respect to the Borel $B_{b_\nu}$ (cf. remark \ref{rem:ComplSlopeDiviGeneral}iv))
}

\prooof
a) We claim that we may just take the triple consisting of the $\ov{P}_{b_\nu}$-torsor $\ov{\MC{P}}_{b_\nu} \coloneqq \MC{I}_0 \times^{I_0(b_\nu)} L^+\ov{P}_{b_\nu}$ viewed as a subtorsor of $\MC{G}$, the restriction $\varphi_{\ov{\MC{P}}_{b_\nu}} \coloneqq \varphi|_{\MC{L}(\MC{I}_0 \times^{I_0(b_\nu)} L^+\ov{P}_{b_\nu})}$ of $\varphi$ to the $L\ov{P}_{b_\nu}$-torsor associated to $\MC{I}_0 \times^{I_0(b_\nu)} L^+\ov{P}_{b_\nu}$ and as $\alpha$ the isomorphism induced by the canonical inclusion $\MC{I}_0 \times^{I_0(b_\nu)} L^+\ov{P}_{b_\nu} \subset \MC{G}$. The only non-trivial part is to show that $\varphi|_{\MC{L}(\MC{I}_0 \times^{I_0(b_\nu)} L^+\ov{P}_{b_\nu})}$ is well-defined. \\
For this consider some fpqc-cover $S' \to S$ together with an isomorphism $\xi: (\MC{G}, \varphi)_{S'} \cong (L^+G, b_\nu\sigma^*)_{S'}$. Then by the previous remark $\xi(\MC{I}_{0 S'}) = I_0(b_\nu)_{S'}$ and hence $\xi((\MC{I}_0 \times^{I_0(b_\nu)} L^+\ov{P}_{b_\nu})_{S'}) = {L^+\ov{P}_{b_\nu}}_{S'}$ and $\xi(\MC{L}(\MC{I}_0 \times^{I_0(b_\nu)} L^+\ov{P}_{b_\nu})_{S'}) = {L\ov{P}_{b_\nu}}_{S'}$. In particular $\varphi|_{\MC{L}(\MC{I}_0 \times^{I_0(b_\nu)} L^+\ov{P}_{b_\nu})}$ is well-defined if and only if $b_\nu\sigma^*$ is well-defined on ${L\ov{P}_{b_\nu}}_{S'}$. But this last statement is obviously true. \\
b) Recall that complete slope divisions are unique by lemma \ref{lem:ComplSlopeDiviUnique}. Thus we check that the tuple $(P_{b_\nu}, s, \nu', (\ov{\MC{P}}_{b_\nu}, \varphi_{\ov{\MC{P}}_{b_\nu}}), \alpha)$ (where $s > 0$ and $\nu' \in X_*(G)$ satisfy $(b_\nu \sigma)^s = z^{\nu'}$) defines indeed a complete slope division. $\nu'$ is dominant in $G$ wrt. $B_{b_\nu}$ by choice of $B_{b_\nu}$ and the Levi subgroup $M_{b_\nu}$ of $P_{b_\nu}$ is by definition the centralizer of $\nu'$ in $G$. Thus it remains to see that $z^{-\nu'}\varphi_{\ov{\MC{P}}_{b_\nu}}^s$ restricts to an isomorphism of $L^+\ov{P}_{b_\nu}$-torsors. We may check this locally. But locally $\varphi_{\ov{\MC{P}}_{b_\nu}}$ is just given by $b_\nu \sigma^*$ and hence $z^{-\nu'}\varphi_{\ov{\MC{P}}_{b_\nu}}^s$ turns out to be the identity.
\exit

\rem{}{
In particular the $I_0(b_\nu)$-torsor $\MC{I}_0$ already contains all the information about the complete slope division.
}

\subsection{Igusa varieties over basic strata}\label{subsec:DefineIgusaBasic}
In the world of abelian varieties and associated $p$-divisible groups, their $p^m$-torsion define finite flat group schemes, called truncated BT-groups. In particular the sheaf of isomorphisms between two of them is representable by schemes. These moduli schemes are (the building blocks) of Igusa varieties, as described in \cite[section 3]{MantoFoliation} or \cite[section 4]{MantoFoliationPEL}. \\
Our aim in this section is to define similar schemes of isomorphisms for local $G$-shtukas. One way to do this, is to define ``$H_n$-truncated local $G$-shtukas`` (depending on a subgroup $H_n \subset L^+G$) as in remark \ref{rem:TruncDefi}iv) and study their isomorphisms by embedding them into isomorphism schemes between finite schemes. However this approach works only in the case of basic completely slope divisible local $G$-shtukas. \\
To get a moduli description for such schemes over any central leaf, we give an ad-hoc definition of $H_n$-truncated isomorphisms. Although it can be (and is) written down in great generality, we do not claim that it gives sensible objects in general. Rather it is tailored to the situation of completely slope divisible local $G$-shtukas. Whenever these moduli problems are representable, the scheme representing them is called an Igusa variety, analogous to the situation in mixed characteristic. \\
As the existence of Igusa varieties contains a few subtle points, we only deal with the general case in the next section (cf. theorem \ref{thm:TruncRepLeaf}), and stick to the basic case for now.  \vspace{3mm} \\
We fix some semi-standard parabolic subgroup $P \subset G$ and a sequence $\ldots \subset H_{d+1} \subset H_d \subset \ldots \subset H_0 \subset L^+P \subset L^+G$ of open subgroups of $L^+P$ such that each $H_i$ is normal in $H_0$. The prototypical example for such subgroups is
\[H_d \coloneqq K_d \cap L^+P = \{g \in L^+P \; | \; g \equiv 1 \bmod z^d\}.\]
In fact in all applications the subgroups $H_d$ are either $K_d \cap L^+P_{b_\nu}$ or the Iwahori-type subgroups $I_d(b_\nu)$, which were defined in the previous section \ref{subsec:IwahoriStructure}. \\
Note that in all cases the quotients $H_0/H_d$ exist as linear algebraic groups and any element in $H_0/H_d(S)$ (for some scheme $S$) can be lifted to an element in $H_0(S)$ by \cite[corollary 2.2]{HaVi2}.

\defi{}{}{
a) A local $G$-shtuka with $H_0$-structure is a quadruple $(\MC{G}, \MC{H}, \varphi, \beta)$ consisting of a local $G$-shtuka $(\MC{G}, \varphi)$, an $H_0$-torsor $\MC{H}$ and an isomorphism $\beta: \MC{H} \times^{H_0} L^+G \cong \MC{G}$. Most of the time we will abbreviate the notation and simply write $(\MC{H}, \varphi)$ assuming wlog. that $\MC{G} \coloneqq \MC{H} \times^{H_0} L^+G$ and $\beta = \id$. \\
b) A quasi-isogeny between two local $G$-shtukas with $H_0$-structure is a quasi-isogeny between the local $G$-shtukas obtained after forgetting the $H_0$-structure. \\
c) An isomorphism between two local $G$-shtukas with $H_0$-structure is an isomorphism between usual local $G$-shtukas which descends to an isomorphism between the respective $H_0$-torsors. 
}

\rem{}{
Note that we do not require that the Frobenius-isomorphism $\varphi$ comes from an isomorphism of $LP$-torsors associated to the $L^+P$-torsor $\MC{H} \times^{H_0} L^+P$. 
}
$\left. \right.$ \\
Let $S$ be any scheme (or more generally DM-stack) over $\Sp E$ and $\MC{H}$ and $\MC{H}'$ be $H_0$-torsors over $S$. Denote for any integer $d > 0$ by $\MC{H}_d$ respectively $\MC{H}'_d$ the associated $H_0/H_d$-torsor coming from the canonical morphism $H_0 \to H_0/H_d$. We will use the following stacks over $S$:
\begin{description}
 \item[$\quad \Isom(\MC{H}, \MC{H}')$] the stack of isomorphisms $\alpha: \MC{H}_{S'} \to \MC{H}'_{S'}$ of $H_0$-torsors over $S$-schemes $S'$.
 \item[$\quad \Aut^d(\MC{H})$] the substack of $\Aut(\MC{H}) = \Isom(\MC{H}, \MC{H})$ of isomorphisms which induce the identity on the associated $H_0/H_d$-torsor $\MC{H}_d$.
 \item[$\quad \Isom^d(\MC{H}, \MC{H}')$] the stack of isomorphisms $\alpha: \MC{H}_{d \, S'} \to \MC{H}'_{d \, S'}$ of $H_0/H_d$-torsors over $S$-schemes $S'$.
\end{description}
Then $\Aut^d(\MC{H})$ and $\Aut^d(\MC{H}')$ act faithfully on $\Isom(\MC{H}, \MC{H}')$. Note that both actions have the same orbits. Indeed for all $\alpha \in \Isom(\MC{H}, \MC{H}')(S')$ and $\varphi \in \Aut^d(\MC{H})(S')$ we have
\[\alpha \cdot \varphi = \alpha \circ \varphi = (\alpha \circ \varphi \circ \alpha^{-1}) \circ \alpha = (\alpha \circ \varphi \circ \alpha^{-1}) \cdot \alpha \in \Aut^d(\MC{H}')(S') \cdot \alpha.\]
Consider now the canonical morphism
\[\Isom(\MC{H}, \MC{H}') \to \Isom^d(\MC{H}, \MC{H}')\]
It is surjective (as over a trivializing \'etale cover it is simply given by the projection $H_0 \to H_0/H_n$) and has as fibers the $\Aut^d(\MC{H})$-orbits. Hence we have canonical isomorphisms
\[\Isom^d(\MC{H}, \MC{H}') \cong \Aut^d(\MC{H}') \backslash \Isom(\MC{H}, \MC{H}') = \Isom(\MC{H}, \MC{H}') / \Aut^d(\MC{H}) \eqqcolon \Isom(\MC{H}, \MC{H}')[d]\]
(where we view the quotients as quotients of \'etale sheaves). In particular any $\alpha_d \in \Isom^d(\MC{H}, \MC{H}')(S)$ may be represented after an \'etale cover $S' \to S$ by an element $\alpha \in \Isom(\MC{H}, \MC{H}')(S')$ such that for the two projections $pr_1, pr_2: S' \times_S S' \to S'$ we have 
\[pr_1^*(\alpha) \in \Aut^d(\MC{H}')(S' \times_S S') \cdot pr_2^*(\alpha) \cdot \Aut^d(\MC{H})(S' \times_S S').\]

\defi{\label{def:TruncDefiModuli}}{}{
A $H_n$-truncated isomorphism between two local $G$-shtukas with $H_0$-structure $(\MC{H}, \varphi)$ and $(\MC{H}', \varphi')$ over $S$ is an element $\alpha_d \in \Isom^d(\MC{H}, \MC{H}')(S)$ such that there is a representative $\alpha: \MC{H} \to \MC{H}'$ over some pro-\'etale cover $S' \to S$, which is compatible with the Frobenius-isomorphisms, i.e. it satisfies
\[\alpha \circ \varphi \circ \sigma^*\alpha^{-1} \in \Aut^d(\MC{H})(S') \cdot \varphi' \cdot \Aut^d(\sigma^*\MC{H})(S').\]
Then define the functor $\op{Ig}_{(\MC{H}, \varphi), (\MC{H}', \varphi')}^{H_d}$ on the \'etale site $S_{\acute{e}t}$ by
\[\op{Ig}_{(\MC{H}, \varphi), (\MC{H}', \varphi')}^{H_d}(T) \coloneqq \{H_d \op{\!-truncated \; isomorphism \;} \alpha: (\MC{H}, \varphi)_T \to (\MC{H}', \varphi')_T\}\]
for any $(T \to S) \in S_{\acute{e}t}$
}

\rem{\label{rem:TruncDefi}}{
i) We require the compatibility with the Frobenius-isomorphisms only for one representative $\alpha$ of $\alpha_d$. But one easily sees that this implies the compatibility for all representatives over any sufficiently large pro-\'etale cover. \\
ii) If $\alpha_d \in \Isom^d(\MC{H}, \MC{H}')(S)$ is an $H_d$-truncated isomorphism, then so is $\alpha_n^{-1} \in \Isom(\MC{H}'_d, \MC{H}_d)(S)$. Normality of the subgroups $H_d$ gives that compositions of $H_d$-truncated isomorphisms are again $H_d$-truncated. \\
iii) By definition of $H_d$-truncated isomorphisms, $\op{Ig}_{(\MC{H}, \varphi), (\MC{H}', \varphi')}^{H_d}$ is obviously a sheaf for the \'etale topology. \\
iv) Truncated isomorphisms can be seen as isomorphisms of objects, which may be called $H_d$-truncated local $G$-shtukas with $H_0$-structure: Define them to be triples
\[(\MC{H}_d, [\MC{H}], \Aut^d(\MC{H}) \cdot \varphi \cdot \Aut^d(\sigma^*\MC{H}))\]
consisting of an $H_0/H_d$-torsor $\MC{H}_d$ over $S$, an isomorphism class of lifts $[\MC{H}]$ of $\MC{H}_d$ to an $H_0$-torsor and a double coset $\Aut^d(\MC{H}) \cdot \varphi \cdot \Aut^d(\sigma^*\MC{H})$ associated to a Frobenius-isomorphism $\varphi: \sigma^*\MC{LG} \to \MC{LG}$. The last constituent is well-defined as changing the lift $\MC{H}$ within its isomorphism class only changes the Frobenius-isomorphism $\varphi$ by an element in $\Aut^d(\MC{H})$ and the same applies to the choice of $\sigma^*\MC{H}$ within its isomorphism class of liftings of $\sigma^*\MC{H}_d$. Furthermore note that $[\MC{H}]$ may not be determined by $\MC{H}_d$, which can be seen from the cohomological description of deformations as worked out in \cite[proof of proposition 1]{Heinloth2010}) in the case of $H_d = K_d$. \\
The natural notion of an isomorphism 
\[\alpha: (\MC{H}_d, [\MC{H}], \Aut^d(\MC{H}) \cdot \varphi \cdot \Aut^d(\sigma^*\MC{H})) \to (\MC{H}'_d, [\MC{H}'], \Aut^d(\MC{H}') \cdot \varphi' \cdot \Aut^d(\sigma^*\MC{H}'))\] between two $H_n$-truncated local $G$-shtukas with $H_0$-structure would be an element $\alpha \in \Isom(\MC{H}_d, \MC{H}'_d)$ such that \'etale locally any representatives $\tilde{\MC{H}} \in [\MC{H}]$, $\tilde{\MC{H}}' \in [\MC{H}']$, $\tilde{\alpha} \in \Isom(\tilde{\MC{H}}, \tilde{\MC{H}}')$, $\tilde{\varphi} \in \Isom(\sigma^*\MC{L}\tilde{\MC{H}}, \MC{L}\tilde{\MC{H}})$ and $\tilde{\varphi}' \in \Isom(\sigma^*\MC{L}\tilde{\MC{H}}', \MC{L}\tilde{\MC{H}}')$ satisfy after a further pro-\'etale cover
\[\tilde{\alpha}^{-1} \circ \tilde{\varphi}' \circ \sigma^*\tilde{\alpha} \in \Aut^d(\tilde{\MC{H}}) \cdot \tilde{\varphi} \cdot \Aut^d(\sigma^*\tilde{\MC{H}}).\]
This condition is equivalent to requiring that it holds (pro-\'etale locally) for one choice of representatives. Therefore assuming that both truncated local $G$-shtukas came from fixed local $G$-shtukas with $H_0$-structure $(\MC{H}, \varphi)$ and $(\MC{H}', \varphi')$, this implies that the set of isomorphisms between the $H_d$-truncated local $G$-shtukas with $H_0$-structure can be canonically identified with the set of with $H_d$-truncated isomorphisms between $(\MC{H}, \varphi)$ and $(\MC{H}', \varphi')$. \\
v) If we assume that the Frobenius-morphisms are already defined over the integral level, i.e. $\varphi \in \Isom(\sigma^*\MC{H}, \MC{H})(S)$ and $\varphi' \in \Isom(\sigma^*\MC{H}', \MC{H}')(S)$, an $H_d$-truncated isomorphism $\alpha_d: (\MC{H}, \varphi) \to (\MC{H}', \varphi')$ can be described in terms of $H_0/H_d$-torsors: 
Let $\MC{H}_d$ and $\MC{H}'_d$ be the $H_0/H_d$-torsors associated to $\MC{H}$ and $\MC{H}'$, then $\varphi$ and $\varphi'$ induce Frobenius-isomorphisms $\varphi_d$ and $\varphi'_d$ of the $H_0/H_d$-torsors. Thus under the isomorphism $\Isom(\MC{H}, \MC{H}')[d] \cong \Isom(\MC{H}_d, \MC{H}'_d)$ a truncated isomorphism $\alpha$ is nothing else than an element $\alpha_n \in \Isom(\MC{H}_d, \MC{H}'_d)(S)$ satisfying $\alpha_d^{-1} \circ \varphi'_d \circ \sigma^*\alpha_d = \varphi_d$ as an equality of morphisms between $H_0/H_d$-torsors. Note that this last condition is satisfied pro-\'etale locally on $S$ if and only if it is satisfied over $S$, hence in this case we can omit passing to suitable covers.
}

Before actually proving the representability, let us first consider the easiest case, namely when the local $G$-shtukas are already trivial.

\lem{\label{prop:TruncTechnicalLemma}}{}{
Let $b_\nu$ be a fundamental alcove, $(\MC{H}, \varphi)$ a local $G$-shtuka with $I_0(b_\nu)$-structure over a scheme $S$ and 
\[\alpha_d: (\MC{H}, \varphi) \to (I_0(b_\nu), b_\nu\sigma^*)\] 
a $I_d(b_\nu)$-truncated isomorphism which admits a representative $\alpha: \MC{H} \to I_0(b_\nu)$ over an \'etale cover $S' \to S$. Then $\Aut^d(I_0(b_\nu))(S') = I_d(b_\nu)(S')$ and
\[\alpha \circ \varphi \circ \sigma^*\alpha^{-1} \in \id_{I_0(b_\nu)} \cdot b_\nu \sigma^* \cdot I_d(b_\nu)(S')\]
without passing to a further pro-\'etale cover of $S'$.
}

\prooof
$\Aut^d(I_0(b_\nu))(S') = I_d(b_\nu)(S')$ is obvious. As $(b_\nu \sigma^*)^{-1} I_d(b_\nu) (b_\nu \sigma^*) = \phi_{b_\nu}^{-1}(I_d(b_\nu)) \subset I_d(b_\nu)$ we have for each scheme $S''$ the equality 
\begin{align*}
 \Aut^d(I_0(b_\nu))(S'') \cdot b_\nu \sigma^* \cdot \Aut^d(\sigma^*I_0(b_\nu))(S'') & = I_d(b_\nu)(S'') \cdot b_\nu\sigma^* \cdot I_d(b_\nu)(S'') \\
  & = b_\nu \sigma^* \cdot I_d(b_\nu)(S'') \\
  & = b_\nu \sigma^* \cdot \Aut^d(\sigma^*I_0(b_\nu))(S'').
\end{align*}
Thus the element $(b_\nu\sigma^*)^{-1} \circ \alpha \circ \varphi \circ \sigma^*\alpha^{-1}$ defines an element over $S'$ which lies after passing to a suitable pro-\'etale cover in the sheaf $\Aut^d(\sigma^*I_0(b_\nu))$. Hence it lies already in $\Aut^d(\sigma^*I_0(b_\nu))(S')$. \exit

\lem{\label{lem:TruncRepFundamental}}{}{
Let $b_\nu$ be a fundamental alcove. Then for the local $G$-shtuka with $I_0(b_\nu)$-structure $(I_0(b_\nu), b_\nu \sigma^*)$ over $\Sp E$ the functor $\op{Ig}_{(I_0(b_\nu), b_\nu\sigma^*), (I_0(b_\nu), b_\nu\sigma^*)}^{I_d(b_\nu)}$ of $I_d(b_\nu)$-truncated automorphisms is representable by a finite disjoint union of points isomorphic to $\Sp E$. In particular it is a finite \'etale Galois-cover.
}

\prooof
Let $S'$ be any $E$-scheme. Then any element $\alpha_n \in \op{Ig}_{I_0(b_\nu), I_0(b_\nu)}^{I_d(b_\nu)}(S')$ can be represented by an isomorphism $\alpha$ of $I_0(b_\nu)$-torsors over $S'$ (without passing to any cover). Thus by the previous lemma it satisfies $\alpha^{-1} \circ b_\nu \sigma^* \circ \sigma^*\alpha \in b_\nu \sigma^* \cdot I_d(b_\nu)(S')$. Denoting representatives of elements $g \in I_0(b_\nu)/I_d(b_\nu)$ by $\tilde{g} \in I_0(b_\nu)$ we get
\begin{align*}
 \op{Ig}_{(I_0(b_\nu), b_\nu\sigma^*), (I_0(b_\nu), b_\nu\sigma^*)}^{I_d(b_\nu)}(S') & = \{g \in I_0(b_\nu)/I_d(b_\nu) \,(S') \,|\,  (b_\nu\sigma)^{-1} \cdot \tilde{g}^{-1} \cdot (b_\nu\sigma) \cdot \tilde{g} \in I_d(b_\nu)(S') \} \\
   & = \{g \in I_0(b_\nu)/I_d(b_\nu) \,(S') \,|\, \phi_{b_{\nu}}^{-1}(g^{-1}) \cdot g = 1 \in I_0(b_\nu)/I_d(b_\nu)(S')\} \\
   & = \ker(L_{b_\nu}: I_0(b_\nu)/I_d(b_\nu) \to I_0(b_\nu)/I_d(b_\nu)) (S')
\end{align*}
where $L_{b_\nu}: I_0(b_\nu)/I_d(b_\nu) \to I_0(b_\nu)/I_d(b_\nu)$, $g \mapsto \phi_{b_{\nu}}^{-1}(g^{-1}) \cdot g$. Hence $\op{Ig}_{(I_0(b_\nu), b_\nu\sigma^*), (I_0(b_\nu), b_\nu\sigma^*)}^{I_d(b_\nu)}$ is represented by $\ker(L_{b_\nu}: I_0(b_\nu)/I_d(b_\nu) \to I_0(b_\nu)/I_d(b_\nu))$. But $I_0(b_\nu)/I_d(b_\nu)$ is a linear algebraic group and by the Theorem of Lang-Steinberg applied to $L_{b_\nu}$ its kernel is indeed finite \'etale over $E$. \\
Note now that the $\op{Ig}_{(I_0(b_\nu), b_\nu\sigma^*), (I_0(b_\nu), b_\nu\sigma^*)}^{I_d(b_\nu)}$ is a group scheme (by composition of isomorphisms) and consider its identity component. As it is finite \'etale over $\Sp E$, it has to be the spectrum of a field. Furthermore it contains an $E$-valued point, namely the identity. Hence the identity component is isomorphic to $\Sp E$ and by the group scheme structure, $\op{Ig}_{(I_0(b_\nu), b_\nu\sigma^*), (I_0(b_\nu), b_\nu\sigma^*)}^{I_d(b_\nu)}$ is a finite disjoint union of points isomorphic to $\Sp E$. \exit 

\rem{}{
Note in particular that we do not yet assume $b_\nu$ to be basic in this lemma. 
}

\prop{\label{prop:IgusaReprCentral}}{}{
Let $(\MC{G}, \varphi)$ be a local $G$-shtuka over a connected $E$-scheme $S$, which is completely slope divisible with $s = 1$ and has a basic Newton point $\nu$. Then $\op{Ig}_{(\MC{G}, \varphi), (L^+G, b_\nu \sigma^*)}^{K_d}$ is representable by a finite \'etale cover over $S$.
}

\prooof
Being completely slope divisible for $s = 1$ implies that the fundamental alcove $b_\nu = z^\nu$ is central and that we may write $\varphi = z^\nu \tilde{\varphi}$ for some isomorphism $\tilde{\varphi}: \sigma^{*s}\MC{G} \to \MC{G}$. Because $z^\nu$ is central, any $\alpha: (\MC{G}, \varphi)_{S'} \to (L^+G, b_\nu \sigma^*)_{S'}$ (over some scheme $S' \to S$) also defines an isomorphism $\alpha: (\MC{G}, \tilde{\varphi})_{S'} \to (L^+G, \sigma^{*s})_{S'}$ of local $G$-shtukas and conversely. Hence we get an isomorphism
\[\op{Ig}_{(\MC{G}, \varphi), (L^+G, b_\nu \sigma^*)_S}^{K_d} \cong \op{Ig}_{(\MC{G}, \tilde{\varphi}), (L^+G, \sigma^*)_S}^{K_d}\]
Now we are in the \'etale setting and by remark \ref{rem:TruncDefi}v) we have an isomorphism
\[\op{Ig}_{(\MC{G}, \tilde{\varphi}), (L^+G, \sigma^*)_S}^{K_d} \cong \Aut_S((\MC{G}_d, \tilde{\varphi}), (L^+G/K_d, \sigma^*)_S)\]
Furthermore the discussion in section \ref{subsec:TateFunctor} yields an isomorphism
\[\Aut_S((\MC{G}_d, \tilde{\varphi}), (L^+G/K_d, \sigma^*)_S) \cong \Aut_S(\MC{G}_d^{\tilde{\varphi}}, (L^+G/K_d)_S^{\sigma^*})\]
where the right-hand side is representable as an isomorphism scheme between two finite schemes. \\
We now show that it is a finite \'etale cover with Galois group $\Aut_E((L^+G, b_\nu \sigma^*), (L^+G, b_\nu \sigma^*))$. It suffices to check this fpqc-locally, i.e. on completions of $S$ along geometric points $x$. By Cohen's structure theorem on complete local rings (cf. \cite[theorem 9]{CohenStructureThm}) we may endow $S^{\wedge x}$ with a structure of a $\ACFq$-scheme. Thus by \cite[proposition 8.1]{HaVi} we may trivialize the restriction of $(\MC{G}, \varphi)$ to $S^{\wedge x}$, hence get an isomorphism between $\op{Ig}_{(\MC{G}, \varphi), (L^+G, b_\nu \sigma^*)}^{K_d}$ and $\op{Ig}_{(L^+G, b_\nu \sigma^*), (L^+G, b_\nu \sigma^*)}^{K_d}$ over it. $b_\nu$ being central implies $I_0(b_\nu) = L^+G$ and $I_d(b_\nu) = K_d$, hence finite \'etaleness follows from lemma \ref{lem:TruncRepFundamental}. \exit

\prop{\label{prop:IgusaReprBasic}}{}{
Let $(\MC{G}, \varphi)$ be a completely slope divisible local $G$-shtuka over a connected $E$-scheme $S$, with basic Newton point $\nu$. Let $b_\nu$ be a fundamental alcove. \\
a) $(\MC{G}, \varphi)$ admits a canonical $I_0(b_\nu)$-structure $\MC{I}_0$. \\
b) $\op{Ig}_{(\MC{I}_0, \varphi), (I_0(b_\nu), b_\nu \sigma^*)}^{I_d(b_\nu)}$ is representable by a finite \'etale cover over $S$.
}

\prooof
a) By proposition \ref{prop:I_0Structure} it suffices to show that the local $G$-shtuka is strongly completely slope divisible, i.e. admits a trivialization after a fpqc-cover. Assume wlog. that $\MC{G} = L^+G$ is trivial. We first prove that $(\MC{G}, \varphi^s)$ admits such a cover, where $s$ is the integer defined by the complete slope division. $(\MC{G}, \varphi^s)$ satisfies the properties of the previous proposition (using the $q^s$-Frobenius for $\sigma$), hence admits finite \'etale covers by Igusa varieties. As they parameterize isomorphisms $\MC{G}_d \to L^+G/K_d$ respecting $\tilde{\varphi^s}$, their inverse limit $\tilde{S} \coloneqq \varprojlim_d \op{Ig}_{(\MC{G}, \varphi^s), (L^+G, (b_\nu \sigma^*)^s)_S}^{K_d}$ parameterizes honest isomorphisms $(\MC{G}, \varphi^s) \to (L^+G, (b_\nu \sigma^*)^s)_S$. It is a pro-finite \'etale cover, hence fpqc. \\
Call the universal isomorphism $\alpha: (\MC{G}, \varphi^s)_{\tilde{S}} \to (L^+G, (b_\nu \sigma^*)^s)_{\tilde{S}}$. Then $\alpha \varphi \sigma^* \alpha^{-1}: \sigma^*L^+G \to L^+G$ satisfies by construction $(\alpha \varphi \sigma^* \alpha^{-1})^s = (b_\nu \sigma^*)^s$. However the elements $b' \in L^+G$ satisfying $(b'\sigma^*)^s = (b_\nu \sigma^*)^s$ are discrete, because they satisfy $b'\sigma^* (b_\nu\sigma^*)^s = (b' \sigma^*)^{s+1} = (b_\nu\sigma^*)^s \cdot \sigma^{* s+1}b'$. In particular the locus $\tilde{S}^0$ in $\tilde{S}$ where $\alpha$ defines an isomorphism $(\MC{G}, \varphi) \to (L^+G, b_\nu \sigma^*)$ is open and closed in $\tilde{S}$. In particular $\tilde{S}^0$ is still fpqc over $S$. It remains to see surjectivity of $\tilde{S}^0 \to S$. But over geometric points, there exists some isomorphism $\alpha: (\MC{G}, \varphi) \to (L^+G, b_\nu \sigma^*)$, giving a point in $\tilde{S}$ by its moduli interpretation, which then lies in the subset $\tilde{S}^0$. \\
This proves the existence of a canonical $I_0(b_\nu)$-structure $\MC{I}_0$ on $(\MC{G}, \varphi)$. \\
b) We proceed in a very similar way: There is a canonical embedding
\[\op{Ig}_{(\MC{I}_0, \varphi), (I_0(b_\nu), b_\nu \sigma^*)}^{I_d(b_\nu)} \subset \Isom^d(\MC{I}_0, I_0(b_\nu)) \times_{\Isom^d(\MC{G}, L^+G)} \op{Ig}_{(\MC{G}, \varphi^s), (L^+G, (b_\nu \sigma^*)^s)_S}^{K_d}.\]
On the first factor it forgets the compatibility with $\varphi$, and on the second factor it forgets the compatibility with the $I_0(b_\nu)$-structures and remembers only the compatibility with $\varphi^s$. We claim that this is a closed immersion. Indeed locally the right-hand side is given by $\{g \in I_0(b_\nu)/I_d(b_\nu) \,|\, (g \varphi \sigma^*g^{-1})^s \in K_d (b_\nu \sigma^*)^s K_d\}$, while the left-hand side is $\{g \in I_0(b_\nu)/I_d(b_\nu) \,|\, g \varphi \sigma^*g^{-1} \in I_d(b_\nu) b_\nu I_d(b_\nu)\}$. But the double coset $I_d(b_\nu) b_\nu I_d(b_\nu)$ is closed in $L_s^{-1}(K_d \cdot (b_\nu \sigma^*)^s \cdot K_d)$, where $L_s: L^+G \to L^+G$, $g \mapsto (g\sigma^*)^s$. So the Igusa variety $\op{Ig}_{(\MC{I}_0, \varphi), (I_0(b_\nu), b_\nu \sigma^*)}^{I_d(b_\nu)}$ is indeed representable as a closed immersion in the fiber product above. \\
The proof that it is finite \'etale is verbatim the same as in the previous proposition. \exit

\subsection{General Igusa varieties over perfect central leaves}\label{subsec:DefineIgusaGeneral}
Unfortunately Igusa varieties over non-basic central leaves do not exist as finite \'etale covers, as the existence of such Igusa varieties implies the triviality of any deformation of the associated local $G$-shtukas. However this problem vanishes over perfect base schemes, because there the local $G$-shtuka is induced from a basic local $M$-shtuka for some Levi subgroup $M$, as shown in proposition \ref{prop:SplitShtukaPerfect}. This allows us to define Igusa varieties over the perfection of arbitrary central leaves in theorem \ref{thm:TruncRepLeaf}.

\lem{\label{lem:SplitShtukaFrobenius}}{}{
Let $(\MC{G}, \varphi)$ be a completely slope divisible local $G$-shtuka over a quasi-compact scheme $S$, with constant Newton point $\nu$. Denote its basic constituent by $(\MC{M}, \varphi_{\MC{M}})$. Then for all $d \geq 0$ there exists some $N_d \gg 0$ such that there is a $K_d$-truncated isomorphism
\[\alpha_d: \sigma^{*N_d}(\MC{G}, \varphi) \to \sigma^{*N_d}(\MC{M}, \varphi_{\MC{M}}).\]
Here $\sigma^{*N_d}(\MC{M}, \varphi_{\MC{M}})$ denotes the local $G$-shtuka with $L^+G$-torsor $\sigma^{*N_d}\MC{M} \times^{\prod L^+M_i} L^+G$ and the Frobenius-isomorphism induced by $\sigma^{*N_d}\varphi_{\MC{M}}: \sigma^{*N_d+1}\MC{LM} \to \sigma^{*N_d}\MC{LM}$.
}

\prooof
Using the canonical $L^+\ov{P}$-structure on completely slope divisible local $G$-shtukas, it suffices to prove the statement only for local $\ov{P}$-shtukas $(\ov{\MC{P}}, \varphi)$. So choose a trivialization $(\ov{\MC{P}}, \varphi) \cong (L^+\ov{P}, \varphi)$ over some \'etale cover. Then $\varphi_{\MC{M}} \in LM$ and by construction $\varphi \cdot \varphi_{\MC{M}}^{-1} \in LR_u\ov{P}$ lies in the unipotent radical. 
We claim that for sufficiently large $N$, $z^{-N \nu'} \varphi \varphi_{\MC{M}}^{-1} z^{N \nu'}$ lies in $LR_u(\ov{P}) \cap K_d$. It suffices to check this on root subgroups $U_\alpha \in \ov{P}$. There conjugation with $z^{-\nu'}$ equals multiplication with $z^{\langle \alpha, -\nu'\rangle}$, where $\langle \alpha, -\nu'\rangle > 0$. Hence any element in $LU_\alpha$ can be conjugated by $z^{-N \nu'}$ into $LU_\alpha \cap K_d$ for sufficiently large $N$. \\
Moreover $z^{-N \nu'} \varphi \varphi_{\MC{M}}^{-1} z^{N \nu'} = z^{-N \nu'} \varphi z^{N \nu'} \cdot \varphi_{\MC{M}}^{-1}$, because $\varphi_{\MC{M}}^{-1}$ lies in $LM$, which is the centralizer of $\nu'$. So the claim above implies that $id: (\ov{\MC{P}}, z^{-N \nu'} \varphi z^{N \nu'}) \to (\MC{M}, \varphi_{\MC{M}})$ is a $K_d$-truncated isomorphism. \\
But then we can just set $N_d = N \cdot s$ and define $\alpha_d$ as
\[\alpha_d: \sigma^{*N_d}(\ov{\MC{P}}, \varphi) \xrightarrow{z^{-N\nu'} \varphi^{Ns}} (\ov{\MC{P}}, z^{-N\nu'} \varphi z^{N\nu'}) \xrightarrow{id} (\MC{M}, \varphi_{\MC{M}}) \xrightarrow{z^{N\nu'} \varphi_{\MC{M}}^{-Ns}} \sigma^{*N_d}(\MC{M}, \varphi_{\MC{M}})\]
because $z^{-N\nu'} \varphi^{Ns}$ is an isomorphism and so are its component $(z^{-N\nu'} \varphi^{Ns})_{\MC{M}} = z^{N\nu'} \varphi_{\MC{M}}^{-Ns}$. As our definition of $\alpha_d$ does not depend on choices made locally, i.e. only on the choice of $N_d$, all these $\alpha_d$ glue together to a $K_d$-truncated isomorphism already over $S$. \exit  

\rem{}{
Analogous to \cite[lemma 4.1]{MantoFoliation}, the constant $N_d$ can be made explicit, though we will never need this. In fact, if the local $G$-shtuka is bounded by some $\mu$, then $N_d$ can be chosen to depend only on $d$ and $\mu$.
}

\lem{\label{lem:SplitShtukaCompatible}}{}{
The $K_d$-truncated isomorphisms $\alpha_d$ defined in the previous lemma do not depend on the choice of $N_d$ in the following sense: The two $K_d$-truncated isomorphisms
\[\sigma^{*s}\alpha_d = \left(z^{N\nu'} \varphi_{\MC{M}}^{-Ns}\right) \circ \left(z^{-N\nu'} \varphi^{Ns}\right): \sigma^{*(N+1)s}(\MC{G}, \varphi) \to \sigma^{*(N+1)s}(\MC{M}, \varphi_{\MC{M}}) \]
and 
\[\alpha_d' = \left(z^{(N+1)\nu'} \varphi_{\MC{M}}^{-(N+1)s}\right) \circ \left(z^{-(N+1)\nu'} \varphi^{(N+1)s}\right): \sigma^{*(N+1)s}(\MC{G}, \varphi) \to \sigma^{*(N+1)s}(\MC{M}, \varphi_{\MC{M}}) \]
coincide.
}

\prooof
Using that $z^{\nu'}$ commutes with $\varphi_{\MC{M}}$, one computes:
\begin{align*}
  \alpha_d' \circ \sigma^{*s}\alpha_d^{-1} & = \left(z^{(N+1)\nu'} \varphi_{\MC{M}}^{-(N+1)s}\right) \circ \left(z^{-(N+1)\nu'} \varphi^{(N+1)s}\right) \circ \left(z^{-N\nu'} \varphi^{Ns}\right)^{-1} \circ \left(z^{N\nu'} \varphi_{\MC{M}}^{-Ns}\right)^{-1} \\
 & = \left(z^{(N+1)\nu'} \varphi_{\MC{M}}^{-(N+1)s}\right) \circ \left(z^{-(N+1)\nu'} \varphi^s z^{N\nu'}\right) \circ \left(z^{N\nu'} \varphi_{\MC{M}}^{-Ns}\right)^{-1} \\
 & = \left(z^{N\nu'} \varphi_{\MC{M}}^{-Ns}\right) \circ \left(z^{\nu'} \varphi_{\MC{M}}^{-s}\right) \circ \left(z^{-(N+1)\nu'} \varphi^s z^{N\nu'}\right) \circ \left(z^{N\nu'} \varphi_{\MC{M}}^{-Ns}\right)^{-1} \\
 & = \left(z^{N\nu'} \varphi_{\MC{M}}^{-Ns}\right) \circ \left(z^{-N\nu'} \varphi_{\MC{M}}^{-s} \varphi^s z^{N\nu'}\right) \circ \left(z^{N\nu'} \varphi_{\MC{M}}^{-Ns}\right)^{-1} 
\end{align*}
But $\varphi_{\MC{M}}^{-s} \varphi^s$ is an element in $L^+R_u\ov{P}$, so after conjugation with $z^{N\nu'}$ one obtains an element in $K_d$. Then further conjugation with the isomorphism $z^{-N\nu'} \varphi_{\MC{M}}^{-Ns}$ preserves this property. Hence $\alpha_d' \circ \sigma^{*s}\alpha_d^{-1}$ is the identity in the category of $K_d$-truncated isomorphisms, i.e. $\alpha_d' = \sigma^{*s}\alpha_d$  \exit

\defi{}{}{
A scheme $S$ over $\Fq$ is perfect, if the absolute Frobenius is surjective on $S$. 
The perfection $S^\sharp$ of a scheme $S$ over $\Fq$ is given by $\varprojlim_{\sigma} S$.
}

\rem{}{
Basic properties of the perfection functor are given in appendix $A$ of \cite{ZhuAffineGrass}.
}

\prop{\label{prop:SplitShtukaPerfect}}{}{
Let $(\MC{G}, \varphi)$ be a completely slope divisible local $G$-shtuka over a perfect scheme $S$, with constant Newton point $\nu$. Denote its basic constituent by $(\MC{M}, \varphi_{\MC{M}})$. Then there is a canonical isomorphism
\[\alpha: (\MC{G}, \varphi) \to (\MC{M}, \varphi_{\MC{M}}).\]
}

\prooof
As we are over a perfect scheme, the Frobenius defines a canonical isomorphism $(\MC{G}, \varphi) \cong (\sigma^*\MC{G}, \sigma^*\varphi)$. Thus we may consider the sequence $(\sigma^{*-N_d}\alpha_d)_d$ of truncated isomorphisms from $(\MC{G}, \varphi)$ to $(\MC{M}, \varphi_{\MC{M}})$, where the $\alpha_d$ are the $K_d$-truncated isomorphisms defined in lemma \ref{lem:SplitShtukaFrobenius}. By lemma \ref{lem:SplitShtukaCompatible}, $\sigma^{*-N_{d+1}}\alpha_{d+1}$ coincides with $\sigma^{*-N_d}\alpha_d$ when both are considered as $K_d$-truncated isomorphisms. Hence their limit defines a honest isomorphism $\alpha: (\MC{G}, \varphi) \to (\MC{M}, \varphi_{\MC{M}})$ as claimed. \exit

\prop{\label{prop:IgusaRepComplSlope}}{}{
Let $(\MC{G}, \varphi)$ be a completely slope divisible local $G$-shtuka over a perfect $E$-scheme $S$, with Newton point $\nu$. Let $b_\nu$ be a fundamental alcove. \\
a) $(\MC{G}, \varphi)$ admits a canonical $I_0(b_\nu)$-structure $\MC{I}_0$. \\
b) $\op{Ig}_{(\MC{I}_0, \varphi), (I_0(b_\nu), b_\nu \sigma^*)}^{I_d(b_\nu)}$ is representable by a finite \'etale cover of $S$.
}

\prooof
Over the perfect scheme $S$, choose an isomorphism $(\MC{G}, \varphi) \cong (\MC{M}, \varphi_{\MC{M}})$ to its basic constituent. Then the basic constituent can be trivialized after a pro-finite \'etale cover by the existence of Igusa varieties in the basic case. Hence $(\MC{G}, \varphi)$ is strongly completely slope divisible and by proposition \ref{prop:I_0Structure} admits an $I_0(b_\nu)$-structure $\MC{I}_0$. \\
For the second part note that it suffices to prove this after an \'etale cover. By proposition \ref{prop:IgusaReprBasic}, the Igusa variety for the basic constituent $(\MC{M}, \varphi_{\MC{M}})$ defines a finite \'etale cover over which there exist truncated isomorphisms $(\MC{M}, \varphi_{\MC{M}}) \to (L^+M, b_{\nu}\sigma^*)$, respecting the Iwahori structures for $M$ on both sides. Then 
\[(\MC{G}, \varphi) \cong (\MC{M}, \varphi_{\MC{M}}) \to (L^+M, b_{\nu} \sigma^*) \cong (L^+G, b_\nu \sigma^*) \]
is an $I_d(b_\nu)$-truncated isomorphism, because $I_d(b_\nu)$ contains the the Iwahori-type subgroup for the Levi subgroup $M$.
Then composition with this truncated isomorphism defines an isomorphism 
\[\op{Ig}_{(\MC{I}_0, \varphi), (I_0(b_\nu), b_\nu \sigma^*)}^{I_d(b_\nu)} \cong \op{Ig}_{(I_0(b_\nu), b_\nu \sigma^*), (I_0(b_\nu), b_\nu \sigma^*)}^{I_d(b_\nu)},\]
where the right-hand side is representable by lemma \ref{lem:TruncRepFundamental}. It is finite \'etale by the very same lemma. \exit

\rem{}{
Actually one can show that $\op{Ig}_{(\MC{I}_0, \varphi), (I_0(b_\nu), b_\nu \sigma^*)}^{I_d(b_\nu)}$ equals the Igusa varieties of the basic constituent. The main difficulty to prove this is to show that any $I_d(b_\nu)$-truncated isomorphism $\alpha$ of local $L^+G$-shtukas has a representative respecting the basic constituent, which is unique up to changes in the Iwahori-type subgroups $I_d(b_{\nu i})$ of the Levi $M$. 
One can see this in the following way: By proposition \ref{prop:IgusaTowerLift} one may lift $\alpha$ after a pro-finite \'etale cover to an actual isomorphism. This isomorphism respects then the basic constituent and any two such lifts differ only by an element in the Iwahori-type sungroup of $M$. This defines the desired representative after a pro-finite \'etale cover, which by uniqueness descends to the representative over the original basis.
}

\thm{\label{thm:TruncRepLeaf}}{}{
The universal local $G$-shtuka $(\MC{G}^{univ}_{c_i, U}, \varphi^{univ}_{c_i, U})$ over the perfection $\MC{C}_U^{(\nu_i) \sharp}$ of the central leaf admits a canonical $I_0(b_\nu)$-structure $\MC{I}^{univ}_{c_i, U}$. \\
Moreover 
\[\op{Ig}_{c_i, U}^{d \sharp} \coloneqq \op{Ig}_{(\MC{I}^{univ}_{c_i, U}, \varphi^{univ}_{c_i, U}), (I_0(b_{\nu_i}), b_{\nu_i} \sigma^*)_{\MC{C}_U^{\sharp}}}^{I_d(b_{\nu_i})}\]
is relatively representable by a finite \'etale cover of $\MC{C}_U^{(\nu_i) \sharp}$. It is called the Igusa variety of level $d$ associated to the characteristic point $c_i$.
}

\prooof
By theorem \ref{thm:UnivComplSlopeDiv} the universal local $G$-shtuka over the central leaf $\MC{C}_U^{(\nu_i)}$ is completely slope divisible, hence the same holds after pullback to the perfection. Thus all properties follow directly from the previous proposition. \exit

\cor{\label{cor:IgusaDefinitionNonPerfect}}{}{
There exists a finite \'etale cover $\op{Ig}_{c_i, U}^d$ over $\MC{C}_U^{(\nu_i)}$, whose perfection is $\op{Ig}_{c_i, U}^{d \sharp}$.
}

\prooof
This follows immediately from the general result \cite[proposition A.4]{ZhuAffineGrass} on perfections. \exit

\warn{}{
$\op{Ig}_{c_i, U}^d$ has no longer any moduli description. In particular it does not represent $I_d(b_\nu)$-truncated isomorphisms $(\MC{I}^{univ}_{c_i, U}, \varphi^{univ}_{c_i, U}) \to (I_0(b_{\nu_i}), b_{\nu_i} \sigma^*)_{\MC{C}_U}$ over the central leaf.
}

\lem{}{}{
Let $L_{b_{\nu_i}}: I_0(b_{\nu_i}) \to I_0(b_{\nu_i})$, $g \mapsto \phi_{b_{\nu}}^{-1}(g^{-1}) \cdot g$. Then for every $d \geq 0$ the short exact sequence of linear algebraic groups over $E$
\[0 \to I_d(b_{\nu_i})/I_{d+1}(b_{\nu_i}) \to I_0(b_{\nu_i})/I_{d+1}(b_{\nu_i}) \to I_0(b_{\nu_i})/I_d(b_{\nu_i}) \to 0\]
induces a short exact sequence of finite group schemes over $\Sp E$
\[0 \to \ker(L_{b_{\nu_i}} \op{on} I_d(b_{\nu_i})/I_{d+1}(b_{\nu_i})) \to \ker(L_{b_{\nu_i}} \op{on} I_0(b_{\nu_i})/I_{d+1}(b_{\nu_i})) \to \ker(L_{b_{\nu_i}} \op{on} I_0(b_{\nu_i})/I_d(b_{\nu_i})) \to 0.\]
}

\prooof
The assertion follows directly from the snake lemma using the surjectivity of $L_{b_{\nu_i}}$ on the group $I_d(b_{\nu_i})/I_{d+1}(b_{\nu_i})$. For this note that while the category of all group schemes is not abelian, the proof of the snake lemma still works in this situation by normality of the subgroups in consideration. \exit

\prop{\label{prop:IgusaTowerLift}}{}{
For every $d \geq 0$ the canonical morphism
\[\op{Ig}^{d+1 \sharp}_{c_i, U} \to \op{Ig}^{d \sharp}_{c_i, U}\]
is a finite \'etale cover with Galois group $\ker(L_{b_{\nu_i}} \,\op{on}\, I_d(b_{\nu_i})/I_{d+1}(b_{\nu_i}))$. 
}

\prooof
Any $I_{d+1}(b_{\nu_i})$-truncated isomorphism defines an $I_d(b_{\nu_i})$-truncated isomorphism: View it as a morphism between $I_0(b_{\nu_i})/I_{d+1}(b_{\nu_i})$-torsors and take the induced morphism between $I_0(b_{\nu_i})/I_d(b_{\nu_i})$-torsors. This defines a canonical morphism
\[\op{Ig}^{d+1 \sharp}_{c_i, U} \to \op{Ig}^{d \sharp}_{c_i, U}.\]
This morphism is compatible with the projection to the central leaf, which is finite \'etale by theorem \ref{thm:TruncRepLeaf}. Hence the morphism between Igusa varieties is finite \'etale as well and it suffices to check surjectivity. As all stacks are reduced, we may do so over geometric points of $\MC{C}_U^{\sharp}$. There both local $G_{c_i}$-shtukas are trivialized and by lemma \ref{lem:TruncRepFundamental} the Igusa varieties are given by $\ker(L_{b_{\nu_i}} \,\op{on}\, I_0(b_{\nu_i})/I_{d+1}(b_{\nu_i}))$ respectively $\ker(L_{b_{\nu_i}} \,\op{on}\, I_0(b_{\nu_i})/I_d(b_{\nu_i})$. Hence the transition map is surjective by the previous lemma. From there we get the description of the Galois-group as well. \exit
$\left. \right.$ \vspace{3mm} \\
Thus we get a tower $(\op{Ig}^{d \sharp}_{c_i, U})$ of Igusa varieties with finite \'etale transition maps. In particular $\op{Ig}^{\infty \sharp}_{c_i, U} = \varprojlim_d \op{Ig}^{d \sharp}_{c_i, U}$ exists as a DM-stack, although not of finite type over $\MC{C}_U^{(\nu_i) \sharp}$. 

\prop{\label{prop:IgusaTowerModuli}}{}{
$\op{Ig}^{\infty \sharp}_{c_i, U}$ is the moduli space representing the functor
\[S \mapsto \{\op{isomorphisms \,} \alpha: (\MC{G}^{univ}_{c_i, U}, \varphi^{univ}_{c_i, U}) \times_{\MC{C}_U^{(\nu_i) \sharp}} S \to (L^+G_{c_i}, b_{\nu_i}\sigma^*) \times_E S \}\]
on the category of schemes $S$ over $\MC{C}_U^{(\nu_i) \sharp}$.
}

\prooof
In other words, we have to check that on schemes $S$ over $\MC{C}_U^{(\nu_i) \sharp}$, $\op{Ig}^{\infty \sharp}_{c_i, U}(S)$ identifies with the set of isomorphisms defined in the statement. \\
Assume we have an isomorphism $\alpha: (\MC{G}^{univ}_{c_i, U}, \varphi^{univ}_{c_i, U}) \times_{\MC{C}_U^{(\nu_i) \sharp}} S \to (L^+G_{c_i}, b_{\nu_i}\sigma^*) \times_E S$. Then by remark \ref{rem:I_0Canonical}, it respects the $I_0(b_{\nu_i})$-structures on both sides and hence gives for each $d \geq 0$ an $I_d(b_{\nu_i})$-truncated isomorphism $\alpha_d \in \op{Ig}_{(\MC{I}^{univ}_{c_i, U}, \varphi^{univ}_{c_i, U}), (I_0(b_{\nu_i}), b_{\nu_i} \sigma^*)}^{I_d(b_{\nu_i})}(S) = \op{Ig}^{d \sharp}_{c_i, U}(S)$. Obviously they define a point in $\op{Ig}^{\infty \sharp}_{c_i, U}(S) = \varprojlim_d \op{Ig}^{d \sharp}_{c_i, U}(S)$. \\
Conversely assume we are given an element in $\op{Ig}^{\infty \sharp}_{c_i, U}(S)$, i.e. a system of $I_d(b_{\nu_i})$-truncated isomorphisms $\alpha_d$. We first claim that these glue to a morphism between $L^+G_{c_i}$-torsors: By definition of truncated isomorphisms we may view 
\[\alpha_d \in \op{Isom}((\MC{I}^{univ}_{c_i, U} \times_{\MC{C}_U^{(\nu_i) \sharp}} S) \times^{I_0(b_{\nu_i})} I_0(b_{\nu_i})/I_d(b_{\nu_i}), I_0(b_{\nu_i})/I_d(b_{\nu_i}) \times_E S)\]
as an isomorphism between $I_0(b_{\nu_i})/I_d(b_{\nu_i})$-torsors. To describe the gluing procedure precisely, note that $\alpha_d$ also gives an isomorphism between $I_0(b_{\nu_i})/K_d$-torsors and then on induced $L^+G_{c_i}/K_d$-torsors: 
\[\alpha_d: (\MC{G}^{univ}_{c_i, U} \times_{\MC{C}_U^{(\nu_i) \sharp}} S) \times^{L^+G_{c_i}} L^+G_{c_i}/K_d \to L^+G_{c_i}/K_d \times_{\Fq} S.\]
This can be viewed as an isomorphism of corresponding $G_{c_i}$-torsors over $S \times_E \Sp E[[z]]/(z^d)$. There it is obvious that they glue to a morphism of $G$-torsors $\alpha_\infty = \varprojlim_d \alpha_d$ over $S \times_E \Spf E[[z]]$. But translating backwards, this is nothing else than a map 
\[\alpha_\infty: \MC{G}^{univ}_{c_i, U} \times_{\MC{C}_U^{(\nu_i) \sharp}} S \to L^+G_{c_i} \times_{\Fq} S.\]
Next we check compatibility with Frobenius-isomorphisms: By assumption (together with remark \ref{rem:TruncDefi}ii)) there is for each $d$ a pro-\'etale cover $S'_d \to S$ such that
\[\alpha_\infty^{-1} \cdot \varphi^{univ}_{c_i, U} \cdot \sigma^*\alpha_\infty \in I_d(b_{\nu_i})(S'_d) \cdot b_{\nu_i}\sigma^* \cdot I_d(b_{\nu_i})(S'_d)\]
using $\alpha_\infty$ as a representative of $\alpha_d$. But by lemma \ref{prop:TruncTechnicalLemma}, this holds already for $S'_d = S$. Thus using
\[\bigcap_{d \geq 0} I_d(b_{\nu_i})(S) \cdot b_{\nu_i}\sigma^* \cdot I_d(b_{\nu_i})(S) = \bigcap_{d \geq 0} b_{\nu_i}\sigma^* \cdot I_d(b_{\nu_i})(S) = b_{\nu_i}\sigma^*\]
we get indeed $\alpha_\infty^{-1} \cdot \varphi^{univ}_{c_i, U} \cdot \sigma^*\alpha_\infty = b_{\nu_i}\sigma^*$. \exit

\defi{\label{def:IgusaGlobal}}{}{
Fix for every characteristic place $c_i$ a positive integer $d_i$. Then the global Igusa variety is
\[\op{Ig}^{(d_i) \sharp}_U \coloneqq \op{Ig}^{d_1 \sharp}_{c_1, U} \times_{\MC{C}_U^{(\nu_i) \sharp}} \op{Ig}^{d_2 \sharp}_{c_2, U} \times_{\MC{C}_U^{(\nu_i) \sharp}} \ldots \times_{\MC{C}_U^{(\nu_i) \sharp}} \op{Ig}^{d_n \sharp}_{c_n, U}.\]
}

As a direct consequence of the statements above, $\op{Ig}^{(d_i) \sharp}_U \to \MC{C}_U^{(\nu_i) \sharp}$ is a finite \'etale cover. Furthermore the limit for growing tuples $(d_i)$ represents the functor 
\[S \mapsto \{\op{for \; each \,} c_i \op{\, an \; isomorphism \,} \alpha_i: (\MC{G}^{univ}_{c_i, U}, \varphi^{univ}_{c_i, U}) \times_{\MC{C}_U^{(\nu_i) \sharp}} S \to (L^+G_{c_i}, b_{\nu_i}\sigma^*) \times_E S \}\]
on the category of schemes over $\MC{C}_U^{(\nu_i) \sharp}$.

\section{A covering of the Newton strata}\label{sec:FiniteCover}
We construct now a finite morphism
\[\pi_{(d_i)}: \prod_i \MB{M}_{b_{\nu_i}}^{\preceq \mu_i, \theta_i \sharp} \times \op{Ig}^{(d_i) \sharp}_U \to \MC{N}^{(\nu_i) \sharp}_U\]
which turns into an \'etale morphism when restricted to an open subset of the source. Here $\MB{M}_{b_{\nu_i}}^{\preceq \mu_i, \theta_i \sharp}$ are perfections of closed subsets of (reduced subschemes underlying) Rapoport-Zink spaces and $\op{Ig}^{(d_i) \sharp}_U$ is the global Igusa variety defined in \ref{def:IgusaGlobal}. \\ 
We briefly review the main idea behind this construction: \\
Points in the global Igusa variety come automatically with a global $G$-shtuka $(\MS{G}, \varphi, \psi)$. This global $G$-shtuka will be modified at each characteristic place $c_i$: A lift of the truncated isomorphism coming from the Igusa variety and the quasi-isogeny defined by the Rapoport-Zink space gives a quasi-isogeny starting with the associated local $G$-shtuka $\MF{L}_{c_i}(\MS{G}, \varphi)$. This allows us to change the global $G$-shtuka by keeping all the structure away from $c_i$ the same and replacing it at a formal neighborhood of $c_i$ according to the quasi-isogeny. This defines a new global $G$-shtuka, giving the image point in $\MC{N}^{(\nu_i) \sharp}_{U}$. \\
In section \ref{subsec:Uniform} we present the general construction of the modification of global $G$-shtukas. This relies heavily on the work of Hartl and Rad \cite{HarRad1}, who call the resulting morphism (in a slightly different setup) the uniformization morphism. In section \ref{subsec:CoverExistence} we apply this construction in our setting to get the morphism above by checking that it is indeed well-defined. \\
Following the approach of Mantovan in the case of mixed characteristic, we show that $\dot{\pi}_{(d_i)}$ is quasi-finite in proposition \ref{prop:CoverQuasiFinite} and satisfies the valuation criterion for properness in proposition \ref{prop:CoverProper}. This almost implies theorem \ref{thm:CoveringMorphismProperties}, stating that $\pi_{(d_i)}$ is a finite morphism and (for sufficiently large $\theta_i$) surjective. Moreover we show that restricting $\pi_{(d_i)}$ to an open subspace of the source turns it into an \'etale morphism $\dot{\pi}_{(d_i)}$.

\subsection{The uniformization morphism}\label{subsec:Uniform}
In this section we define the uniformization morphism: Let (for simplicity during this introduction) $(\MS{G}_0, \varphi_0, \psi_0) \in \nabla_{(c_i)}^{\mmu}\MC{H}^1_U(C, G)(S_0)$ be a global $G$-shtuka over an $E$-scheme $S_0$ with characteristic places $(c_i)$ and bounded by $\mmu = (\mu_i)_i$. Let $(\MC{G}_{c_i}, \varphi_{c_i}) \coloneqq \MF{L}_{c_i}(\MS{G}_0, \varphi_0)$ be the associated local $G$-shtukas. Assume that they are constant and equal to a decent local $G$-shtuka, i.e. that we are given an identification $(\MC{G}_{c_i}, \varphi_{c_i}) \cong ({L^+G_{c_i}}_E, b_{c_i}\sigma^*)_{S_0}$, with $({L^+G_{c_i}}_E, b_{c_i}\sigma^*)$ being a decent local $G$-shtuka over $\Sp E$ as in definition \ref{def:RapoZinkDecent}. Then one has a canonical morphism of DM-stacks over $\Sp E$
\[S_0 \times_{E} \prod\nolimits_i \MB{M}_{b_{c_i}}^{\preceq \mu_i} \to \nabla_n^{\mmu}\MC{H}^1_U(C, G)\]
where $\MB{M}_{b_{c_i}}^{\preceq \mu_i}$ is the reduced fiber of the formal scheme $\MC{M}_{b_{c_i}}^{\preceq \mu_i}$ (c.f. \ref{Thm:RapoZinkBounded}). \\
This morphism was first constructed in \cite[section 5]{HarRad1} (although only in the case where the global $G$-shtuka is constant on $S_0$) by changing $(\MC{G}_0, \varphi_0, \psi_0)$ at every characteristic place according to the respective point in $\MB{M}_{b_{c_i}}^{\preceq \mu_i}$. Unfortunately we were not able to understand either the descent argument in lemma \cite[4.23]{HarRad1} (c.f. \ref{Rem:UnifProblem}i)) or the one in their key lemma \cite[5.1]{HarRad1}. So let us provide more elaborate arguments, following the very same approach. \\
Note that the construction below is given for any global $G$-shtuka $(\MS{G}_0, \varphi_0, \psi_0)$ with characteristic places only in a formal neighborhood of the $(c_i)$. As the results of this section are only needed in full generality in the next section \ref{sec:FormalLifting}, the reader may wish to stick to the situation considered in this introduction and ignore the slight variation of characteristic places during a first read. \\
As the construction is not restricted to characteristic places, let us fix any place $v \in C$ with residue field $\kappa(v)$, though in applications it is one of the $c_i$. Let $\M{A}_v^{int} \cong \kappa(v)[[z]]$ be the ring of integral adeles at $v$, $\M{A}_v \cong \kappa(v)((z))$ its quotient field and $G_v = Res_{\kappa(v)/\Fq}(G)$. Moreover fix a local coordinate $\zeta$ at the place $v$ and consider the category $Nilp_{\ENil}$. Our first aim is to identify certain classes of local $G_v$-shtukas with certain pairs consisting of a $G$-torsor over $\Sp \M{A}_v^{int}$ and a Frobenius-isomorphism on the generic fiber.

\defi{}{}{
Let $R$ be an $\ENil$-algebra with $\zeta$ locally nilpotent. A torsor (for some group) over a (formal) scheme $S$ over $\Sp R$ is called \textit{nice} if it trivializes over a finite \'etale cover of $S$ coming by base change from a finite \'etale cover of $\Sp R$. \\
A local $G_v$-shtuka over $\Sp R$ is called nice, if the underlying $L^+G_v$-torsor is nice.
}

\rem{}{
Here is the reason why nice torsors are helpful in the study of arbitrary global $G$-shtukas: If $\MS{G}$ is a $G$-torsor over $C \times_{\Fq} S$, then there is exists a Zariski-open cover $\{V_i\}_i$ of $S$ such that $\MS{G}|_{C \times_{\Fq} V_i}$ is nice in the formal neighborhood of the fixed place $v$. This will be shown in lemma \ref{lem:UnifLem1} and \ref{lem:UnifLem2}. In particular any local shtuka associated to a global $G$-shtuka is Zariski-locally nice.
}

\lem{\label{Lem:UnifSemiTorsor}}{}{
Let $R$ be an $\ENil$-algebra with $\zeta$ locally nilpotent. Then we have equivalences of categories \\
a) \vspace{-5mm}
\begin{align*}
 \{\textnormal{nice} \, L^+G_v\textnormal{-torsors over} \, \Sp R\} & \cong \{\textnormal{nice formal} \, G\textnormal{-torsors over} \, \Spf (R \widehat{\otimes}_{\Fq} \M{A}_v^{int})\} \\
 & \cong \{\textnormal{nice} \, G\textnormal{-torsors over} \, \Sp (R \widehat{\otimes}_{\Fq} \M{A}_v^{int})\}
\end{align*}
where the morphisms are bijective morphisms which are equivariant for the respective group action. \\
b) \[\{\textnormal{nice} \, LG_v\textnormal{-torsors over} \, \Sp R\} \cong \{\textnormal{nice} \, G\textnormal{-torsors over} \, \Sp (R \widehat{\otimes}_{\Fq} \M{A}_v)\} \]
where the morphisms are isomorphisms of the respective torsors. \\
If $R_1 \to R_2$ is a morphism of $\ENil$-algebras, then all equivalences above commute with taking pullbacks of the torsors along the induced morphism on the base schemes.
}

\prooof
a) We have already seen in lemma \ref{lem:LoGloRes}a) that $L^+G_v$-torsors over $\Sp R$ correspond to $L^+G$-torsors over $\Sp (R \otimes_{\Fq} \kappa(v))$. In \cite[proposition 2.2 a)]{HaVi} a bijection was constructed between $L^+G$-torsors over $S = \Sp (R \otimes_{\Fq} \kappa(v))$ and formal $G$-torsors over $\Spf (R \otimes_{\Fq} \kappa(v))[[z]] \cong \Spf (R \widehat{\otimes}_{\Fq} \M{A}_v^{int})$ which trivialize over an \'etale cover induced by an \'etale cover of $S$. This construction obviously restricts to a bijection between nice torsors on both sides. Furthermore one can easily see, that the construction in \cite{HaVi} extends to morphisms and hence indeed defines the equivalence of categories on the left-hand side. \\
Let us now tackle the second comparison: As $R \widehat{\otimes}_{\Fq} \M{A}_v^{int} \cong (R \otimes_{\Fq} \kappa(c_i))[[z]]$ we may replace $R$ by the finite \'etale algebra $R \otimes_{\Fq} \kappa(c_i)$ and hence reduce to $G$-torsors over $\Spf R[[z]]$ and $\Sp R[[z]]$. Now giving a nice $G$-torsor over $\Spf R[[z]]$ is the same as giving a finite \'etale cover $\Spf R'[[z]] \to \Spf R[[z]]$ coming from an \'etale morphism $R \to R'$ and a descent datum 
\begin{align*}
 h \in \varprojlim G(R'[[z]]/(z^n) \otimes_{R[[z]]/(z^n)} R'[[z]]/(z^n)) & \cong G(\varprojlim R'[[z]]/(z^n) \otimes_{R[[z]]/(z^n)} R'[[z]]/(z^n)) \\
 & = G(R'[[z]] \widehat{\otimes}_{R[[z]]} R'[[z]])
\end{align*}
for the trivial $G$-torsor on $\Spf R'[[z]]$ (see also \cite[proof of proposition 2.2 a)]{HaVi}). By finiteness of the cover we may view $h \in G(R'[[z]] \otimes_{R[[z]]} R'[[z]])$ and get therefore a unique descent datum for the trivial $G$-torsor on the finite \'etale cover $\Sp R'[[z]] \to \Sp R[[z]]$. But this determines a unique object in the category on the right-hand side. Note that this construction does not depend on the actual choice of the finite \'etale cover. Thus this defines an essential bijection between the objects in both categories. \\
Consider now an isomorphism $\alpha: \MC{G}_1 \to \MC{G}_2$ between nice formal $G$-shtukas. If they trivialize over $\Spf R'[[z]] \to \Spf R[[z]]$ with descent data given by $h_1, h_2 \in G(R'[[z]] \widehat{\otimes}_{R[[z]]} R'[[z]])$ then $\alpha$ can be represented by a unique element $g \in G(R'[[z]])$ satisfying $pr_1(g) = h_2^{-1} \cdot pr_2(g) \cdot h_1$ where $pr_i: G(R'[[z]]) \to G(R'[[z]] \widehat{\otimes}_{R[[z]]} R'[[z]])$ is given by the inclusion of $R'[[z]]$ in the $i$th factor. But the same data determines morphisms between the associated nice $G$-torsors over $\Sp R[[z]]$, with the only difference that we now consider $pr_1(g) = h_2^{-1} \cdot pr_2(g) \cdot h_1$ as an equality in $G(R'[[z]] \otimes_{R[[z]]} R'[[z]])$. Hence we get a canonical bijection between the sets of morphisms. \\
b) Let $R \to R'$ be a finite \'etale morphism. Then $LG$-torsors over $\Sp R$ which trivialize over $\Sp R'$ are given by a descent datum in $G((R' \otimes_R R')((z)))$. On the other hand $G$-torsors over $\Sp R((z))$ trivializing over $\Sp R'((z))$ are given by a descent datum in $G(R'((z)) \otimes_{R((z))} R'((z)))$. But by finiteness of $R'$, there is a canonical isomorphism $(R' \otimes_R R')((z)) \cong R'((z)) \otimes_{R((z))} R'((z))$. Therefore we have a bijection between the objects on both sides. For morphisms a very similar argument as in part a) works. \\
For the compatibility with pull-backs, note first that if a $L^+G_v$-torsor over $\Sp R_1$ trivializes over a finite \'etale cover $\Sp R_1' \to \Sp R_1$, then its pullback to $\Sp R_2$ trivializes over the finite \'etale cover $\Sp R_2 \otimes_{R_1} R_1' \to \Sp R_2$. The same holds for $G$-torsors over $\Sp (R \widehat{\otimes}_{\Fq} \M{A}_v^{int})$. Hence the compatibility of the constructions above and the pullback functor follows from the obvious facts that it holds for trivial torsors and that the descent datum for the pullback to $R_2$ coincides with the pullback of the descent datum for the torsor over $R_1$. \exit

\rem{\label{Rem:UnifProblem}}{
i) The functor from nice $G$-torsors over $\Sp (R \widehat{\otimes}_{\Fq} \M{A}_v^{int})$ to nice formal $G$-torsors over $\Spf (R \widehat{\otimes}_{\Fq} \M{A}_v^{int})$ is in fact the pullback functor. As pulling back is not restricted to nice $G$-torsors one may ask, whether this defines an equivalence of categories
\[\{\textnormal{formal} \, G\textnormal{-torsors over} \, \Spf R[[z]]\} \stackrel{?}{\cong} \{G\textnormal{-torsors over} \, \Sp R[[z]]\}\]
(assuming for simplicity $\kappa(v) \cong \Fq$). Nevertheless we were not able to prove this due to two problems: \\
$\alpha$) The functor $\widehat{\otimes}_{\Fq} \Fq[[z]]$ associating to $R$ the ring $R[[z]]$ does not commute with localizations. In particular assume we have a $G$-torsor $\MC{G}$ over $\Spf R[[z]]$ and an open formal subscheme $\Spf R[f^{-1}][[z]]$ (for some $f \in R$) such that $\MC{G}$ restricts to a nice $G$-torsor over $\Spf R[f^{-1}][[z]]$. Then applying the equivalence above gives a $G$-torsor over $\Sp R[f^{-1}][[z]]$. But this is not an open subscheme of $\Sp R[[z]]$ and hence we may not glue such $G$-torsors. \\
$\beta$) Assume there is a functor $\{\textnormal{formal} \, G\textnormal{-torsors over} \, \Spf R[[z]]\} \to \{G\textnormal{-torsors over} \, \Sp R[[z]]\}$ which commutes with pullbacks to $\Spf R[f^{-1}][[z]]$ resp. $\Sp R[f^{-1}][[z]]$ for any $f \in R$. Then any $G$-torsor in the image has the property that there is a set of elements $\{f_i\}_i \in R$ such that the morphisms $\Sp R[f_i^{-1}][[z]] \to \Sp R[[z]]$ are jointly surjective and the pullback to each $\Sp R[f_i^{-1}][[z]]$ is nice. But it is not known to us, whether every $G$-torsor over $\Sp R[[z]]$ has this property. Nevertheless note the partial result \ref{lem:UnifLem1} below. \\
ii) The lemma implies, that we have a bijection between nice $L^+G$-torsors over $\Spf R[[t]]$ and nice $L^+G$-torsors over $\Sp R[[t]]$. To see this view the $L^+G$-torsors on both sides as an injective limit of $G$-torsors over $\Spf (R[[z]]/(z^n))[[t]]$ respectively $\Sp (R[[z]]/(z^n))[[t]]$ and use the correspondence above. As in \cite[proposition 3.16]{HaVi} this induces an equivalence of categories
\begin{align*}
 & \{\textnormal{nice formal local} \, G\textnormal{\!-shtukas over} \, \Spf R[[t]] \, \textnormal{bounded by} \, \mu\} \cong \\
 & \hspace{1cm} \cong \{\textnormal{nice local} \, G\textnormal{\!-shtukas over} \, \Sp R[[t]] \, \textnormal{bounded by} \, \mu\}
\end{align*}
where $\mu$ is as usual a $\Gamma$-invariant dominant cocharacter.
Nevertheless the same problems as stated in part i) of this remark refrain us from omitting the condition `nice'. Therefore we do not know that \cite[proposition 3.16]{HaVi} as stated there is actually correct, although this seems likely.
}

\defi{}{}{
a) Let $R$ be an $\ENil$-algebra with $\zeta$ locally nilpotent. A semi-local $G_v$-shtuka over $R$ is a pair $(\MC{G}, \varphi)$ consisting of a $G$-torsor $\MC{G}$ over $\Sp (R \widehat{\otimes}_{\Fq} \M{A}_v^{int})$ and a $G$-equivariant morphism $\varphi: \sigma^*\MC{G}|_\eta \to \MC{G}|_\eta$ over the generic fiber $\eta = \Sp (R \widehat{\otimes}_{\Fq} \M{A}_v)$.\\
b) A semi-local $G_v$-shtuka is called nice if the underlying $G$-torsor is nice.
}

\rem{}{
i) Note that the notion of a (nice) semi-local $G$-shtuka behaves very badly, e.g. one cannot glue semi-local $G$-shtukas over an open cover of $\Sp R$.\\
ii) Pulling back global $G$-shtukas via the map $\Sp (R \widehat{\otimes}_{\Fq} \M{A}_v^{int}) \to C \times_{\Fq} \Sp R$ defines a semi-local $G_v$-shtuka. Here we use the pullback along $\Sp (R \widehat{\otimes}_{\Fq} \M{A}_v) \to (C \setminus \{c_i\}_i) \times_{\Fq} \Sp R$ to get the local Frobenius-morphism.
}

\prop{\label{Prop:UnifShtukaPowerseries}}{}{
Let $R$ be an $E$-algebra. Then we have an equivalence of categories
\[\{\textnormal{nice local} \, G_v\textnormal{-shtukas over} \, \Sp R\} \cong \{\textnormal{nice semi-local} \, G_v\textnormal{-shtukas over} \, R\}\]
Furthermore it commutes with taking the pullback along $\Sp R_1 \to \Sp R_2$ respectively $\Sp R_1[[z]] \to \Sp R_2[[z]]$ for any morphism $R_2 \to R_1$ of $E$-algebras.
}

\prooof
We have already seen this on the level of torsors. 
Next note that a $LG_v$-torsor over an affine base $\Sp R$ is nothing else than a $G_v$-torsor over $\Sp (R \widehat{\otimes}_{\Fq} \M{A}_v)$. Furthermore associating the $LG_v$-torsor to a nice $L^+G_v$-torsor over $\Sp R$ corresponds exactly to restricting a $G_v$-torsor over $\Sp (R \widehat{\otimes}_{\Fq} \M{A}_v^{int})$ to its generic fiber.
As these identifications were given by equivalences of categories, this defines a bijection between the set of possible Frobenius-isomorphisms on both sides. In other words, we have a bijection between nice local $G_v$-shtukas over $\Sp R$ and nice semi-local $G_v$-shtukas over $R$. \\
That the morphism sets coincide follows in a very similar way after noticing that the compatibility conditions with the Frobenius-morphism coincide on both sides. \\
The compatibility with taking pullbacks follows from the corresponding statement in the previous lemma. \exit 

\rem{}{
The bijection of nice $LG_v$-torsors over $\Sp R$ and nice $G_v$-torsors over $\Sp (R \widehat{\otimes}_{\Fq} \M{A}_v)$ was only used for convenience. In fact we may pass to some trivializing \'etale cover, where the bijection on Frobenius-isomorphisms is clear. Then one easily checks that this bijection preserves the property to commute with the descent data, which induces the bijection already over $R$.
}

\nota{}{
From now on fix $n$ distinct $E$-valued points $c_i$ on the curve $C$ together with local coordinates $\zeta_i$. Let $Nilp_{E[[\zeta_1, \ldots, \zeta_n]]}$ be the category of schemes $S$ over $\Spf E[[\zeta_1, \ldots, \zeta_n]]$, such that the ideal $(\zeta_1, \ldots, \zeta_n)$ is locally nilpotent on $S$. Equivalently it is the category of schemes $S$ over $C^n \setminus \Delta$, that factor over the formal completion of $C^n \setminus \Delta$ along $(c_i)_i$. We automatically consider schemes $S \in Nilp_{E[[\zeta_1, \ldots, \zeta_n]]}$ as elements of $Nilp_{E[[\zeta_j]]}$ via the canonical projection $\Spf E[[\zeta_1, \ldots, \zeta_n]] \to \Spf E[[\zeta_j]]$. \\
Define for any scheme $S$ respectively affine scheme $\Sp R$ in $Nilp_{E[[\zeta_1, \ldots, \zeta_n]]}$ the categories
\begin{description}[]
 \item[] $\nabla_{(\hat{c}_i)} \MC{H}^1_U(C, G)(S)$ with objects global $G$-shtukas with $U$-level structure, whose characteristic places are given by the projections of $S \to C^n \setminus \Delta$ to the various components $C$. Hence $\nabla_{(\hat{c}_i)} \MC{H}^1_U(C, G)(S)$ can be identified with the fiber of the stack $\nabla_n \MC{H}^1_U(C, G) \to C^n \setminus \Delta$ over $S$.
 \item[] $\nabla^{\mmu}_{(\hat{c}_i)} \MC{H}^1_U(C, G)(S)$ similarly as $\nabla_{(\hat{c}_i)} \MC{H}^1_U(C, G)(S)$, but for bounded global $G$-shtukas.
 \item[] $\nabla_{(\hat{c}_i)} \MC{H}^1_U(V, G)(S)$ with the very same definition as for the category $\nabla_{(\hat{c}_i)} \MC{H}^1_U(C, G)(S)$, but $C$ now replaced by an open subscheme $V \subset C$. This requires $U \subset G(\M{A}^{int (c_i)}_V)$, where $\M{A}^{int (c_i)}_V$ is the ring of integral adeles of $V \setminus \{c_i\}_i$ and the condition on the compatibility of characteristic places with $S$ only applies to places $c_i$, that actually lie inside $V$. 
 \item[] $SemiSht_{G_{c_j}}(\Sp R)$ the category of semi-local $G_{c_j}$-shtukas over $R$.
 \item[] $QSemiSht_{G_{c_j}}(\Sp R)$ with objects $(\MC{G}, \varphi)$ consisting of a $G_{c_j}$-torsor $\MC{G}$ over $\Sp R((z))$ and an isomorphism $\varphi: \sigma^*\MC{G} \to \MC{G}$.
 \item[] $Sht_{G_{c_j}}(S)$ the category of local $G_{c_j}$-shtukas over $S$.
 \item[] $QSht_{G_{c_j}}(S)$ with objects $(\MC{LG}, \varphi)$ consisting of a $LG_{c_j}$-torsor $\MC{LG}$ over $S$ and an isomorphism $\varphi: \sigma^*\MC{LG} \to \MC{LG}$.
\end{description}
In all cases the morphisms are given by quasi-isogenies.
}

\prop{\label{Prop:UnifPreGluing}}{}{
Fix one of the characteristic places $c_j$ as above and let $\Sp R \in Nilp_{E[[\zeta_1, \ldots, \zeta_n]]}$ be an affine scheme. Then there is a fiber product diagram of categories
\[
 \begin{xy} \xymatrix{
   \nabla_{(\hat{c}_i)}\MC{H}^1_U(C, G)(\Sp R) \ar[r] \ar[d] & \nabla_{(\hat{c}_i)_{i \neq j}}\MC{H}^1_U(C \setminus \{c_j\}, G)(\Sp R) \ar[d] \\
   SemiSht_{G_{c_j}}(\Sp R) \ar[r] & QSemiSht_{G_{c_j}}(\Sp R)
 } \end{xy} 
\]
where the vertical functors are given by pulling back objects and morphisms on a formal neighborhood the $j$th characteristic place and the horizontal functors are given by restriction to the open subset where $c_j$ is missing. Moreover the fiber product diagram is natural in $R$.
}

\prooof
As already mentioned, we get a canonical morphism $\Sp R \to C^n \setminus \Delta \xrightarrow{pr_j} C$. We denote the image of its graph morphism by $Z \subset C \times_{\Fq} \Sp R$ and the formal completion of $C \times_{\Fq} \Sp R$ along $Z$ by $\widehat{Z}$. Then $\Sp (R \widehat{\otimes}_{\Fq} \M{A}_{c_j}^{int})$ is isomorphic to $\widehat{Z}$ and the vertical morphism on the left is just the restriction to $\widehat{Z}$. Furthermore $(C \setminus \{c_j\}) \times_{\Fq} \Sp R$ identifies with $(C \times_{\Fq} \Sp R) \setminus Z$, because the $j$th characteristic place varies only in a formal neighborhood. \\
Hence by \cite{BeauLaszlo} (or even earlier \cite[proposition 4.2]{FerrandRaynaud}) we have a fiber product diagram of categories of quasi-coherent sheaves
\[
 \begin{xy} \xymatrix{
   QCoh(C \times_{\Fq} \Sp R) \ar[r] \ar[d] & QCoh((C \setminus \{c_j\}) \times_{\Fq} \Sp R) \ar[d] \\
   QCoh(\Sp (R \widehat{\otimes}_{\Fq} \M{A}_{c_j}^{int})) \ar[r] & QCoh(\Sp (R \widehat{\otimes}_{\Fq} \M{A}_{c_j}))
 } \end{xy} 
\]
which is natural in $R$. \\
Assume now $\MS{G}^\circ$ is a $G$-torsor on $(C \setminus \{c_j\}) \times_{\Fq} \Sp R$, $\MC{G}$ is a $G$-torsor on $\Sp (R \widehat{\otimes}_{\Fq} \M{A}_{c_j}^{int})$ and $\alpha: \MS{G}^\circ|_{\Sp (R \widehat{\otimes}_{\Fq} \M{A}_{c_j})} \to \MC{G}|_{\Sp (R \widehat{\otimes}_{\Fq} \M{A}_{c_j})}$ is an isomorphisms between their restrictions to $\Sp (R \widehat{\otimes}_{\Fq} \M{A}_{c_j})$. Let $\MS{G}$ be the quasi-coherent sheaf obtained by gluing $\MS{G}^\circ$ and $\MC{G}$ along $\alpha$. Then the $G$-actions on $\MS{G}^\circ$ and $\MC{G}$ glue to a $G$-action on $\MS{G}$. 
As $((C \setminus \{c_j\}) \times_{\Fq} \Sp R) \coprod \Sp (R \widehat{\otimes}_{\Fq} \M{A}_{c_j}^{int}) \to C \times_{\Fq} \Sp R$ is a fpqc-cover, this $G$-bundle $\MS{G}$ is in fact a $G$-torsor for the fpqc-topology (and hence for the \'etale topology). Thus we have a fiber product diagram
\[
 \begin{xy} \xymatrix{
   \{G\op{\!-torsors \; on \,} C \times_{\Fq} \Sp R\} \ar[r] \ar[d] & \{G\op{\!-torsors \; on \,} (C \setminus \{c_j\}) \times_{\Fq} \Sp R\} \ar[d] \\
   \{G\op{\!-torsors \; on \,} \Sp (R \widehat{\otimes}_{\Fq} \M{A}_{c_j}^{int})\} \ar[r] & \{G\op{\!-torsors \; on \,} \Sp (R \widehat{\otimes}_{\Fq} \M{A}_{c_j})\}
 } \end{xy} 
\]
Choose now elements $(\MS{G}^\circ, \varphi^\circ, \psi^\circ) \in \nabla_{(\hat{c}_i)_{i \neq j}}\MC{H}^1_U(C \setminus \{c_j\}, G)(\Sp R)$, $(\MC{G}, \varphi) \in SemiSht_{G_{c_j}}(\Sp R)$ and an isomorphism $\alpha: \MS{G}|_{\Sp (R \widehat{\otimes}_{\Fq} \M{A}_{c_j})} \to \MC{G}|_{\Sp (R \widehat{\otimes}_{\Fq} \M{A}_{c_j})}$ compatible with the isomorphisms $\varphi^\circ$ and $\varphi$. We have to show that these objects define a unique global $G$-shtuka. We already know that $\MS{G}^\circ$ and $\MC{G}$ glue to a unique $G$-torsor $\MS{G}$ over $C \times_{\Fq} \Sp R$. 
For the Frobenius-isomorphism notice that the one in $(\MC{G}, \varphi) \in SemiSht_{G_{c_j}}(\Sp R)$ is already determined by its image in $QSemiSht_{G_{c_j}}(\Sp R)$. Thus simply setting 
\[\varphi = \varphi^\circ: \sigma^*\MS{G}|_{(C \times_{\Fq} \Sp R) \setminus Z} \to \MS{G}|_{(C \times_{\Fq} \Sp R) \setminus Z}\]
(using the canonical $\MS{G}|_{(C \times_{\Fq} \Sp R) \setminus Z} \cong \MS{G}^\circ|_{(C \times_{\Fq} \Sp R) \setminus Z}$) defines a Frobenius-morphism which restricts to the given $\varphi$ on the generic fiber of the formal neighborhood of $c_j$ (at least after applying the canonical identification $\alpha$).
The $U$-level structure $\psi^\circ$ gives a $U$-level structure $\psi = \psi^{\circ}$ on $\MS{G}$ as $\MS{G}$ coincides with $\MS{G}^\circ$ on the formal neighborhood of the support of $U$. \\
Hence we have indeed the stated fiber product diagram. \exit

\thm{\label{Thm:UnifGluingOne}}{}{
For every $S \in Nilp_{E[[\zeta_1, \ldots, \zeta_n]]}$ there is a canonical fiber product diagram of categories
\[
 \begin{xy} \xymatrix{
   \nabla_{(\hat{c}_i)}\MC{H}^1_U(C, G)(S) \ar[r] \ar_{\MF{L}_{c_j}}[d] & \nabla_{(\hat{c}_i)_{i \neq j}}\MC{H}^1_U(C \setminus \{c_j\}, G)(S) \ar^{\MF{L}_{c_j}}[d] \\
   Sht_{G_{c_j}}(S) \ar[r] & QSht_{G_{c_j}}(S)
 } \end{xy} 
\]
It is natural in $S$.
}

\prooof
We construct a functor
\[\nabla_{(\hat{c}_i)_{i \neq j}}\MC{H}^1_U(C \setminus \{c_j\}, G)(S) \times_{QSht_{G_{c_j}}(S)} Sht_{G_{c_j}}(S) \longrightarrow \nabla_{(\hat{c}_i)}\MC{H}^1_U(C, G)(S)\]
which is an inverse to the obvious functor in the other direction. \\
For this consider any elements $(\MS{G}^\circ, \psi^\circ, \varphi^\circ) \in \nabla_{(\hat{c}_i)_{i \neq j}}\MC{H}^1_U(C \setminus \{c_j\}, G)(S)$, $(\MC{G}, \varphi) \in Sht_{G_{c_j}}(S)$ and an isomorphism $\alpha: \MF{L}_{c_j}(\MS{G}^\circ) \to \MC{LG}$ compatible with the isomorphisms $\varphi^\circ$ and $\varphi$. Choose now a basis for the topology $\MF{B}$ of affine open subschemes of $S$ such that for any $W \in \MF{B}$ the torsor $\MC{G}|_W$ is nice (with respect to the base scheme $W$). Consider now any such $W = \Sp R \in \MF{B}$. 
By proposition \ref{Prop:UnifShtukaPowerseries} $(\MC{G}, \varphi)$ defines a unique element $(\MC{G}', \varphi')$ in $SemiSht_{G_{c_j}}(\Sp R)$ and $\alpha$ defines a unique isomorphism $\alpha': \MC{G}^\circ|_{\Sp (R \widehat{\otimes}_{\Fq} \M{A}_{c_j}^{int})} \to \MC{G}'|_{\Sp (R \widehat{\otimes}_{\Fq} \M{A}_{c_j}^{int})}$. Thus the triple $((\MC{G}^\circ, \varphi^\circ, \psi^\circ), (\MC{G}', \varphi'), \alpha')$ defines a uniquely determined global $G$-shtuka $(\MS{G}_W, \psi_W, \varphi_W) \in \nabla_n\MC{H}^1_U(C, G)(W)$ by proposition \ref{Prop:UnifPreGluing}. \\
If $W_2 \subset W_1$ is an inclusion of elements in $\MF{B}$ then $(\MS{G}_{W_1}, \varphi_{W_1}, \psi_{W_1})|_{W_2}$ and $(\MS{G}_{W_2}, \varphi_{W_2}, \psi_{W_2})$ define the same elements in $\nabla_{(\hat{c}_i)_{i \neq j}}\MC{H}^1_U(C \setminus \{c_j\}, G)(W_2)$ and $SemiSht_{G_{c_j}}(W_2)$ (using the compatibility with pullbacks in \ref{Prop:UnifShtukaPowerseries}), hence there is a uniquely determined isomorphisms $(\MS{G}_{W_1}, \varphi_{W_1}, \psi_{W_1})|_{W_2} \cong (\MS{G}_{W_2}, \varphi_{W_2}, \psi_{W_2})$. Thus we may glue the global $G$-shtukas $(\MS{G}_W, \varphi_W, \psi_W)_{W \in \MF{B}}$ to a global $G$-shtuka $(\MS{G}, \varphi, \psi) \in \nabla_{(\hat{c}_i)}\MC{H}^1_U(C, G)(S)$. It is clear that this construction does not depend on the choice of the basis $\MF{B}$. \\
It remains to be shown, that these functors are inverses. If one starts with an element in the fiber product and takes the image in $\nabla_{(\hat{c}_i)}\MC{H}^1_U(C, G)(S)$, then it is obvious from the construction that the restriction to $(C \setminus \{c_j\}) \times_{\Fq} S$ gives back the original element in $\nabla_{(\hat{c}_i)_{i \neq j}}\MC{H}^1_U(C \setminus \{c_j\}, G)(S)$. We have to see that $\MF{L}_{c_j}$ returns the original local $G_{c_j}$-shtuka as well. This can be done locally on the elements $W = \Sp R \in \MF{B}$. There proposition \ref{Prop:UnifPreGluing} and proposition \ref{Prop:UnifShtukaPowerseries} imply that we may obtain the original local $G_{c_j}$-shtuka via
\[\nabla_{(\hat{c}_i)}\MC{H}^1_U(C, G)(W) \to SemiSht_{G_{c_j}}(W) \cong Sht_{G_{c_j}}(W)\]
where the first functor is the pullback to $\Sp (R \widehat{\times}_{\Fq} \Sp \M{A}_{c_j}^{int})$ and the second functor consists of pulling further back to the formal completion of $\Sp (R \widehat{\times}_{\Fq} \Sp \M{A}_{c_j}^{int})$ and translation from $G$-torsors to $L^+G_{c_j}$-torsors. But this is by construction nothing else than the functor $\MF{L}_{c_j}$. \\
Conversely start with a global $G$-shtuka $(\MS{G}_0, \varphi_0, \psi_0) \in \nabla_{(\hat{c}_i)}\MC{H}^1_U(C, G)(S)$ with restriction $(\MS{G}^\circ, \varphi^\circ, \psi^\circ) \in \nabla_{(\hat{c}_i)_{i \neq j}}\MC{H}^1_U(C \setminus \{c_j\}, G)(S)$ and associated local $G_{c_j}$-shtuka $(\MC{G}, \varphi)$. Let $(\MS{G}, \varphi, \psi)$ be the global $G$-shtuka obtained from gluing $(\MS{G}^\circ, \varphi^\circ, \psi^\circ)$ and $(\MC{G}, \varphi)$. By the following lemma \ref{lem:UnifLem1} one can choose a basis of the topology $\MF{B}$ of affine open subschemes of $S$ such that the restriction of $(\MC{G}_0, \psi_0, \varphi_0)$ to $\Sp (R \widehat{\otimes}_{\Fq} \M{A}_{c_j}^{int})$ is nice for every element $\Sp R \in \MF{B}$. 
Consider now one such element $\Sp R \in \MF{B}$. Then by choice on $\Sp R$ there is a canonical isomorphism between the semi-local $G$-shtukas $(\MC{G}_0, \psi_0, \varphi_0)|_{\Sp (R \widehat{\otimes}_{\Fq} \M{A}_{c_j}^{int})}$ and the one associated to the nice local $G_{c_j}$-shtuka $(\MC{G}, \varphi)|_{\Sp R}$. This gives a canonical isomorphism $(\MS{G}_0, \varphi_0, \psi_0)|_{C \times_{\Fq} \Sp R} \cong (\MS{G}, \varphi, \psi)|_{C \times_{\Fq} \Sp R}$. As they are compatible with restrictions to other subschemes in $\MF{B}$, they glue to an isomorphism $(\MS{G}_0, \varphi_0, \psi_0) \cong (\MS{G}, \varphi, \psi)$. \\
The same construction can be done for morphisms in the respective categories, i.e. for quasi-isogenies of global respective local $G$-shtukas. \exit

\lem{\label{lem:UnifLem2}}{}{
Let $R$ be a $\Fq$-algebra. \\
a) The pullback along $\Sp R = \Sp R[[z]]/(z) \to \Sp R[[z]]$ gives an isomorphism  
\[\pi_1(\Sp R[[z]]) \cong \pi_1(\Sp R)\]
b) Any $G$-torsor $\MC{G}$ over $\Sp (R \widehat{\otimes}_{\Fq} \kappa(c_j)[[z]])$ which admits a trivialization after a finite \'etale cover, also admits one induced from a finite \'etale cover of $R$. 
}

\prooof
a) This is surely well-known, though we were not able to find a reference. Surjectivity is obvious as the inverse map associates to every finite \'etale morphism $R \to R'_1$ the finite \'etale morphism $R[[z]] \to R'_1[[z]]$. \\
Conversely let $R[[z]] \to R'$ be any finite \'etale morphism and define $R_n = R[[z]]/(z^n)$ and $R'_n = R' \otimes_{R[[z]]} R_n$ for every $n \geq 1$. We have to show that $R' \cong R_1'[[z]]$. 
But by \cite[IX, proposition 1.7]{GroSGAI} $\pi_1(\Sp R_n) \cong \pi_1(\Sp R)$. Applied to the finite \'etale map $R_n \to R'_n$ we get isomorphisms $R'_n \cong R_{0, n}' \otimes_R R_n$ for each $n$, where $R_{0, n}'$ is finite \'etale over $R$. Then $R_{0, n+1}' \otimes_R R_n \cong R'_{n + 1} \otimes_{R_{n + 1}} R_n \cong R'_n \cong R_{0, n}' \otimes_R R_n$ and all rings $R_{0, n}'$ are isomorphic to the pullback $R_1' = R_{0, 1}'$. Then we have an isomorphism of $R$-algebras 
\[R' = \varprojlim R'_n \cong \varprojlim (R_1' \otimes_R R_n) \cong R_1' \otimes_R \varprojlim R_n = R_1' \otimes_R R[[z]]\]
using that $R'$ equals its $z$-adic completion as it is finitely presented over $R[[z]]$. Note that in our case the inverse limit does commute with the tensor product, because $R_1'$ is finite over $R$ and hence both sides equal $R_1'{}^{\M{N}}$ (as $R$-modules). \\
b) The $G$-torsor $\MC{G}$ corresponds to a $G_{c_j}$-torsor over $R[[z]]$ and any trivializing cover of $R[[z]]$ for the $G_{c_j}$-torsor gives also a trivializing cover for $\MC{G}$. But by part a) any finite \'etale cover of $R[[z]]$ comes from a finite \'etale extension of $R$. \exit

\lem{\label{lem:UnifLem1}}{}{
Let $\Sp A$ be an affine curve over $\Fq$ with a regular point $v$. Fix a local coordinate $\zeta$ of $\Sp A$ at $v$ and an $\ENil$-algebra $R$ with $\zeta$ locally nilpotent in $R$. Denote the completion of $\Sp A \times_{\Fq} \Sp R$ along $\{v\} \times_{\Fq} \Sp R$ by $(\Sp A \times_{\Fq} \Sp R)^{\wedge v}$. \\
Then for any $G$-torsor $\MS{G}$ on $\Sp A \times_{\Fq} \Sp R$ and any point $x \in \Sp R$ there is an open neighborhood $x \in \Sp R_0 \subset \Sp R$ such that the pullback of $\MS{G}$ to $(\Sp A \times_{\Fq} \Sp R_0)^{\wedge v}$ is nice with respect to $\Sp R_0$, i.e. trivializes over a finite \'etale cover coming from one over $\Sp R_0$. 
}

\rem{}{
Note that $(\Sp A \times_{\Fq} \Sp R)^{\wedge v}$ equals the completion along the graph of the canonical morphism $\Sp R \to A$ as constructed at the beginning of the proof of proposition \ref{Prop:UnifPreGluing}.
}

\prooof
We may choose an open subset $U \subset \Sp R \times_{\Fq} \Sp A$ containing the point $(x, v)$ and a finite \'etale cover $U' \to U$ such that $\MS{G}$ trivializes over $U'$. We may assume wlog. that $U$ is of the form $D(f)$ for some element $f \in A \otimes_{\Fq} R$. By restricting to an open subset $\Sp R_0 \subset \Sp R$ containing $x$ we may assume that 
\[\{v\} \times_{\Fq} \Sp R_0 \subset U = D(f) \subset \Sp A \times_{\Fq} \Sp R.\]
Then $f$ is a unit over $\{v\} \times_{\Fq} \Sp R_0$, hence as well over the completion $(\Sp A \times_{\Fq} \Sp R_0)^{\wedge v}$. Thus the inclusion $(\Sp A \times_{\Fq} \Sp R)^{\wedge v} \to \Sp A \times_{\Fq} \Sp R$ factors over $\Sp (A \otimes_{\Fq} R_0)[f^{-1}] = U \times_{\Sp R} \Sp R_0$. Now we just observe that by definition the pullback of $\MS{G}$ to $U \times_{\Sp R} \Sp R_0$ trivializes over the finite \'etale cover $U' \times_{\Sp R} \Sp R_0$ and hence the pullback to $(\Sp A \times_{\Fq} \Sp R_0)^{\wedge v}$ trivializes over the finite \'etale cover $U' \times_{\Sp A \times \Sp R} (\Sp A \times_{\Fq} \Sp R_0)^{\wedge v}$. \\
Finally note that by regularity of the point $v$, there exists an isomorphism $(\Sp A \times_{\Fq} \Sp R_0)^{\wedge v} \cong \Sp (\kappa(v) \otimes_{\Fq} R_0)[[z]]$. Hence by lemma \ref{lem:UnifLem2}b) the trivializing cover $U' \times_{\Sp A \times \Sp R} (\Sp A \times_{\Fq} \Sp R_0)^{\wedge v}$ already comes from a finite \'etale cover of $R_0$, hence $\MS{G}|_{(\Sp A \times_{\Fq} \Sp R_0)^{\wedge v}}$ is indeed nice. \exit

\rem{}{
Similarly for any non-characteristic place $v \in C$, we get a canonical fiber product diagram of categories
\[
 \begin{xy} \xymatrix{
   \nabla_{(\hat{c}_i)}\MC{H}^1_U(C, G)(S) \ar[r] \ar_{\MF{L}_v}[d] & \nabla_{(\hat{c}_i)} \MC{H}^1_U(C \setminus \{v\}, G)(S) \ar^{\MF{L}_v}[d] \\
   \acute{E}tSht_{G_v}(S) \ar[r] & QSht_{G_v}(S)
 } \end{xy} 
\]
}

\cor{\label{cor:UnifGluing}}{}{
Let $\mmu = (\mu_i)_i$ a bound. Then for each $S \in Nilp_{E[[\zeta_1, \ldots, \zeta_n]]}$ there is a canonical fiber product diagram of categories
\[
 \begin{xy} \xymatrix{
   \nabla^{\mmu}_{(\hat{c}_i)} \MC{H}^1_U(C, G)(S) \ar[r] \ar_{\MF{L}}[d] & \nabla_0 \MC{H}^1_U(C \setminus \{c_i\}_i, G)(S) \ar^{\MF{L}}[d] \\
   \prod_i Sht_{G_{c_i}}^{\preceq \mu_i}(S) \ar[r] & \prod_i QSht_{G_{c_i}}(S)
 } \end{xy} 
\]
which is natural is $S$.
}

\prooof
Without the bounds, this follows from the repeated use of theorem \ref{Thm:UnifGluingOne} and the naturality assertion on the base scheme. That we may restrict to bounded subsets follows directly from the definition of a bounded global $G$-shtuka.  \exit

\rem{}{
i) All categories above define stacks over $\Spf E[[\zeta_1, \ldots, \zeta_n]]$. The stack $\nabla^{\mmu}_{(\hat{c}_i)} \MC{H}^1_U(C, G)$ is even representable by the formal completion of $\nabla^{\mmu}_n \MC{H}^1_U(C, G)$ along $\nabla^{\mmu}_{(c_i)} \MC{H}^1_U(C, G)$. However the other stacks seem not to admit any useful representablility result. \\
ii) Note that by general nonsense (i.e. any Artin-stack given by a functor from schemes to groupoids extends uniquely to a functor on all Artin-stacks), we still get a cartesian diagram of categories when considering $S$ over $\Spf E[[\zeta_1, \ldots, \zeta_n]]$ which are only DM-stacks (or even Artin stacks). In other words, the modification at characteristic places can be carried out even for families over DM-stacks.
}

\thm{\label{thm:UnifMorphism}}{}{
Let $S$ be any DM-stack over $\Spf E[[\zeta_1, \ldots, \zeta_n]]$ such that $(\zeta_1, \ldots, \zeta_n)$ is locally nilpotent on $S$. Fix a global $G$-shtuka with $U$-level structure $(\MS{G}_0, \varphi_0, \psi_0) \in \nabla_n\MC{H}^1_U(C, G)(S)$. For each characteristic place let $(\MC{G}_{c_i}, \varphi_{c_i}) \coloneqq \MF{L}_{c_i}(\MS{G}_0, \varphi_0)$ be the associated local $G$-shtukas and assume there exists an isomorphism $\beta_i: (\MC{G}_{c_i}, \varphi_{c_i}) \cong ({L^+G_{c_i}}_E, b_{c_i}\sigma^*)_S$ to a decent local $G_{c_i}$-shtuka $({L^+G_{c_i}}_E, b_{c_i}\sigma^*)$ over $\Sp E$ (as in definition \ref{def:RapoZinkDecent}). \\
Then one has a morphism of formal DM-stacks over $\Spf E[[\zeta_1, \ldots, \zeta_n]]$ (depending on the choice of $\beta_i$)
\[S \times_{\Spf E[[\zeta_1, \ldots, \zeta_n]]} \prod_i \MC{M}_{b_{c_i}}^{\preceq \mu_i} \to \nabla_n^{\mmu}\MC{H}^1_U(C, G) \times_{C^n \setminus \Delta} \widehat{C^n \setminus \Delta}^{(c_i)_i},\]
where $\widehat{C^n \setminus \Delta}^{(c_i)_i}$ denotes the formal completion of $C^n \setminus \Delta$ at $(c_i)_i$.
}

\prooof
Write $\MC{M}_{b_{c_i}}^{\preceq \mu_i} = \varinjlim_r M_{r, i}$ for some $M_{r, i} \in Nilp_{E[[\zeta_i]]}$. By pulling back the universal object of $M_{r, i}$ to $S'_r \coloneqq S \times_{\Spf E[[\zeta_1, \ldots, \zeta_n]]} \prod_i M_{r, i} \in Nilp_{E[[\zeta_1, \ldots, \zeta_n]]}$ via the canonical projection one gets for each $i$ a local $G_{c_i}$-shtuka $(\MC{G}'_{c_i}, \varphi'_{c_i})$ over $S'_r$ and a quasi-isogeny $\alpha_i: (\MC{G}'_{c_i}, \varphi'_{c_i}) \to ({L^+G_{c_i}}_E, b_{c_i}\sigma^*)_{S'_r}$. Therefore we may consider
\begin{itemize}
 \item The element $(\MS{G}_0, \varphi_0, \psi_0)_{S'_r} \in \nabla_0 \MC{H}^1_U(C \setminus \{c_i\}_i, G)(S'_r)$ defined by pulling back the global $G$-shtuka $(\MS{G}_0, \varphi_0, \psi_0)|_{(C \setminus \{c_i\}_i) \times_{\Fq} S}$ via $(C \setminus \{c_i\}_i) \times_{\Fq} S'_r \to (C \setminus \{c_i\}_i) \times_{\Fq} S$. 
 \item $(\MC{G}'_{c_i}, \varphi'_{c_i}) \in Sht_{G_{c_i}}^{\preceq \mu_i}(S'_r)$
 \item $\alpha_i^{-1} \circ \beta_i: \MC{L}\MF{L}_{c_i}(\MS{G}_0, \varphi_0)_{S'_r} = (\MC{L}\MC{G}_{c_i}, \varphi_{c_i})_{S'_r} \to ({LG_{c_i}}_E, b_{c_i}\sigma^*)_{S'_r} \to (\MC{L}\MC{G}'_{c_i}, \varphi'_{c_i})$ in $QSht_{G_{c_i}}(S'_r)$
\end{itemize}
This defines by corollary \ref{cor:UnifGluing} a global $G$-shtuka $(\MC{G}, \psi, \varphi) \in \nabla^{\mmu}_{(\hat{c}_i)} \MC{H}^1_U(C, G)(S'_r)$, i.e. a morphism
\[S \times_{\Spf E[[\zeta_1, \ldots, \zeta_n]]} \prod_i M_{r, i} = S'_r \to \nabla_n^{\mmu}\MC{H}^1_U(C, G) \times_{C^n \setminus \Delta} \widehat{C^n \setminus \Delta}^{(c_i)_i}\] 
As these morphisms are obviously compatible for varying $r$, we get the desired morphism of the theorem. \exit

\rem{}{
This morphism automatically induces a map between the perfections, i.e.
\[S \times_{\Spf E[[\zeta_1, \ldots, \zeta_n]]^\sharp} \prod_i \MC{M}_{b_{c_i}}^{\preceq \mu_i \sharp} \to \nabla_n^{\mmu}\MC{H}^1_U(C, G)^\sharp \times_{(C^n \setminus \Delta)^\sharp} \widehat{C^n \setminus \Delta}^{(c_i)_i \sharp},\]
for perfect DM-stacks $S$.
}

\subsection{Construction of the covering morphism}\label{subsec:CoverExistence}
Recall that in section \ref{sec:Igusa} we constructed DM-stacks $\op{Ig}^{(d_i) \sharp}_U$ parametrizing tuples consisting of
\begin{itemize}
 \item a global $G$-shtuka $(\MS{G}, \varphi, \psi)$ whose associated local $G$-shtukas at the characteristic places lie in fixed isomorphism classes,
 \item for each characteristic place $c_i$ an $I_{d_i}(b_{\nu_i})$-truncated isomorphism $\alpha_{d_i}: \MF{L}_{c_i}(\MS{G}, \varphi) \to (L^+G_{c_i}, b_{\nu_i}\sigma^*)$ between local $G_{c_i}$-shtukas with canonical $I_0(b_{\nu_i})$-structures.
\end{itemize}
Recall furthermore that we saw already in theorem \ref{Thm:RapoZinkBounded} the existence of Rapoport-Zink spaces, i.e. of a formal scheme $\MC{M}_{b_{\nu_i}}^{\preceq \mu_i \sharp}$ parametrizing
\begin{itemize}
 \item a local $G_{c_i}$-shtuka $(\MC{G}_i', \varphi_i')$ bounded by some $\mu_i$
 \item together with a quasi-isogeny $\beta_i: (\MC{G}_i', \varphi_i') \to (L^+G_{c_i}, b_{\nu_i}\sigma^*)$.
\end{itemize}
The underlying reduced subscheme will be denoted by $\MB{M}_{b_{\nu_i}}^{\preceq \mu_i \sharp}$. Using certain closed subspaces $\MB{M}_{b_{\nu_i}}^{\preceq \mu_i, \theta_i \sharp}$ (for finite $\Gamma$-invariant subsets of $X_*(T)_{\rm dom}$) defined in \ref{def:CoverDefSubspace}, we will now combine these spaces to construct (for sufficiently large $(d_i)$) a morphism
\[\pi_{(d_i)}: \prod_i \MB{M}_{b_{\nu_i}}^{\preceq \mu_i, \theta_i \sharp} \times_E \op{Ig}^{(d_i) \sharp}_U \to \MC{N}^{(\nu_i) \sharp}_U\]
into the corresponding Newton stratum of the special fiber of the perfection of $\nabla_n^{\mmu}\MC{H}^1_U(C, G)$. Here ``sufficiently large $(d_i)$'' is meant with respect to the partial order on $\M{Z}^n$ defined by $(d_i)_i \leq (d_i')_i$ if and only if $d_i \leq d_i'$ for each $i$. \\
In order to apply theorem \ref{thm:UnifMorphism} to define $\pi_{(d_i)}$, we have to lift the truncated isomorphism on the Igusa variety to an actual isomorphism of local $G_{c_i}$-shtukas. We have already seen that this can be done after a pro-finite \'etale cover. Thus our main objective in this section is to see, that the image of the uniformization morphism is independent of the choice of this lift.

\rem{}{
Until section \ref{sec:FormalLifting} we will work only over the special fiber. So all results of the previous section are applied only in the case $\zeta_1 = \ldots = \zeta_n = 0$. For more general versions of the morphism $\pi_{(d_i)}$ see sections \ref{subsec:FormalInfiniteCover} to \ref{subsec:AdicCoverMorphism}.
}

\const{\label{const:CoverMorphism}}{}{
Let us start by defining a morphism of DM-stacks over $\Sp E$
\[\pi_{(\infty_i)}: \prod_i \MB{M}_{b_{\nu_i}}^{\preceq \mu_i \sharp} \times_E \op{Ig}^{(\infty_i) \sharp}_U \to \MC{N}^{(\nu_i) \sharp}_U\]
where $\op{Ig}^{(\infty_i) \sharp}_U = \varprojlim_{(d_i)} \op{Ig}^{(d_i) \sharp}_U$. \\
We do this by describing the image of $S$-valued points for perfect schemes $S$, using the very same construction as in theorem \ref{thm:UnifMorphism}. As we frequently use the precise definition, let us describe it again: Pick $((\MC{G}_i', \varphi_i'), \beta_i) \in \MB{M}_{b_{\nu_i}}^{\preceq \mu_i \sharp}(S)$ (for each $i$) and $((\MS{G}, \varphi, \psi), (\alpha_{\infty_i})) \in \op{Ig}^{(\infty_i) \sharp}_U(S)$, where we keep the notation as above, but have actual isomorphisms $\alpha_{\infty_i}: \MF{L}_{c_i}(\MS{G}, \varphi) \to (L^+G_{c_i}, b_{\nu_i}\sigma^*)$ now by \ref{prop:IgusaTowerModuli}. This allows us to consider the objects
\begin{itemize}
 \item $(\MS{G}, \varphi, \psi)|_{(C \setminus \{c_i\}_i) \times_{\Fq} S} \in \nabla_0 \MC{H}^1_U(C \setminus \{c_i\}_i, G)(S)$, where $C \setminus \{c_i\}_i$ denotes again the open complement of all characteristic places.
 \item $((\MC{G}_i', \varphi_i')_i) \in \prod_i Sht_{G}^{\preceq \mu_i}(S)$
 \item for each $c_i$ an isomorphism $(\MC{L}\MC{G}_i', \varphi_i') \xrightarrow{\beta_i} (L^+G_{c_i}, b_{\nu_i}\sigma^*)_S \xrightarrow{\alpha_{\infty_i}^{-1}} \MF{L}_{c_i}((\MS{G}, \varphi)|_{(C \setminus \{c_i\}_i) \times_{\Fq} S})$
\end{itemize}
Thus we are in the situation of corollary \ref{cor:UnifGluing} and we get a global $G$-shtuka $(\stackrel{\sim}{\MS{G}}, \stackrel{\sim}{\varphi}, \psi)$ which is the modification of $(\MS{G}, \varphi, \psi)$ at the characteristic places. By construction its associated local $G_{c_i}$-shtuka at a characteristic place $c_i$ is isomorphic to $(\MC{G}_i', \varphi_i')$, hence bounded by $\mu_i$ (by definition of the Rapoport-Zink space). In particular $(\stackrel{\sim}{\MS{G}}, \stackrel{\sim}{\varphi}, \psi) \in \MB{X}^{\mmu}_U(S)$. Furthermore the local $G$-shtukas at characteristic points $c_i$ of $(\MS{G}, \varphi, \psi)$ and $(\stackrel{\sim}{\MS{G}}, \psi, \stackrel{\sim}{\varphi})$ are quasi-isogenous (by construction), hence both global $G$-shtukas lie in the same Newton stratum, i.e. in $\MC{N}^{(\nu_i)}_U(S)$. As $S$ was assumed to be perfect, it lifts uniquely to a point in $\MC{N}^{(\nu_i) \sharp}_U(S)$. This gives the desired image of the chosen points under $\pi_{(\infty_i)}$.
}

It remains to see that we may descend $\pi_{(\infty_i)}$ to some finite level, at least on quasi-compact subschemes of the source. To make this precise, let us introduce specific subschemes of the (reduced subscheme underlying) Rapoport-Zink spaces $\MC{M}_{b_{\nu_i}}^{\preceq \mu_i}$, and therefore as well of their perfections $\MC{M}_{b_{\nu_i}}^{\preceq \mu_i \sharp}$. As we wish them to cover the Rapoport-Zink space, we need a weaker definition of bounds for local $G_{c_i}$-shtukas:

\defi{\label{def:CoverWeakBound}}{}{
Let $\theta \subset X_*(T)_{\rm dom}$ a finite $\Gamma$-invariant subset of dominant cocharacters for the group $G_{c_i}$. \\
a) Let $\MC{G}$ and $\MC{G}'$ be two $L^+G_{c_i}$-torsors over a DM-stack $S \in Nilp_{\ENil}$ and $\alpha: \MC{LG} \to \MC{LG}'$ a morphism between the associated $LG_{c_i}$-torsors. Then $\alpha$ is weakly bounded by $\theta$ if there exists a finite (scheme theoretically) surjective morphism $Z \to S$ such that $G_{c_i}$ splits over $Z$ and $Z$ can be covered by open substacks $Z_\vartheta$ (parametrized by $\vartheta \in \theta$), such that the restriction of $\alpha$ to $Z_\vartheta$ is bounded by $\vartheta$. \\
b) A quasi-isogeny $\alpha: (\MC{G}, \varphi) \to (\MC{G}', \varphi')$ is weakly bounded by $\theta$ if the morphism $\alpha: \MC{LG} \to \MC{LG}'$ between the associated $LG$-torsors is bounded by $\theta$. 
}

\rem{\label{rem:CoverWeakComment}}{
i) Part a) of the definition should be thought of as requiring that for every $\ACFq[[t]]/(t^n)$-valued point $x$ of $S$ (for some $n$) there exists some element in $\theta$ such that $\alpha|_x$ is bounded by it, though this property is slightly stronger than being weakly bounded. Note that it is essential that we consider only thickenings of geometric points in \textit{one} direction. \\
ii) The splitting of $G_{c_i}$ over $Z$ does involve choices. But any two such splittings differ locally by an element of $\Gamma$. Moreover if $\alpha$ is (locally) bounded by some $\vartheta$, its pullback along the morphism induced by some $\gamma \in \Gamma$ is bounded by $\gamma.\vartheta$. Thus by $\Gamma$-invariance of $\theta$, the boundedness condition does not depend on the actual choice of the splitting. \\
iii) Combining proposition \ref{prop:CoverWeakProperty}b) and theorem \ref{thm:RapoZinkUnbounded} it follows that \ref{def:CoverWeakBound} defines in fact (for each $\theta$) a bound $\hat{Z}$ in the sense of \cite[definition 4.5]{HarRad1}. It can even be seen that any bound given by some $\hat{Z}$ equals, at least on reduced schemes (or DM-stacks), a weak bound defined by some $\theta$ (though bounds of the form $\hat{Z}$ likely allow for much more flexibility concerning the nilpotent structure).
}

\lem{\label{lem:CoverWeakFpqcLocal}}{}{
It suffices to check weak boundedness by $\theta$ fpqc-locally on $S$.
}

\prooof
Consider a finite field extension $E'$ such that $G_{c_i}$ splits over $E'$. Let $\theta = \{\vartheta_1, \ldots, \vartheta_k\} \subset X_*(T)_{\rm dom}$. Proposition \ref{prop:BoundClosed} applied to the pullback of $\alpha$ to $S \times_E \Sp E'$ and each of the coweights $\vartheta_i$ defines closed substack $\tilde{Z}_{\vartheta_i} \subset S \times_E \Sp E'$ such that $\alpha|_{\tilde{Z}_{\vartheta_i}}$ is bounded by $\vartheta_i$. Hence it remains to see that we may take
\[\tilde{Z} \coloneqq \coprod\nolimits_i \tilde{Z}_{\vartheta_i} \to S \times_E \Sp E' \to S\]
as the finite surjective morphism used above. Finiteness is obvious and surjectivity may be checked fpqc-locally. Thus we may assume that there is a finite surjective morphism $Z \to S$ as in definition \ref{def:CoverWeakBound}a). Then $Z \times_E \Sp E' \to S$ is again such a morphism and we may consider the closed subschemes $Z_{\vartheta_i} \subset Z \times_E \Sp E'$ given as the locus where $\alpha$ is bounded by $\vartheta_i$. It is clear that $\coprod_i Z_{\vartheta_i} \to Z \times_E \Sp E'$ is finite surjective. But by choice of the $\tilde{Z}_{\vartheta_i}$, each morphism $Z_{\vartheta_i} \to S \times_E \Sp E'$ factors over $\tilde{Z}_{\vartheta_i}$, i.e. $\coprod_i Z_{\vartheta_i} \to S \times_E \Sp E'$ factors over $\coprod_i \tilde{Z}_{\vartheta_i}$. Hence $\coprod_i \tilde{Z}_{\vartheta_i} \to S \times_E \Sp E'$ and then $\coprod_i \tilde{Z}_{\vartheta_i} \to S$ are surjective. \exit

\prop{\label{prop:CoverWeakProperty}}{}{
Let $\alpha: (\MC{G}, \varphi) \to (\MC{G}', \varphi')$ be a quasi-isogeny between local $G_{c_i}$-shtukas over a DM-stack $S$ over $E$. \\
a) If $\alpha$ is weakly bounded by $\theta = \{\vartheta_1, \ldots, \vartheta_k\}$, then $\alpha^{-1}$ is weakly bounded by $-\theta \coloneqq \{-\vartheta_1, \ldots, -\vartheta_k\}$. \\
b) If $\alpha$ is weakly bounded by $\theta$ and $\alpha': (\MC{G}', \varphi') \to (\MC{G}'', \varphi'')$ is weakly bounded by some $\theta'$, then the composition $\alpha' \circ \alpha$ is weakly bounded by $\theta \oplus \theta' = \{\vartheta + \vartheta' \,|\, \vartheta \in \theta, \vartheta' \in \theta'\}$. \\
c) Fix a finite $\Gamma$-invariant subset $\theta \subset X_*(T)_{\rm dom}$. Then the locus where $\alpha$ is weakly bounded by $\theta$ is representable by a closed immersion into $S$. \\
d) If $S$ is reduced, noetherian and quasi-compact, then there exists some finite $\Gamma$-invariant subset $\theta \subset X_*(T)_{\rm dom}$ such that $\alpha$ is weakly bounded by $\theta$.
}

\prooof
a) This follows directly from proposition \ref{Prop:BoundInverse}. \\
b) Choose finite surjective morphisms $Z \to S$ resp. $Z' \to S$ suitable for checking the weak boundedness of $\alpha$ resp. $\alpha'$. Then the assertion follows from lemma \ref{Lem:BoundTrivia}b) after pulling back along the finite surjective morphism $Z \times_S Z' \to S$. \\
c) Let $E'/E$ be a finite field extension such that $G_{c_i}$ splits over $E$. By the previous lemma the locus of weak boundedness is a fpqc-sheaf on $S$, hence it suffices to show representability on the fpqc-cover $S \times_E \Sp E'$. Let $\theta = \{\vartheta_1, \ldots \vartheta_k\}$. Then by proposition \ref{prop:BoundClosed} there is a closed immersion $Z_i \to S \times_E \Sp E'$ representing the locus where $\alpha$ is bounded by $\vartheta_i$. Let $Z \subset S \times_E \Sp E'$ be the scheme-theoretic image of the finite morphism $\coprod_i Z_i \to S \times_E \Sp E'$. By construction of $Z$ it is clear that it represents the locus of weak boundedness for $\theta$ on the scheme $S \times_E \Sp E'$. \\
d) By passing to some \'etale cover we may assume that $S$ is a scheme. By part c) it suffices to show, that $\alpha$ can be weakly bounded on each of the finitely many generic points. Thus assume wlog. that $S = \Sp k$ for some field $k$. Choose now $E'$ as in c). Then $\Gamma$ acts transitively on the connected components of $S \times_E \Sp E'$, which are all isomorphic to some field extension of $E'$. Fix one of these components $x \in S \times_E \Sp E'$. Then by lemma \ref{Lem:BoundFields}b), the restriction of $\alpha$ to $x$ is bounded by some $\vartheta \in X_*(T)_{\rm dom}$ (for the base field $E'$). Hence by $\Gamma$-action, we can bound $\alpha$ Zariski-locally on $S \times_E \Sp E'$ by elements in $\theta \coloneqq \Gamma \vartheta$ and $\alpha$ is indeed weakly bounded by $\theta$ over $S$. \exit

\defi{\label{def:CoverDefSubspace}}{}{
Let $\theta_i \subset X_*(T)_{\rm dom}$ be a finite $\Gamma$-invariant subset of dominant coweights. Then define $\MB{M}_{b_{\nu_i}}^{\preceq \mu_i, \theta_i} \subset \MB{M}_{b_{\nu_i}}^{\preceq \mu_i}$ to be the closed subscheme consisting of all points where the universal quasi-isogeny $\beta_i^{univ}$ is weakly bounded by $\theta_i$. We endow $\MB{M}_{b_{\nu_i}}^{\preceq \mu_i, \theta_i}$ with the reduced subscheme structure, even though weak boundedness by $\theta_i$ already defines it as a closed subscheme of $\MB{M}_{b_{\nu_i}}^{\preceq \mu_i}$ which usually has non-trivial nilpotent structure. \\
Let $\MB{M}_{b_{\nu_i}}^{\preceq \mu_i, \theta_i \sharp}$ be the perfection of $\MB{M}_{b_{\nu_i}}^{\preceq \mu_i, \theta_i}$ or equivalently the locus of weak boundedness by $\theta_i$ inside $\MB{M}_{b_{\nu_i}}^{\preceq \mu_i \sharp}$. 
}

\rem{\label{rem:CoverDefSubspaceProper}}{
i) Consider any quasi-compact reduced open subscheme $Z \subset \MB{M}_{b_{\nu_i}}^{\preceq \mu_i}$. Then by \ref{prop:CoverWeakProperty}d) the universal quasi-isogeny over $Z$ can be weakly bounded by some $\theta_i$. Hence we have
\[\MB{M}_{b_{\nu_i}}^{\preceq \mu_i} = \bigcup_{\theta_i} \MB{M}_{b_{\nu_i}}^{\preceq \mu_i, \theta_i}\]
and the subspaces on the right-hand side define a filtered system of closed subspaces. \\
ii) The proof of \cite[theorem 6.3]{HaVi} shows, that $\MB{M}_{b_{\nu_i}}^{\preceq \mu_i, \theta_i}$ is proper, as it admits a closed immersion into one of the proper spaces called $\MC{M}^n$ in loc. cit. Moreover as in the proof of \cite[corollary 6.5]{HaVi}, it may be shown that each irreducible component of $\MB{M}_{b_{\nu_i}}^{\preceq \mu_i, \theta_i}$ is in fact projective.
}

\prop{\label{prop:CoverMorphDescentPoint}}{}{
Let $S$ be any perfect simply connected scheme over $E$. Then for $(d_i)_i$ sufficiently large (depending on the $\theta_i$) there is a unique factorization
\[\begin{xy}
 \xymatrix {
    \prod_i \MB{M}_{b_{\nu_i}}^{\preceq \mu_i, \theta_i \sharp}(S) \times_E \op{Ig}^{(\infty_i) \sharp}_U(S) \ar^{\pi_{(\infty_i)}}[dr] \ar[dd] & \\
   &  \MC{N}^{(\nu_i) \sharp}_U(S) \\
    \prod_i \MB{M}_{b_{\nu_i}}^{\preceq \mu_i, \theta_i \sharp}(S) \times_E \op{Ig}^{(d_i) \sharp}_U(S) \ar_{\pi_{(d_i)}}[ur] }
\end{xy} \]
of maps between $S$-valued points.
}

\prooof
The vertical map is surjective as it is induced from an (infinite) Galois-cover and $\pi_1(S)$ is trivial by assumption. Hence we have to show, that $\pi_{(\infty_i)}$ is constant on fibers of the vertical map: \\
Consider any element $x \in \prod_i \MB{M}_{b_{\nu_i}}^{\preceq \mu_i, \theta_i \sharp}(S) \times \op{Ig}^{(d_i) \sharp}_U(S)$ (where we start with arbitrary $d_i$'s and specify them later on). As usual we may describe $x$ by $((\MC{G}_i', \varphi_i'), \beta_i) \in \MB{M}_{b_{\nu_i}}^{\preceq \mu_i, \theta_i \sharp}(S)$ and $((\MS{G}, \varphi, \psi), (\alpha_{d_i})) \in \op{Ig}^{(d_i)}_U(S)$. Any two preimages $x_1$ and $x_2$ of $x$ under the Galois-cover are given by the very same datum, except for $(\alpha_{d_i})$ being replaced by a lift $(\alpha_{\infty_i, 1})$ resp. $(\alpha_{\infty_i, 2})$ to actual isomorphisms of local $G_{c_i}$-shtukas. Thus we have to see that the definition of $\pi_{(\infty_i)}$ does not depend on the choice of this lift. 
To simplify notations let us only consider the modification of the global $G$-shtuka at just one characteristic place and denote the resulting global $G$-shtukas by $(\stackrel{\sim}{\MS{G}}_1, \stackrel{\sim}{\varphi}_1, \psi)$ resp. $(\stackrel{\sim}{\MS{G}}_2, \stackrel{\sim}{\varphi}_2, \psi)$. The general case follows immediately by applying this argument for each characteristic place separately. \\
Consider the following diagram
\[\begin{xy}
 \xymatrix {
  (\stackrel{\sim}{\MS{G}}_1, \stackrel{\sim}{\varphi}_1, \psi) \ar@{-->}[dr] \ar@{~>}[ddd] \ar@{~>}[rrr] & & & (\MS{G}, \varphi, \psi) \ar@{=}[dr] \ar@{~>}[dd] &  \\
   & (\stackrel{\sim}{\MS{G}}_2, \stackrel{\sim}{\varphi}_2, \psi) \ar@{~>}[ddd] \ar@{~>}[rrr] & & & (\MS{G}, \varphi, \psi) \ar@{~>}[dd] \\
   & & & \MF{L}_{c_i}(\MS{G}, \varphi) \ar@{=}[dr] &  \\
   (\MC{G}_i', \varphi_i') \ar@{-->}[dr] \ar@{~>}[rr] & & (\MC{L}\MC{G}_i', \varphi_i') \ar_{\alpha_{\infty_i, 1}^{-1} \circ \beta_i}[ur] \ar_(0.4){\beta_i^{-1} \circ \alpha_{\infty_i, 2} \circ \alpha_{\infty_i, 1}^{-1} \circ \beta_i \quad}[dr] & & \MF{L}_{c_i}(\MS{G}, \varphi) \\
   & (\MC{G}_i', \varphi_i') \ar@{~>}[rr] & & (\MC{L}\MC{G}_i', \varphi_i') \ar_{\alpha_{\infty_i, 2}^{-1} \circ \beta_i}[ur] &  }
\end{xy} \]
where the squiggly arrows are given by functors and the dashed arrows are not yet established. \\
We first prove that, after choosing $d_i$ appropriately, the quasi-isogeny $\beta_i^{-1} \circ \alpha_{\infty_i, 2} \circ \alpha_{\infty_i, 1}^{-1} \circ \beta_i$ defines in fact an automorphism of the local $G$-shtuka $(\MC{G}_i', \varphi_i')$: Note first that $\alpha_{\infty_i, 1}$ and $\alpha_{\infty_i, 2}$ lie in one $I_{d_i}(b_{\nu_i})$-double coset (assuming for simplicity that we have chosen trivializations for the $G_{c_i}$-torsors in question), because they both lift the $I_{d_i}(b_{\nu_i})$-truncated isomorphism $\alpha_{d_i}$. Hence by lemma \ref{prop:TruncTechnicalLemma} (or fpqc-locally by the very definition) we find that $\alpha_{\infty_i, 2} \circ \alpha_{\infty_i, 1}^{-1} \in I_{d_i}(b_{\nu_i})(S) \subset K_{d_i}(S)$. \\
Furthermore by definition of $\MB{M}_{b_{\nu_i}}^{\preceq \mu_i, \theta_i \sharp}$, $\beta_i$ is weakly bounded by $\theta_i$. Thus the following lemma ensures that for sufficiently large $d_i$, we have 
\[\beta_i^{-1} \circ \alpha_{\infty_i, 2} \circ \alpha_{\infty_i, 1}^{-1} \circ \beta_i \in \beta_i^{-1} \cdot K_{d_i}(S) \cdot \beta_i \subset L^+G_{c_i}(S)\]
Now the automorphism $\beta_i^{-1} \circ \alpha_{\infty_i, 2} \circ \alpha_{\infty_i, 1}^{-1} \circ \beta_i: (\MC{G}_i', \varphi_i') \to (\MC{G}_i', \varphi_i')$ defines together with the identity $id: (\MS{G}, \varphi, \psi)|_{C \setminus \{c_i\}_i} \to (\MS{G}, \varphi, \psi)|_{C \setminus \{c_i\}_i}$ an isomorphism
\[(\stackrel{\sim}{\MS{G}}_1, \stackrel{\sim}{\varphi}_1, \psi) \to (\stackrel{\sim}{\MS{G}}_2, \stackrel{\sim}{\varphi}_2, \psi)\]
by the fiber product diagram of categories in corollary \ref{cor:UnifGluing}. Hence we get as desired
\[\pi_{(\infty_i)}(x_1) = \pi_{(\infty_i)}(x_2)\] \exit

\lem{}{}{
Let $\theta_i$ be a $\Gamma$-orbit of dominant cocharacters. Then for sufficiently large $d_i$ (depending only on $\theta_i$) the following holds: For any $L^+G_{c_i}$-torsor $\MC{G}$ over a scheme (or even DM-stack) $S$ and any 
\begin{itemize}
 \item $\alpha \in \Aut^{d_i}(\MC{G})$, i.e. a morphism $\alpha: \MC{G} \to \MC{G}$ which induces the identity on $\MC{G} \times^{L^+G_{c_i}} L^+G_{c_i}/K_{d_i}$
 \item and $\beta_i \in \Aut(\MC{L}\MC{G})$ (between associated $LG_{c_i}$-torsors) which is bounded by $\theta_i$,
\end{itemize}
the isomorphism of $LG_{c_i}$-torsors
\[\beta_i^{-1} \circ \alpha \circ \beta_i: \MC{L}\MC{G} \to \MC{L}\MC{G}\]
actually restricts to an automorphism of the $L^+G_{c_i}$-torsor $\MC{G}$.
}

\prooof
Let us first prove the lemma for $G_{c_i} = GL_n$: \\
As the $d_i$ defined below is independent of $S$ and $\MC{G}$, we may pass to some \'etale cover of $S$ and assume wlog. that $\MC{G} = {L^+GL_n}_S$ is the trivial torsor. Then we may represent $\alpha \in K_{N_i}(S)$ and $\beta_i \in LGL_n(S) = GL_n(S((z)))$. As $\det(\beta_i^{-1} \circ \alpha \circ \beta_i) = \det(\alpha) \in (S[[z]])^\times$, we have to see that we may choose $d_i$ in such a fashion, that $\beta_i^{-1} \circ \alpha \circ \beta_i \in GL_n(S((z)))$ has coefficients in $S[[z]]$. 
As $GL_n$ is split, $\Gamma$ acts trivially on $X_*(T)$ and $\theta_i = \{\vartheta_i\}$ contains only one element. Thus $\beta_i$ is bounded by $\vartheta_i$ over $S$ and by proposition \ref{Prop:BoundInverseGL} the morphism $\beta_i^{-1}$ is bounded by $(-\vartheta_i)_{\rm dom}$. Hence using \cite[lemma 4.3]{HaVi} or equivalently lemma \ref{Lem:BoundGLEasy} (both for $i = 1$), there are constants $d_i'$ and $d_i''$ (depending only on $\theta_i$) such that both $z^{d_i'}\beta_i$ and $z^{d_i''} \beta_i^{-1}$ have coefficients in $S[[z]]$. Thus take any (positive) $d_i \geq d_i' + d_i''$. Then writing $\alpha \in K_{d_i}(S)$ as $k = 1 + z^{d_i' + d_i''}\alpha'$, where $\alpha' \in \op{Mat}_{n \times n}(S[[z]])$, we get
\[\beta_i^{-1} \circ \alpha \circ \beta_i = 1 + z^{d_i''} \beta_i^{-1} \cdot \alpha' \cdot z^{d_i'}\beta_i\]
which obviously has coefficients in $S[[z]]$. Hence the lemma is shown for $G_{c_i} = GL_n$. \\
We come now to the general case: It suffices to check the assertion after pullback along a finite surjective morphism. Hence we may assume that $G_{c_i}$ splits over $S$ (and we fix such a splitting). Moreover as $\theta_i$ is a finite set, we may assume that $\beta_i$ is bounded by a single element $\vartheta_i \in \theta_i \subset X_*(T)_{\rm dom}$. Consider then the representation $\rho: G_{c_i} \to GL \coloneqq GL(\bigoplus_{\lambda} V_{G_{c_i}}(\lambda))$, where $V_{G_{c_i}}(\lambda)$ is the Weyl module with highest weight $\lambda$ and the direct sum runs over a finite generating system in $X^*(T)_{\rm dom}$. By \cite[proposition 3.14]{HaVi} the map $\rho$ is a closed immersion. 
Furthermore lemma \ref{Lem:BoundNiceBasis} (except for property ii)) and lemma \ref{Lem:BoundAfterRepr}a) hold even for $\bigoplus_{\lambda} V_{G_{c_i}}(\lambda)$ instead of $V_{G_{c_i}}(\lambda)$ (using the very same proofs). Thus there is a maximal torus and a Borel $T_{GL} \subset B_{GL} \subset GL$ containing the images of the corresponding subgroups of $G_{c_i}$ such that ${\beta_i}_{GL}: \MC{L}\MC{G} \times^{LG_{c_i}} LGL \to \MC{L}\MC{G} \times^{LG_{c_i}} LGL$ is bounded by $\rho \circ \theta_i \in X_*(T_{GL})$. 
Furthermore it is obvious that any $\alpha \in \Aut^{d_i}(\MC{G})$ gives an automorphism $\alpha_{GL}: \MC{G} \times^{L^+G_{c_i}} L^+GL \to \MC{G} \times^{L^+G_{c_i}} L^+GL$, which actually lies in $\Aut^{d_i}(\MC{G} \times^{L^+G_{c_i}} L^+GL)$. Thus by the first part of this proof, we may choose some constant $d_i$ (depending only on $\rho \circ \vartheta_i$) such that for any $\alpha$ and $\beta_i$ as in the statement, we get an isomorphism of $L^+GL$-torsors
\[(\beta_i^{-1} \circ \alpha \circ \beta_i)_{GL} = {\beta_i}_{GL}^{-1} \circ \alpha_{GL} \circ {\beta_i}_{GL} \in \Aut(\MC{G} \times^{L^+G_{c_i}} L^+GL).\]
Using the faithfulness of $\rho$, this implies that already 
\[\beta_i^{-1} \circ \alpha \circ \beta_i \in \Aut(\MC{G})\]
as desired. \exit

\rem{}{
This lemma is a (slight) generalization of the following well-known statement: \\
\textit{
Let $\theta_i$ be any finite $\Gamma$-invariant set of dominant cocharacters. Then there is a constant $d_i$ such that for any algebraically closed field $k$ (over $\Fq$), any $g \in K_{d_i}(k)$ and any $b \in LG_{c_i}(k)$ such that its Hodge point satisfies $\mu(b) \preceq \vartheta$ (using the partial order in $X_*(T)$) for some $\vartheta \in \theta_i$, we have:
\[b^{-1} \cdot g \cdot b \in L^+G_{c_i}(k)\]
} \\
As (over an algebraically closed field) being bounded by $\theta_i$ is actually equivalent to satisfying $\mu(b) \preceq \vartheta$ for some $\vartheta \in \theta_i$ (cf. lemma \ref{Lem:BoundFields}), this statement is indeed a special case of the lemma above.
}

\prop{\label{prop:CoverMorphDescent}}{}{
For $(d_i)_i$ sufficiently large (depending only on $\theta_i$), $\pi_{(\infty_i)}$ factors over a morphism
\[\pi_{(d_i)}: \prod_i \MB{M}_{b_{\nu_i}}^{\preceq \mu_i, \theta_i \sharp} \times_E \op{Ig}^{(d_i) \sharp}_U \to \MC{N}^{(\nu_i) \sharp}_U\]
}

\prooof
The map
\[\prod_i \MB{M}_{b_{\nu_i}}^{\preceq \mu_i, \theta_i \sharp} \times_E \op{Ig}^{(\infty_i) \sharp}_U \to \prod_i \MB{M}_{b_{\nu_i}}^{\preceq \mu_i, \theta_i \sharp} \times_E \op{Ig}^{(d_i) \sharp}_U\]
is an infinite Galois-cover. Hence it suffices to check the factorization of $\pi_{(\infty_i)}$ for formal neighborhoods of geometric points. But this was done in the previous proposition. \exit

\lem{\label{lem:CoverInfiniteSurjective}}{}{
The morphism $\pi_{(\infty_i)}$ is surjective. In fact for sufficiently large $\theta_i$ (depending on $\nu_i$ and the fixed cocharacters $\mu_i$), even the morphism 
\[\pi_{(d_i)}: \prod_i \MB{M}_{b_{\nu_i}}^{\preceq \mu_i, \theta_i \sharp} \times_E \op{Ig}^{(d_i) \sharp}_U \to \MC{N}^{(\nu_i) \sharp}_U\] 
is surjective.
}

\prooof
As all stacks are reduced, it suffices to see this on geometric points. In particular it suffices to show surjectivity of 
\[\pi_{(\infty_i)}: \prod_i \MB{M}_{b_{\nu_i}}^{\preceq \mu_i, \theta_i \sharp} \times_E \op{Ig}^{(\infty_i) \sharp}_U \to \MC{N}^{(\nu_i) \sharp}_U\] 
for sufficiently large $\theta_i$. Thus consider a global $G$-shtuka $(\MS{G}_0, \varphi_0, \psi_0)$ in the Newton stratum. By definition there is for each characteristic place $c_i$ some quasi-isogeny $\beta_i: \MF{L}_{c_i}(\MS{G}_0, \varphi_0) \to (L^+G_{c_i}, b_{\nu_i}\sigma^*)$. By \cite[theorem 1.4]{RapoZinkBTbuilding} this $\beta_i$ may be chosen to be weakly bounded by $\theta_i$ for suitably chosen $\theta_i$. Then we may use corollary \ref{cor:UnifGluing} again to glue $(\MS{G}_0, \varphi_0, \psi_0)|_{C \setminus \{c_i\}_i}$ along the $\beta_i$ to the local $G_{c_i}$-shtukas $(L^+G_{c_i}, b_{\nu_i}\sigma^*)$. 
In this way we obtain a new global $G$-shtuka $(\MS{G}, \varphi, \psi)$ with trivializations $\alpha_{\infty_i} = \id: \MF{L}_{c_i}(\MS{G}, \varphi) \to (L^+G_{c_i}, b_{\nu_i}\sigma^*)$. It is now clear from these constructions that the points defined by $(\MF{L}_{c_i}(\MS{G}_0, \varphi_0), \beta_i)$ inside $\MB{M}_{b_{\nu_i}}^{\preceq \mu_i, \theta_i \sharp}$ and by $(\MS{G}, \varphi, \psi, (\alpha_{\infty_i})_i)$ inside $\op{Ig}^{(\infty_i) \sharp}_U$ define a preimage of the given point under $\pi_{(\infty_i)}$. \exit

\rem{}{
We give a detailed description of the fibers of $\pi_{(\infty_i)}$ in proposition \ref{prop:JActionFibers}.
}

\subsection{Action of quasi-isogenies on Igusa varieties}\label{subsec:JAction}
We define and study the action of the group of self-quasi-isogenies of the local $G_{c_i}$-shtukas $(L^+G_{c_i}, b_{\nu_i} \sigma^*)$ on Igusa varieties. The corresponding actions on Shimura varieties and their Igusa varieties are defined in \cite[sections 3.4 and 4.3.1]{MantoFoliation} or \cite[lemma 5]{MantoFoliationPEL}. We apply this action to show that $\pi_{(d_i)}$ is quasi-finite.

\defi{\label{def:JGroup}}{}{
Let $b_{\nu_i} \in LG_{c_i}(E)$ be (as usual) a fundamental alcove defined over the reflex field $E$. Fix an integer $s \geq 1$ such that $b_{\nu_i}$ is decent for $s$, i.e. $(b_{\nu_i}\sigma)^s = z^{s\nu_i'}$ for the Newton point $\nu_i'$. Let $E' = E \cdot \M{F}_{q^s}$. Then define
\[J_i = QIsog((L^+G_{c_i}, b_{\nu_i} \sigma^*)_{E'}, (L^+G_{c_i}, b_{\nu_i} \sigma^*)_{E'}) = \{g \in LG_{c_i}(E') \;|\; g^{-1} b_{\nu_i} \sigma(g) = b_{\nu_i}\}\]
We endow $J_i$ with the $z$-adic topology via its canonical inclusion into $LG_{c_i}(E')$.
}

\rem{}{
Consider the functor on noetherian $\Fq((z))$-algebras $R$ given by
\[J_i(R) = \{g \in G_{c_i}(R \otimes_{\Fq((z))} E'((z))) \;|\; g^{-1} b_{\nu_i} \sigma(g) = b_{\nu_i}\}\]
where $\sigma$ acts trivially on $R$. By literally the same proof as in \cite[proposition 1.12]{RapoZinkSpaces}, this defines a smooth affine group scheme over $\Fq((z))$. Its $\Fq((z))$-valued points identify with $J_i$.
}

\lem{}{}{
If $E''/E'$ is any field extension then 
\[J_i = QIsog((L^+G_{c_i}, b_{\nu_i} \sigma^*)_{E''}, (L^+G_{c_i}, b_{\nu_i} \sigma^*)_{E''})\]
}

\prooof
As above the right-hand side equals $\{g \in LG_{c_i}(E'') \;|\; g^{-1} b_{\nu_i} \sigma(g) = b_{\nu_i}\}$. By the decency equation such an element satisfies $g^{-1} z^{s\nu_i'} \sigma^s(g) = z^{s\nu_i'}$ as well. Now it is easy to see, that $g$ has to lie in the centralizer of $\nu_i'$ (using that $z$ is $\sigma$-invariant) and hence satisfies $g = \sigma^s(g)$. In particular $g \in LG_{c_i}(E_s)$ and the quasi-isogeny associated to $g$ is already defined over $E'$, i.e. lies in $J_i$. The other inclusion is trivial. \exit

\prop{}{}{
There is a canonical action of $\prod_i J_i$ on the Igusa variety $\op{Ig}^{(\infty_i) \sharp}_U \times_E \Sp E'$.
}

\prooof
Let us first explain how one single factor $J_j$ acts: Let $S$ be any perfect DM-stack over $E'$ and $\gamma_j \in J_j$ viewed by base change as a quasi-isogeny over $S$. Consider a $S$-valued point in $\op{Ig}^{(\infty_i) \sharp}_U \times_E \Sp E'$, i.e. a global $G$-shtuka $(\MS{G}, \varphi, \psi)$ over $S$ and for each characteristic place $c_i$ an isomorphism $\alpha_{(\infty_i)}: \MF{L}_{c_i}(\MS{G}, \varphi) \to (L^+G_{c_i}, b_{\nu_i} \sigma^*)_S$. Then $\alpha_{(\infty_j)}^{-1} \circ \gamma_j \circ \alpha_{(\infty_j)}$ defines a self-quasi-isogeny of $\MF{L}_{c_j}(\MS{G}, \varphi)$.
Thus using theorem \ref{Thm:UnifGluingOne} to glue the global $G$-shtuka $(\MS{G}, \varphi, \psi)$ and the local $G_{c_j}$-shtuka $\MF{L}_{c_j}(\MS{G}, \varphi)$ along the (non-trivial) quasi-isogeny $\alpha_{(\infty_j)}^{-1} \circ \gamma_j \circ \alpha_{(\infty_j)}$ one obtains a new global $G$-shtuka $(\MS{G}', \varphi', \psi')$. As by construction $\MF{L}_{c_j}(\MS{G}', \varphi') = \MF{L}_{c_j}(\MS{G}, \varphi)$, we have a canonical trivialization $\alpha_{(\infty_j)}: \MF{L}_{c_j}(\MS{G}', \varphi') \to (L^+G_{c_j}, b_{\nu_j} \sigma^*)_S$. For all other places $c_i \neq c_j$ the local $G_{c_i}$-shtuka does not change and we may simply keep the trivializations $\alpha_{(\infty_i)}$. Thus the quadrupel $(\MS{G}', \varphi', \psi', (\alpha_{(\infty_i)})_i)$ defines a point in $\op{Ig}^{(\infty_i) \sharp}_U \times_E \Sp E'$ and we set
\[\gamma_j.(\MS{G}, \varphi, \psi, (\alpha_{(\infty_i)})_i) \coloneqq (\MS{G}', \varphi', \psi', (\alpha_{(\infty_i)})_i).\]
As this construction is functorial in $S$, we get indeed the desired $J_j$-action. \\
All the actions of the groups $J_j$ commute, because changing a global $G$-shtuka at different places is independent of the order. Hence they induce an action of $\prod_i J_i$. \exit

\prop{\label{prop:JActionDescent}}{}{
Fix some element $(\gamma_i)_i \in \prod_i J_i$. Then there is a tuple $(d_{i, \gamma_i})_i$ depending only on $(\gamma_i)_i$, such that there exists for every tuple $(d_i)_i$ a unique morphism
\[(\gamma_i)_i: \op{Ig}^{(d_i + d_{i, \gamma_i}) \sharp}_U \times_E \Sp E' \to \op{Ig}^{(d_i) \sharp}_U \times_E \Sp E'\]
such that the diagram
\[\begin{xy}
 \xymatrix {
  \op{Ig}^{(\infty_i) \sharp}_U \times_E \Sp E' \ar^-{(\gamma_i)_i}[r] \ar[d] & \op{Ig}^{(\infty_i) \sharp}_U \times_E \Sp E' \ar[d]  \\
  \op{Ig}^{(d_i + d_{i, \gamma_i}) \sharp}_U \times_E \Sp E' \ar^-{(\gamma_i)_i}[r] & \op{Ig}^{(d_i) \sharp}_U \times_E \Sp E' }
\end{xy} \]
commutes. Here the upper horizontal morphism is given by the $\prod_i J_i$-action defined in the previous proposition and the vertical maps are the canonical projections. 
}

\prooof
The quasi-isogeny $\gamma_i$ is bounded by some $\Gamma$-invariant finite subset of cocharacters. For such elements, we showed in the last section that there is some $d_{i, \gamma_i}$ such that $\gamma_i \cdot I_{b_{\nu_i}}(d_{i, \gamma_i}) \cdot \gamma_i^{-1} \subset I_{b_{\nu_i}}(0)$. This automatically implies for each $d_i \geq 0$
\[\gamma_i \cdot I_{b_{\nu_i}}(d_i + d_{i, \gamma_i}) \cdot \gamma_i^{-1} \subset I_{b_{\nu_i}}(d_i).\]
As in proposition \ref{prop:CoverMorphDescent} it suffices to check the assertion on points with values in perfect simply connected schemes $S$. Thus let $S$ be such a scheme over $E'$ and consider now two preimages $(\MS{G}, \varphi, \psi, (\alpha_{(\infty_i) 1})_i)$ and $(\MS{G}, \varphi, \psi, (\alpha_{(\infty_i) 2})_i)$ in $\op{Ig}^{(\infty_i) \sharp}_U$ of one $S$-valued point in $\op{Ig}^{(d_i + d_{i, \gamma_i}) \sharp}_U \times_E \Sp E'$. Denote their images under $(\gamma_i)_i$ by $(\MS{G}'_1, \varphi'_1, \psi'_1, (\alpha_{(\infty_i) 1})_i)$ respectively $(\MS{G}'_2, \varphi'_2, \psi'_2, (\alpha_{(\infty_i) 2})_i)$. \\
We proceed as in the proof of proposition \ref{prop:CoverMorphDescentPoint}: As $\alpha_{(\infty_i) 1}$ and $\alpha_{(\infty_i) 2}$ differ only by an element in $I_{b_{\nu_i}}(d_i + d_{i, \gamma_i})$ the quasi-isogeny 
\begin{align*}
 \delta_i & \coloneqq (\alpha_{(\infty_i) 1}^{-1} \gamma_i \alpha_{(\infty_i) 1})^{-1} \circ (\alpha_{(\infty_i) 2}^{-1} \gamma_i \alpha_{(\infty_i) 2}) \\ 
 & \, = \alpha_{(\infty_i) 1}^{-1} \circ \gamma_i^{-1} (\alpha_{(\infty_i) 1}\alpha_{(\infty_i) 2}^{-1}) \gamma_i \circ \alpha_{(\infty_i) 2}:  \MF{L}_{c_i}(\MS{G}'_1, \varphi'_1) \to \MF{L}_{c_i}(\MS{G}'_2, \varphi'_2)
\end{align*}
is by choice of $d_{i, \gamma_i}$ an isomorphism. Hence we may glue the identity on $(\MS{G}, \varphi, \psi)|_{(C \setminus \{c_i\}_i) \times_{\Fq} S}$ and the isomorphisms $\delta_i$ (for each $i$) to an isomorphism
\[(\MS{G}'_1, \varphi'_1, \psi'_1) \cong (\MS{G}'_2, \varphi'_2, \psi'_2).\]
We have to see that under this isomorphism, the trivializations $\alpha_{(\infty_i) 1}$ and $\alpha_{(\infty_i) 2}$ define at each characteristic place $c_i$ the same $I_{b_{\nu_i}}(d_i)$-truncation class. For this note that we have the commutative diagram
\[\begin{xy}
 \xymatrix {
  & \MF{L}_{c_i}(\MS{G}, \varphi) \ar^{\alpha_{(\infty_i) 2}}[rr] \ar^(.3){\alpha_{(\infty_i) 2}^{-1} \gamma_i \alpha_{(\infty_i) 2}}[dd]  && (L^+G_{c_i}, b_{\nu_i} \sigma^*)_S \ar^{\gamma_i}[dd] \\
  \MF{L}_{c_i}(\MS{G}, \varphi) \ar^(.3){\alpha_{(\infty_i) 1}}[rr] \ar_{\alpha_{(\infty_i) 1}^{-1} \gamma_i \alpha_{(\infty_i) 1}}[dd] \ar@{=}[ur] && (L^+G_{c_i}, b_{\nu_i} \sigma^*)_S \ar^(.3){\gamma_i}[dd] \ar_{\quad \alpha_{(\infty_i) 2} \alpha_{(\infty_i) 1}^{-1}}[ur] \\
  & \MF{L}_{c_i}(\MS{G}'_2, \varphi'_2) \ar^(.3){\alpha_{(\infty_i) 2}}[rr]  && (L^+G_{c_i}, b_{\nu_i} \sigma^*)_S \\
  \MF{L}_{c_i}(\MS{G}'_1, \varphi'_1) \ar^{\alpha_{(\infty_i) 1}}[rr] \ar^{\delta_i}[ur] && (L^+G_{c_i}, b_{\nu_i} \sigma^*)_S \ar_{\qquad \gamma_i \alpha_{(\infty_i) 2} \alpha_{(\infty_i) 1}^{-1} \gamma_i^{-1}}[ur] }
\end{xy} \]
But as $\alpha_{(\infty_i) 1}$ and $\alpha_{(\infty_i) 2}$ are both lifts of one $I_{b_{\nu_i}}(d_i + d_{i, \gamma})$-truncation class, we have by choice of $d_{i, \gamma_i}$ that $\gamma_i \alpha_{(\infty_i) 2} \alpha_{(\infty_i) 1}^{-1} \gamma_i^{-1}$ induces the identity modulo $I_{b_{\nu_i}}(d_i)$. This shows that $(\MS{G}'_1, \varphi'_1, \psi'_1, (\alpha_{(\infty_i) 1})_i)$ and $(\MS{G}'_2, \varphi'_2, \psi'_2, (\alpha_{(\infty_i) 2})_i)$ indeed define the same element in $\op{Ig}^{(d_i) \sharp}_U \times_E \Sp E'$. \exit

\lem{\label{lem:JActionIgusaFinite}}{}{
With the notations of the previous proposition, the morphism
\[(\gamma_i)_i: \op{Ig}^{(d_i + d_{i, \gamma_i}) \sharp}_U \times_E \Sp E' \to \op{Ig}^{(d_i) \sharp}_U \times_E \Sp E'\]
is finite \'etale.
}

\prooof
Observe that we have a commutative diagram
\[\begin{xy}
 \xymatrix @C=0.5pc {
  \op{Ig}^{(d_i + d_{i, \gamma_i}) \sharp}_U \times_E \Sp E' \ar^-{\gamma_i}[rr] \ar[dr] && \op{Ig}^{(d_i) \sharp}_U \times_E \Sp E' \ar[dl]  \\
  & \MC{C}_U^{\sharp} \times_E \Sp E' & }
\end{xy} \]
where the projections to the central leaf are finite \'etale. Hence the lemma follows from the cancellation property of finite \'etale morphisms. \exit

As for Rapoport-Zink spaces of $p$-divisible groups, we have

\lem{\label{lem:JActionRZ}}{}{
a) There is a canonical (left) action on $J_i$ on the Rapoport-Zink space $\MC{M}_{b_{\nu_i}}^{\preceq \mu_i} \times_E \Sp E'$ and on its underlying reduced fiber $\MB{M}_{b_{\nu_i}}^{\preceq \mu_i} \times_E \Sp E'$. \\
b) For every $\gamma_i \in J_i$ there exists a finite $\Gamma$-invariant subset $\theta_{i, \gamma_i} \subset X_*(T)$ such that the action of a) restricts to a morphism
\[\gamma_i: \MB{M}_{b_{\nu_i}}^{\preceq \mu_i, \theta_i} \times_E \Sp E' \to \MB{M}_{b_{\nu_i}}^{\preceq \mu_i, \theta_i \oplus \theta_{i, \gamma_i}} \times_E \Sp E'\]
Here $\theta_i \oplus \theta_{i, \gamma_i} = \{\vartheta + \vartheta' \;|\; \vartheta \in \theta_i, \vartheta' \in \theta_{i, \gamma_i}\}$ as in \ref{prop:CoverWeakProperty}b). \\
c) In the situation of part b), the morphism $\gamma_i$ is a closed immersion. \\
All properties hold as well after taking the perfection.
}

\prooof
a) The action is given by postcomposition of the universal quasi-isogeny with the self-quasi-isogeny defined by the element in $J_i$. \\
b) Choose $\theta_{i, \gamma_i}$ such that $\gamma_i$ is bounded by it. By proposition \ref{prop:CoverWeakProperty}d) such a subset exists over any field. Then the assertion follows directly from \ref{prop:CoverWeakProperty}b). \\
c) Consider the commutative diagram 
\[\begin{xy}
 \xymatrix {
  \MB{M}_{b_{\nu_i}}^{\preceq \mu_i, \theta_i} \times_E \Sp E' \ar^-{\gamma_i}[r] \ar^{j_{\theta_i, \infty_i}}@{_{(}->}[d] & \MB{M}_{b_{\nu_i}}^{\preceq \mu_i, \theta_i + \theta_{i, \gamma_i}} \times_E \Sp E' \ar^{j_{\theta_i + \theta_{i, \gamma_i}, \infty_i}}@{_{(}->}[d]  \\
  \MB{M}_{b_{\nu_i}}^{\preceq \mu_i} \times_E \Sp E' \ar^{\gamma_i}[r] & \MB{M}_{b_{\nu_i}}^{\preceq \mu_i} \times_E \Sp E' }
\end{xy} \]
As $\gamma_i$ is an automorphisms of $\MB{M}_{b_{\nu_i}}^{\preceq \mu_i} \times_E \Sp E'$ we see that the compositions $j_{\theta_i + \theta_{i, \gamma_i}, \infty_i} \circ \gamma_i = \gamma_i \circ j_{\theta_i, \infty_i}$ are closed immersions. As $j_{\theta_i + \theta_{i, \gamma_i}, \infty_i}$ is a closed immersion as well, we see that $\gamma_i: \MB{M}_{b_{\nu_i}}^{\preceq \mu_i, \theta_i} \times_E \Sp E' \to \MB{M}_{b_{\nu_i}}^{\preceq \mu_i, \theta_i + \theta_{i, \gamma_i}} \times_E \Sp E'$ is just another closed immersion. \\
All assertions are stable under taking perfections, hence the statements hold as well for the $J_i$-action on $\MC{M}_{b_{\nu_i}}^{\preceq \mu_i \sharp} \times_E \Sp E'$. \exit

\thm{\label{thm:JActionEquivariant}}{}{
The morphism 
\[\pi_{(\infty_i)}: \prod_i \MB{M}_{b_{\nu_i}}^{\preceq \mu_i \sharp} \times_E \op{Ig}^{(\infty_i) \sharp}_U \times_E \Sp E' \to \MC{N}^{(\nu_i) \sharp}_U \times_E \Sp E'\]
constructed in \ref{const:CoverMorphism} is equivariant for the diagonal action of $\prod_i J_i$ on the source and the trivial action on the target. \\
Moreover for every element $(\gamma_i)_i \in \prod_i J_i$ the diagram
\[\begin{xy}
 \xymatrix @C=2.4pc {
  \prod_i \MB{M}_{b_{\nu_i}}^{\preceq \mu_i, \theta_i \sharp} \times_E \op{Ig}^{(d_i) \sharp}_U \times_E \Sp E' \ar^-{(\gamma_i)_i}[rr] \ar_{\pi_{(d_i)}}[dr] && \prod_i \MB{M}_{b_{\nu_i}}^{\preceq \mu_i, \theta_i' \sharp} \times_E \op{Ig}^{(d_i') \sharp}_U \times_E \Sp E' \ar^{\pi_{(d_i')}}[dl] \\
  &  \MC{N}^{(\nu_i) \sharp}_U \times_E \Sp E' &  }
\end{xy} \]
commutes, whenever $\theta_i, \theta_i', d_i$ and $d_i'$ are chosen in such a way that all morphisms are well-defined.
}

\prooof
Fix a perfect DM-stack $S$ over $E'$, an element $(\gamma_i)_i \in \prod_i J_i$ viewed by base-change as a quasi-isogeny over $S$, a point $x \in \op{Ig}^{(\infty_i) \sharp}_U(S)$ given by $(\MS{G}, \varphi, \psi, (\alpha_{\infty_i})_i)$ and for each $i$ a point $y_i \in \MB{M}_{b_{\nu_i}}^{\preceq \mu_i \sharp}(S)$ given by $(\MC{G}_i', \varphi_i', \beta_i)$. To check equivariance of $\pi_{(\infty_i)}$ it suffices to check that for each characteristic place $c_i$ the quasi-isogenies used to ``glue'' $(\MS{G}, \varphi, \psi)$ and $(\MC{G}_i', \varphi_i')$ respectively $(\gamma_i)_i.(\MS{G}, \varphi, \psi)$ and $\gamma_i. (\MC{G}_i', \varphi_i')$ coincide. \\
$(\gamma_i)_i.(\MS{G}, \varphi, \psi)$ is constructed using (locally at $c_i$) the quasi-isogeny
\[\alpha_{\infty_i}^{-1} \circ \gamma \circ \alpha_{\infty_i}: \MF{L}_{c_i}(\MS{G}, \varphi) \to \MF{L}_{c_i}((\gamma_i)_i.(\MS{G}, \varphi))\]
Moreover we glue $(\gamma_i)_i.(\MS{G}, \varphi, \psi)$ and $\gamma_i.(\MC{G}_i', \varphi_i')$ at $c_i$ via
\[\gamma_i. (\MC{G}_i', \varphi_i') \xrightarrow{\gamma \circ \beta_i} (L^+G_{c_i}, b_{\nu_i}\sigma^*)_S \xrightarrow{\alpha_{\infty_i}^{-1}} \MF{L}_{c_i}((\gamma_i)_i.(\MS{G}, \varphi)).\]
Thus we obtain the global $G$-shtuka (locally around $c_i$) defining $\pi_{(\infty_i)}((\gamma_i.y_i)_i, (\gamma_i)_i.x)$ by changing $(\MS{G}, \varphi, \psi)$ along the quasi-isogeny 
\[\alpha_{\infty_i}^{-1} \circ \beta_i: \gamma_i. (\MC{G}_i', \varphi_i') \xrightarrow{\gamma \circ \beta_i} (L^+G_{c_i}, b_{\nu_i}\sigma^*)_S \xrightarrow{\alpha_{\infty_i}^{-1}} \MF{L}_{c_i}((\gamma_i)_i.(\MS{G}, \varphi)) \xrightarrow{(\alpha_{\infty_i}^{-1} \circ \gamma \circ \alpha_{\infty_i})^{-1}} \MF{L}_{c_i}(\MS{G}, \varphi)\]
But this is precisely the quasi-isogeny used to get the global $G$-shtuka defining $\pi_{(\infty_i)}((y_i)_i, x)$ (again only locally around $c_i$). Hence we get indeed the desired equivariance. \\
The second assertion obviously follows from the compatibility statements in the previous propositions. \exit

\prop{\label{prop:JActionFibers}}{}{
Let $\ov{\M{F}}$ be any algebraically closed field of characteristic $p$. Fix $x \in \MC{N}^{(\nu_i) \sharp}_U(\ov{\M{F}})$ a geometric point. Then its fiber $\pi_{(\infty_i)}^{-1}(x)$ is a torsor for the $z$-adic group $\prod_i J_i$. 
}

\prooof
By the previous theorem, $\prod_i J_i$ indeed acts on fibers of geometric points, which are non-empty by lemma \ref{lem:CoverInfiniteSurjective}. \\
We first show that the action is simple on $\op{Ig}^{(\infty_i) \sharp}_U(\ov{\M{F}})$ and therefore simple on the whole source as well: Assume there is some $(\gamma_i)_i \in \prod_i J_i$ and an element $(\MS{G}, \varphi, \psi, (\alpha_{\infty_i})_i)$ in the Igusa variety together with an isomorphism
\[\zeta: (\gamma_i)_i.(\MS{G}, \varphi, \psi, (\alpha_{\infty_i})_i) \to (\MS{G}, \varphi, \psi, (\alpha_{\infty_i})_i) \]
Considered locally around a characteristic place $c_i$, we get the commutative diagram
\[\begin{xy}
 \xymatrix {
  (\gamma_i)_i.(\MS{G}, \varphi, \psi) \ar^{\zeta}[dr] \ar@{~>}[ddd] \ar@{~>}[rrr] & & & (\MS{G}, \varphi, \psi) \ar^{\zeta}[dr] \ar@{~>}[dd] &  \\
   & (\MS{G}, \varphi, \psi) \ar@{~>}[ddd] \ar@{~>}[rrr] & & & (\MS{G}, \varphi, \psi) \ar@{~>}[dd] \\
   & & & \MF{L}_{c_i}(\MS{G}, \varphi) \ar^{\zeta}[dr] &  \\
   \MF{L}_{c_i}(\MS{G}, \varphi) \ar^{\zeta}[dr] \ar@{~>}[rr] & & \MF{L}_{c_i}(\MS{G}, \varphi) \ar_{\quad \alpha_{\infty_i, i}^{-1} \circ \gamma_i \circ \alpha_{\infty_i, i}}[ur] \ar^{\zeta}[dr] & & \MF{L}_{c_i}(\MS{G}, \varphi) \\
   & \MF{L}_{c_i}(\MS{G}, \varphi) \ar@{~>}[rr] & & \MF{L}_{c_i}(\MS{G}, \varphi) \ar@{=}[ur] &  }
\end{xy} \]
where the squiggly arrows are again given by functors. Hence the lower right square implies $\gamma_i = \id$, which implies simplicity of the action. \\
It suffices now to show transitivity of the action on fibers. Let $(\MS{G}, \varphi, \psi)$ be the global $G$-shtuka corresponding to the point $x$. Pick two preimages under $\pi_{(\infty_i)}$ corresponding to $(\MS{G}_j, \varphi_j, \psi_j, (\alpha_{\infty_i j})_i) \in \op{Ig}^{(\infty_i) \sharp}_U(\ov{\M{F}})$ and $(\MC{G}_{i j}, \varphi_{i j}, \beta_{i j}) \in \MB{M}_{b_{\nu_i}}^{\preceq \mu_i}(\ov{\M{F}})$ for $j = 1, 2$. Then by construction
\[(\MC{G}_{i 1}, \varphi_{i 1}) = \MF{L}_{c_i}(\MS{G}, \varphi) = (\MC{G}_{i 2}, \varphi_{i 2})\]
and we may define $(\gamma_i)_i \in \prod_i J_i$ via
\[\gamma_i = \beta_{i 2} \circ \beta_{i 1}^{-1} \in QIsog((L^+G_{c_i}, b_{\nu_i} \sigma^*)_{\ov{\M{F}}}, (L^+G_{c_i}, b_{\nu_i} \sigma^*)_{\ov{\M{F}}}) = J_i.\]
Then it is clear that $\gamma_i.(\MC{G}_{i 1}, \varphi_{i 1}, \beta_{i 1}) = (\MC{G}_{i 2}, \varphi_{i 2}, \beta_{i 2})$ for each $i$ and we are left to show
\[(\gamma_i)_i.(\MS{G}_1, \varphi_1, \psi_1, (\alpha_{\infty_i 1})_i) \cong (\MS{G}_2, \varphi_2, \psi_2, (\alpha_{\infty_i 2})_i).\]
By construction $(\MS{G}, \varphi, \psi)$ is obtained by changing $(\MS{G}_1, \varphi_1, \psi_1)$ at every characteristic place along the quasi-isogeny 
\[\alpha_{\infty_i 1}^{-1} \circ \beta_{i 1}: \MF{L}_{c_i}(\MS{G}, \varphi) = (\MC{G}_{i 1}, \varphi_{i 1}) \to \MF{L}_{c_i}(\MS{G}_1, \varphi_1)\]
and similarly for $(\MS{G}_2, \varphi_2, \psi_2)$. Hence we obtain $(\MS{G}_2, \varphi_2, \psi_2)$ if we change $(\MS{G}_1, \varphi_1, \psi_1)$ at every characteristic place along
\[(\alpha_{\infty_i 2}^{-1} \circ \beta_{i 2}) \circ (\alpha_{\infty_i 1}^{-1} \circ \beta_{i 1})^{-1} = \alpha_{\infty_i 2}^{-1} \circ \gamma_i \circ \alpha_{\infty_i 1}: \MF{L}_{c_i}(\MS{G}_1, \varphi_1) \to \MF{L}_{c_i}(\MS{G}_2, \varphi_2)\]
Thus we can define an isomorphism $(\MS{G}_2, \varphi_2, \psi_2, (\alpha_{\infty_i 2})_i) \cong (\gamma_i)_i.(\MS{G}_1, \varphi_1, \psi_1, (\alpha_{\infty_i 1})_i)$ by gluing the identity on 
\[(\MS{G}_1, \varphi_1, \psi_1)|_{(C \setminus \{c_i\}_i) \times_{\Fq} \ov{\M{F}}} = (\MS{G}, \varphi, \psi)|_{(C \setminus \{c_i\}_i) \times_{\Fq} \ov{\M{F}}} = (\MS{G}_2, \varphi_2, \psi_2)|_{(C \setminus \{c_i\}_i) \times_{\Fq} \ov{\M{F}}}\] 
and for each $c_i$ the isomorphism 
\[\alpha_{\infty_i 1}^{-1} \circ \alpha_{\infty_i 2}: \MF{L}_{c_i}(\MS{G}_2, \varphi_2) \to \MF{L}_{c_i}((\gamma_i)_i.(\MS{G}_1, \varphi_1)).\]
This shows that the element $(\gamma_i)_i$ identifies not only the points on Rapoport-Zink spaces but also on Igusa-varieties. \exit

\rem{\label{rem:FiberTopologyProFinite}}{
The $z$-adic topology on $\pi_{(\infty_i)}^{-1}(x)$ induced by the $z$-adic topology on $\prod_i J_i$ coincides with the coarsest topology such that for all sufficiently large $d_i$ and $\theta_i$ the canonical map 
\[\pi_{(\infty_i)}^{-1}(x) \cap \left(\prod_i \MB{M}_{b_{\nu_i}}^{\preceq \mu_i, \theta_i \sharp} \times_E \op{Ig}^{(\infty_i) \sharp}_U\right) \to \pi_{(d_i)}^{-1}(x) \cap \left(\prod_i \MB{M}_{b_{\nu_i}}^{\preceq \mu_i, \theta_i \sharp} \times_E \op{Ig}^{(d_i) \sharp}_U\right)\]
is continuous for the discrete topology on the target (which is a finite set by the next lemma).
}

\lem{}{}{
Fix an algebraically closed field $\ov{\M{F}}$ of characteristic $p$ and let $\theta_i$ be a finite $\Gamma$-invariant subset of $X_*(T)$. Then each $J_i$-orbit of $\ov{\M{F}}$-valued points in $\MB{M}_{b_{\nu_i}}^{\preceq \mu_i \sharp}$ intersects $\MB{M}_{b_{\nu_i}}^{\preceq \mu_i, \theta_i \sharp}$ in only finitely many points.
}

\prooof
Fix one orbit $J_i.x \subset \MB{M}_{b_{\nu_i}}^{\preceq \mu_i \sharp}$ with wlog. $x \in \MB{M}_{b_{\nu_i}}^{\preceq \mu_i, \theta_i \sharp}(\ov{\M{F}})$ corresponding to $(\MC{G}_i, \varphi_i, \beta_i)$. To simplify notations we fix a trivialization of $\MC{G}_i$ which allows us to view $\beta_i \in LG(\ov{\M{F}})$. Moreover we fix a faithful representation $\rho: G_{c_i} \to GL(V)$ (over $\ov{\M{F}}$). \\
Define the set
\[\Delta \coloneqq \{\gamma_i \in J_i \;|\; \gamma_i.x \in \MB{M}_{b_{\nu_i}}^{\preceq \mu_i, \theta_i \sharp}\} = \{\gamma_i \in J_i \;|\; \gamma_i \circ \beta_i \, \textnormal{weakly bounded by} \, \theta_i\} \subset J_i\]
where we view $\gamma_i$ in the second definition as a quasi-isogeny over $\ov{\M{F}}$. Using $\rho$ it is contained in the set
\[\Delta' \coloneqq \{\gamma_i \in J_i \;|\; \rho(\gamma_i \circ \beta_i) \, \textnormal{weakly bounded by} \, \rho(\theta_i)\}\]
But for quasi-isogenies between local $GL(V)$-shtukas we have in lemma \ref{Lem:BoundGLEasy} worked out explicit conditions for boundedness. Note for this that $GL(V)$ splits and $\rho(\gamma_i \circ \beta_i)$ being weakly bounded by $\rho(\theta_i)$ is equivalent to the existence of one $\vartheta_i \in \rho(\theta_i)$ such that $\rho(\gamma_i \circ \beta_i)$ is bounded by $\vartheta_i$. In particular these explicit conditions imply that $\Delta'$ is ($z$-adically) closed in $J_i$ and that there are constants $N_i$ and $N_i'$ (depending on $\theta_i$) such that for each $\gamma_i \in \Delta'$ every coefficient of the matrix $\rho(\gamma_i \circ \beta_i)$ lies in $z^{-N_i} \ov{\M{F}}[[z]]$ and every coefficient of $\rho(\gamma_i \circ \beta_i)^{-1}$ lies in $z^{-N_i'} \ov{\M{F}}[[z]]$. Thus $\Delta'$ is compact as well. \\
Consider now the quotients (by right group actions)
\[\Delta/(\beta_i \Aut_{\ov{\M{F}}}(\MC{G}_i, \varphi_i) \beta_i^{-1}) \;\subset\; \Delta'/(\beta_i \Aut_{\ov{\M{F}}}(\MC{G}_i, \varphi_i) \beta_i^{-1}) \;\subset\; J_i/(\beta_i \Aut_{\ov{\M{F}}}(\MC{G}_i, \varphi_i) \beta_i^{-1}).\]
The group $\beta_i \Aut_{\ov{\M{F}}}(\MC{G}_i, \varphi_i) \beta_i^{-1}$ is open in $J_i$. Together with $J_i$ being defined over a field $E'$ with finite residue field and compactness of $\Delta'$, the quotient in the middle is finite. Hence so is the quotient on the left. But one easily sees that $\beta_i \Aut_{\ov{\M{F}}}(\MC{G}_i, \varphi_i) \beta_i^{-1}$ is nothing else than the stabilizer of $x$. This gives a bijection (of sets)
\[J_i.x \cap \MB{M}_{b_{\nu_i}}^{\preceq \mu_i, \theta_i \sharp}(\ov{\M{F}}) \cong \Delta/(\beta_i \Aut_{\ov{\M{F}}}(\MC{G}_i, \varphi_i) \beta_i^{-1})\]
which proves the desired finiteness result. \exit

\prop{\label{prop:CoverQuasiFinite}}{}{
Assume that $(\theta_i)_i$ and $(d_i)_i$ are chosen in such a way that 
\[\pi_{(d_i)}: \prod_i \MB{M}_{b_{\nu_i}}^{\preceq \mu_i, \theta_i \sharp} \times_E \op{Ig}^{(d_i) \sharp}_U \to \MC{N}^{(\nu_i) \sharp}_U\]
exists. Then $\pi_{(d_i)}$ is quasi-finite. 
}

\prooof
Fix any geometric point $x \in \MC{N}^{(\nu_i) \sharp}_U(\ov{\M{F}})$ and consider its preimage $\pi_{(d_i)}^{-1}(x)$. On $\ov{\M{F}}$-valued points we have the $\prod_i J_i$-action and combining the propositions \ref{prop:JActionDescent} and \ref{prop:JActionFibers}, any two preimages of $x$ are translates via some element in $\prod_i J_i$. Denote the projection on the first factor by $pr_1: \prod_i \MB{M}_{b_{\nu_i}}^{\preceq \mu_i, \theta_i \sharp} \times_E \op{Ig}^{(d_i) \sharp}_U \to \prod_i \MB{M}_{b_{\nu_i}}^{\preceq \mu_i, \theta_i \sharp}$. Then the previous lemma implies the finiteness of the set
\[pr_1(\pi_{(d_i)}^{-1}(x)) \subset \prod_i \MB{M}_{b_{\nu_i}}^{\preceq \mu_i, \theta_i \sharp}(\ov{\M{F}})\]
Fix one element $y \in pr_1(\pi_{(d_i)}^{-1}(x))$. Then it suffices to show the finiteness of the set
\[\pi_{(d_i)}^{-1}(x) \cap pr_1^{-1}(y).\]
If the $i$th component of $y$ corresponds to $(\MC{G}_i, \varphi_i, \beta_i)$, then any two elements in this intersection are translates under the stabilizer of $y$, which equals
\[\prod_i (\beta_i \Aut_{\ov{\M{F}}}(\MC{G}_i, \varphi_i) \beta_i^{-1}) \subset \prod_i J_i.\]
Next observe that the group $\prod_i \Aut^{d_i}(L^+G_{c_i}, b_{\nu_i}\sigma^*) \subset \prod_i J_i$ of automorphisms, which for each $i$ induce the identity modulo $I_{d_i}(b_{\nu_i})$, acts trivially on $\op{Ig}^{(d_i) \sharp}_U$: Indeed for such an element we change the global $G$-shtuka at each characteristic place via some isomorphism. This induces an isomorphisms between the global $G$-shtukas over the whole curve $C$. Moreover the trivializations at the characteristic places get changed only by an element in $I_n(b_{\nu_i})$, hence remain in the same $I_n(b_{\nu_i})$-truncation class. \\
Therefore the subgroup $\prod_i (\Aut^{d_i}(L^+G_{c_i}, b_{\nu_i}\sigma^*) \cap \beta_i \Aut_{\ov{\M{F}}}(\MC{G}_i, \varphi_i) \beta_i^{-1})$ acts trivially on $\pi_{(d_i)}^{-1}(x) \cap pr_1^{-1}(y)$. Thus $\pi_{(d_i)}^{-1}(x) \cap pr_1^{-1}(y)$ is one orbit under the group
\[\prod_i \left(\beta_i \Aut_{\ov{\M{F}}}(\MC{G}_i, \varphi_i) \beta_i^{-1}/ (\Aut^{d_i}(L^+G_{c_i}, b_{\nu_i}\sigma^*) \cap \beta_i \Aut_{\ov{\M{F}}}(\MC{G}_i, \varphi_i) \beta_i^{-1})\right)\]
which is finite as the quotient group is open and $\Aut_{\ov{\M{F}}}(\MC{G}_i, \varphi_i)$ is isomorphic to an automorphism group of a local $G$-shtuka defined over a finite field, and hence can be defined over a finite field itself. \exit

\subsection{Finiteness of the cover}\label{subsec:FinitenessCover}
We saw in the last statement that $\pi_{(d_i)}$ is quasi-finite. The next aim is to prove properness, which implies that $\pi_{(d_i)}$ is finite and for sufficiently large $(d_i)_i$ surjective.

\prop{\label{prop:CoverProper}}{}{
$\pi_{(d_i)}$ is satisfies the valuation criterion for properness (whenever well-defined).
}

\prooof
Let $R$ be any (necessarily perfect) valuation ring with function field $K$ and consider a commutative diagram
\[\begin{xy}
 \xymatrix @C=3pc {
  \Sp K \ar^-{g}[r] \ar^{\eta}[d] & \prod_i \MB{M}_{b_{\nu_i}}^{\preceq \mu_i, \theta_i \sharp} \times_E \op{Ig}^{(d_i) \sharp}_U \ar^{\pi_{(d_i)}}[d]  \\
  \Sp R \ar^{f}[r] \ar^{h}[ur] & \MC{N}^{(\nu_i) \sharp}_U }
\end{xy} \]
where $f$ is defined by a global $G$-shtuka $(\MS{G}_0, \varphi_0, \psi_0)$ over $\Sp R$ and $g$ is defined by an element $(\MS{G}, \varphi, \psi, (\alpha_{d_i})_i) \in \op{Ig}^{(d_i) \sharp}_U(\Sp K)$ and for each $c_i$ a tuple $(\MC{G}_i, \varphi_i, \beta_i) \in \MB{M}_{b_{\nu_i}}^{\preceq \mu_i, \theta_i \sharp}(\Sp K)$. We lift each $\alpha_{d_i}$ to a trivialization $\alpha_{\infty_i}$. Moreover we will always identify $(\MS{G}_0, \varphi_0, \psi_0)_K$ over the generic fiber with the global $G$-shtuka constructed out of the objects given by the points in the source of $\pi_{(d_i)}$. For example we will identify the restriction of $(\MS{G}_0, \varphi_0, \psi_0)_K$ to $C \setminus \{c_i\}_i$ with the restriction of $(\MS{G}, \varphi, \psi)$ to the same subscheme. \\
\textbf{Claim 1:} Each $(\MC{G}_i, \varphi_i, \beta_i)$ admits a model $(\widetilde{\MC{G}_i}, \widetilde{\varphi_i}, \widetilde{\beta_i})$ over $\Sp R$. \\
By construction we identify $\MF{L}_{c_i}(\MS{G}_0, \varphi_0)_K = (\MC{G}_i, \varphi_i)$. Thus the local $G$-shtuka $(\widetilde{\MC{G}_i}, \widetilde{\varphi_i}) \coloneqq \MF{L}_{c_i}(\MS{G}_0, \varphi_0)$ over $\Sp R$ has generic fiber $(\MC{G}_i, \varphi_i)$. By Tate's theorem \ref{thm:ThmTateGShtuka} we may extend the quasi-isogeny $\beta_i$ which is a priori only defined over $K$ to a quasi-isogeny
\[\widetilde{\beta_i}: (\widetilde{\MC{G}_i}, \widetilde{\varphi_i}) \to (L^+G_{c_i}, b_{\nu_i}\sigma^*)_R\]
over $\Sp R$. 
As boundedness by $\theta_i$ is a closed condition, $\widetilde{\beta_i}$ is indeed bounded by $\theta_i$ and $(\widetilde{\MC{G}_i}, \widetilde{\varphi_i}, \widetilde{\beta_i})$ defines a point in $\MB{M}_{b_{\nu_i}}^{\preceq \mu_i, \theta_i \sharp}$. \\
\textbf{Claim 2:} $(\MS{G}, \varphi, \psi)$ admits a model $(\widetilde{\MS{G}}, \widetilde{\varphi}, \widetilde{\psi})$ over $\Sp R$. \\
For each $i$ the element $(\MF{L}_{c_i}(\MS{G}, \varphi), \alpha_{\infty_i})$ is a local $G$-shtuka over $\Sp K$ together with a quasi-isogeny (and even an isomorphism) to $(L^+G_{c_i}, b_{\nu_i}\sigma^*)$. Hence it defines an element in $\MB{M}_{b_{\nu_i}}^{\preceq \mu_i \sharp}$. By corollary \ref{cor:RapoZinkProper}, i.e. properness of irreducible components of $\MB{M}_{b_{\nu_i}}^{\preceq \mu_i \sharp}$, it extends uniquely to a local $G$-shtuka $(\widehat{\MC{G}_i}, \widehat{\varphi_i})$ over $\Sp R$ and a quasi-isogeny $\widetilde{\alpha_{\infty_i}}: (\widehat{\MC{G}}, \widehat{\varphi}) \to (L^+G_{c_i}, b_{\nu_i}\sigma^*)_R$ (which we will prove to be an isomorphism later on). Next observe that we have a quasi-isogeny 
\[\widetilde{\alpha_{\infty_i}}^{-1} \circ \widetilde{\beta_i}: \MF{L}_{c_i}(\MS{G}_0, \varphi_0) = (\widetilde{\MC{G}_i}, \widetilde{\varphi_i}) \to (\widehat{\MC{G}_i}, \widehat{\varphi_i})\]
defined over $\Sp R$. Hence we may glue $(\MS{G}_0, \varphi_0, \psi_0)|_{(C \setminus \{c_i\}_i \times_E \Sp R}$ along these quasi-isogenies to the local $G_{c_i}$-shtukas $(\widehat{\MC{G}_i}, \widehat{\varphi_i})$. This gives a global $G$-shtuka $(\widetilde{\MS{G}}, \widetilde{\varphi}, \widetilde{\psi})$ over $\Sp R$. Over the generic point it is given by modification of $(\MS{G}_0, \varphi_0, \psi_0)$ along $\alpha_{\infty_i}^{-1} \circ \beta_i$, i.e. precisely reversing the construction defining $\pi_{(d_i)}$. Thus over the generic point $(\widetilde{\MS{G}}, \widetilde{\varphi}, \widetilde{\psi})$ equals the given element $(\MS{G}, \varphi, \psi)$. \\
\textbf{Claim 3:} $(\widetilde{\MS{G}}, \widetilde{\varphi}, \widetilde{\psi})$ defines a point in the central leaf $\MC{C}^{(\nu_i) \sharp}_U$. \\
By construction of $(\widetilde{\MS{G}}, \widetilde{\varphi}, \widetilde{\psi})$ as a modification of $(\MS{G}_0, \varphi_0, \psi_0)$ at the characteristic places, it defines a $R$-valued point in $\MC{N}^{(\nu_i) \sharp}_U$. We have seen that its generic point lies in the central leaf. As the central leaf is closed by \ref{prop:CentralLeafProperties}a), this implies claim 3. \\
\textbf{Claim 4:} $\widetilde{\alpha_{\infty_i}}$ is an isomorphism of local $G_{c_i}$-shtukas. \\
$\op{Ig}^{(\infty_i) \sharp}_U \to \MC{C}^{(\nu_i) \sharp}_U$ is a pro-finite \'etale cover by theorem \ref{thm:TruncRepLeaf}. In particular it satisfies the valuation criterion and we may lift the $R$-valued point defined by $(\widetilde{\MS{G}}, \widetilde{\varphi}, \widetilde{\psi})$ to a $R$-valued point $(\widetilde{\MS{G}}, \widetilde{\varphi}, \widetilde{\psi}, (\alpha'_{\infty_i})_i)$ in $\op{Ig}^{(\infty_i) \sharp}_U$ whose generic point equals $(\MS{G}, \varphi, \psi, (\alpha_{\infty_i})_i)$. 
Now observe that both $\alpha'_{\infty_i}$ and $\widetilde{\alpha_{\infty_i}}$ define quasi-isogenies from $(\widehat{\MC{G}}, \widehat{\varphi}) = \MF{L}_{c_i}(\widetilde{\MS{G}}, \widetilde{\varphi})$ to $(L^+G_{c_i}, b_{\nu_i}\sigma^*)_R$ which coincide on the generic point. Hence they have to coincide everywhere. In particular is $\widetilde{\alpha_{\infty_i}}$ an isomorphism. \\
Let $\widetilde{\alpha_{d_i}}$ be the $I_{d_i}(b_{\nu_i})$-truncation class of $\widetilde{\alpha_{\infty_i}}$. Then $(\widetilde{\MS{G}}, \widetilde{\varphi}, \widetilde{\psi}, (\widetilde{\alpha_{d_i}})_i)$ defines a $R$-valued point in $\op{Ig}^{(d_i) \sharp}_U$ and extends the tuple $(\MS{G}, \varphi, \psi, (\alpha_{d_i})_i)$. Therefore we may consider the morphism $h: \Sp R \to \prod_i \MB{M}_{b_{\nu_i}}^{\preceq \mu_i, \theta_i \sharp} \times_E \op{Ig}^{(d_i) \sharp}_U$ defined by $(\widetilde{\MC{G}_i}, \widetilde{\varphi_i}, \widetilde{\beta_i})$ (for each $i$) and $(\widetilde{\MS{G}}, \widetilde{\varphi}, \widetilde{\psi}, (\widetilde{\alpha_{d_i}})_i)$. \\
Note that lifting $(\widetilde{\MS{G}}, \widetilde{\varphi}, \widetilde{\psi})$ to infinite level is only done for convenience and only a lift to $\op{Ig}^{(d_i) \sharp}_U$ is actually necessary. \\
\textbf{Claim 5:} The morphism $h$ makes the diagram above commutative. \\
$g = h \circ \eta$ is precisely the statement, that the shtukas defining $h$ extend the given ones on the generic fiber. $f = \pi_{(d_i)} \circ h$ rephrases as the condition that the modification of $(\widetilde{\MS{G}}, \widetilde{\varphi}, \widetilde{\psi})$ along $\widetilde{\alpha_{\infty_i}}^{-1} \circ \widetilde{\beta_i}$ (for each $i$) is isomorphic to $(\MS{G}_0, \varphi_0, \psi_0)$. But $(\widetilde{\MS{G}}, \widetilde{\varphi}, \widetilde{\psi})$ was constructed precisely to have this property. \exit

\prop{\label{thm:CoverFinite}}{}{
Assume that $(\theta_i)_i$ and $(d_i)_i$ are chosen in such a way that 
\[\pi_{(d_i)}: \prod_i \MB{M}_{b_{\nu_i}}^{\preceq \mu_i, \theta_i \sharp} \times_E \op{Ig}^{(d_i) \sharp}_U \to \MC{N}^{(\nu_i) \sharp}_U\]
exists. Then $\pi_{(d_i)}$ is the perfection of a finite morphism. Moreover for sufficiently large $\theta_i$ (and thus $d_i$) it is surjective.
}

\prooof
Newton strata $\MC{N}^{(\nu_i)}_U$, Rapoport-Zink spaces $\MB{M}_{b_{\nu_i}}^{\preceq \mu_i, \theta_i}$ and Igusa varieties $\op{Ig}^{(d_i)}_U$ are DM-stacks of finite type. Hence the composition of $\pi_{(d_i)}$ with the projection to $\MC{N}^{(\nu_i)}_U$ has to factor through a finite totally ramified cover of $\prod_i \MB{M}_{b_{\nu_i}}^{\preceq \mu_i, \theta_i} \times_E \op{Ig}^{(d_i)}_U$. Denote this morphism by
\[\pi_{(d_i)}: \prod_i \sigma^{*N}\MB{M}_{b_{\nu_i}}^{\preceq \mu_i, \theta_i} \times_E \sigma^{*N}\op{Ig}^{(d_i)}_U \to \MC{N}^{(\nu_i)}_U\]
Propositions \ref{prop:CoverQuasiFinite} and \ref{prop:CoverProper} imply that it is quasi-finite and satisfies the valuation criterion for properness (because one may check these after a cover). Moreover it is of finite type, hence finite. \\
The surjectivity assertion was already shown in \ref{lem:CoverInfiniteSurjective}. \exit

\rem{}{
i) We show that $\pi_{(d_i)}$ itself is finite in theorem \ref{thm:CoveringMorphismProperties}b). \\
ii) Up to Frobenius-pullbacks of the Rapoport-Zink space, 
\[\pi_{(d_i)}: \prod_i \sigma^{*N}\MB{M}_{b_{\nu_i}}^{\preceq \mu_i, \theta_i} \times_E \sigma^{*N}\op{Ig}^{(d_i)}_U \to \MC{N}^{(\nu_i)}_U\]
is precisely the analogue of the covering morphism constructed by Mantovan in \cite{MantoFoliation} and \cite{MantoFoliationPEL}.
}

\subsection{An \'etale version of the covering morphism}\label{subsec:EtaleCoveringMorphism}
So far we defined the finite morphism $\pi_{(d_i)}$. Nevertheless it will be more convenient to have an \'etale and quasi-finite version of this morphism. Fortunately this is easily constructed by restricting $\pi_{(d_i)}$ to some open subset of the source.

\defi{}{}{
a) Let $\MB{M}_{b_{\nu_i}}^{\circ \preceq \mu_i, \theta_i}$ be the largest open subspace of $\MB{M}_{b_{\nu_i}}^{\preceq \mu_i}$, which is contained in $\MB{M}_{b_{\nu_i}}^{\preceq \mu_i, \theta_i}$. \\
Let $\MB{M}_{b_{\nu_i}}^{\circ \preceq \mu_i, \theta_i \sharp}$ be its perfection, which coincides with the largest open subspace of $\MB{M}_{b_{\nu_i}}^{\preceq \mu_i \sharp}$ contained in $\MB{M}_{b_{\nu_i}}^{\preceq \mu_i, \theta_i \sharp}$. \\
b) Let $\dot{\pi}_{(d_i)}: \prod_i \MB{M}_{b_{\nu_i}}^{\circ \preceq \mu_i, \theta_i \sharp} \times_E \op{Ig}^{(d_i) \sharp}_U \to \MC{N}^{(\nu_i) \sharp}_U$ be the restriction of $\pi_{(d_i)}$.
}

\lem{\label{lem:EtaleCoveringMcircProperties}}{}{
For every $\theta_i$ there exists some bound $\theta_i'$ such that
\[\MB{M}_{b_{\nu_i}}^{\circ \preceq \mu_i, \theta_i} \subset \MB{M}_{b_{\nu_i}}^{\preceq \mu_i, \theta_i} \subset \MB{M}_{b_{\nu_i}}^{\circ \preceq \mu_i, \theta_i'}\]
}

\prooof
The first inclusion is immediate from the definition. For the second inclusion recall that the space $\MB{M}_{b_{\nu_i}}^{\preceq \mu_i, \theta_i}$ is quasi-compact by remark \ref{rem:CoverDefSubspaceProper}ii), hence contained in a finite union of irreducible components $Z_i \subset \MB{M}_{b_{\nu_i}}^{\preceq \mu_i}$. As $\MB{M}_{b_{\nu_i}}^{\preceq \mu_i}$ is locally of finite type, each irreducible component intersects only finitely many others. So the union $Z'_i$ of all irreducible components of $\MB{M}_{b_{\nu_i}}^{\preceq \mu_i}$ intersecting $\MB{M}_{b_{\nu_i}}^{\preceq \mu_i}$ is still quasi-compact. Thus we may choose $\theta_i'$ such that $Z'_i \subset \MB{M}_{b_{\nu_i}}^{\preceq \mu_i, \theta_i'}$. Then by our choices
\[\MB{M}_{b_{\nu_i}}^{\preceq \mu_i, \theta_i} \subset \MB{M}_{b_{\nu_i}}^{\circ \preceq \mu_i, \theta_i'}.\]
\exit

\prop{\label{thm:EtaleCoveringMorphismProperties}}{}{
Assume that $(\theta_i)_i$ and $(d_i)_i$ are chosen in such a way that 
\[\dot{\pi}_{(d_i)}: \prod_i \MB{M}_{b_{\nu_i}}^{\circ \preceq \mu_i, \theta_i \sharp} \times_E \op{Ig}^{(d_i) \sharp}_U \to \MC{N}^{(\nu_i) \sharp}_U\]
exists. Then $\dot{\pi}_{(d_i)}$ is formally \'etale. 
}

\prooof
Consider a closed immersion $T_0 \subset T$ of affine schemes given by a square-zero ideal and the test diagram
\[\begin{xy}
 \xymatrix @C=3pc {
  T_0 \ar^-{(f_0, g_0)}[r] \ar[d] & \prod_i \MB{M}_{b_{\nu_i}}^{\circ \preceq \mu_i, \theta_i \sharp} \times \op{Ig}^{(d_i) \sharp}_U \ar^-{\dot{\pi}_{(d_i)}}[d] \\
  T \ar_-{h}[r] \ar@{-->}^-{(\tilde{f}, \tilde{g})}[ur] & \MC{N}^{(\nu_i) \sharp}_U
  }
\end{xy} \]
As maps from an affine space $\Sp A$ into perfect schemes factor over the $\Sp A'$, where $A' \subset A$ is the maximal subring on which the Frobenius is surjective, we may wlog. assume that $T_0$ and $T$ are perfect themselves.
Let now $(\MS{G}, \varphi, \psi)$ (over $T$) be the pullback of the universal global $G$-shtuka along $h$, $(\MC{G}_i, \varphi_i, \beta_i)$ (over $T_0$) be the pullback of the universal family on the Rapoport-Zink space $\MB{M}_{b_{\nu_i}}^{\circ \preceq \mu_i, \theta_i \sharp}$ along $f_0$ and $((\MS{G}', \varphi', \psi'), (\alpha'_{d_i}))$ (over $T_0$) the pullback of the universal family over $\op{Ig}^{(\infty_i) \sharp}_U$ along $g_0$. 
Then by assumption $\MF{L}_{c_i}(\MS{G}, \varphi)|_{T_0} = (\MC{G}_i, \varphi_i)$. Thus we may extend $\beta_i$ by proposition \ref{Prop:BoundRigid} to a quasi-isogeny $\tilde{\beta}_i: \MF{L}_{c_i}(\MS{G}, \varphi) \to (L^+G_{c_i}, b_{\nu_i}\sigma^*)$ over $T$. Then $(\MF{L}_{c_i}(\MS{G}, \varphi), \tilde{\beta}_i)$ defines a morphism $T \to \MB{M}_{b_{\nu_i}}^{\circ \preceq \mu_i, \theta_i}$, which by perfectness of $T$ lifts to $\tilde{f}: T \to \MB{M}_{b_{\nu_i}}^{\circ \preceq \mu_i, \theta_i \sharp}$. By construction $\tilde{f}$ extends $f_0$.
Now the modification of $(\MS{G}, \varphi, \psi)$ along the quasi-isogenies $\tilde{\beta}_i$ on each characteristic place, defines a new global $G$-shtuka $(\tilde{\MS{G}}, \tilde{\varphi}, \psi)$ over $T$ together with canonical $I_{d_i}(b_{\nu_i})$-truncated trivializations $\tilde{\alpha}_{d_i}: \MF{L}_{c_i}(\tilde{\MS{G}}, \tilde{\varphi}) \to \MF{L}_{c_i}(\MS{G}, \varphi)$. This defines a morphism $\tilde{g}: T \to \op{Ig}^{(d_i) \sharp}_U$. \\
We have to see that both triangles in the test diagram commute. The lower one is obvious from the construction of $\tilde{g}$. For the upper one, note that $(\MS{G}', \varphi', \psi')$ can be identified with the modification of $(\MS{G}, \varphi, \psi)|_{T_0}$ along the $\beta_i^{-1} \circ \alpha'_{d_i}$ (for any choice of fpqc-locally representing $\alpha'_{d_i}$ by a honest isomorphism). Thus there is a canonical isomorphism $(\tilde{\MS{G}}, \tilde{\varphi}, \psi)|_{T_0} \cong (\MS{G}', \varphi', \psi')$, which is by construction compatible with the $I_{d_i}(b_{\nu_i})$-truncated trivializations $\tilde{\alpha}_{d_i}$ and $\alpha'_{d_i}$. This shows $\tilde{g}|_{T_0} = g_0$. The equality $\tilde{f}|_{T_0} = f_0$ is clear by construction of $\tilde{f}$. \\
Finally we have to see that $(\tilde{f}, \tilde{g})$ is unique. By definition of $\dot{\pi}_{(d_i)}$, the pullback of the universal local $G$-shtuka along $\tilde{f}$ has to coincide with $\MF{L}_{c_i}(\MS{G}, \varphi)$, hence is uniquely determined. As extensions of quasi-isogenies along nilpotent thickenings are unique, we get uniqueness of our choice of $\tilde{f}$. Uniqueness of $\tilde{g}$ follows now from the definition of $\dot{\pi}_{(d_i)}$, because modification of global $G$-shtukas along local quasi-isogenies is an invertible operation. \exit

\lem{\label{lem:EtalePerfections}}{}{
Let $f: X \to Y$ be a morphism between schemes of finite type over $\Fq$. Assume that the perfection $f^\sharp: X^\sharp \to Y^\sharp$ is formally \'etale. \\
a) $f$ factors over a totally ramified cover $\widetilde{Y} \to Y$ such that $\tilde{f}: X \to \widetilde{Y}$ is \'etale. \\
b) $f^\sharp$ is \'etale. 
}

\prooof
a) Define $\widetilde{Y}$ as the maximal totally ramified extension such that $f$ factors over it. In other words adjoin to the structure sheaf of $Y$ all $p$-power roots of elements, whose image in the structure sheaf of $X$ have $p$-power roots of the same order, and then define $\widetilde{Y}$ as the relative spectrum over $Y$. It is clear that $f$ factors over this $\widetilde{Y}$ and that $\widetilde{Y}$ is totally ramified over $Y$. Moreover $\widetilde{Y}^\sharp = Y^\sharp$ and $\tilde{f}^\sharp = f^\sharp$. \\
To see that $\tilde{f}$ is \'etale, note that satisfying the valuation criterion for \'etaleness implies that $f^\sharp$ defines isomorphisms on completions on local rings. Hence for a geometric point $x \in X$ with image $y = \tilde{f} \in \widetilde{Y}$, the ring $\widehat{\MC{O}_{X, x}}$ is a totally ramified extension of $\widehat{\MC{O}_{\widetilde{Y}, y}}$, because they define the same ring after perfection. So by construction of $\widetilde{Y}$, they have to be equal.
But on DM-stacks of finite type, this implies \'etaleness. \\
b) is a direct consequence of \cite[lemma A.6]{ZhuAffineGrass}: $\tilde{f}$ is \'etale, hence so is $\tilde{f}^\sharp$ which coincides with $f^\sharp$. \exit

\thm{\label{thm:CoveringMorphismProperties}}{}{
Assume that $(\theta_i)_i$ and $(d_i)_i$ are chosen in such a way, that all morphisms exist. \\
a) The morphism 
\[\dot{\pi}_{(d_i)}: \prod_i \MB{M}_{b_{\nu_i}}^{\circ \preceq \mu_i, \theta_i \sharp} \times_E \op{Ig}^{(d_i) \sharp}_U \to \MC{N}^{(\nu_i) \sharp}_U\]
is \'etale and quasi-finite. \\
b) The morphism  
\[\pi_{(d_i)}: \prod_i \MB{M}_{b_{\nu_i}}^{\preceq \mu_i, \theta_i \sharp} \times_E \op{Ig}^{(d_i) \sharp}_U \to \MC{N}^{(\nu_i) \sharp}_U\]
finite. \\
Moreover for all sufficiently large $(\theta_i)_i$ and $(d_i)_i$ both morphisms are surjective.
}

\prooof
Proposition \ref{thm:EtaleCoveringMorphismProperties} implies that we can apply the previous lemma to $\dot{\pi}_{(d_i)}$, proving \'etaleness. $\dot{\pi}_{(d_i)}$ is quasi-finite as a restriction of the quasi-finite morphism $\pi_{(d_i)}$, cf. proposition \ref{prop:CoverQuasiFinite}. \\
We already know that $\pi_{(d_i)}$ is quasi-finite by proposition \ref{prop:CoverQuasiFinite} and satisfies the valuation criterion for properness by proposition \ref{prop:CoverProper}. So it remains to see that it is of finite type, which will follow from the same property for $\dot{\pi}_{(d_i)}$: As the transition morphisms between Igusa varieties are finite \'etale, it suffices to check this for all sufficiently large $(d_i)_i$. 
Now choose $\theta_i'$ as in lemma \ref{lem:EtaleCoveringMcircProperties} such that we have a closed immersion
\[\MB{M}_{b_{\nu_i}}^{\preceq \mu_i, \theta_i} \subset \MB{M}_{b_{\nu_i}}^{\circ \preceq \mu_i, \theta_i'}.\]
Hence the same inclusion of their perfections is again a closed immersion and 
\[\pi_{(d_i)}: \prod_i \MB{M}_{b_{\nu_i}}^{\preceq \mu_i, \theta_i \sharp} \times_E \op{Ig}^{(d_i) \sharp}_U \hookrightarrow \prod_i \MB{M}_{b_{\nu_i}}^{\circ \preceq \mu_i, \theta_i' \sharp} \times_E \op{Ig}^{(d_i) \sharp}_U \xrightarrow{\dot{\pi}_{(d_i)}} \MC{N}^{(\nu_i) \sharp}_U\]
is a composition of morphisms of finite type, assuming that $(d_i)_i$ is large enough for $\dot{\pi}_{(d_i)}$ to exist. So $\pi_{(d_i)}$ is indeed of finite type.  \\
Surjectivity of $\pi_{(d_i)}$ was already shown in proposition \ref{thm:CoverFinite}. Surjectivity of $\dot{\pi}_{(d_i)}$ follows now, because $\dot{\pi}_{(d_i)}$ extends the surjective morphism $\pi_{(d_i)}$ for all sufficiently large $(\theta_i')_i$. \exit

\prop{\label{prop:JActionEtaleVersion}}{}{
Let $(\gamma_i)_i \in \prod_i J_i$. Then there are $d_{i, \gamma_i}$ and $\theta_{i, \gamma_i}$ such that there exists an \'etale morphism
\[(\gamma_i)_i: \prod_i \MB{M}_{b_{\nu_i}}^{\circ \preceq \mu_i, \theta_i \sharp} \times_E \op{Ig}^{(d_i + d_{i, \gamma_i}) \sharp}_U \times_E E' \to \prod_i \MB{M}_{b_{\nu_i}}^{\circ \preceq \mu_i, \theta_i \oplus \theta_{i, \gamma_i} \sharp} \times_E \op{Ig}^{(d_i) \sharp}_U \times_E E'\]
induced from the diagonal action of $(\gamma_i)_i$ on $\prod_i \MB{M}_{b_{\nu_i}}^{\preceq \mu_i \sharp} \times_E \op{Ig}^{(\infty_i) \sharp}_U$. 
}

\prooof
As $(\gamma_i)_i$ acts diagonally, we can handle Rapoport-Zink spaces and Igusa varieties separately. For Igusa varieties this is just proposition \ref{prop:JActionDescent} and lemma \ref{lem:JActionIgusaFinite}. In the case of Rapoport-Zink spaces, first note that $\gamma_i(\MB{M}_{b_{\nu_i}}^{\circ \preceq \mu_i, \theta_i} \sharp) \subset \MB{M}_{b_{\nu_i}}^{\preceq \mu_i \sharp}$ is an open subscheme, as the image of an open subscheme under some automorphism. Moreover lemma \ref{lem:JActionRZ} shows the existence of some $\theta_{i, \gamma_i}$ such that 
\[\gamma_i: \MB{M}_{b_{\nu_i}}^{\preceq \mu_i, \theta_i \sharp} \times_E \Sp E' \to \MB{M}_{b_{\nu_i}}^{\preceq \mu_i, \theta_i \oplus \theta_{i, \gamma_i} \sharp} \times_E \Sp E'\]
is well-defined. Thus $\gamma_i(\MB{M}_{b_{\nu_i}}^{\circ \preceq \mu_i, \theta_i \sharp})$ is an open subscheme of the Rapoport-Zink space, which lies in $\MB{M}_{b_{\nu_i}}^{\preceq \mu_i, \theta_i \oplus \theta_{i, \gamma_i} \sharp} \times_E \Sp E'$. Hence it is contained in $\MB{M}_{b_{\nu_i}}^{\circ \preceq \mu_i, \theta_i \oplus \theta_{i, \gamma_i} \sharp} \times_E \Sp E'$. But this is nothing else that the existence of some $\theta_i$ such that
\[\gamma_i: \MB{M}_{b_{\nu_i}}^{\circ \preceq \mu_i, \theta_i \sharp} \times_E \Sp E' \to \MB{M}_{b_{\nu_i}}^{\circ \preceq \mu_i, \theta_i \oplus \theta_{i, \gamma_i} \sharp} \times_E \Sp E'\]
exists. Moreover it is clearly an open immersion, hence \'etale. \exit

\rem{}{
Alternatively one can deduce \'etaleness of $(\gamma_i)_i$ by considering the commutative diagram
\[\begin{xy}
 \xymatrix @C=0.25pc {
  \prod_i \MB{M}_{b_{\nu_i}}^{\circ \preceq \mu_i, \theta_i \sharp} \times_E \op{Ig}^{(d_i + d_{i, \gamma_i} \sharp)}_U \times_E \Sp E' \ar^-{\gamma_i}[rr] \ar_-{\dot{\pi}_{(d_i + d_{i, \gamma_i})}}[dr] && \prod_i \MB{M}_{b_{\nu_i}}^{\circ \preceq \mu_i, \theta_i \oplus \theta_{i, \gamma_i} \sharp} \times_E \op{Ig}^{(d_i) \sharp}_U \times_E \Sp E' \ar^-{\dot{\pi}_{(d_i)}}[dl]  \\
  & \MC{N}^{(\nu_i) \sharp}_U \times_E \Sp E' & }
\end{xy} \]
and using that the morphisms $\dot{\pi}_{(d_i)}$ and $\dot{\pi}_{(d_i + d_{i, \gamma_i})}$ are \'etale by the previous theorem.
}

\lem{\label{lem:JActionFiberProduct}}{}{
In the situation of the previous proposition, fix some $d_i' \geq d_i$. Then the commutative square
\[\begin{xy}
 \xymatrix @C=2.5pc {
  \prod_i \MB{M}_{b_{\nu_i}}^{\circ \preceq \mu_i, \theta_i \sharp} \times_E \op{Ig}^{(d_i' + d_{i, \gamma_i}) \sharp}_U \times_E \Sp E' \ar^-{\gamma_i}[r] \ar[d] & \prod_i \MB{M}_{b_{\nu_i}}^{\circ \preceq \mu_i, \theta_i \oplus \theta_{i, \gamma_i} \sharp} \times_E \op{Ig}^{(d_i') \sharp}_U \times_E \Sp E' \ar[d]  \\
  \prod_i \MB{M}_{b_{\nu_i}}^{\circ \preceq \mu_i, \theta_i \sharp} \times_E \op{Ig}^{(d_i + d_{i, \gamma_i}) \sharp}_U \times_E \Sp E' \ar^-{\gamma_i}[r] & \prod_i \MB{M}_{b_{\nu_i}}^{\circ \preceq \mu_i, \theta_i \oplus \theta_{i, \gamma_i} \sharp} \times_E \op{Ig}^{(d_i) \sharp}_U \times_E \Sp E'  }
\end{xy} \]
is a fiber product diagram.
}

\prooof
Denote the fiber product space by $X$. Then the universal property gives a morphism
\[f: \prod_i \MB{M}_{b_{\nu_i}}^{\circ \preceq \mu_i, \theta_i \sharp} \times_E \op{Ig}^{(d_i' + d_{i, \gamma_i}) \sharp}_U \times_E E' \to X\]
But both spaces are finite \'etale Galois-covers over $\prod_i \MB{M}_{b_{\nu_i}}^{\circ \preceq \mu_i, \theta_i \sharp} \times_E \op{Ig}^{(d_i + d_{i, \gamma_i}) \sharp}_U \times_E E'$ of the same degree $\prod_i [I_{d_i' + d_{i, \gamma_i}}(b_{\nu_i}):I_{d_i'}(b_{\nu_i})] = \prod_i [I_{d_i + d_{i, \gamma_i}}(b_{\nu_i}):I_{d_i}(b_{\nu_i})]$. Hence the morphism $f$ is a finite \'etale Galois-cover of degree $1$, i.e. an isomorphism. \exit 

\rem{}{
The same result obviously holds as well for the closed spaces $\prod_i \MB{M}_{b_{\nu_i}}^{\preceq \mu_i, \theta_i \sharp}$ instead of $\prod_i \MB{M}_{b_{\nu_i}}^{\circ \preceq \mu_i, \theta_i \sharp}$
}

\subsection{Dimensions of arbitrary leaves}\label{subsec:DimensionLeaf}
As a first application of this product decomposition, we show that every leaf, i.e. the locus in one Newton stratum with fixed isomorphism class of associated local $G$-shtukas, has the same dimension. 

\defi{}{}{
Fix for each characteristic place $c_i$ a local $G_{c_i}$-shtuka $(\MC{G}_i, \varphi_i)$ over a finite field $E'$. Then define the set
\[\MC{C}_U^{(\MC{G}_i, \varphi_i)_i} = \{x \in \MB{X}^{\mmu}_U \;|\; \MF{L}_{c_i}(\MS{G}^{univ}, \varphi^{univ}) \times_{\MB{X}^{\mmu}_U} \Sp k \cong (\MC{G}_i, \varphi_i) \times_{E'} \Sp k \,\textnormal{for each}\, c_i\},\]
where $\ov{x}: \Sp k \to \MB{X}^{\mmu}_U$ is any geometric point with image $x$. 
}

As for central leaves, $\MC{C}_U^{(\MC{G}_i, \varphi_i)_i}$ is a closed subset in the corresponding Newton stratum, hence locally closed in $\MB{X}^{\mmu}_U$. We endow it with the reduced substack structure, making it into a DM-stack over $E'$ and call it a \textit{leaf} in $\MB{X}^{\mmu}_U$.

\prop{}{}{
Fix a Newton stratum $\MC{N}^{(\nu_i)}_U$. Then every leaf in $\MC{N}^{(\nu_i)}_U$ has the same dimension.
}

\prooof
We prove that the leaf corresponding to local $G_{c_i}$-shtukas $(\MC{G}_i, \varphi_i)$ has the same dimension as our fixed central leaf. As perfections preserve dimensions, it suffices to compare the dimensions of $\MC{C}_U^{(\MC{G}_i, \varphi_i)_i \sharp}$ and $\MC{C}_U^{(\nu_i) \sharp}$. \\
Choose for each $i$ a quasi-isogeny $\beta_i: (\MC{G}_i, \varphi_i) \to (L^+G_{c_i}, b_{\nu_i}\sigma^*)$, which exists at least after base-change to $\ACFq$. This defines a set of points $y_i \in \MB{M}_{b_{\nu_i}}^{\preceq \mu_i \sharp}(\Sp \ACFq)$. Choose bounds $\theta_i$ sufficiently large such that $y_i \in \MB{M}_{b_{\nu_i}}^{\preceq \mu_i, \theta_i \sharp}$ for each $i$ and fix a tuple $(d_i)_i$ such that $\pi_{(d_i)}$ exists. We claim now that the morphism
\[\pi_{(d_i)}: (y_i)_i \times_E \op{Ig}^{(d_i) \sharp}_U \to \MC{C}_U^{(\MC{G}_i, \varphi_i)_i \sharp}\]
is well-defined, finite and surjective. Well-defined is obvious from the definition and finiteness follows from theorem \ref{thm:CoveringMorphismProperties}b) by precomposition with the closed immersion defined by the points $y_i$. To see the surjectivity, we proceed as in lemma \ref{lem:CoverInfiniteSurjective}: Fix any geometric point $x \in \MC{C}_U^{(\MC{G}_i, \varphi_i)_i \sharp}$ corresponding to the global $G$-shtuka $(\MS{G}_0, \varphi_0, \psi_0)$. Then by definition of the leaf, we may consider the quasi-isogenies
\[\MF{L}_{c_i}(\MS{G}_0, \varphi_0) \xrightarrow{\sim} (\MC{G}_i, \varphi_i) \xrightarrow{\beta_i} (L^+G_{c_i}, b_{\nu_i}\sigma^*)\]
which define again the points $y_i$ in the Rapoport-Zink space. Moreover the modification of $(\MS{G}_0, \varphi_0, \psi_0)$ along these quasi-isogenies define a global $G$-shtuka in $\op{Ig}^{(d_i) \sharp}_U$ (using the identity as the trivialization of the associated local $G_{c_i}$-shtukas), which lies in the preimage of $x$ under $\pi_{(d_i)}$. This shows surjectivity. \\
In particular we have
\[\dim \MC{C}_U^{(\MC{G}_i, \varphi_i)_i \sharp} = \dim \op{Ig}^{(d_i) \sharp}_U = \dim \MC{C}^{(\nu_i) \sharp}_U.\] \exit

\section{Lifting to formal schemes and adic spaces}\label{sec:FormalLifting}
So far we worked only over the special fiber corresponding to the fixed characteristic places $(c_i)_i$. We now wish to extend $\pi_{(d_i)}$ and $\dot{\pi}_{(d_i)}$ to morphisms of formal schemes, whose underlying reduced subschemes are the ones considered in section \ref{sec:FiniteCover} or at least reasonably close to them. By applying Huber's generic fiber functor, we obtain similar covering morphisms for associated adic spaces.

\subsection{Central leaves as formal schemes}\label{sec:FormalIgusaVarieties}
We construct now the central leaf, not as a reduced subscheme of the special fiber of the moduli space of global $G$-shtukas, but as a formal scheme. A priori we could just take the formal completion of the moduli space $\nabla^{\mmu}_n\MC{H}^1_U(C, G)$ along the central leaf $\MC{C}^{(\nu_i)}_U$. However there are disadvantages in using the resulting formal scheme: First of all the construction of the extended covering morphism could not be made using the universal families as done in theorem \ref{thm:InfiniteFormalCoveringMorphism}, but only using the equivalence of deformations of global and local $G$-shtukas, cf. remark \ref{rem:InfiniteFormalCoveringAlternativeDefinition}. Secondly it would turn out, that the resulting morphism is not adic, having therefore none of the expected properties like being finite or \'etale. 
Thus we define the (formal) central leaf $\MF{C}^{(\nu_i)}_U$ (loosely speaking) as the locus in $\nabla^{\mmu}_n\MC{H}^1_U(C, G)$, where the universal global $G$-shtuka has associated local $G$-shtukas fpqc-locally isomorphic to $(L^+G_{c_i}, b_{\nu_i} \sigma^*)$, cf. \ref{def:FormalCentralLeafDefinition}. \vspace{5mm} \\
From now on we fix local coordinates $\zeta_i$ for $\{c_1\} \times \ldots \{c_{i-1}\} \times C \times \{c_{i+1}\} \times \ldots \times \{c_n\} \subset C^n \setminus \Delta$ at the point $c_i$ once and for all. This way we identify $\Spf E[[\zeta_1, \ldots, \zeta_n]]$ with the formal completion of $C^n \setminus \Delta$ at the point $(c_i)_i$.

\defi{}{}{
a) Let $\MF{X}^{\mmu}_U$ be the formal completion of $\nabla^{\mmu}_n\MC{H}^1_U(C, G)$ along the special fiber $\MB{X}^{\mmu}_U = \nabla^{\mmu}_{(c_i)}\MC{H}^1_U(C, G)$. It is by construction a $(\zeta_1, \ldots, \zeta_n)$-adic formal scheme over $\Spf E[[\zeta_1, \ldots, \zeta_n]] \subset C^n \setminus \Delta$. \\
Let $\MF{X}^{\mmu \sharp}_U$ be its perfection over $\Spf E[[\zeta_1, \ldots, \zeta_n]]^\sharp$. \\
b) Let $\MF{N}^{(\nu_i)}_U \subset \MF{X}^{\mmu}_U$ be the formal completion of $\nabla^{\mmu}_n\MC{H}^1_U(C, G)$ along $\MC{N}^{(\nu_i)}_U$. The formal scheme $\MF{N}^{(\nu_i)}_U$ is not $(\zeta_1, \ldots, \zeta_n)$-adic nor is the inclusion $\MF{N}^{(\nu_i)}_U \subset \MF{X}^{\mmu}_U$ an adic morphism. \\
Let $\MF{N}^{(\nu_i) \sharp}_U$ be its perfection.
}

\rem{}{
i) The formal scheme $\MF{N}^{(\nu_i)}_U$ represents the locus in $\MF{X}^{\mmu}_U$ that admits a fpqc-cover over which for each characteristic place the local $G$-shtuka associated to the universal global one is quasi-isogenous to $(L^+G_{c_i}, b_{\nu_i} \sigma^*)$. Indeed: It is obvious, that the reduced subscheme underlying such a locus has to be contained in $\MC{N}^{(\nu_i)}_U$. Hence the locus of existence of such a fpqc-local quasi-isogeny has to be contained in $\MF{N}^{(\nu_i)}_U$. Conversely such a fpqc-cover over $\MF{N}^{(\nu_i)}_U$ is constructed in the proof of theorem \ref{thm:FormalCentralLeafExistence} (and called $Y$ there). Actually we show the existence of such a cover only modulo some power of an ideal of definition, but the constructions are easily adapted. \\
ii) Alternatively $\MF{X}^{\mmu \sharp}_U$ can be defined as the formal completion of the perfection $\nabla^{\mmu}_n\MC{H}^1_U(C, G)^\sharp$ along the special fiber $\MB{X}^{\mmu \sharp}_U$. The same applies to the Newton stratum $\MF{N}^{(\nu_i) \sharp}_U$.
}

We have already seen that $\op{Ig}^{(\infty_i) \sharp}_U \to \MC{C}^{(\nu_i) \sharp}_U \to \MC{C}^{(\nu_i)}_U$ is a fpqc-cover over which the local $G$-shtukas associated to the universal global $G$-shtuka admit a (canonical) isomorphism to $(L^+G_{c_i}, b_{\nu_i} \sigma^*)$. Therefore it is natural to make the following

\defi{\label{def:FormalCentralLeafDefinition}}{}{
$\MF{C}^{(\nu_i)}_U$ is the locus in $\MF{X}^{\mmu}_U$, where some fpqc-cover exists such that for each characteristic place $c_i$ the following holds: The local $G$-shtuka associated to the universal global $G$-shtuka at this characteristic place $c_i$ admits (over this cover) an isomorphism to the constant local $G$-shtuka $(L^+G_{c_i}, b_{\nu_i} \sigma^*)$. \\
To be more precise, $\MF{C}^{(\nu_i)}_U$ is the fpqc-sheaf on $\MF{X}^{\mmu}_U$, which has on a test scheme $S \to \MF{X}^{\mmu}_U$ exactly one element if some fpqc-cover of $S$ exists with the properties above. $\MF{C}^{(\nu_i)}_U(S)$ is empty otherwise. \\
Let $\MF{C}^{(\nu_i) \sharp}_U$ be the perfection of $\MF{C}^{(\nu_i)}_U$. 
}

Before we can start proving representability results, we need the following general lemma, which is surely nothing new.

\lem{\label{lem:AlgClosedFpqcCover}}{}{
Let $X$ be any quasi-compact noetherian scheme over a field $k$. Then there exists a fpqc-cover $Y \to X$ such that each connected component of the underlying reduced subscheme $Y^{\red}$ is irreducible, normal and has an algebraically closed function field.
}

\prooof
Assume first that $X$ is irreducible with function field $k(X)$. Then we can just take as $Y$ the normalization of $X$ in the algebraic closure $\ov{k(X)}$. \\
In the general case we may assume wlog that the reduced subscheme $X^{\red}$ is irreducible. Then after passing to some open cover, we may embed $X$ into $X^{\red} \times_k \Spf k[[t_1, \ldots, t_n]]$ for some $n \gg 0$. Choose now a fpqc-cover $Y^{\red} \to X^{\red}$ which is normal and has algebraic closed function field. Then $Y^{\red} \times_k \Spf k[[t_1, \ldots, t_n]] \to X^{\red} \times_k \Spf k[[t_1, \ldots, t_n]]$ is still a fpqc-cover with this property. We define now
\[Y = X \times_{ X^{\red} \times \Spf k[[t_1, \ldots, t_n]]} (Y^{\red} \times_k \Spf k[[t_1, \ldots, t_n]])\]
It is obviously a fpqc-cover and the reduced subscheme of $Y$ equally $Y^{\red}$, hence is normal and has as well an algebraically closed function field. \exit

\rem{}{
Note that although $X$ is noetherian, $Y$ will no longer have this property.
}

\thm{\label{thm:FormalCentralLeafExistence}}{}{
$\MF{C}^{(\nu_i)}_U$ exists as a formal scheme over $\Spf E[[\zeta_1, \ldots, \zeta_n]]$. We call it the (formal) central leaf.
}

\prooof
Obviously $\MF{C}^{(\nu_i)}_U$ lies inside $\MF{N}^{(\nu_i)}_U$, which equals the locus, where the associated local $G$-shtukas are quasi-isogenous to $(L^+G_{c_i}, b_{\nu_i} \sigma^*)$. Choose an ideal of definition $I$ for $\MF{N}^{(\nu_i)}_U$. Then we have to show that for each $d \geq 1$ the locus inside $\MF{N}^{(\nu_i)}_U/I^d$, where the associated local $G$-shtukas are isomorphic to $(L^+G_{c_i}, b_{\nu_i} \sigma^*)$, is representable. But it suffices to do this fpqc-locally. Hence choose any $Y \to \MF{N}^{(\nu_i)}_U/I^d$ as in lemma \ref{lem:AlgClosedFpqcCover} and assume wlog. that $Y^{\red}$ is irreducible. We denote by $(\MC{G}_i^{univ}, \varphi_i^{univ})$ the local $G$-shtuka over $Y$ associated to the universal global $G$-shtuka at the characteristic place $c_i$. As there are only finitely many characteristic places, it suffices to deal with each characteristic place separately.
\\
\textbf{Claim 1:} There exists a $J_{b_{\nu_i}}$-torsor of quasi-isogenies $\alpha: (\MC{G}_i^{univ}, \varphi_i^{univ}) \to (L^+G_{c_i}, b_{\nu_i} \sigma^*)$. \\
As the generic point $\eta$ has an algebraically closed function field, we can choose some quasi-isogeny $\alpha_\eta$ over it. By Tate's theorem of extending quasi-isogenies \ref{thm:ThmTateGShtuka} and normality of $Y^{\red}$ it extends to a quasi-isogeny $\alpha_{Y^{\red}}$ over $Y^{\red}$. But now rigidity of quasi-isogenies \ref{Prop:BoundRigid} for nilpotent thickenings implies, that $\alpha_{Y^{\red}}$ extends uniquely to a quasi-isogeny $\alpha$ over all of $Y$. The claimed $J_{b_{\nu_i}}$-torsor of quasi-isogenies comes now by post-composing $\alpha$ with self-isogenies of $(L^+G_{c_i}, b_{\nu_i} \sigma^*)$, which are parametrized by $J_{b_{\nu_i}}$. \\
\textbf{Claim 2:} There are only finitely many $J_{b_{\nu_i}} \cap L^+G_{c_i}$-cosets of these quasi-isogenies, that restrict to an isomorphism on some geometric point of $Y$. \\
Fix one quasi-isogeny $\alpha$. The goal is to bound all $\gamma \in J_{b_{\nu_i}}$ such that $\gamma \circ \alpha$ restricts to an isomorphism at some point. For this note first that one can bound $\alpha$ by some $\mu'$:
Over the generic point $\alpha_\eta$ is bounded by some $\mu$ (as $\eta$ has algebraically closed function field). As the locus of boundedness is closed, $\alpha_{Y^{\red}}$ is bounded by $\mu$ as well. Now the $d$-th power of the ideal defining $Y^{\red} \subset Y$ vanishes, hence by proposition \ref{Prop:BoundRigid} there is some other bound $\mu'$ such that all of $\alpha$ is bounded by it. \\
Consider now $\gamma \in J_{b_{\nu_i}}$ and assume the quasi-isogeny $\gamma \circ \alpha$ is an isomorphism at some point $y \in Y$. Then $\gamma$ (at this point $y$) is the composition of an isomorphism and $\alpha_y^{-1}$, which is bounded by $\mu'^{-1}$. Hence $\gamma$ itself is bounded by $\mu'^{-1}$. But there are only finitely many choices of $J_{b_{\nu_i}} \cap L^+G_{c_i}$-cosets of self-isogenies $\gamma$, which are bounded by $\mu'^{-1}$. \\
\textbf{Claim 3:} The locus $Z$ where at least one of the quasi-isogenies $\alpha$ is an isomorphism, is representable by a closed immersion. \\
Being an isomorphism is equivalent to being bounded by the zero cocharacter. Hence the locus where one fixed quasi-isogeny $\alpha$ is an isomorphism is representable by a closed immersion. Next observe that all quasi-isogenies in $J_{b_{\nu_i}} \cap L^+G_{c_i}$ are isomorphisms, hence this locus coincides for all elements in one $J_{b_{\nu_i}} \cap L^+G_{c_i}$-coset. Thus claim $2$ implies that the locus, where at least one quasi-isogeny is an isomorphism, is a finite union of closed immersions, hence representable by a closed immersion. \\
\textbf{Claim 4:} If $S \to Y$ is any fqpc-morphism such that exists an isomorphism $\beta: (\MC{G}_i^{univ}, \varphi_i^{univ}) \to (L^+G_{c_i}, b_{\nu_i} \sigma^*)$ over $S$, then $S$ already factors through $Z$. \\
We even claim that the isomorphism $\beta$ over $S$ actually comes by pull-back from one of the quasi-isogenies $\alpha$ defined above. To show this, pick any quasi-isogeny $\alpha$ over $Y$ and pull it back to $S$. Then $\beta \circ \alpha^{-1}$ defines a self-quasi-isogeny of $(L^+G_{c_i}, b_{\nu_i} \sigma^*)$ over $S$. As the set $J_{b_{\nu_i}}$ of self-quasi-isogenies is totally disconnected, $\beta \circ \alpha^{-1}$ has to be constantly equal to one $\gamma \in J_{b_{\nu_i}}$. But this shows that $\beta$ is simply the pullback of the quasi-isogeny $\gamma \circ \alpha$ from $Y$ to $S$. \\
Now claim $4$ implies, that the formal scheme $Z$ we constructed in claim $3$ equals $(\MF{C}^{(\nu_i)}_U/I^d) \times_{\MF{N}^{(\nu_i)}_U/I^d} Y$. This proves the theorem. \exit

\prop{\label{prop:FormalCentralLeafAdic}}{}{
$\MF{C}^{(\nu_i)}_U$ is $(\zeta_1, \ldots, \zeta_n)$-adic.
}

\prooof
It suffices to see that the special fiber $\MF{C}^{(\nu_i)}_U/(\zeta_1, \ldots, \zeta_n) \coloneqq \MF{C}^{(\nu_i)}_U \times_{\Spf E[[\zeta_1, \ldots, \zeta_n]]} \Sp E$ exists as a scheme. We already know by theorem \ref{thm:FormalCentralLeafExistence} that it exists as a formal scheme. So it suffices to show that for every algebraically closed field $k$ a morphism $f: \Sp k[[t]] \to \MB{X}^{\mmu}_U$ factors through $\MF{C}^{(\nu_i)}_U/(\zeta_1, \ldots, \zeta_n)$ if its restriction to $\Spf k[[t]]$ does. 
As all $L^+G$-torsors over $\Sp k[[t]]$ are trivial, we may apply \cite[proposition 3.16]{HaVi} (using that all torsors are nice in the sense of section \ref{subsec:Uniform}) to establish an equivalence of categories between bounded local $G$-shtukas over $\Spf k[[t]]$ and over $\Sp k[[t]]$. The very same arguments show as well, that we have an equivalence between bounded quasi-isogenies over $\Spf k[[t]]$ and over $\Sp k[[t]]$. In particular taking the zero cocharacter as a bound, we see that an isomorphism over $\Spf k[[t]]$ extends uniquely to a quasi-isogeny over $\Sp k[[t]]$, where it is again an isomorphism. Hence $f$ factors indeed though $\MF{C}^{(\nu_i)}_U/(\zeta_1, \ldots, \zeta_n)$, if its restriction to $\Spf k[[t]]$ does. \exit

\lem{}{}{
$\MF{C}^{(\nu_i) \sharp}_U$ is the locus in $\MF{N}^{(\nu_i) \sharp}_U$, where some fpqc-cover exists such that the local $G$-shtukas over all $c_i$ can be trivialized in the sense of definition \ref{def:FormalCentralLeafDefinition}.
}

\prooof
This is immediate as $\MF{C}^{(\nu_i) \sharp}_U \to \MF{C}^{(\nu_i)}_U$ is a fpqc-cover itself. \exit

\cor{\label{cor:FormalIgusaVarAdic}}{}{
For all tuples $(d_i)_i$ the sheaf $\MF{Ig}^{(d_i) \sharp}_U$ over $\MF{C}^{(\nu_i) \sharp}_U$ defined by
\[\MF{Ig}^{(d_i)}_U(T) =  \begin{Bmatrix} I_{d_i}(b_{\nu_i})\textnormal{-truncated isomorphisms between the associated local} \; G\textnormal{-shtukas} \\ \textnormal{and} \; (L^+G_{c_i}, b_{\nu_i} \sigma^*) \; \textnormal{for all characteristic places} \; c_i \end{Bmatrix} \]
for all schemes $T$ over $\MF{C}^{(\nu_i) \sharp}_U$, is representable by a $(\zeta_1, \ldots, \zeta_n)$-adic formal scheme, which is finite \'etale over $\MF{C}^{(\nu_i) \sharp}_U$. We call the $\MF{Ig}^{(d_i) \sharp}_U$ (formal) Igusa varieties. \\
Similarly $\MF{Ig}^{(\infty_i) \sharp}_U$ parametrizing trivializations of the whole associated local $G$-shtukas is representable by a $(\zeta_1, \ldots, \zeta_n)$-adic formal scheme. It is isomorphic to $\varprojlim_{(d_i)} \MF{Ig}^{(d_i) \sharp}_U$.
}

\prooof
All assertions are immediate from section \ref{sec:Igusa}, when considered over $\MF{C}^{(\nu_i) \sharp}_U/(\zeta_1, \ldots, \zeta_n)^r$ for any $r \geq 1$. Obviously the Igusa varieties are compatible for varying $r$ and the corollary follows for the stated formal schemes. \exit

\rem{}{
The reduced subscheme underlying $\MF{C}^{(\nu_i) \sharp}_U$ is, as already remarked above, equal to $\MC{C}^{(\nu_i) \sharp}_U$. Moreover we have isomorphisms $\op{Ig}^{(d_i) \sharp}_U = \MF{Ig}^{(d_i) \sharp}_U \times_{\MF{C}^{(\nu_i) \sharp}_U} \MC{C}^{(\nu_i) \sharp}_U$. As all Igusa varieties (whether formal or not) are finite \'etale over the respective basis, we may characterize $\MF{Ig}^{(d_i) \sharp}_U$ as the unique finite \'etale formal scheme over $\MF{C}^{(\nu_i) \sharp}_U$ with reduced subscheme $\op{Ig}^{(d_i) \sharp}_U$.
}

\subsection{The infinite covering morphism of formal schemes}\label{subsec:FormalInfiniteCover}
We first deal with extending the morphism
\[\pi_{(\infty_i)}: \prod_i \MB{M}_{b_{\nu_i}}^{\preceq \mu_i \sharp} \times_E \op{Ig}^{(\infty_i) \sharp}_U \to \MC{N}^{(\nu_i) \sharp}_U.\]

\thm{\label{thm:InfiniteFormalCoveringMorphism}}{}{
There exists a canonical morphism of formal schemes over $\Spf E[[\zeta_1, \ldots, \zeta_n]]^\sharp$
\[\widehat{\pi}_{(\infty_i)}: \prod_i \MC{M}_{b_{\nu_i}}^{\preceq \mu_i \sharp} \times_{\Spf E[[\zeta_1, \ldots, \zeta_n]]^\sharp} \MF{Ig}^{(\infty_i)\sharp}_U \to \MF{N}^{(\nu_i) \sharp}_U\]
restricting to $\pi_{(\infty_i)}$ on reduced subschemes. $\widehat{\pi}_{(\infty_i)}$ is formally \'etale.
}

\prooof
The construction is word for word the same as on the special fiber, cf. \ref{const:CoverMorphism}.  Note however that we now need theorem \ref{thm:UnifMorphism} in its full generality. That $\widehat{\pi}_{(\infty_i)}$ restricts to $\pi_{(\infty_i)}$ is obvious. \\
To show formal \'etaleness consider a closed immersion $T_0 \subset T$ of affine schemes given by a square-zero ideal and the test diagram
\[\begin{xy}
 \xymatrix @C=3pc {
  T_0 \ar^-{(f_0, g_0)}[r] \ar[d] & \prod_i \MC{M}_{b_{\nu_i}}^{\preceq \mu_i \sharp} \times \MF{Ig}^{(\infty_i) \sharp}_U \ar^-{\widehat{\pi}_{(\infty_i)}}[d] \\
  T \ar_-{h}[r] \ar@{-->}^-{(\tilde{f}, \tilde{g})}[ur] & \MF{N}^{(\nu_i) \sharp}_U
  }
\end{xy} \]
As maps from an affine space $\Sp A$ into perfect schemes factor over the $\Sp A'$, where $A' \subset A$ is the maximal subring on which the Frobenius is surjective, we may wlog. assume that $T_0$ and $T$ are perfect themselves.
Let now $(\MS{G}, \varphi, \psi)$ (over $T$) be the pullback of the universal global $G$-shtuka along $h$, $(\MC{G}_i, \varphi_i, \beta_i)$ (over $T_0$) be the pullback of the universal family on the Rapoport-Zink space $\MC{M}_{b_{\nu_i}}^{\preceq \mu_i \sharp}$ along $f_0$ and $((\MS{G}', \varphi', \psi'), (\alpha'_{\infty_i}))$ (over $T_0$) the pullback of the universal family over $\MF{Ig}^{(\infty_i) \sharp}_U$ along $g_0$. 
Then by assumption $\MF{L}_{c_i}(\MS{G}, \varphi)|_{T_0} = (\MC{G}_i, \varphi_i)$. Thus we may extend $\beta_i$ by proposition \ref{Prop:BoundRigid} to a quasi-isogeny $\tilde{\beta}_i: \MF{L}_{c_i}(\MS{G}, \varphi) \to (L^+G_{c_i}, b_{\nu_i}\sigma^*)$ over $T$. Then $T$ being a perfect scheme, $(\MF{L}_{c_i}(\MS{G}, \varphi), \tilde{\beta}_i)$ defines an extension $\tilde{f}: T \to \MC{M}_{b_{\nu_i}}^{\preceq \mu_i \sharp}$ of $f_0$.
Now the modification of $(\MS{G}, \varphi, \psi)$ along the quasi-isogenies $\tilde{\beta}_i$ on each characteristic place, defines a new global $G$-shtuka $(\tilde{\MS{G}}, \tilde{\varphi}, \psi)$ over $T$ together with canonical trivializations $\tilde{\alpha}_{\infty_i}: \MF{L}_{c_i}(\tilde{\MS{G}}, \tilde{\varphi}) \to \MF{L}_{c_i}(\MS{G}, \varphi)$. This defines a morphism $\tilde{g}: T \to \MF{Ig}^{(\infty_i) \sharp}_U$. \\
We have to see that both triangles in the test diagram commute. The lower one is obvious from the construction of $\tilde{g}$. For the upper one, note that $(\MS{G}', \varphi', \psi')$ can be identified with the modification of $(\MS{G}, \varphi, \psi)|_{T_0}$ along the $\beta_i^{-1} \circ \alpha'_{\infty_i}$. Thus there is a canonical isomorphism $(\tilde{\MS{G}}, \tilde{\varphi}, \psi)|_{T_0} \cong (\MS{G}', \varphi', \psi')$, which is by construction compatible with the trivializations $\tilde{\alpha}_{\infty_i}$ and $\alpha'_{\infty_i}$. This shows $\tilde{g}|_{T_0} = g_0$. The equality $\tilde{f}|_{T_0} = f_0$ is clear by construction of $\tilde{f}$. \\
Finally we have to see that $(\tilde{f}, \tilde{g})$ is unique. By definition of $\widehat{\pi}_{(\infty_i)}$, the pullback of the universal local $G$-shtuka along $\tilde{f}$ has to coincide with $\MF{L}_{c_i}(\MS{G}, \varphi)$, hence is uniquely determined. As extensions of quasi-isogenies along nilpotent thickenings are unique, we get uniqueness of our choice of $\tilde{f}$. Uniqueness of $\tilde{g}$ follows now from the definition of $\widehat{\pi}_{(\infty_i)}$, because modification of global $G$-shtukas along local quasi-isogenies is an invertible operation. \exit

\lem{\label{lem:InfiniteFormalCoveringMorphismCharacterize}}{}{
Let $Y$ be any formal scheme. Then the diagram 
\[\begin{xy}
 \xymatrix @C=3pc {
  Y \ar^-{f}[r] \ar_-{g}[dr] & \prod_i \MC{M}_{b_{\nu_i}}^{\preceq \mu_i \sharp} \times \MF{Ig}^{(\infty_i) \sharp}_U \ar^-{\widehat{\pi}_{(\infty_i)}}[d] \\
   & \MF{N}^{(\nu_i) \sharp}_U
  }
\end{xy} \]
commutes if and only if the following two conditions are satisfied:
\begin{itemize}
 \item The diagram commutes, when restricted to $Y \times_{\prod \MC{M}_{b_{\nu_i}}^{\preceq \mu_i} \sharp} \prod_i \MB{M}_{b_{\nu_i}}^{\preceq \mu_i \sharp}$.
 \item For each characteristic place, the universal local $G$-shtuka obtained from pulling back the universal family over $\MC{M}_{b_{\nu_i}}^{\preceq \mu_i \sharp}$ along $f$ coincides with the local $G$-shtuka associated to the global $G$-shtuka obtained from pulling back the universal family over $\MF{N}^{(\nu_i) \sharp}_U$ along $g$.
\end{itemize}
}

\prooof
Using the first condition, the assertion reduces to comparing deformations of global $G$-shtukas. By \cite[theorem 5.10]{HarRad1} such deformations are determined by the deformations of associated local $G$-shtukas at the characteristic places. Hence the second condition ensures, that the pullback of the universal family via $g$ coincides with the pullback along $f \circ \widehat{\pi}_{(\infty_i)}$. The lemma follows. \exit

\rem{\label{rem:InfiniteFormalCoveringAlternativeDefinition}}{
Of course we could have defined $\widehat{\pi}_{(\infty_i)}$, by taking $\pi_{(\infty_i)}$ on reduced subschemes and then extending everything via the equivalence of deformations \cite[theorem 5.10]{HarRad1}. However in the way presented above, one retains a very good control of $\widehat{\pi}_{(\infty_i)}$ over the whole formal scheme, and not only over its underlying reduced subscheme.
}

\prop{\label{prop:InfiniteFormalCoveringMorphismAdic}}{}{
$\widehat{\pi}_{(\infty_i)}$ is an adic morphism.
}

\prooof
Let $I$ be an ideal of definition for $\MF{N}^{(\nu_i)}_U$, which then by construction is an ideal of definition for $\MF{N}^{(\nu_i) \sharp}_U$ as well. Now the intersection of all powers of $I$ is zero, so the same holds after pulling it back along $\widehat{\pi}_{(\infty_i)}$. Thus it suffices to show that $Y = \MF{N}^{(\nu_i) \sharp}_U/I \times_{\MF{N}^{(\nu_i) \sharp}_U} \left(\prod_i \MC{M}_{b_{\nu_i}}^{\preceq \mu_i \sharp} \times_{\Spf E[[\zeta_1, \ldots, \zeta_n]]^\sharp} \MF{Ig}^{(\infty_i) \sharp}_U\right)$ exists as a scheme. As $(\zeta_1, \ldots, \zeta_n) \subset I$, we may get the more convenient form 
\[Y = \MF{N}^{(\nu_i) \sharp}_U/I \times_{\MF{N}^{(\nu_i) \sharp}_U/(\zeta_1, \ldots, \zeta_n)} \left(\prod_i \MC{M}_{b_{\nu_i}}^{\preceq \mu_i \sharp}/(\zeta_1, \ldots, \zeta_n) \times \MF{Ig}^{(\infty_i) \sharp}_U/(\zeta_1, \ldots, \zeta_n)\right)\]
As it obviously exists as a formal scheme, we proceed similarly to proposition \ref{prop:FormalCentralLeafAdic}, i.e. we show that for any perfect valuation ring $A$, any morphism $f: \Spf A \to Y$ extends to a morphism $f: \Sp A \to Y$. First of all $\MF{N}^{(\nu_i) \sharp}_U/I$ is a scheme, hence $pr_1 \circ f: \Spf A \to Y \to \MF{N}^{(\nu_i)\sharp}_U/I$ extends over $\Sp A$. Using corollary \ref{cor:FormalIgusaVarAdic} one can extend $pr_2 \circ f: \Spf A \to Y \to \MF{Ig}^{(\infty_i )\sharp}_U/(\zeta_1, \ldots, \zeta_n)$ by the same reasoning. \\
Moreover using the extended version of $pr_1 \circ f$, one gets a global $G$-shtuka over $\Sp A$. Its associated local $G$-shtukas are by construction algebraizations of the local $G$-shtukas defined by pulling back the universal ones from $\prod_i \MC{M}_{b_{\nu_i}}^{\preceq \mu_i \sharp}$. Now the quasi-isogenies defined by $pr_2 \circ f: \Spf A \to \MC{M}_{b_{\nu_i}}^{\preceq \mu_i \sharp}$ are bounded, hence can be algebraized by the argument given in \cite[proposition 3.16]{HaVi}. This defines a morphism $g: \Sp A \to \prod_i \MC{M}_{b_{\nu_i}}^{\preceq \mu_i \sharp}$, which factors over $\prod_i \MC{M}_{b_{\nu_i}}^{\preceq \mu_i \sharp}/(\zeta_1, \ldots, \zeta_n)$, because the global $G$-shtukas inducing the local ones are defined over $\Spf E[[\zeta_1, \ldots, \zeta_n]]^\sharp/(\zeta_1, \ldots, \zeta_n)$.
Then the triple $(pr_1 \circ f, g, pr_2 \circ f): \Sp A \to Y$ extends by construction $f$ as desired. \exit

\subsection{An \'etale covering morphism of formal schemes}
Our next objective is to define a generalization of $\dot{\pi}_{(d_i)}$ to formal schemes. It will be adic, \'etale and surjective (in a strong sense).

\defi{}{}{
Fix a weak bound $(\theta_i)_i$ and recall the open immersion $\MB{M}_{b_{\nu_i}}^{\circ \preceq \mu_i, \theta_i \sharp} \hookrightarrow \MB{M}_{b_{\nu_i}}^{\preceq \mu_i \sharp}$. As \'etale morphisms extend uniquely to nilpotent thickenings, there is a unique open formal subscheme $\MC{M}_{b_{\nu_i}}^{\circ \preceq \mu_i, \theta_i \sharp} \hookrightarrow \MC{M}_{b_{\nu_i}}^{\preceq \mu_i \sharp}$, whose reduced subscheme equals $\MB{M}_{b_{\nu_i}}^{\circ \preceq \mu_i, \theta_i \sharp}$. \\
We view $\prod_i \MC{M}_{b_{\nu_i}}^{\circ \preceq \mu_i, \theta_i \sharp}$ as a (non-adic) formal scheme over $\Spf E[[\zeta_1, \ldots, \zeta_n]]^\sharp$ in the usual way.
}

\prop{\label{prop:FormalCoveringMorphEtaleConstruction}}{}{
Fix a weak bound $(\theta_i)_i$. Then for all sufficiently large tuples $(d_i)_i$ (depending only on $(\theta_i)_i$ as in proposition \ref{prop:CoverMorphDescent}), there exists a canonical morphism of formal schemes over $\Spf E[[\zeta_1, \ldots, \zeta_n]]^\sharp$
\[\widehat{\dot{\pi}}_{(d_i)}: \prod_i \MC{M}_{b_{\nu_i}}^{\circ \preceq \mu_i, \theta_i \sharp} \times_{\Spf E[[\zeta_1, \ldots, \zeta_n]]^\sharp} \MF{Ig}^{(d_i) \sharp}_U \to \MF{N}^{(\nu_i) \sharp}_U.\]
restricting to $\dot{\pi}_{(d_i)}$ on reduced subschemes. Moreover $\widehat{\pi}_{(\infty_i)}|_{\prod \MC{M}_{b_{\nu_i}}^{\circ \preceq \mu_i, \theta_i \sharp} \times \MF{Ig}^{(\infty_i) \sharp}_U}$ factors through $\widehat{\dot{\pi}}_{(d_i)}$.
}

\prooof
Proposition \ref{prop:CoverMorphDescentPoint} still holds (with the very same proof) for $S$-valued points, where $S$ is now a perfect formal scheme over $\Spf E[[\zeta_1, \ldots, \zeta_n]]^\sharp$ with simply connected underlying reduced subscheme. In particular we get as in the reduced case (cf. proposition \ref{prop:CoverMorphDescent}) a factorization of $\widehat{\pi}_{(\infty_i)}|_{\prod \MB{M}_{b_{\nu_i}}^{\circ \preceq \mu_i, \theta_i \sharp} \times \MF{Ig}^{(\infty_i) \sharp}_U}$ over $\prod_i \MB{M}_{b_{\nu_i}}^{\circ \preceq \mu_i, \theta_i \sharp} \times_{\Spf E[[\zeta_1, \ldots, \zeta_n]]^\sharp} \MF{Ig}^{(d_i) \sharp}_U$. 
Thus we obtain a morphism
\[\widehat{\dot{\pi}}_{(d_i)}: \prod_i \MB{M}_{b_{\nu_i}}^{\circ \preceq \mu_i, \theta_i \sharp} \times_{\Spf E[[\zeta_1, \ldots, \zeta_n]]^\sharp} \MF{Ig}^{(d_i) \sharp}_U \to \MF{N}^{(\nu_i) \sharp}_U\]
coinciding with $\dot{\pi}_{(d_i)}$ over the underlying reduced subscheme. \\
Now let $(\MS{G}_0, \varphi_0, \psi_0)$ be the pullback of the universal global $G$-shtuka to $\prod_i \MB{M}_{b_{\nu_i}}^{\circ \preceq \mu_i, \theta_i \sharp} \times_{\Spf E[[\zeta_1, \ldots, \zeta_n]]^\sharp} \MF{Ig}^{(d_i) \sharp}_U$ and denote by $(\MC{G}_i, \varphi_i)$ (for each characteristic place $c_i$) the universal local $G$-shtuka defined by the Rapoport-Zink space $\MC{M}_{b_{\nu_i}}^{\circ \preceq \mu_i, \theta_i \sharp}$. Then $(\MC{G}_i, \varphi_i)$ coincides with $\MF{L}_{c_i}(\MS{G}_0, \varphi_0)$ where both are defined. 
By \cite[theorem 5.10]{HarRad1} deformations of global $G$-shtukas are uniquely determined by the deformations of their associated local $G$-shtukas. Hence the $(\MC{G}_i, \varphi_i)$ for each characteristic place define a unique deformation $(\MS{G}, \varphi, \psi)$ of $(\MS{G}_0, \varphi_0, \psi_0)$ to $\prod_i \MC{M}_{b_{\nu_i}}^{\circ \preceq \mu_i, \theta_i \sharp} \times_{\Spf E[[\zeta_1, \ldots, \zeta_n]]^\sharp} \MF{Ig}^{(d_i) \sharp}_U$. This family $(\MS{G}, \varphi, \psi)$ (on a perfect basis) defines the morphism $\widehat{\dot{\pi}}_{(d_i)}$. \\
It remains to see that $\widehat{\pi}_{(\infty_i)}$ factors through $\widehat{\dot{\pi}}_{(d_i)}$. This follows essentially from lemma \ref{lem:InfiniteFormalCoveringMorphismCharacterize}, but for the morphism $\widehat{\dot{\pi}}_{(d_i)}$ instead of $\widehat{\pi}_{(\infty_i)}$. We nevertheless give the argument for sake of completeness: 
Denote by $(\MS{G}', \varphi', \psi')$ the global $G$-shtuka over $\prod_i\MC{M}_{b_{\nu_i}}^{\circ \preceq \mu_i, \theta_i \sharp} \times \MF{Ig}^{(\infty_i) \sharp}_U$ obtained by pulling back the universal family along $\widehat{\pi}_{(\infty_i)}$. If we pull back $(\MS{G}, \varphi, \psi)$ along the canonical morphism $r_{\infty_i}: \prod_i \MC{M}_{b_{\nu_i}}^{\circ \preceq \mu_i, \theta_i \sharp} \times \MF{Ig}^{(\infty_i) \sharp}_U \to \prod_i \MC{M}_{b_{\nu_i}}^{\circ \preceq \mu_i, \theta_i \sharp} \times \MF{Ig}^{(d_i) \sharp}_U$, we obtain another global $G$-shtuka, still denoted by $(\MS{G}, \varphi, \psi)$. We have to see that $(\MS{G}', \varphi', \psi')$ and $(\MS{G}, \varphi, \psi)$ are isomorphic. 
Indeed over $\prod_i \MB{M}_{b_{\nu_i}}^{\circ \preceq \mu_i, \theta_i \sharp} \times \MF{Ig}^{(\infty_i) \sharp}_U$ this follows from the construction of $\widehat{\dot{\pi}}_{(d_i)}$. Moreover the associated local $G$-shtukas are both canonically isomorphic to $(\MC{G}_i, \varphi_i)$. Hence the local description of deformations of global $G$-shtukas \cite[theorem 5.10]{HarRad1} lets us conclude $(\MS{G}', \varphi', \psi') \cong (\MS{G}, \varphi, \psi)$. \exit

\prop{\label{prop:FormalCoveringMorphEtaleEtale}}{}{
$\widehat{\dot{\pi}}_{(d_i)}$ is adic and \'etale.
}

\prooof
We first prove that $\widehat{\dot{\pi}}_{(d_i)}$ is formally \'etale. This follows in the very same way as for $\widehat{\pi}_{(\infty_i)}$, cf. theorem \ref{thm:InfiniteFormalCoveringMorphism}. One only has to notice, that by definition any morphism $S \to \MC{M}_{b_{\nu_i}}^{\preceq \mu_i \sharp}$ of formal schemes, whose underlying reduced subscheme maps into $\MB{M}_{b_{\nu_i}}^{\circ \preceq \mu_i, \theta_i \sharp}$ actually factors through $\MC{M}_{b_{\nu_i}}^{\circ \preceq \mu_i, \theta_i \sharp}$. \\
That $\widehat{\dot{\pi}}_{(d_i)}$ is adic follows now in precisely the same way as in proposition \ref{prop:InfiniteFormalCoveringMorphismAdic}. That it is \'etale is again a consequence of lemma \ref{lem:EtalePerfections}, which generalizes to adic morphisms between formal schemes: Take some power $I^n$ of the ideal of definition on the target and consider the morphism between the closed subscheme defined by $\widehat{\dot{\pi}}_{(d_i)}^*(I^n)$ of the source and the closed subscheme defined by $I^n$ on the target. By the lemma this is \'etale, hence so is the limit morphism $\widehat{\dot{\pi}}_{(d_i)}$. \exit

\rem{}{
Instead copying the proof of proposition \ref{prop:InfiniteFormalCoveringMorphismAdic}, one can argue as follows: As a restriction of an adic morphism to an open formal subscheme, $\widehat{\pi}_{(\infty_i)}|_{\prod \MC{M}_{b_{\nu_i}}^{\circ \preceq \mu_i, \theta_i \sharp} \times \MF{Ig}^{(\infty_i) \sharp}_U}$ is adic. It is a composition of $\widehat{\dot{\pi}}_{(d_i)}$ and an adic surjective (and pro-\'etale) morphism. Hence $\widehat{\dot{\pi}}_{(d_i)}$ is adic as well.
}

\prop{\label{prop:FormalCoveringMorphEtaleSurjective}}{}{
For all sufficiently large tuples $(\theta_i)_i$ and $(d_i)_i$, the morphism $\widehat{\dot{\pi}}_{(d_i)}$ satisfies the following surjectivity assertion: Let $S$ be any perfection of a quasi-compact noetherian formal scheme over $\Spf E[[\zeta_1, \ldots, \zeta_n]]^\sharp$ and $f: S \to \MF{N}^{(\nu_i) \sharp}_U$ any morphism of formal schemes. Then there exists an \'etale cover $S' \to S$ and a morphism $\tilde{f}: S' \to \prod_i \MC{M}_{b_{\nu_i}}^{\circ \preceq \mu_i, \theta_i \sharp} \times_{\Spf E[[\zeta_1, \ldots, \zeta_n]]^\sharp} \MF{Ig}^{(d_i) \sharp}_U$ such that
\[\begin{xy}
 \xymatrix @C=3pc {
  S' \ar^-{\tilde{f}}[r] \ar[d] & \prod_i \MC{M}_{b_{\nu_i}}^{\circ \preceq \mu_i, \theta_i \sharp} \times_{\Spf E[[\zeta_1, \ldots, \zeta_n]]^\sharp} \MF{Ig}^{(d_i) \sharp}_U \ar^-{\widehat{\dot{\pi}}_{(d_i)}}[d] \\
  S \ar_-{f}[r] & \MF{N}^{(\nu_i) \sharp}_U
  }
\end{xy} \]
commutes.
}

\prooof
We define 
\[S' = S \times_{\MF{N}^{(\nu_i) \sharp}_U} \left(\prod_i \MC{M}_{b_{\nu_i}}^{\circ \preceq \mu_i, \theta_i \sharp} \times_{\Spf E[[\zeta_1, \ldots, \zeta_n]]^\sharp} \MF{Ig}^{(d_i) \sharp}_U\right)\]
Then $\tilde{f}$ is nothing else than the projection onto the second factor and the previous proposition implies, that $S' \to S$ is \'etale. Surjectivity for \'etale morphisms can be checked on underlying reduced subschemes. Now $\dot{\pi}_{(d_i)}$ is surjective, which was shown in lemma \ref{lem:CoverInfiniteSurjective} by considering geometric points. Hence $S'_{\red} \to S_{\red}$ is surjective as the base-change of $\dot{\pi}_{(d_i)}$ along $f_{\red}$. \exit

\rem{}{
The formal schemes $S$ to test surjectivity above need not be adic or admissible. In particular $S$ could as well be any (usual) scheme.
}

\subsection{A finite covering morphism of formal schemes}
We define now a generalization of the finite morphism $\pi_{(d_i)}$ to formal schemes. This morphism $\widehat{\pi}_{(d_i)}$ will turn out to be adic, finite and surjective (at least for sufficiently large $(d_i)_i$). But the underlying reduced subschemes of the source will not coincide with $\prod \MB{M}_{b_{\nu_i}}^{\preceq \mu_i, \theta_i \sharp} \times_E \op{Ig}^{(d_i) \sharp}_U$ but rather with a closed subscheme of it.

\defi{}{}{
Consider the coherent sheaf of ideals in $\MC{O}_{\prod \MC{M}_{b_{\nu_i}}^{\preceq \mu_i \sharp}}$ of all functions vanishing on $\prod_i \MC{M}_{b_{\nu_i}}^{\circ \preceq \mu_i, \theta_i \sharp}$. By \cite[section 10.14]{GroEGAI} this ideal sheaf corresponds to a closed formal subscheme $\prod_i \MC{M}_{b_{\nu_i}}^{\preceq \mu_i, \theta_i \sharp} \subset \prod_i \MC{M}_{b_{\nu_i}}^{\preceq \mu_i \sharp}$ such that the inclusion morphism is adic and of finite type. \\
The underlying reduced subscheme is the closure of $\prod_i \MB{M}_{b_{\nu_i}}^{\circ \preceq \mu_i, \theta_i \sharp}$ in $\prod_i \MB{M}_{b_{\nu_i}}^{\preceq \mu_i \sharp}$.
}

\prop{\label{prop:FormalCoveringMorphFinite}}{}{
Fix a weak bound $(\theta_i)_i$. Then for all sufficiently large tuples $(d_i)_i$ (depending only on $(\theta_i)_i$ as in proposition \ref{prop:CoverMorphDescent}), there exists a canonical morphism of formal schemes over $\Spf E[[\zeta_1, \ldots, \zeta_n]]^\sharp$
\[\widehat{\pi}_{(d_i)}: \prod_i \MC{M}_{b_{\nu_i}}^{\preceq \mu_i, \theta_i \sharp} \times_{\Spf E[[\zeta_1, \ldots, \zeta_n]]^\sharp} \MF{Ig}^{(d_i) \sharp}_U \to \MF{N}^{(\nu_i) \sharp}_U.\]
coinciding with $\pi_{(d_i)}$ on the underlying reduced subscheme. \\
a) $\widehat{\pi}_{(\infty_i)}|_{\prod \MC{M}_{b_{\nu_i}}^{\preceq \mu_i, \theta_i \sharp} \times \MF{Ig}^{(\infty_i) \sharp}_U}$ factors through $\widehat{\pi}_{(d_i)}$. \\
b) $\widehat{\pi}_{(d_i)}$ is an adic and finite morphism. \\
c) For sufficiently large $(\theta_i)_i$ and $(d_i)_i$, $\widehat{\pi}_{(d_i)}$ satisfies the surjectivity property stated in proposition \ref{prop:FormalCoveringMorphEtaleSurjective}, i.e. the same as for $\widehat{\dot{\pi}}_{(d_i)}$.
}

\prooof
The underlying reduced subscheme of $\MC{M}_{b_{\nu_i}}^{\preceq \mu_i, \theta_i \sharp}$ lies in $\MB{M}_{b_{\nu_i}}^{\preceq \mu_i, \theta_i \sharp}$. Thus the construction of $\widehat{\pi}_{(d_i)}$ is the very same as for $\widehat{\dot{\pi}}_{(d_i)}$. For property a) we may use the proof of proposition \ref{prop:FormalCoveringMorphEtaleConstruction} as well. \\
Choose now some $(\theta_i')_i$ such that $\prod_i \MB{M}_{b_{\nu_i}}^{\preceq \mu_i, \theta_i \sharp} \subset \prod_i \MB{M}_{b_{\nu_i}}^{\circ \preceq \mu_i, \theta_i' \sharp}$ on reduced subschemes (cf. proof of theorem \ref{thm:CoveringMorphismProperties} for the construction of such $(\theta_i')_i$). Then $\prod_i \MC{M}_{b_{\nu_i}}^{\preceq \mu_i, \theta_i \sharp}$ is contained in the formal completion of $\prod_i \MB{M}_{b_{\nu_i}}^{\preceq \mu_i, \theta_i \sharp}$ and hence in $\prod_i \MC{M}_{b_{\nu_i}}^{\circ \preceq \mu_i, \theta_i' \sharp}$. The inclusion $\prod_i \MC{M}_{b_{\nu_i}}^{\preceq \mu_i, \theta_i \sharp} \subset \prod_i \MC{M}_{b_{\nu_i}}^{\circ \preceq \mu_i, \theta_i' \sharp}$ is adic, because both formal schemes are adic formal subschemes of $\prod_i \MC{M}_{b_{\nu_i}}^{\preceq \mu_i \sharp}$. 
Fix now some $(d_i')_i \geq (d_i)_i$ such that there is a morphism
\[\widehat{\dot{\pi}}_{(d_i')}: \prod_i \MC{M}_{b_{\nu_i}}^{\circ \preceq \mu_i, \theta_i' \sharp} \times_{\Spf E[[\zeta_1, \ldots, \zeta_n]]^\sharp} \MF{Ig}^{(d_i') \sharp}_U \to \MF{N}^{(\nu_i) \sharp}_U.\]
Then $\widehat{\pi}_{(d_i')}$ is adic as the restriction of an adic morphism to an adic closed formal subscheme. As $\MF{Ig}^{(d_i') \sharp}_U \to \MF{Ig}^{(d_i) \sharp}_U$ is a  finite \'etale cover, it follows that $\widehat{\pi}_{(d_i)}$ is adic as well. With the same arguments it follows as well, that $\widehat{\pi}_{(d_i)}$ is locally of finite type. \\
So to show finiteness, it suffices to check this on reduced subschemes. But there $\widehat{\pi}_{(d_i)}$ is just the restriction of $\pi_{(d_i)}$ to a closed subscheme, hence finite by theorem \ref{thm:CoveringMorphismProperties}b). This proves part b). \\
To get the surjectivity assertion, assume that $(\theta_i)_i$, $(d_i)_i$ are sufficiently large such that $\widehat{\dot{\pi}}_{(d_i)}$ satisfies the surjectivity assertion of proposition \ref{prop:FormalCoveringMorphEtaleSurjective}. Then $\widehat{\dot{\pi}}_{(d_i)}$ is just the restriction of $\widehat{\pi}_{(d_i)}$ to $\prod_i \MC{M}_{b_{\nu_i}}^{\circ \preceq \mu_i, \theta_i \sharp} \times_{\Spf E[[\zeta_1, \ldots, \zeta_n]]^\sharp} \MF{Ig}^{(d_i) \sharp}_U$, hence we may just take the lift constructed in \ref{prop:FormalCoveringMorphEtaleSurjective} and view it as a morphism to $\prod_i \MC{M}_{b_{\nu_i}}^{\preceq \mu_i, \theta_i \sharp} \times_{\Spf E[[\zeta_1, \ldots, \zeta_n]]^\sharp} \MF{Ig}^{(d_i) \sharp}_U$. \exit

\rem{}{
We just briefly sketch another way to obtain a finite covering morphism on formal schemes. Consider for this the locus $\widetilde{\MC{M}}_{b_{\nu_i}}^{\preceq \mu_i, \theta_i \sharp}$ in the Rapoport-Zink space $\MC{M}_{b_{\nu_i}}^{\preceq \mu_i \sharp}$, where the universal quasi-isogeny is weakly bounded by $\theta_i$. It is not hard to see that this defines an adic formal scheme locally of finite type over $\Spf E[[\zeta_i]]^\sharp$. Then the construction \ref{const:CoverMorphism} immediately yields a morphism of formal schemes
\[\widehat{\tilde{\pi}}_{(d_i)}: \prod_i \widetilde{\MC{M}}_{b_{\nu_i}}^{\preceq \mu_i, \theta_i \sharp} \times_{\Spf E[[\zeta_1, \ldots, \zeta_n]]^\sharp} \MF{Ig}^{(d_i) \sharp}_U \to \MF{N}^{(\nu_i) \sharp}_U,\]
that is adic and finite. However it will never be surjective, i.e. its scheme-theoretic image will always be a proper closed formal subscheme of $\MF{N}^{(\nu_i) \sharp}_U$, which is adic over $\Spf E[[\zeta_1, \ldots, \zeta_n]]^\sharp$. Moreover it will not have good properties with respect to passing to the generic fiber as done in the next section: The associated adic spaces $\widetilde{\MC{M}}_{b_{\nu_i}}^{\preceq \mu_i, \theta_i \sharp \; an}$ will be rather small compared to $\MC{M}_{b_{\nu_i}}^{\preceq \mu_i, \theta_i \sharp \; an}$. In particular the $\widetilde{\MC{M}}_{b_{\nu_i}}^{\preceq \mu_i, \theta_i \sharp \; an}$ for varying $\theta_i$ do not cover the space $\MC{M}_{b_{\nu_i}}^{\preceq \mu_i \sharp \; an}$. 
}

\subsection{Morphisms on Huber's adic spaces}\label{subsec:AdicCoverMorphism}
We deal now with the structure of the generic fiber of the formal schemes considered above. As already mentioned, the spaces are rarely adic over $\Spf E[[\zeta_1, \ldots, \zeta_n]]^\sharp$, which enforces us to use Huber's theory in its most general version. So let us recall Huber's generic fiber functor and some of its properties here for the reader's convenience:

\thm{\label{thm:DefAssociatedAdicSpace}}{\cite[propositions 4.1 and 4.2]{HuberAdicSpaces1994}}{
There exists a functor (called $t$ in \cite{HuberAdicSpaces1994})
\[(-)^{an}: \{\textnormal{locally noetherian formal schemes}\} \to \{\textnormal{adic spaces}\}\]
together with morphisms of locally ringed topological spaces $\pi_X: X^{an} \to X$ (for any locally noetherian formal scheme $X$), having the following universal property: 
If $X$ is any locally noetherian formal scheme, $T$ any adic space and $f: T \to X$ any morphism of locally ringed topological spaces, then there exists a unique morphism $\tilde{f}: T \to X^{an}$ of adic spaces with $f = \pi_X \circ \tilde{f}$. \\
Let $f: X \to X'$ be a morphism of locally noetherian formal schemes and $f^{an}: X^{an} \to X'^{an}$ the corresponding morphism of adic spaces. Then 
\begin{itemize}
 \item $f$ is adic if and only if $f^{an}$ is adic.
 \item $f$ is locally of finite type if and only if $f^{an}$ is locally of finite type.
\end{itemize}
}

\rem{}{
i) The analytification of a formal scheme is in general not analytic in the sense of adic spaces (cf. \cite[after theorem 3.5]{HuberAdicSpaces1994}). But as remarked in \cite[section 1.9]{HuberEtaleCohom}, passing to the open subspace of analytic points does not lose much information, e.g. one still has some (slightly more restrictive) universal property with respect to morphisms to the original formal scheme. For more details see corollary \ref{cor:CoveringMorphismAnalytic}. \\
ii) Note that the analytification functor extends to perfections $X^\sharp$ of a locally noetherian formal schemes $X$ by defining $(X^\sharp)^{an}$ as the perfection of $X^{an}$. This way, we get adic spaces for all formal schemes considered above, except for the infinite Igusa variety $\MF{Ig}^{(\infty_i)}_U$.
}

\rem{}{
In the analogous world of Shimura varieties in mixed characteristic, one usually views the associated generic fiber as a rigid analytic space. So let us discuss the relative merits of adic spaces versus rigid analytic spaces: \\
In \cite{RaynaudRigid} Raynaud defined a functor
\[(-)^{rig}: \begin{Bmatrix} \textnormal{formal schemes topologically} \\ \textnormal{of finite type over a DVR} \; R\end{Bmatrix} \to \begin{Bmatrix} \textnormal{quasi-compact rigid analytic spaces} \\ \textnormal{topologically of finite type over} \, \Sp R \end{Bmatrix}\]
So two problems occur here: First of all we are over the base ring $E[[\zeta_1, \ldots, \zeta_n]]^\sharp$, which is not a DVR and secondly our formal schemes are usually not topologically of finite type (or even adic) over $\Spf E[[\zeta_1, \ldots, \zeta_n]]^\sharp$. 
The first problem was at least partially solved by Bosch and L\"utkebohmert in \cite{BoschLutkeRigid}, where they followed Raynaud's approach to define the category of rigid analytic spaces as the localization of the category of admissible formal schemes (in the sense that it is an adic formal schemes locally of finite presentation, such that the ideal of definition does not contain any nilpotent elements) by admissible blow-ups. In this way they obtain a functor
\[(-)^{rig}: \begin{Bmatrix} \textnormal{admissible quasi-compact} \\ \textnormal{formal schemes over} \; S \end{Bmatrix} \to \begin{Bmatrix} \textnormal{quasi-separated rigid analytic spaces} \\ \textnormal{locally of finite presentation over} \, S^{rig} \end{Bmatrix}\]
where $S$ is any noetherian formal scheme, whose structure sheaf is generated by some coherent ideal. On the other hand Berthelot \cite{BerthelotRigidGeometry} generalized Raynaud's construction to obtain a functor
\[(-)^{rig}: \begin{Bmatrix} \textnormal{locally noetherian} \\ \textnormal{formal schemes over a DVR} \; R \end{Bmatrix} \to \begin{Bmatrix} \textnormal{rigid analytic} \\ \textnormal{spaces over} \, \Sp R \end{Bmatrix}.\]
Although there are several technical problems, it seems possible to construct a functor from formal schemes over $\Spf E[[\zeta_1, \ldots, \zeta_n]]^\sharp$ to rigid analytic spaces over $\Spf E[[\zeta_1, \ldots, \zeta_n]]^{\sharp rig}$, by combining both constructions. However we do not know any reference, where this is worked out. \\
Note in particular that many of the following lemmas are well-known for the functor from formal schemes to rigid analytic spaces (at least for Raynaud's original setup) on the one hand, and for the functor from rigid analytic spaces to adic spaces (cf. the various statements in \cite[section 1]{HuberEtaleCohom}) on the other hand. So generalizing the arguments presented there to our situation, should give alternative proofs of these facts. \\
Overall we felt, that using Huber's already existing functor and then proving certain properties by hand, is the easier and more accessible way. 
}

\lem{\label{lem:AdicSpaceFiberProduct}}{}{
If $X$, $Y$ and $Z$ are (perfections of) locally noetherian formal schemes and $f: X \to Z$, $g: Y \to Z$ any morphisms, then the fiber product $X^{an} \times_{Z^{an}} Y^{an}$ exists in the category of adic spaces. In fact it is equal to $(X \times_Z Y)^{an}$. 
}

\rem{}{
We do not assume that $f$ or $g$ are adic or even locally of finite type. Hence we cannot simply refer to \cite[proposition 1.2.2]{HuberEtaleCohom} for the existence of $X^{an} \times_{Z^{an}} Y^{an}$. 
}

\prooof
$X \times_Z Y$ is again locally noetherian (or a perfection thereof), hence $(X \times_Z Y)^{an}$ exists as an adic space. The analytification of the canonical projections $pr_1: X \times_Z Y \to X$ respectively $pr_2: X \times_Z Y \to Y$ induce maps of adic spaces $pr_1^{an}: (X \times_Z Y)^{an} \to X^{an}$ and $pr_2^{an}: (X \times_Z Y)^{an} \to Y^{an}$. We check that the triple $((X \times_Z Y)^{an}, pr_1^{an}, pr_2^{an})$ has the universal property of a fiber product in the category of adic spaces: 
Let $T$ be any adic space and $a: T \to X^{an}$, $b: T \to Y^{an}$ two morphisms of adic spaces satisfying $f^{an} \circ a = g^{an} \circ b: T \to Z^{an}$. We obtain morphisms of locally ringed spaces $\pi_X \circ a: T \to X$ and $\pi_y \circ b: T \to Y$, which coincide when mapped to $Z$. Now observe that the fiber product of locally ringed spaces coincides with the fiber product for (formal) schemes, cf. \cite[Errata to EGAI: proposition 1.8.1]{GroEGAII} (or for a more detailed proof \cite[theorem 8]{GillamLocRingedSpaces}). Thus we obtain a morphism of locally ringed spaces $(\pi_X \circ a, \pi_Y \circ b): T \to X \times_Z Y$. 
By the universal property of the analytification functor, we get $(\pi_X \circ a, \pi_Y \circ b)^{\sim}: T \to (X \times_Z Y)^{an}$. One easily checks that $pr_1^{an} \circ (\pi_X \circ a, \pi_Y \circ b)^{\sim} = a$ and similarly for the second projection. Indeed it suffices to prove this in the category of locally ringed spaces after composing it with $\pi_X: X^{an} \to X$. Then equality holds by construction. \\
Uniqueness is just an exercise in shifting around formal symbols: Assume we have any map $(a, b): T \to (X \times_Z Y)^{an}$ fitting in the fiber product diagram of adic spaces. Then observe
\[pr_1 \circ \pi_{X \times_Z Y} \circ (a, b) = \pi_X \circ pr_1^{an} \circ (a, b) = \pi_X \circ a: T \to X\]
and similarly for the other projection. In particular $\pi_{X \times_Z Y} \circ (a, b): T \to X \times_Z Y$ is uniquely determined. By the universal property of the analytification, this in turn determines $(a, b)$ uniquely. \exit

\lem{\label{lem:AdicSpaceOpenImmersion}}{}{
If $f: U \to X$ is an open embedding of (perfections of) locally noetherian formal schemes, then its analytification $f^{an}: X^{an} \to Y^{an}$ is an open embedding of adic spaces.
}

\prooof
This is an immediate consequence of the construction of the analytification functor. \exit

\lem{\label{lem:AdicSpaceEtale}}{}{
a) If $f: X \to Y$ is an \'etale morphism between (perfections of) locally noetherian formal schemes, then $f^{an}: X^{an} \to Y^{an}$ is an \'etale morphism of adic spaces (as defined in \cite[definition 1.6.5]{HuberEtaleCohom}). \\
b) If $f: X \to Y$ is a finite morphism between (perfections of) locally noetherian formal schemes, then $f^{an}: X^{an} \to Y^{an}$ is a finite morphism of adic spaces (as defined in \cite[1.4.4]{HuberEtaleCohom}).
}

\prooof
a) This is just \cite[lemma 3.5.1i)]{HuberEtaleCohom}. \\
b) We may check this locally on $Y$ resp. $Y^{an}$. So assume that $Y = \Spf A$ for some adic ring $A$. Then by \cite[proposition 4.8.1]{GroEGAIII} there exists a finite adic $A$-algebra $B$ with $X = \Spf B$. So $f^{an}: X^{an} = \Spa (B, B) \to Y^{an} = \Spa (A, A)$ is the morphism induced by the morphism of affinoid rings $(A, A) \to (B, B)$. But this is finite in the sense of \cite[1.4.2]{HuberEtaleCohom}, which gives the desired result. \vspace{5mm} \exit

We now apply all this theory to Rapoport-Zink spaces, Igusa varieties, Newton strata and the covering morphisms as constructed in the previous section. Throughout the remaining part of this section we abbreviate $\Spa E[[\zeta_1, \ldots, \zeta_n]]^\sharp \coloneqq \Spf E[[\zeta_1, \ldots, \zeta_n]]^{\sharp an} = \Spa (E[[\zeta_1, \ldots, \zeta_n]]^\sharp, E[[\zeta_1, \ldots, \zeta_n]]^\sharp)$. 

\thm{\label{thm:GenericCoveringMorphism}}{}{
Fix $(\theta_i)_i$ and $(d_i)_i$ sufficiently large as in proposition \ref{prop:CoverMorphDescent}. Then there is an \'etale map
\[\widehat{\dot{\pi}}_{(d_i)}^{an}: \prod_i \MC{M}_{b_{\nu_i}}^{\circ \preceq \mu_i, \theta_i \sharp \; an} \times_{\Spa E[[\zeta_1, \ldots, \zeta_n]]^\sharp} \MF{Ig}^{(d_i) \sharp \; an}_U \to \MF{N}^{(\nu_i) \sharp \; an}_U\]
and a finite map
\[\widehat{\pi}_{(d_i)}^{an}: \prod_i \MC{M}_{b_{\nu_i}}^{\preceq \mu_i, \theta_i \sharp \; an} \times_{\Spa E[[\zeta_1, \ldots, \zeta_n]]^\sharp} \MF{Ig}^{(d_i) \sharp \; an}_U \to \MF{N}^{(\nu_i) \sharp \; an}_U\]
of adic spaces over $\Spa E[[\zeta_1, \ldots, \zeta_n]]^\sharp$. They are compatible for varying $(\theta_i)_i$ and $(d_i)_i$.
}

\prooof
Apply Huber's generic fiber functor to the morphism constructed in proposition \ref{prop:FormalCoveringMorphEtaleConstruction} respectively \ref{prop:FormalCoveringMorphFinite}. Then use lemma \ref{lem:AdicSpaceFiberProduct} to get the above description of the source spaces and use lemma \ref{lem:AdicSpaceEtale} to carry over the property \'etale respectively finite. \exit

\rem{}{
To get a generic version of the infinite covering morphism $\widehat{\pi}_{(\infty_i)}$ is far more difficult (if not impossible). The main problem is that $\MF{Ig}^{(\infty_i) \sharp}_U$ is not locally noetherian (or even a perfection thereof), so Huber's functor cannot be applied. The alternative approach to define it via the limit $\varprojlim_{(d_i)} \widehat{\dot{\pi}}_{(d_i)}^{an}$ respectively the limit $\varprojlim_{(d_i)} \widehat{\pi}_{(d_i)}^{an}$ causes problems as well, because infinite inverse limits of adic spaces do not need to exist in the category of adic spaces, even if the transition morphisms are finite \'etale as in our situation.
}

Now that we have a good understanding of single adic Newton strata $\MF{N}^{(\nu_i) \sharp \; an}_U$, it remains to discuss their relationship with $\MF{X}^{\mmu \sharp \; an}_U$, the adic space associated to (the perfection of the formal completion of the special fiber of) the moduli space of all bounded global $G$-shtukas. For convenience we will switch to the non-perfect setting in order to have noetherian adic spaces. Note however that most assertions are stable under passing to the perfection anyway.

\lem{\label{lem:AdicStrataAffine}}{}{
Let $A$ be a noetherian $E[[\zeta_1, \ldots, \zeta_n]]$-algebra, which is complete for the ideal $I_0 = (\zeta_1, \ldots, \zeta_n)$. Let $S \subset A$ be a multiplicative subset and $I \subset A$ some radical ideal containing $I_0$. Let $\widehat{A[S^{-1}]}{}^{I}$ be the $I$-adic completion of the localization of $A$ at $S$. 
Then there exists a morphism of adic spaces $\delta^{an}: \Spf \widehat{A[S^{-1}]}{}^{I an} \to \Spf A^{an}$ over $\Spf E[[\zeta_1, \ldots, \zeta_n]]^{an}$.
Moreover:\\
a) $\Spf \widehat{A[S^{-1}]}{}^{I}$ is the formal completion of $\Spf A$ along $\Sp A[S^{-1}]/I \subset \Sp A/I_0$. \\
b) $\delta^{an}$ is injective (as a map of topological spaces). \\
c) The image of $\delta^{an}$ consists of all points in $\Spf A^{an}$ corresponding to a continuous valuation $v: A \to \Gamma \cup \{0\}$ satisfying 
\begin{enumerate}
 \item $v(s) = 1$ for all $s \in S$ and
 \item for each $x \in I$ there is some positive real number $\varepsilon_x$ such that $v(x) < 1 - \varepsilon_x$.
\end{enumerate}
}

\prooof
By definition of the rings there exists a natural morphism $\delta: \Spf \widehat{A[S^{-1}]}{}^{I} \to \Spf A$ of formal schemes over $\Spf E[[\zeta_1, \ldots, \zeta_n]]$. Then define $\delta^{an}$ as the analytification of this morphism. Part a) of the lemma is immediate from the construction. \\
For part b) consider any continuous valuation $v: A \to \Gamma \cup \{0\}$. Then there is at most one extension of this valuation to the ring $A[S^{-1}]$ and by continuity at most one extension to the completion $\widehat{A[S^{-1}]}{}^{I_0}$. Now the elements in $\widehat{A[S^{-1}]}{}^{I_0}$ are dense in $\widehat{A[S^{-1}]}{}^I$ (especially since we change from the $I_0$-adic topology to the coarser $I$-adic one). Hence there is at most one valuation $\tilde{v}: \widehat{A[S^{-1}]}{}^I \to \Gamma \cup \{0\}$ extending $v$, i.e. at most one preimage of $v$ under $\delta^{an}$. Note that this argument actually works for any valuation spectrum, and not only for adic spaces. \\
To prove assertion c) we first show that the conditions 1. and 2. are necessary. If $v$ lies in the image, then it extends to a continuous valuation $\tilde{v}: \widehat{A[S^{-1}]}{}^I \to \Gamma \cup \{0\}$ satisfying $\tilde{v}(x) \leq 1$ for all $x \in \widehat{A[S^{-1}]}{}^I$. From this last condition it follows immediately that all units in $\widehat{A[S^{-1}]}{}^I$ have valuation $1$. In particular this holds for the elements in $S$. Moreover continuity for the $I$-adic topology implies that all elements $x \in I$ are topologically nilpotent, i.e. satisfy $\tilde{v}(x) < 1 - \varepsilon_x$ for some $\varepsilon_x$. \\
Conversely, assume now that $v \in \Spf A^{an} = \Spa(A, A)$ is a valuation satisfying conditions 1. and 2.. Then we have to see that all extensions in the proof of part b) actually exist. The extension $v'$ to $A[S^{-1}]$ and hence to $\widehat{A[S^{-1}]}{}^{I_0}$ exists because $v(s) \neq 0$ for all $s \in S_{\nu_i}$. However only the full condition $v(s) = 1$ implies that $v'(x) \leq 1$ for all $x \in \widehat{A[S^{-1}]}{}^{I_0}$. Thus we obtain an element $v' \in \Spa (\widehat{A[S^{-1}]}{}^{I_0}, \widehat{A[S^{-1}]}{}^{I_0})$.
Next we have to pass from the $I_0$-adic topology to the $I$-adic one. But a valuation $v'$ of $\widehat{A[S^{-1}]}{}^{I_0}$ is continuous for the $I$-adic topology if (and only if) there exists one positive real number $\varepsilon$ such that $v'(x) < 1 - \varepsilon$ for all $x \in I$. As $A$ was assumed to be noetherian, $I$ is finitely generated and hence this condition is actually equivalent to condition 2.. 
Now any continuous valuation extends automatically to the completion of the underlying ring (with respect to the chosen topology). Thus $v'$ extends to a continuous valuation $\tilde{v}: \widehat{A[S^{-1}]}{}^I \to \Gamma \cup \{0\}$. Moreover $\tilde{v}(x) \leq 1$ for all $x \in \widehat{A[S^{-1}]}{}^I$, because this holds over the dense subset $\widehat{A[S^{-1}]}{}^{I_0}$. Thus we get a preimage $\tilde{v} \in \Spf \widehat{A[S^{-1}]}{}^{I an} = \Spa(\widehat{A[S^{-1}]}{}^I, \widehat{A[S^{-1}]}{}^I)$ of $v$ under $\delta^{an}$ as claimed.  \exit

\prop{\label{prop:AdicNewtonStrataImage}}{}{
a) There is a canonical morphism $\delta_{\nu_i}: \MF{N}^{(\nu_i) \; an}_U \to \MF{X}^{\mmu \; an}_U$ over $\Spa E[[\zeta_1, \ldots, \zeta_n]]$. \\
b) $\delta_{\nu_i}$ is injective (as a map of topological spaces). \\
c) The images of $\delta_{\nu_i}$ for varying $(\nu_i)_i$ are disjoint. \\
d) The set of morphisms $\{\delta_{\nu_i}\}_{(\nu_i)_i}$ is jointly surjective on rank-$1$-points of $\MF{X}^{\mmu \; an}_U$.
}

\prooof
a) By definition we have a canonical (non-adic) inclusion of formal schemes $\MF{N}^{(\nu_i)}_U \to \MF{X}^{\mmu}_U$. Then $\delta_{\nu_i}$ is just the analytification of this inclusion. \\
The next two assertions are shown by reduction to the affine situation of the previous lemma. As $\MC{N}^{(\nu_i)}_U \subset \MB{X}^{\mmu}_U$ is locally closed, there are open affine subsets $V_j = \Sp A_{0, j} \subset \MB{X}^{\mmu}_U$ covering $\MC{N}^{(\nu_i)}_U$ such that $V_j \cap \MC{N}^{(\nu_i)}_U = \Sp A_{0, j}[S_{0, j}^{-1}]/I_{0, j}$ for some multiplicative subset $S_{0, j}$ and an ideal $I_{0, j}$. Then passing to formal completions, one gets open formal subschemes $\MF{V}_j = \Spf A_j \subset \MF{X}^{\mmu}_U$ with underlying reduced fibers $V_j$ such that $\MF{N}^{(\nu_i)}_U \cap \MF{V}_j \to \MF{X}^{\mmu}_U \cap \MF{V}_j$ is given by $\Spf \widehat{A_j[S_j^{-1}]}{}^{I_j} \to \Spf A_j$ for any lift $S_j$ of $S_{0, j}$ and $I_j$ the preimage of $I_{0, j}$ under the canonical morphism $A_j \to A_{0, j}$. \\
b) Consider such an open cover $\{V_j\}_j$ as described above. Then $\MF{V}_j$ being open in $\MF{X}^{\mmu}_U$ implies
\[\delta_{\nu_i}^{-1}(\MF{V}_j^{an}) = \delta_{\nu_i}^{-1} \pi_{\MF{X}^{\mmu}_U}^{-1} (V_j) = \pi_{\MF{N}^{(\nu_i)}_U}^{-1} (V_j \cap \MC{N}^{(\nu_i)}_U) = (\MF{N}^{(\nu_i)}_U \cap \MF{V}_j)^{an}\]
and lemma \ref{lem:AdicStrataAffine}b) gives the injectivity of 
\[\delta_{\nu_i}: \delta_{\nu_i}^{-1}(\MF{V}_j^{an}) = (\MF{N}^{(\nu_i)}_U \cap \MF{V}_j)^{an} \to \MF{V}_j^{an}.\]
As the $\MF{V}_j^{an}$ cover the image of $\delta_{\nu_i}$ by construction, the assertion follows. \\
c) It suffices to show that the images of $\pi_{\MF{X}^{\mmu}_U} \circ \delta_{\nu_i}: \MF{N}^{(\nu_i) \; an}_U \to \MB{X}^{\mmu}_U$ are disjoint. But
\[\pi_{\MF{X}^{\mmu}_U} \circ \delta_{\nu_i}(\MF{N}^{(\nu_i) \; an}_U) = \pi_{\MF{N}^{(\nu_i)}_U}(\MF{N}^{(\nu_i) \; an}_U) = \MC{N}^{(\nu_i)}_U \subset \MB{X}^{\mmu}_U\]
and the different Newton strata $\MC{N}^{(\nu_i)}_U$ in the special fiber are by construction disjoint. \\
d) Let $y \in \MF{X}^{\mmu \; an}_U$ be a point given by a rank-$1$-valuation. Then $y$ corresponds to a morphism $y: \Spa(K, K^+) \to \MF{X}^{\mmu \; an}_U$, where $(K, K^+)$ is a complete affinoid field. From the assumption on the rank, it follows that $K^+ = K^{\circ}$ is the subalgebra of power-bounded elements. In particular the underlying topological space of $\Spa(K, K^+)$ has only one single point. 
Now $y$ gives a morphism of locally ringed spaces $\pi_{\MF{X}^{\mmu}_U} \circ y: \Spa(K, K^+) \to \MF{X}^{\mmu}_U$, which maps into some point $|y| \in |\MF{X}^{\mmu}_U| = |\MB{X}^{\mmu}_U|$, where $|\cdot|$ denotes the underlying topological space. As the Newton strata cover $\MB{X}^{\mmu}_U$, there is some $\nu_i$ with $|y| \in |\MC{N}^{(\nu_i)}_U|$. Now $\MF{N}^{(\nu_i)}_U$ is the formal completion of $\MF{X}^{\mmu}_U$ in $\MC{N}^{(\nu_i)}_U$. Hence any morphism to $\MF{X}^{\mmu}_U$, which topologically maps into $|\MC{N}^{(\nu_i)}_U|$, actually factors through $\MF{N}^{(\nu_i)}_U$. 
In particular we obtain a morphism $\pi_{\MF{X}^{\mmu}_U} \circ y: \Spa(K, K^+) \to \MF{N}^{(\nu_i)}_U$. By the universal property of analytification, there is as well a morphism of adic spaces $(\pi_{\MF{X}^{\mmu}_U} \circ y)^{an}: \Spa(K, K^+) \to \MF{N}^{(\nu_i) \; an}_U$. Then by construction, $(\pi_{\MF{X}^{\mmu}_U} \circ y)^{an}$ defines a preimage of $y$ under $\delta_{\nu_i}$. \\
Note that for arbitrary $y \in \MF{X}^{\mmu \; an}_U$ the argument fails, because the image of $\pi_{\MF{X}^{\mmu}_U} \circ y$ consists no longer of one point and hence may spread out (topologically) over several Newton strata. \exit

\rem{}{
The failure of the $\delta_{\nu_i}(\MF{N}^{(\nu_i) \; an}_U)$ to be locally closed or to be jointly surjective has a rather simple explanation: In the case of schemes, a section can either be a unit or topologically nilpotent in a formal neighborhood of a point and such distinctions determine, in which stratum (respectively formal completion of it) the point lies. However in the case of adic spaces, there is a third possibility (which appears only when considering valuations of rank at least $2$) for the section, namely being power-bounded but not a unit in the valuation ring. 
}

\lem{\label{lem:AdicNewtonStrataPartition}}{}{
There exists a partition of $\MF{X}^{\mmu \; an}_U$ (over $\Spa E[[\zeta_1, \ldots, \zeta_n]]$) into locally closed subspaces $\MF{N}^{(\nu_i) \; an, strat}_U$, such that each $\MF{N}^{(\nu_i) \; an, strat}_U$ contains $\MF{N}^{(\nu_i) \; an}_U$ as an open dense subset. Moreover for all $(\nu_i)_i$ the union
\[\bigcup_{(\nu_i') \preceq (\nu_i)} \MF{N}^{(\nu_i') \; an, strat}_U\]
is closed.
}

\prooof
Define $\MF{N}^{(\nu_i) \; an, strat}_U \coloneqq \pi_{\MF{X}^{\mmu}_U}^{-1}(\MC{N}^{(\nu_i)}_U)$ as the preimage of the Newton stratum under the specialization morphism. Then the $\MF{N}^{(\nu_i) \; an, strat}_U$ are obviously locally closed, form a partition of $\MF{X}^{\mmu \; an}_U$ and $\bigcup_{(\nu_i') \preceq (\nu_i)} \MF{N}^{(\nu_i') \; an, strat}_U$ is closed as a preimage of a closed set.
Moreover by definition $\MF{N}^{(\nu_i) \; an}_U \subset \MF{N}^{(\nu_i) \; an, strat}_U$. Both form a partition of the rank-$1$-points, hence these points coincide for both $\MF{N}^{(\nu_i) \; an}_U$ and $\MF{N}^{(\nu_i) \; an, strat}_U$. But rank-$1$-points are dense in $\MF{N}^{(\nu_i) \; an, strat}_U$, because specialization morphisms collapse the ``cloud'' of higher-rank specializations around a point of rank $1$. Thus $\MF{N}^{(\nu_i) \; an}_U$ is dense in $\MF{N}^{(\nu_i) \; an, strat}_U$. \\
It remains to see that $\MF{N}^{(\nu_i) \; an}_U$ is indeed open in $\MF{N}^{(\nu_i) \; an, strat}_U$. This can be checked locally and by choosing a cover $\{V_j\}_j$ as in proposition \ref{prop:AdicNewtonStrataImage}, one is reduced to the affine situation as in lemma \ref{lem:AdicStrataAffine}. In other words for a locally closed formal subscheme $\Spf \widehat{A[S^{-1}]}{}^{I} \subset \Spf A$ with underlying reduced fiber $\Sp A[S^{-1}]/I$, we have to see that $\Spf \widehat{A[S^{-1}]}{}^{I \; an}$ is open in $\pi_{\Spf A}^{-1}(\Sp A[S^{-1}]/I)$. But $\pi_{\Spf A}^{-1}(\Sp A[S^{-1}]/I)$ is (essentially by definition) the set of continuous valuations $v: A \to \Gamma \cup \{0\}$ satisfying
\begin{enumerate}
 \item $v(s) = 1$ for all $s \in S$ and
 \item $v(x) < 1$ for all $x \in I$.
\end{enumerate}
Thus $\pi_{\Spf A}^{-1}(\Sp A[S^{-1}]/I) \setminus \Spf \widehat{A[S^{-1}]}{}^{I \; an}$ is given by the condition $v(x) > 1 - \varepsilon$ for all $x \in I$ and all real numbers $\varepsilon > 0$. As $I$ is finitely generated, this is a closed condition for a fixed $\varepsilon$. Hence the desired complement is closed as an infinite intersection of closed subspaces for varying $\varepsilon$. \exit

\rem{}{
i) The $\MF{N}^{(\nu_i) \; an, strat}_U$ likely do not form a stratification in general, as closures of $\MF{N}^{(\nu_i) \; an, strat}_U$ should be far smaller than unions of other strata, due to the finer topology. \\
ii) Note that in general the subsets $\MF{N}^{(\nu_i) \; an, strat}_U$ are only locally closed subsets as topological spaces and not the intersection of an open adic subspace with a closed adic subspace. In particular $\MF{N}^{(\nu_i) \; an, strat}_U$ has in general no canonical structure as an adic space and geometric points $\Spa(K, K^+)$ of $\MF{X}^{\mmu \; an}_U$ will not necessarily factor over one of the $\MF{N}^{(\nu_i) \; an, strat}_U$ (even as maps of topological spaces).
}

For technical reasons, analytic adic spaces have in general better properties with respect to their cohomology. So let us give a version of the covering morphisms in the setup of analytic adic spaces. As there is little difference between both settings, we will not distinguish notationally between them.

\cor{\label{cor:CoveringMorphismAnalytic}}{}{
Fix $(\theta_i)_i$ and $(d_i)_i$ as in proposition \ref{prop:CoverMorphDescent}. \\
a) There is an \'etale map between analytic adic spaces
\[\widehat{\dot{\pi}}_{(d_i)}^{an}: \prod_i \MC{M}_{b_{\nu_i}}^{\circ \preceq \mu_i, \theta_i \; an} \times_{\Spa E[[\zeta_1, \ldots, \zeta_n]]} \MF{Ig}^{(d_i) \; an}_U \to \MF{N}^{(\nu_i) \; an}_U,\]
which are compatible for varying $(\theta_i)_i$ and $(d_i)_i$. \\
b) There is a finite map between analytic adic spaces
\[\widehat{\pi}_{(d_i)}^{an}: \prod_i \MC{M}_{b_{\nu_i}}^{\preceq \mu_i, \theta_i \; an} \times_{\Spa E[[\zeta_1, \ldots, \zeta_n]]} \MF{Ig}^{(d_i) \; an}_U \to \MF{N}^{(\nu_i) \; an}_U,\]
which are compatible for varying $(\theta_i)_i$ and $(d_i)_i$. \\
c) The canonical maps $\delta_{\nu_i}: \MF{N}^{(\nu_i) \; an}_U \to \MF{X}^{\mmu \; an}_U$ are injective and have open image. For varying $(\nu_i)$ their images are disjoint and cover all rank-$1$ points in $\MF{X}^{\mmu \; an}_U$.
}

\prooof
All analytic adic spaces are constructed by taking the locus of analytic points in the respective adic space. Then \cite[proposition 1.9.1]{HuberEtaleCohom} states that these analytic adic spaces still have a good universal property, hence all morphisms appearing in the corollary are well-defined. As restrictions of \'etale respective finite morphisms, they retain these properties. This proves part a) and b). \\
All assertions in part c) are immediate from proposition \ref{prop:AdicNewtonStrataImage}, except for the images being open. As in the previous statements, we may pass to the affine situation of lemma \ref{lem:AdicStrataAffine} and have to see that $\Spf \widehat{A[S^{-1}]}{}^{I \; an} \subset \Spf A^{an}$ is open in the analytic setting. As localization at $S$ defines an open subset anyway (regardless whether taking arbitrary adic spaces or only analytic ones), we may wlog. just consider $\Spf \widehat{A}{}^{I \; an} \subset \Spf A^{an}$. \\
Fix now any analytic point in $\Spf \widehat{A}{}^{I \; an} \subset \Spf A^{an}$ with corresponding valuation $v: A \to \Gamma \cup \{0\}$. Having the $(\zeta_1, \ldots, \zeta_n)$-adic topology on $A$ implies by the definition of analytic points, that there is at least one $\zeta_i$ with $v(\zeta_i) > 0$. As $I$ is finitely generated, we can find by lemma \ref{lem:AdicStrataAffine}c) some positive integer $r \gg 0$ such that $v(x^r) \leq v(\zeta_i)$ for all $x \in I$. Thus 
\[\{v \in \Spf A^{an} \;|\; v(x^r) \leq v(\zeta_i) \neq 0  \; \forall \, x \in I\}\] 
is an open (rational) subspace of (the subset of analytic points of) the adic space $\Spf A^{an}$, which contains $v$. It lies itself inside $\Spf \widehat{A}{}^{I \; an}$, because $\zeta_i$ is topologically nilpotent or equivalently $v(\zeta_i) < 1 - \varepsilon$ for some positive real number $\varepsilon$. \exit

\section{Representations in cohomology}\label{sec:Cohomology}
Our goal is to study the cohomology group
\[\sum\nolimits_i (-1)^i H^i_c\left(\ov{\nabla^{\mmu}_n\MC{H}^1(C, G)}_\eta, \M{Q}_\ell\right),\] 
i.e. the alternating sum of the cohomology groups of the complex
\[R\Gamma_c\left(\ov{\nabla^{\mmu}_n\MC{H}^1(C, G)}_\eta, \M{Q}_\ell\right) \coloneqq \varinjlim_{U \subset G(\M{A}^{int c_i})} \varprojlim_r R\Gamma_c\left(\ov{\nabla^{\mmu}_n\MC{H}_U^1(C, G)}_\eta, \M{Z}/\ell^r\M{Z}\right) \otimes_{\M{Z}_\ell} \M{Q}_\ell\]
as a representation of $G(\M{A}^{c_i}) \times \Gamma_{E'}$. For the definition of the $G(\M{A}^{c_i})$-action, see \ref{const:ActionGenericCohomology}. \\
The main result is theorem \ref{thm:FinalFormula}, which expresses this representation in terms of cohomology groups of Igusa varieties and Rapoport-Zink-spaces. \\
Nevertheless the first step is to analyze the cohomology of the special fiber in detail, which will be done in the next two sections. This allows us to deal in section \ref{sec:TorsionCohomology} with the representation in the generic fiber for torsion sheaves. Finally one has to deal with the limits over these torsion sheaves, which is done in section \ref{sec:LimitsOfMath}.



\subsection{Torsion sheaves with group actions on the special fiber}\label{sec:TorsionSheaves}
Before we can deal with even the cohomology of the special fiber, we need some information about the sheaves we are to deal with. More precisely given a torsion sheaf $\MC{L}$ on $\MC{N}^{(\nu_i) \sharp}_U$, we construct sheaves $\MC{F}^{(\theta_i, d_i)}$ by first pulling back and then pushing forward along the morphism $\dot{\pi}_{(d_i)}$. Note that it is essential to use the \'etale version of the covering morphism here, though we will give an overview at the end of this section what one can still do if one uses $\pi_{(d_i)}$ instead. The direct limit of the sheaves $\MC{F}^{(\theta_i, d_i)}$ has then an action of $\prod_i J_i$. In fact $\MC{F} = \varinjlim_{(\theta_i, d_i)} \MC{F}^{(\theta_i, d_i)}$ is nothing else than the compactly induced $\prod_i J_i$-representation obtained from $\MC{L}$ (with trivial action), at least on geometric points. \\
For Shimura varieties of PEL-type the corresponding statements can be found in \cite[section 5.2]{MantoFoliation}.

Throughout this section we fix a sheaf $\MC{L}$ of abelian torsion groups over $\MC{N}^{(\nu_i) \sharp}_U$, whose torsion order is prime to $p$. Note that by \cite[proposition A.4]{ZhuAffineGrass} this is equivalent to giving such a sheaf on  $\MC{N}^{(\nu_i)}_U$. For each tuple of pairs $(\theta_i, d_i)_i$, such that $\dot{\pi}_{(d_i)}$ exists, define
\[\MC{F}^{(\theta_i, d_i)} \coloneqq \dot{\pi}_{(d_i) \,!} \dot{\pi}_{(d_i)}^* (\MC{L})\]
which is again a sheaf in abelian torsion groups with torsion order prime to $p$. 
\vspace{2mm} \\
As we will deal frequently with the transition morphisms between various Igusa varieties and open subschemes of Rapoport-Zink spaces, we fix the following

\nota{}{
The canonical finite \'etale projection morphism between Igusa varieties for $(d_i') \geq (d_i)$ is denoted by
\[r_{(d_i' - d_i)}: \op{Ig}^{(d_i') \sharp}_U \to \op{Ig}^{(d_i) \sharp}_U.\]
The canonical open immersion between subspaces of Rapoport-Zink spaces for $\theta_i'$ sufficiently large with respect to $\theta_i$ is denoted by
\[\iota_{\theta_i, \theta_i'}: \MB{M}_{b_{\nu_i}}^{\circ \preceq \mu_i, \theta_i \sharp} \to \MB{M}_{b_{\nu_i}}^{\circ \preceq \mu_i, \theta_i' \sharp}\]
We will not distinguish (notationally) between such morphisms and the induced morphism on product spaces $\prod_i \MB{M}_{b_{\nu_i}}^{\circ \preceq \mu_i, \theta_i \sharp} \times_E \op{Ig}^{(d_i) \sharp}_U$.
}

\lem{}{\label{lem:SheafTransitionMaps}}{
a) If $(d_i') \geq (d_i)$, then there is a canonical morphism $r_{(d_i' - d_i)}^*: \MC{F}^{(\theta_i, d_i)} \to \MC{F}^{(\theta_i, d_i')}$. \\
b) If $\theta_i$ and $\theta_i'$ are two weak bounds such that $\MB{M}_{b_{\nu_i}}^{\circ \preceq \mu_i, \theta_i} \subset \MB{M}_{b_{\nu_i}}^{\circ \preceq \mu_i, \theta_i'}$ for all $i$, then there is a canonical morphism $\iota_{\theta_i, \theta_i' \,!}: \MC{F}^{(\theta_i, d_i)} \to \MC{F}^{(\theta_i', d_i)}$. In particular this happens if for every $\vartheta_i \in \theta_i$ there exists a $\vartheta_i' \in \theta_i'$ with $\vartheta_i \preceq \vartheta_i'$. 
}

\prooof
a) We may factor $\dot{\pi}_{(d_i')} = \dot{\pi}_{(d_i)} \circ r_{(d_i' - d_i)}$. Then finiteness of $r_{(d_i' - d_i)}$ implies $r_{(d_i' - d_i) \,!} = r_{(d_i' - d_i) \,*}$, and we may use adjunction of $(r_{(d_i' - d_i) \,*}, r_{(d_i' - d_i)}^*)$ to define a canonical morphism
\[r_{(d_i' - d_i)}^*: \MC{F}^{(\theta_i, d_i)} = \dot{\pi}_{(d_i) \,!} \dot{\pi}_{(d_i)}^* (\MC{L}) \to \dot{\pi}_{(d_i) \,!} r_{(d_i' - d_i) \,!} r_{(d_i' - d_i)}^* \dot{\pi}_{(d_i)}^* (\MC{L}) = \dot{\pi}_{(d_i') \,!} \dot{\pi}_{(d_i')}^* (\MC{L}) = \MC{F}^{(\theta_i, d_i')}\]
b) We proceed as in a): Note first $\dot{\pi}_{(d_i)} = \dot{\pi}_{(d_i)} \circ \iota_{\theta_i, \theta_i'}$. As $\iota_{\theta_i, \theta_i'}$ is an open immersion, we may use adjunction of $(\iota_{\theta_i, \theta_i' \,!}, \iota_{\theta_i, \theta_i'}^!)$ to get
\[\iota_{\theta_i, \theta_i' \,!}: \MC{F}^{(\theta_i, d_i)} = \dot{\pi}_{(d_i) \,!} \dot{\pi}_{(d_i)}^* (\MC{L}) = \dot{\pi}_{(d_i) \,!} \iota_{\theta_i, \theta_i' \,!} \iota_{\theta_i, \theta_i'}^* \dot{\pi}_{(d_i)}^* (\MC{L}) \to \dot{\pi}_{(d_i) \,!} \dot{\pi}_{(d_i)}^* (\MC{L}) = \MC{F}^{(\theta_i', d_i)}\]
\exit

Giving all these transition morphisms, we can form the sheaf of abelian torsion groups
\[\MC{F} = \varinjlim_{(\theta_i, d_i)} \MC{F}^{(\theta_i, d_i)}.\]

\lem{}{\label{lem:SheafGlobalProjection}}{
There exists, up to choice of normalization, a canonical morphism $\dot{\pi}_!: \MC{F} \to \MC{L}$.
}

\prooof
As $\dot{\pi}_{(d_i)}$ is \'etale, adjunction of $(\dot{\pi}_{(d_i) \,!}, \dot{\pi}_{(d_i)}^!)$ gives a canonical morphism
\[\dot{\pi}_{(d_i) \,!}: \MC{F}^{(\theta_i, d_i)} \to \MC{L}\]
As $r_{(d_i' - d_i)}$ is a finite \'etale cover of degree $\prod_i [I_{d_i'}(b_{\nu_i}):I_{d_i}(b_{\nu_i})]$ we get equalities
\[\dot{\pi}_{(d_i') \,!} \circ r_{(d_i' - d_i)}^* = \prod_i [I_{d_i'}(b_{\nu_i}):I_{d_i}(b_{\nu_i})] \cdot \dot{\pi}_{(d_i) \,!} \quad {\rm and} \quad \dot{\pi}_{(d_i) \,!} \circ \iota_{\theta_i, \theta_i' \,!} = \dot{\pi}_{(d_i) \,!}\]
Note that $[I_{d_i'}(b_{\nu_i}):I_{d_i}(b_{\nu_i})]$ is a $p$-power for $d_i \geq 1$, hence invertible on $\MC{L}$. Thus we can set
\[\dot{\pi}_! = \varinjlim_{(\theta_i, d_i)} \prod_i [I_{d_i}(b_{\nu_i}):I_1(b_{\nu_i})]^{-1} \dot{\pi}_{(d_i) \,!}: \MC{F} \to \MC{L}\]
\exit

For the next proposition note that all constructions above commute with base-change to finite field extensions of $E$. In particular we may base-change to $E'$ as defined in \ref{def:JGroup}, i.e. the composite of $E$ and $\M{F}_{q^s}$, where $s$ is the integer appearing in the decency equation.

\prop{}{\label{prop:SheafJActionDefinition}}{
The sheaf $\MC{F}$ over $\MC{N}^{(\nu_i) \sharp}_U \times_E \Sp E'$ admits a canonical smooth action of the group $\prod_i J_i$, induced by the action of this group on the tower $\left(\prod_i \MB{M}_{b_{\nu_i}}^{\preceq \mu_i, \theta_i \sharp} \times_E \op{Ig}^{(d_i) \sharp}_U \times_E \Sp E'\right)_{(\theta_i, d_i)}$. Moreover if we consider the trivial action of $\prod_i J_i$ on $\MC{L}$, then the morphism $\dot{\pi}_!: \MC{F} \to \MC{L}$ is equivariant for the $\prod_i J_i$-action. 
}

\prooof
Consider any element $(\gamma_i)_i \in \prod_i J_i$. Choose now integers $d_{i, \gamma_i}$, $d_{i, \gamma_{\mathrlap{i}}^{}{}^{\!-1}}$ and $\Gamma$-invariant subsets $\theta_{i, \gamma_i}$, $\theta_{i, \gamma_{\mathrlap{i}}^{}{}^{\!-1}}$ as in proposition \ref{prop:JActionEtaleVersion} (for the respective elements $(\gamma_i)_i$ and $(\gamma_i^{-1})_i$). In particular we have an \'etale morphism
\[(\gamma_i)_i: \prod_i \MB{M}_{b_{\nu_i}}^{\circ \preceq \mu_i, \theta_i} \times_E \op{Ig}_{\mathrlap{U}}{}^{\! (d_i + d_{i, \gamma_i} + d_{i, \gamma_{\mathrlap{i}}^{}{}^{\!-1}})} \times_E \Sp E' \to \prod_i \MB{M}_{b_{\nu_i}}^{\circ \preceq \mu_i, \theta_i \oplus \theta_{i, \gamma_i}} \times_E \op{Ig}_{\mathrlap{U}}{}^{\!(d_i + d_{i, \gamma_{\mathrlap{i}}^{}{}^{\!-1}})} \times_E \Sp E'.\]
It follows immediately from theorem \ref{thm:JActionEquivariant} that $\dot{\pi}_{(d_i + d_{i, \gamma_i} + d_{i, \gamma_{\mathrlap{i}}^{}{}^{\!-1}})} = \dot{\pi}_{(d_i + d_{i, \gamma_{\mathrlap{i}}^{}{}^{\!-1}})} \circ (\gamma_i)_i$ holds.  \\
Thus we can use the adjunction between $(\gamma_i)_{i \,!}$ and $(\gamma_i)_i^* = (\gamma_i)_i^!$ to get a morphism
\begin{align*}
\prod_i [I_{d_{i, \gamma_i}}(b_{\nu_i}):I_1(b_{\nu_i})] \cdot (\gamma_i)_i: \MC{F}^{(\theta_i, d_i)} & \xrightarrow{r_{(d_{i, \gamma_i} + d_{i, \gamma_{\mathrlap{i}}^{}{}^{\!-1}})}^*} \MC{F}^{(\theta_i, d_i + d_{i, \gamma_i} + d_{i, \gamma_{\mathrlap{i}}^{}{}^{\!-1}})} \\
 & = \dot{\pi}_{(d_i + d_{i, \gamma_i} + d_{i, \gamma_{\mathrlap{i}}^{}{}^{\!-1}}) \, !} \dot{\pi}_{(d_i + d_{i, \gamma_i} + d_{i, \gamma_{\mathrlap{i}}^{}{}^{\!-1}})}^* (\MC{L}) \\
 & = \dot{\pi}_{(d_i + d_{i, \gamma_{\mathrlap{i}}^{}{}^{\!-1}}) \, !} (\gamma_i)_{i \,!} (\gamma_i)_i^* \dot{\pi}_{(d_i + d_{i, \gamma_{\mathrlap{i}}^{}{}^{\!-1}})}^* (\MC{L}) \\
 & \xrightarrow{(\gamma_i)_{i \,!}} \dot{\pi}_{(d_i + d_{i, \gamma_{\mathrlap{i}}^{}{}^{\!-1}}) \, !} \dot{\pi}_{(d_i + d_{i, \gamma_{\mathrlap{i}}^{}{}^{\!-1}})}^* (\MC{L}) = \MC{F}^{(\theta_i \oplus \theta_{i, \gamma_i}, d_i + d_{i, \gamma_{\mathrlap{i}}^{}{}^{\!-1}})}
\end{align*}
of sheaves on $\MC{N}^{(\nu_i) \sharp}_U \times_E \Sp E'$. As $\prod_i [I_{d_{i, \gamma_i}}(b_{\nu_i}):I_0(b_{\nu_i})]$ acts invertibly on $\MC{F}^{(\theta_i \oplus \theta_{i, \gamma_i}, d_i)}$, this indeed defines a morphism
\[(\gamma_i)_i: \MC{F}^{(\theta_i, d_i)} \to \MC{F}^{(\theta_i \oplus \theta_{i, \gamma_i}, d_i + d_{i, \gamma_{\mathrlap{i}}^{}{}^{\!-1}})}.\]
\textbf{Claim 1:} $(\gamma_i)_i$ does not depend on the choices of $d_{i, \gamma_i}$ and $\theta_{i, \gamma_i}$. \\
It is straight-forward to see that $(\gamma_i)_i: \MC{F}^{(\theta_i, d_i)} \to \MC{F}^{(\theta_i \oplus \theta_{i, \gamma_i}, d_i + d_{i, \gamma_{\mathrlap{i}}^{}{}^{\!-1}})}$ is actually independent of our choice of $\theta_{i, \gamma_i}$'s in the sense that 
\[(\gamma_i')_i = \iota_{\theta_i \oplus \theta_{i, \gamma_i}, \theta_i \oplus \theta_{i, \gamma_i}' \,!} \circ (\gamma_i)_i: \MC{F}^{(\theta_i, d_i)} \to \MC{F}^{(\theta_i \oplus \theta_{i, \gamma_i}', d_i + d_{i, \gamma_{\mathrlap{i}}^{}{}^{\!-1}})}\]
if $\theta_{i, \gamma_i}$ and $\theta_{i, \gamma_i}'$ are chosen in such a way, that an open immersion $\iota_{\theta_i \oplus \theta_{i, \gamma_i}, \theta_i \oplus \theta_{i, \gamma_i}'}$ exists. Here the use of $(\gamma_i')_i$ does not mean, that we change the element in $\prod_i J_i$; it just indicates that this map is related to the second choice of $\theta_{i, \gamma_i}'$'s. A similar notation is used for varying $d_{i, \gamma_i}$ later on. \\
To get independence of the choice of the $d_{i, \gamma_i}$, we have to consider the following diagram 
\[\begin{xy}
 \xymatrix @C=5pc {
   \MC{F}^{(\theta_i, d_i)} \ar^-{r_{(d_{i, \gamma_i} + d_{i, \gamma_{\mathrlap{i}}^{}{}^{\!-1}})}^*}[r] \ar@/_1.0pc/[dr]_-{r_{(d_{i, \gamma_i}' + d_{i, \gamma_{\mathrlap{i}}^{}{}^{\!-1}})}^*} & \MC{F}^{(\theta_i, d_i + d_{i, \gamma_i} + d_{i, \gamma_{\mathrlap{i}}^{}{}^{\!-1}})} \ar^-{(\gamma_i)_{i \,!}}[dr] \ar_-{r_{(d_{i, \gamma_i}' - d_{i, \gamma_i})}^*}[d] &  \\
  & \MC{F}^{(\theta_i, d_i + d_{i, \gamma_i}' + d_{i, \gamma_{\mathrlap{i}}^{}{}^{\!-1}})} \ar_-{(\gamma_i')_{i \,!}}[r] & \MC{F}^{(\theta_i \oplus \theta_{i, \gamma_i}, d_i + d_{i, \gamma_{\mathrlap{i}}^{}{}^{\!-1}})} 
  }
\end{xy} \]
Then one computes
\begin{align*}
(\gamma_i)_i & = \prod_i [I_{d_{i, \gamma_i}}(b_{\nu_i}):I_1(b_{\nu_i})]^{-1} \cdot  (\gamma_i)_{i \,!} \circ r_{(d_{i, \gamma_i} + d_{i, \gamma_{\mathrlap{i}}^{}{}^{\!-1}})}^* \\
& = \prod_i [I_{d_{i, \gamma_i}'}(b_{\nu_i}):I_1(b_{\nu_i})]^{-1} \cdot (\gamma_i)_{i \,!} \circ r_{(d_{i, \gamma_i}' - d_{i, \gamma_i}) \,!} \circ r_{(d_{i, \gamma_i}' - d_{i, \gamma_i})}^* \circ r_{(d_{i, \gamma_i} + d_{i, \gamma_{\mathrlap{i}}^{}{}^{\!-1}})}^* \\
& = \prod_i [I_{d_{i, \gamma_i}'}(b_{\nu_i}):I_1(b_{\nu_i})]^{-1} \cdot (\gamma_i')_{i \,!} \circ r_{(d_{i, \gamma_i}' + d_{i, \gamma_{\mathrlap{i}}^{}{}^{\!-1}})}^* \\
& = (\gamma_i')_i
\end{align*}
\textbf{Claim 2:} $(\gamma_i)_i$ does not depend on the choice of $d_{i, \gamma_{\mathrlap{i}}^{}{}^{\!-1}}$. \\
Consider another choice $d'{}_{\!\!i, \gamma_{\mathrlap{i}}^{}{}^{\!-1}} \geq d_{i, \gamma_{\mathrlap{i}}^{}{}^{\!-1}}$. As in claim $1$ we denote $(\gamma_i)_i: \MC{F}^{(\theta_i, d_i)} \to \MC{F}^{(\theta_i \oplus \theta_{i, \gamma_i}, d_i + d_{i, \gamma_{\mathrlap{i}}^{}{}^{\!-1}})}$ and $(\gamma_i')_i: \MC{F}^{(\theta_i, d_i)} \to \MC{F}^{(\theta_i \oplus \theta_{i, \gamma_i}, d_i + d'{}_{\!\!i, \gamma_{\mathrlap{i}}^{}{}^{\!-1}})}$. Then we wish to show:
\[(\gamma_i')_i = r_{(d'{}_{\!\!i, \gamma_{\mathrlap{i}}^{}{}^{\!-1}} - d_{i, \gamma_{\mathrlap{i}}^{}{}^{\!-1}})}^* \circ (\gamma_i)_i: \MC{F}^{(\theta_i, d_i)} \to \MC{F}^{(\theta_i \oplus \theta_{i, \gamma_i}, d_i + d'{}_{\!\!i, \gamma_{\mathrlap{i}}^{}{}^{\!-1}})}\]
We can compute directly as morphisms from $\MC{F}^{(\theta_i, d_i)}$ to $\MC{F}^{(\theta_i \oplus \theta_{i, \gamma_i}, d_i + d_{i, \gamma_{\mathrlap{i}}^{}{}^{\!-1}})}$:
\begin{align*}
 & r_{(d'{}_{\!\!i, \gamma_{\mathrlap{i}}^{}{}^{\!-1}} - d_{i, \gamma_{\mathrlap{i}}^{}{}^{\!-1}}) \,!} \circ (\gamma_i')_i = \\
 & \quad = \prod_i [I_{d_{i, \gamma_i}}(b_{\nu_i}):I_1(b_{\nu_i})]^{-1} \cdot  r_{(d'{}_{\!\!i, \gamma_{\mathrlap{i}}^{}{}^{\!-1}} - d_{i, \gamma_{\mathrlap{i}}^{}{}^{\!-1}}) \,!} \circ (\gamma_i')_{i \,!} \circ r_{(d_{i, \gamma_i} + d'{}_{\!\!i, \gamma_{\mathrlap{i}}^{}{}^{\!-1}})}^* \\
 & \quad = \prod_i [I_{d_{i, \gamma_i}}(b_{\nu_i}):I_1(b_{\nu_i})]^{-1} \cdot (\gamma_i)_{i \,!} \circ r_{(d'{}_{\!\!i, \gamma_{\mathrlap{i}}^{}{}^{\!-1}} - d_{i, \gamma_{\mathrlap{i}}^{}{}^{\!-1}}) \,!} \circ r_{(d_{i, \gamma_i} + d'{}_{\!\!i, \gamma_{\mathrlap{i}}^{}{}^{\!-1}})}^* \\
 & \quad = \prod_i [I_{d_{i, \gamma_i}}(b_{\nu_i}):I_1(b_{\nu_i})]^{-1} \cdot \prod_i [I_{d_{i, \gamma_{\mathrlap{i}}^{}{}^{\!-1}}}(b_{\nu_i}):I_{d'{}_{\!\!i, \gamma_{\mathrlap{i}}^{}{}^{\!-1}}}(b_{\nu_i})] \cdot (\gamma_i)_{i \,!} \circ  r_{(d_{i, \gamma_i} + d_{i, \gamma_{\mathrlap{i}}^{}{}^{\!-1}})}^* \\
 & \quad = \prod_i [I_{d_{i, \gamma_{\mathrlap{i}}^{}{}^{\!-1}}}(b_{\nu_i}):I_{d'{}_{\!\!i, \gamma_{\mathrlap{i}}^{}{}^{\!-1}}}(b_{\nu_i})] \cdot (\gamma_i)_i \\
 & \quad = r_{(d'{}_{\!\!i, \gamma_{\mathrlap{i}}^{}{}^{\!-1}} - d_{i, \gamma_{\mathrlap{i}}^{}{}^{\!-1}}) \,!} \circ r_{(d_i + d_{i, \gamma_{\mathrlap{i}}^{}{}^{\!-1}}), (d_i + d'{}_{\!\!i, \gamma_{\mathrlap{i}}^{}{}^{\!-1}})}^* \circ (\gamma_i)_i
\end{align*}
View now elements $s \in \MC{F}^{(\theta_i, d_i)}$ as sections (with proper support) over $\prod_i \MB{M}_{b_{\nu_i}}^{\circ \preceq \mu_i, \theta_i \sharp} \times_E \op{Ig}^{(d_i) \sharp}_U \times_E \Sp E'$ with values in $\dot{\pi}_{(d_i)}^*(\MC{L})$ and similarly for other sheaves like $\MC{F}^{(\theta_i \oplus \theta_{i, \gamma_i}, d_i + d'{}_{\!\!i, \gamma_{\mathrlap{i}}^{}{}^{\!-1}})}$. 
As $r_{(d'{}_{\!\!i, \gamma_{\mathrlap{i}}^{}{}^{\!-1}} - d_{i, \gamma_{\mathrlap{i}}^{}{}^{\!-1}})}$ is a finite \'etale morphism, the previous computation implies that it suffices to show the following: For any geometric point $(x, y) \in \prod_i \MB{M}_{b_{\nu_i}}^{\circ \preceq \mu_i, \theta_i \oplus \theta_{i, \gamma_i} \sharp} \times_E \op{Ig}_{\mathrlap{U}}{}^{\! (d_i + d'{}_{\!\!i, \gamma_{\mathrlap{i}}^{}{}^{\!-1}}) \sharp} \times_E \Sp E'$ with image $(x, \ov{y}) \in \prod_i \MB{M}_{b_{\nu_i}}^{\circ \preceq \mu_i, \theta_i \oplus \theta_{i, \gamma_i} \sharp} \times_E \op{Ig}_{\mathrlap{U}}{}^{\! (d_i + d'{}_{\!\!i, \gamma_{\mathrlap{i}}^{}{}^{\!-1}}) \sharp} \times_E \Sp E' $ and formal neighborhoods $\widehat{(x, y)}$ respectively $\widehat{(x, \ov{y})}$ (which can be canonically identified), we have $(\gamma_i')_i(s)|_{\widehat{(x, y)}} = (\gamma_i)_i(s)|_{\widehat{(x, \ov{y})}}$. 
\[\begin{xy}
 \xymatrix @C=2pc {
  \prod_i \MB{M}_{b_{\nu_i}}^{\circ \preceq \mu_i, \theta_i \sharp} \times_E \op{Ig}_{\mathrlap{U}}{}^{\! (d_i + d_{i, \gamma_i} + d'{}_{\!\!i, \gamma_{\mathrlap{i}}^{}{}^{\!-1}}) \sharp} \times_E \Sp E' \ar[d] \ar^{(\gamma_i')_i}[r] & \prod_i \MB{M}_{b_{\nu_i}}^{\circ \preceq \mu_i, \theta_i \oplus \theta_{i, \gamma_i} \sharp} \times_E \op{Ig}_{\mathrlap{U}}{}^{\! (d_i + d'{}_{\!\!i, \gamma_{\mathrlap{i}}^{}{}^{\!-1}}) \sharp} \times_E \Sp E' \ar[d] \\
  \prod_i \MB{M}_{b_{\nu_i}}^{\circ \preceq \mu_i, \theta_i \sharp} \times_E \op{Ig}_{\mathrlap{U}}{}^{\! (d_i + d_{i, \gamma_i} + d_{i, \gamma_{\mathrlap{i}}^{}{}^{\!-1}}) \sharp} \times_E \Sp E' \ar[d] \ar^{(\gamma_i)_i}[r] & \prod_i \MB{M}_{b_{\nu_i}}^{\circ \preceq \mu_i, \theta_i \oplus \theta_{i, \gamma_i} \sharp} \times_E \op{Ig}_{\mathrlap{U}}{}^{\! (d_i + d_{i, \gamma_{\mathrlap{i}}^{}{}^{\!-1}}) \sharp} \times_E \Sp E' \ar^{(\gamma_i^{-1})_i}[d] \\
  \prod_i \MB{M}_{b_{\nu_i}}^{\circ \preceq \mu_i, \theta_i \sharp} \times_E \op{Ig}^{(d_i) \sharp}_U \times_E \Sp E' \quad \ar^-{\iota_{\theta_i, \theta_i \oplus \theta_{i, \gamma_i} \oplus \theta_{i, \gamma_{\mathrlap{i}}^{}{}^{\!-1}}}}[r] & \prod_i \MB{M}_{\mathrlap{b_{\nu_i}}}^{}{}^{ \circ \preceq \mu_i, \theta_i \oplus \theta_{i, \gamma_i} \oplus \theta_{i, \gamma_{\mathrlap{i}}^{}{}^{\!-1}} \sharp} \times_E \op{Ig}^{(d_i) \sharp}_U \times_E \Sp E' 
  }
\end{xy} \]
If $(x, \ov{y})$ does not lie in the image of the morphism $(\gamma_i)_i$, then $(x, y)$ cannot lie in the image of $(\gamma_i')_i$ and both sections are zero. Otherwise lemma \ref{lem:JActionFiberProduct} shows that the upper square is a fiber product diagram and we may identify as well
\[(\gamma_i')_i^{-1}(\widehat{(x, y)}) = (\gamma_i)_i^{-1}(\widehat{(x, \ov{y})})\]
Hence by construction of the morphisms $(\gamma_i)_i$ and $(\gamma_i')_i$ on sheaves it suffices to see
\[r_{(d_{i, \gamma_i} + d'{}_{\!\!i, \gamma_{\mathrlap{i}}^{}{}^{\!-1}})}^*(s)|_{(\gamma_i')_i^{-1}(\widehat{(x, y)})} = r_{(d_{i, \gamma_i} + d_{i, \gamma_{\mathrlap{i}}^{}{}^{\!-1}})}^*(s)|_{(\gamma_i)_i^{-1}(\widehat{(x, \ov{y})})} \qquad (*)\]
To show this consider $(x_0, y_0) \coloneqq (\gamma_i^{-1})_i(x, \ov{y}) \in \prod_i \MB{M}_{\mathrlap{b_{\nu_i}}}^{}{}^{\circ \preceq \mu_i, \theta_i \oplus \theta_{i, \gamma_i} \oplus \theta_{i, \gamma_{\mathrlap{i}}^{}{}^{\!-1}} \sharp} \times_E \op{Ig}^{(d_i) \sharp}_U \times_E \Sp E'$. As $(x, \ov{y})$ has a preimage under $(\gamma_i)_i$ it follows that $(x_0, y_0)$ actually lies in $\prod_i \MB{M}_{b_{\nu_i}}^{\circ \preceq \mu_i, \theta_i \sharp} \times_E \op{Ig}^{(d_i) \sharp}_U \times_E \Sp E'$. 
But indeed 
\begin{align*}
 (\gamma_i)_i^{-1}(\widehat{(x, \ov{y})}) & \subseteq r_{(d_{i, \gamma_i} + d_{i, \gamma_{\mathrlap{i}}^{}{}^{\!-1}})}^{-1}(\widehat{(x_0, y_0)}) \\
 (\gamma_i')_i^{-1}(\widehat{(x, y)}) & \subseteq r_{(d_{i, \gamma_i} + d'{}_{\!\!i, \gamma_{\mathrlap{i}}^{}{}^{\!-1}})}^{-1}(\widehat{(x_0, y_0)}) 
\end{align*}
where $\widehat{(x_0, y_0)}$ denotes the formal neighborhood of $(x_0, y_0)$ in $\prod_i \MB{M}_{b_{\nu_i}}^{\circ \preceq \mu_i, \theta_i \sharp} \times_E \op{Ig}^{(d_i) \sharp}_U \times_E \Sp E'$. The construction of the $\prod_i J_i$-action on (truncated) Igusa varieties implies this on geometric points, from which the statement for formal neighborhoods follows directly. Hence the two sections of $(*)$ are nothing else than the pullback of $s|_{\widehat{(x_0, y_0)}}$ to the respective fibers and are therefore equal. \\
\textbf{Claim 3:} The morphisms $(\gamma_i)_i$ glue to an endomorphism of $\MC{F}$. \\
For this purpose fix all the $d_{i, \gamma_i}$'s, $d_{i, \gamma_{\mathrlap{i}}^{}{}^{\!-1}}$'s and $\theta_{i, \gamma_i}$'s. Then we have to see that the following two diagrams commute:
\[\begin{xy}
 \xymatrix @C=2.5pc {
  \MC{F}^{(\theta_i, d_i)} \ar^-{r_{(d_{i, \gamma_i} + d_{i, \gamma_{\mathrlap{i}}^{}{}^{\!-1}})}^*}[rr] \ar_-{\iota_{\theta_i, \theta_i' \,!}}[d] && \MC{F}^{(\theta_i, d_i + d_{i, \gamma_i} + d_{i, \gamma_{\mathrlap{i}}^{}{}^{\!-1}})} \ar^-{(\gamma_i)_{i \,!}}[r] \ar_-{\iota_{\theta_i, \theta_i' \,!}}[d] & \MC{F}^{(\theta_i \oplus \theta_{i, \gamma_i}, d_i + d_{i, \gamma_{\mathrlap{i}}^{}{}^{\!-1}})} \ar^-{\iota_{\theta_i \oplus \theta_{i, \gamma_i}, \theta_i' \oplus \theta_{i, \gamma_i} \,!}}[d]   \\
  \MC{F}^{(\theta_i', d_i)} \ar_-{r_{(d_{i, \gamma_i} + d_{i, \gamma_{\mathrlap{i}}^{}{}^{\!-1}})}^*}[rr] && \MC{F}^{(\theta_i', d_i + d_{i, \gamma_i} + d_{i, \gamma_{\mathrlap{i}}^{}{}^{\!-1}})} \ar_-{(\gamma_i)_{i \,!}}[r] & \MC{F}^{(\theta_i' \oplus \theta_{i, \gamma_i}, d_i + d_{i, \gamma_{\mathrlap{i}}^{}{}^{\!-1}})} \\
  \MC{F}^{(\theta_i, d_i)} \ar^-{r_{(d_{i, \gamma_i} + d_{i, \gamma_{\mathrlap{i}}^{}{}^{\!-1}})}^*}[rr] \ar_-{r_{(d_i' - d_i)}^*}[d] && \MC{F}^{(\theta_i, d_i + d_{i, \gamma_i} + d_{i, \gamma_{\mathrlap{i}}^{}{}^{\!-1}})} \ar^-{(\gamma_i)_{i \,!}}[r] \ar_-{r_{(d_i' - d_i)}^*}[d] & \MC{F}^{(\theta_i \oplus \theta_{i, \gamma_i}, d_i + d_{i, \gamma_{\mathrlap{i}}^{}{}^{\!-1}})} \ar^-{r_{(d_i' - d_i)}^*}[d] \\ 
  \MC{F}^{(\theta_i, d_i')} \ar_-{r_{(d_{i, \gamma_i} + d_{i, \gamma_{\mathrlap{i}}^{}{}^{\!-1}})}^*}[rr] && \MC{F}^{(\theta_i, d_i' + d_{i, \gamma_i} + d_{i, \gamma_{\mathrlap{i}}^{}{}^{\!-1}})} \ar_-{(\gamma_i)_{i \,!}}[r] & \MC{F}^{(\theta_i \oplus \theta_{i, \gamma_i}, d_i' + d_{i, \gamma_{\mathrlap{i}}^{}{}^{\!-1}})} 
  }
\end{xy} \]
For the upper diagram the left-hand side obviously commutes. The other square does so, because it already commutes on the level of schemes. The lower diagram was already dealt with in claim $2$ by setting $d_{i, \gamma_{\mathrlap{i}}^{}{}^{\!-1}}' = d_{i, \gamma_{\mathrlap{i}}^{}{}^{\!-1}} + (d_i' - d_i)$. \\
\textbf{Claim 4:} The morphisms $(\gamma_i)_i$ define an action of $\prod_i J_i$ on $\MC{F}$. \\
Consider two elements $(\beta_i)_i, (\gamma_i)_i \in \prod_i J_i$ and choose $d_{i, \beta_i}$, $d_{i, \beta_{\mathrlap{i}}^{}{}^{\!-1}}$, $d_{i, \gamma_i}$, $d_{i, \gamma_{\mathrlap{i}}^{}{}^{\!-1}}$, $\theta_{i, \beta_i}$ and $\theta_{i, \gamma_i}$ as usual. Then we may set $d_{i, \beta_i \gamma_i} = d_{i, \beta_i} + d_{i, \gamma_i}$, $d_{i, \gamma_{\mathrlap{i}}^{}{}^{\!-1} \beta_{\mathrlap{i}}^{}{}^{\!-1}} = d_{i, \beta_{\mathrlap{i}}^{}{}^{\!-1}} + d_{i, \gamma_{\mathrlap{i}}^{}{}^{\!-1}}$ and $\theta_{i, \beta_i \gamma_i} = \theta_{i, \beta_i} \oplus \theta_{i, \gamma_i}$ for the element $(\beta_i \cdot \gamma_i)_i$. We obtain the diagram
\[\begin{xy}
 \xymatrix @C=6pc {
  \MC{F}^{(\theta_i, d_i)} \ar^-{r_{(d_{i, \gamma_i} + d_{i, \gamma_{\mathrlap{i}}^{}{}^{\!-1}})}^*}[r] & \MC{F}^{(\theta_i, d_i + d_{i, \gamma_i} + d_{i, \gamma_{\mathrlap{i}}^{}{}^{\!-1}})} \ar^-{r_{(d_{i, \beta_i} + d_{i, \beta_{\mathrlap{i}}^{}{}^{\!-1}})}^*}[r] \ar_-{(\gamma_i)_{i \,!}}[d] & \MC{F}^{(\theta_i, d_i + d_{i, \beta_i \gamma_i} + d_{i, \gamma_{\mathrlap{i}}^{}{}^{\!-1} \beta_{\mathrlap{i}}^{}{}^{\!-1}})} \ar^-{(\gamma_i)_{i \,!}}[d] \\
  & \MC{F}^{(\theta_i \oplus \theta_{i, \gamma_i}, d_i + d_{i, \gamma_{\mathrlap{i}}^{}{}^{\!-1}})} \ar^-{r_{(d_{i, \beta_i} + d_{i, \beta_{\mathrlap{i}}^{}{}^{\!-1}})}^*}[r] & \MC{F}^{(\theta_i \oplus \theta_{i, \gamma_i}, d_i + d_{i, \beta_i} + d_{i, \gamma_{\mathrlap{i}}^{}{}^{\!-1} \beta_{\mathrlap{i}}^{}{}^{\!-1}})} \ar^-{(\beta_i)_{i \,!}}[d]   \\
  & & \MC{F}^{(\theta_i \oplus \theta_{i, \beta_i \gamma_i}, d_i + d_{i, \gamma_{\mathrlap{i}}^{}{}^{\!-1} \beta_{\mathrlap{i}}^{}{}^{\!-1}})}
  }
\end{xy} \]
and the computation done in claim $2$ allow us to conclude
\begin{align*}
& (\beta_i)_i \circ (\gamma_i)_i \\
& = \prod_i [I_{d_{i, \beta_i} + d_{i, \gamma_i}}(b_{\nu_i}):I_0(b_{\nu_i})]^{-1} \cdot (\beta_i)_{i \,!} \circ r_{(d_{i, \beta_i} + d_{i, \beta_{\mathrlap{i}}^{}{}^{\!-1}})}^* \circ (\gamma_i)_{i \,!} \circ r_{(d_{i, \gamma_i} + d_{i, \gamma_{\mathrlap{i}}^{}{}^{\!-1}})}^* \\
& = \prod_i [I_{d_{i, \beta_i} + d_{i, \gamma_i}}(b_{\nu_i}):I_0(b_{\nu_i})]^{-1} \cdot (\beta_i)_{i \,!} \circ (\gamma_i)_{i \,!} \circ r_{(d_{i, \beta_i} + d_{i, \beta_{\mathrlap{i}}^{}{}^{\!-1}})}^* \circ r_{(d_{i, \gamma_i} + d_{i, \gamma_{\mathrlap{i}}^{}{}^{\!-1}})}^* \\
& = \prod_i [I_{d_{i, \beta_i \gamma_i}}(b_{\nu_i}):I_0(b_{\nu_i})]^{-1} \cdot (\beta_i \cdot \gamma_i)_{i \,!} \circ  r_{(d_{i, \beta_i \gamma_i} + d_{i, \gamma_{\mathrlap{i}}^{}{}^{\!-1} \beta_{\mathrlap{i}}^{}{}^{\!-1}})}^* \\
& = (\beta_i \cdot \gamma_i)_i.
\end{align*}
\textbf{Claim 5:} This action is smooth. \\
We wish to show, that every element in $\MC{F}$ is fixed by some open subgroup of $\prod_i J_i$. Let $s \in \MC{F}$ be any section and choose $(\theta_i, d_i)_i$ such that $s \in \MC{F}^{(\theta_i, d_i)}$. Now we consider the open compact subgroup $\Gamma_{(\theta_i, d_i)} \subset \prod_i J_i$ of all quasi-isogenies, which are actually automorphisms such that
\begin{itemize}
 \item they induce the identity modulo $\prod_i I_{d_i}(b_{\nu_i})$ and
 \item stay automorphisms after conjugation by any quasi-isogeny bounded by $(\theta_i)_i$.
\end{itemize}
The first condition ensures that $\Gamma_{(\theta_i, d_i)}$ acts trivially on $\op{Ig}^{(d_i) \sharp}_U \times_E \Sp E'$ and the second condition implies the same for $\prod_i \MB{M}_{b_{\nu_i}}^{\circ \preceq \mu_i, \theta_i \sharp} \times_E \Sp E'$. In particular we may take $d_{i, \gamma_i} = d_{i, \gamma_{\mathrlap{i}}^{}{}^{\!-1}} = 0$ and $\theta_{i, \gamma_i} = \emptyset$ for all $i$ in the definition of $(\gamma_i)_i \in \Gamma_{(\theta_i, d_i)}$ on sheaves and $\MC{F}^{(\theta_i, d_i)}$ is fixed (point-wise) under $\Gamma_{(\theta_i, d_i)}$. This shows that $s \in \MC{F}$ is indeed a smooth vector. \\
\textbf{Claim 6:} $\dot{\pi}_!: \MC{F} \to \MC{L}$ is equivariant. \\
Using $\dot{\pi}_{(d_i + d_{i, \gamma_{\mathrlap{i}}^{}{}^{\!-1}})} \circ (\gamma_i)_i = \dot{\pi}_{(d_i + d_{i, \gamma_i} + d_{i, \gamma_{\mathrlap{i}}^{}{}^{\!-1}})} = \dot{\pi}_{(d_i)} \circ r_{(d_{i, \gamma_i} + d_{i, \gamma_{\mathrlap{i}}^{}{}^{\!-1}})}$ already on the level of schemes, we get
\begin{align*}
 \dot{\pi}_{(d_i + d_{i, \gamma_{\mathrlap{i}}^{}{}^{\!-1}}) \,!} \circ (\gamma_i)_i & = \prod_i [I_{d_{i, \gamma_i}}(b_{\nu_i}):I_0(b_{\nu_i})]^{-1} \cdot \dot{\pi}_{(d_i + d_{i, \gamma_{\mathrlap{i}}^{}{}^{\!-1}}) \,!} \circ (\gamma_i)_{i \,!} \circ r_{(d_{i, \gamma_i} + d_{i, \gamma_{\mathrlap{i}}^{}{}^{\!-1}})}^* \\
 & = \prod_i [I_{d_{i, \gamma_i}}(b_{\nu_i}):I_0(b_{\nu_i})]^{-1} \cdot \dot{\pi}_{(d_i + d_{i, \gamma_i} + d_{i, \gamma_{\mathrlap{i}}^{}{}^{\!-1}}) \,!} \circ r_{(d_{i, \gamma_i} + d_{i, \gamma_{\mathrlap{i}}^{}{}^{\!-1}})}^* \\
 & = \prod_i [I_{d_{i, \gamma_i}}(b_{\nu_i}):I_0(b_{\nu_i})]^{-1} \cdot \dot{\pi}_{(d_i) \,!} \circ r_{(d_{i, \gamma_i} + d_{i, \gamma_{\mathrlap{i}}^{}{}^{\!-1}}) \,!} \circ r_{(d_{i, \gamma_i} + d_{i, \gamma_{\mathrlap{i}}^{}{}^{\!-1}})}^* \\
 & = \prod_i [I_{d_{i, \gamma_{\mathrlap{i}}^{}{}^{\!-1}}}(b_{\nu_i}):I_0(b_{\nu_i})] \cdot \dot{\pi}_{(d_i) \,!}
\end{align*}
The claim follows now by taking the normalization factors appearing in the definition \ref{lem:SheafGlobalProjection} of $\dot{\pi}_!$ into account. \exit
\vspace{3mm} \\
We wish now to study stalks of $\MC{F}$ over geometric points, where the $\prod_i J_i$-action has a particularly nice description.

\lem{\label{lem:SheafFiberHomSpace}}{}{
Let $x \in \MC{N}^{(\nu_i) \sharp}_U \times_E \Sp E'$ be a geometric point. Then there is a canonical isomorphism
\[\MC{F}_x \cong \varinjlim_{(\theta_i, d_i)} \Hom_{sets}(\dot{\pi}_{(d_i)}^{-1}(x), \MC{L}_x)\]
}

\prooof
By definition of the sheaf $\MC{F}$ we have
\begin{align*}
 \MC{F}_x & = \varinjlim_{(\theta_i, d_i)} \MC{F}^{(\theta_i, d_i)}_x = \varinjlim_{(\theta_i, d_i)} (\dot{\pi}_{(d_i) \,!} \dot{\pi}_{(d_i)}^* \MC{L})_x \subseteq \varinjlim_{(\theta_i, d_i)} (\dot{\pi}_{(d_i) \,*} \dot{\pi}_{(d_i)}^* \MC{L})_x = \varinjlim_{(\theta_i, d_i)} \prod_{y \in \dot{\pi}_{(d_i)}^{-1}(x)} (\dot{\pi}_{(d_i)}^* \MC{L})_y \\
 &  = \varinjlim_{(\theta_i, d_i)} \prod_{y \in \dot{\pi}_{(d_i)}^{-1}(x)} \MC{L}_{\dot{\pi}_{(d_i)}(y)} = \varinjlim_{(\theta_i, d_i)} \prod_{y \in \dot{\pi}_{(d_i)}^{-1}(x)} \MC{L}_x = \varinjlim_{(\theta_i, d_i)} \Hom_{sets}(\dot{\pi}_{(d_i)}^{-1}(x), \MC{L}_x)
\end{align*}
Thus we are left to see that any element in $\varinjlim_{(\theta_i, d_i)} \prod_{y \in \dot{\pi}_{(d_i)}^{-1}(x)} (\dot{\pi}_{(d_i)}^* \MC{L})_y$ actually comes from an element in $\varinjlim_{(\theta_i, d_i)} (\dot{\pi}_{(d_i) \,!} \dot{\pi}_{(d_i)}^* \MC{L})_x$. So fix some $(s_y)_y \in \prod_{y \in \dot{\pi}_{(d_i)}^{-1}(x)} (\dot{\pi}_{(d_i)}^* \MC{L})_y$. As there are only finitely many points $y$, we may find some \'etale local neighborhood $V \to \MC{N}^{(\nu_i) \sharp}_U \times_E \Sp E'$ of $x$ such that we may represent $(s_y)_y$ by a section $s$ over $V \times_{\MC{N}^{(\nu_i) \sharp}_U \times_E \Sp E'} \left(\prod_i \MB{M}_{b_{\nu_i}}^{\circ \preceq \mu_i, \theta_i \sharp} \times_E \op{Ig}^{(d_i) \sharp}_U \times_E \Sp E'\right)$. Unfortunately the map ${\rm supp}(s) \to U$ will not be proper in general. As we did in the proof of theorem \ref{thm:CoveringMorphismProperties}, choose for each $i$ some $d_i'$ and $\theta_i'$ such that 
\[\MB{M}_{b_{\nu_i}}^{\circ \preceq \mu_i, \theta_i \sharp} \subset \MB{M}_{b_{\nu_i}}^{\preceq \mu_i, \theta_i \sharp} \subset \MB{M}_{b_{\nu_i}}^{\circ \preceq \mu_i, \theta_i' \sharp}.\]
and consider $s$ as a section of $\dot{\pi}_{(d_i')}^* \MC{L}$ over $V \times_{\MC{N}^{(\nu_i) \sharp}_U \times_E \Sp E'} \left(\prod_i \MB{M}_{b_{\nu_i}}^{\circ \preceq \mu_i, \theta_i' \sharp} \times_E \op{Ig}^{(d_i') \sharp}_U \times_E \Sp E'\right)$ or equivalently as an element in $\MC{F}^{(\theta_i', d_i')}$ (after applying the corresponding normalization factor). By construction its support lies in $V \times_{\MC{N}^{(\nu_i) \sharp}_U \times_E \Sp E'} \left(\prod_i \MB{M}_{b_{\nu_i}}^{\preceq \mu_i, \theta_i \sharp} \times_E \op{Ig}^{(d_i') \sharp}_U \times_E \Sp E'\right)$, which is finite over $V$ by theorem \ref{thm:CoveringMorphismProperties}b). Hence $s \in (\dot{\pi}_{(d_i') \,!} \dot{\pi}_{(d_i')}^* \MC{L})(V)$ as desired. \exit

\warn{}{
The statement is wrong without taking the direct limit, i.e. the canonical inclusion $\MC{F}^{(\theta_i, d_i)}_x \subset \Hom_{sets}(\dot{\pi}_{(d_i)}^{-1}(x), \MC{L}_x)$ is not a bijection in general.
}

Under the identification given in the lemma, the transition morphisms take the following form: Fix any $s \in \MC{F}^{(\theta_i, d_i)}_x = \Hom_{sets}(\dot{\pi}_{(d_i)}^{-1}(x), \MC{L}_x)$. Then $r_{(d_i' - d_i) *}(s) \in \MC{F}^{(\theta_i, d_i')}_x \subset \Hom_{sets}(\dot{\pi}_{(d_i')}^{-1}(x), \MC{L}_x)$ is given by
\[r_{(d_i' - d_i)}^*(s)(y) = s(r_{(d_i' - d_i)}(y)) \in \MC{L}_x  \qquad \textnormal{for all} \; y\]
and $\iota_{\theta_i, \theta_i' \,!}(s) \in \MC{F}^{(\theta_i', d_i)}_x \subset \Hom_{sets}(\dot{\pi}_{(d_i)}^{-1}(x), \MC{L}_x)$ is given by
\[\iota_{\theta_i, \theta_i' \,!}(s)(y) = \left\{ \begin{array}{cc} s(\iota_{\theta_i, \theta_i'}^{-1}(y)) & {\rm if} \; y \in \prod_i \MB{M}_{b_{\nu_i}}^{\preceq \mu_i, \theta_i \sharp} \times_E \op{Ig}^{(d_i) \sharp}_U \times_E \Sp E' \\ 0  & \rm{otherwise} \end{array} \right. \qquad \textnormal{for all} \; y\]
Moreover we have $\dot{\pi}_{(d_i) \,!}: \MC{F}^{(\theta_i, d_i)}_x \subset \Hom_{sets}(\dot{\pi}_{(d_i)}^{-1}(x), \MC{L}_x) \to \MC{L}_x$ which is given by
\[\dot{\pi}_{(d_i) \,!}(s) = \sum_{y \in \dot{\pi}_{(d_i)}^{-1}(x)} s(y)\]
allowing us to define $\dot{\pi}_!: \MC{F}_x \to \MC{L}_x$ via the injective limit of the morphisms $\prod_i [I_{d_i}(b_{\nu_i}):I_1(b_{\nu_i})]^{-1} \dot{\pi}_{(d_i) \,!}$. \\
Finally we describe the action of $(\gamma_i)_i \in \prod_i J_i$. Choose again $d_{i, \gamma_i}$, $d_{i, \gamma_{\mathrlap{i}}^{}{}^{\!-1}}$ and $\theta_{i, \gamma_i}$ as usual. Then for any $s \in \MC{F}^{(\theta_i, d_i)}_x \subset \Hom_{sets}(\dot{\pi}_{(d_i)}^{-1}(x), \MC{L}_x)$ we have
\[(\gamma_i)_i (s)(y) = \left\{ \begin{array}{cc} s(r_{(d_{i, \gamma_i} + d_{i, \gamma_{\mathrlap{i}}^{}{}^{\!-1}})}(z)) & {\rm if} \; y \; \textnormal{lies in the image of} \; (\gamma_i)_i \\ 0  & \rm{otherwise} \end{array} \right.\]
Here $z \in \prod_i \MB{M}_{b_{\nu_i}}^{\preceq \mu_i, \theta_i \sharp} \times_E \op{Ig}_{\mathrlap{U}}{}^{\! (d_i + d_{i, \gamma_i} + d_{i, \gamma_{\mathrlap{i}}^{}{}^{\!-1}}) \sharp} \times_E \Sp E'$ is any preimage of $y \in \prod_i \MB{M}_{b_{\nu_i}}^{\preceq \mu_i, \theta_i \oplus \theta_{i, \gamma_i} \sharp} \times_E \op{Ig}_{\mathrlap{U}}{}^{\! (d_i + d_{i, \gamma_{\mathrlap{i}}^{}{}^{\!-1}}) \sharp} \times_E \Sp E'$ under $(\gamma_i)_i$.

\prop{\label{prop:SheafFiberRepresentation}}{}{
There is an isomorphism of $\prod_i J_i$-modules
\[\MC{F}_x \cong C^{\infty}_c(\pi_{(\infty_i)}^{-1}(x), \MC{L}_x)\]
}

\prooof
We define $\Theta: \MC{F}_x \to C^{\infty}_c(\pi_{(\infty_i)}^{-1}(x), \MC{L}_x)$ on an element $s \in \MC{F}^{(\theta_i, d_i)}_x \subset \prod_{y \in \dot{\pi}_{(d_i)}^{-1}(x)} \MC{L}_x$ as follows:
\[\Theta(s)(y) = \left\{ \begin{array}{cc} s(r_{(\infty_i), (d_i)}(y)) & {\rm if} \; y \in \prod_i \MB{M}_{b_{\nu_i}}^{\preceq \mu_i, \theta_i \sharp} \times_E \op{Ig}^{(\infty_i) \sharp}_U \times_E \Sp E' \\ 0  & \rm{otherwise} \end{array} \right.\]
It is clear that this definition is independent of the choice of $\MC{F}^{(\theta_i, d_i)}_x$ in which we represent $s$. Moreover its support is the preimage of a finite set of points in $\prod_i \MB{M}_{b_{\nu_i}}^{\preceq \mu_i, \theta_i \sharp} \times_E \op{Ig}^{(d_i) \sharp}_U \times_E \Sp E'$ , hence it is compact (recall the description of the topology on $\pi_{(\infty_i)}^{-1}(x)$ given in remark \ref{rem:FiberTopologyProFinite}). Moreover $\Theta(s)$ is constant for the action of $\Gamma_{(\theta_i, d_i)}$ (as defined in claim $5$ of the proof of proposition \ref{prop:SheafJActionDefinition}) on $\pi_{(\infty_i)}^{-1}(x)$, as $\Gamma_{(\theta_i, d_i)}$ only permutes the points in the fiber of $r_{(\infty_i), (d_i)}$ (except for points not lying over $\prod_i \MB{M}_{b_{\nu_i}}^{\preceq \mu_i, \theta_i \sharp} \times_E \Sp E'$, where $\Theta(s)$ is zero anyway). Thus $\Theta(s)$ is smooth. \\
It is immediate that $\Theta$ is injective. To show surjectivity consider any $f \in C^{\infty}_c(\pi_{(\infty_i)}^{-1}(x), \MC{L}_x)$. As $f$ has compact support there is some $\theta_i$ such that the ${\rm supp}(f) \subset \prod_i \MB{M}_{b_{\nu_i}}^{\preceq \mu_i, \theta_i \sharp} \times_E \op{Ig}^{(\infty_i) \sharp}_U \times_E \Sp E'$. Moreover smoothness implies that $f$ is invariant under some $\Gamma_{(\theta_i, d_i)}$, i.e. $f$ comes by pullback from a function on $\dot{\pi}_{(d_i)}^{-1}(x) \subset \prod_i \MB{M}_{b_{\nu_i}}^{\preceq \mu_i, \theta_i \sharp} \times_E \op{Ig}^{(d_i) \sharp}_U \times_E \Sp E'$. But the previous proposition identifies functions on $\dot{\pi}_{(d_i)}^{-1}(x)$ with some sections of $\MC{F}_x$, which give the desired preimage of $f$ under $\Theta$. \\
It remains to see invariance of this isomorphism under the $\prod_i J_i$-action. But this immediately follows from the description of the action on the left-hand side stated above. \exit

\cor{\label{cor:SheafFiberRepresentation2}}{}{
There exists a (non-canonical) isomorphism of $\prod_i J_i$-modules
\[\MC{F}_x \cong c\op{-}Ind_{\{1\}}^{\prod J_i} (\MC{L}_x)\]
where $c\op{-}Ind$ denotes the compactly induced representation.
}

\prooof
By proposition \ref{prop:JActionFibers} we may (non-canonically) identify $\pi_{(\infty_i)}^{-1}(x)$ and $\prod_i J_i$ (as schemes with $\prod_i J_i$-action). Under this identification $C^{\infty}_c(\pi_{(\infty_i)}^{-1}(x), \MC{L}_x)$ is just the definition of $c\op{-}Ind_{\{1\}}^{\prod J_i} (\MC{L}_x)$. Hence proposition \ref{prop:SheafFiberRepresentation} translates precisely into the statement of the corollary. \exit

\rem{}{
Most of the statements above can be modified to work for $\pi_{(d_i)}$ instead of $\dot{\pi}_{(d_i)}$ as well: In order to define all maps between the sheaves, we have to change the definition of $\MC{F}^{(\theta_i, d_i)}$ into
\[\MC{F}^{(\theta_i, d_i)} \coloneqq \pi_{(d_i) \,*} \pi_{(d_i)}^! (\MC{L})\]
which exists as a sheaf (and not only in the derived category) by \cite[\S XVIII, proposition 3.1.8a)]{GroSGAVI}.
Then lemma \ref{lem:SheafTransitionMaps} holds as well, one can define the morphism $\pi_!: \MC{F} = \varinjlim_{(\theta_i, d_i)} \MC{F}^{(\theta_i, d_i)} \to \MC{L}$ and proposition \ref{prop:SheafJActionDefinition} is still valid, using essentially the same construction of the $\prod_i J_i$-action. One only has to change some of the arguments in the proofs. \\
When analyzing the fibers at geometric points, the main problem is no longer the identification of $\prod_{y \in \dot{\pi}_{(d_i)}^{-1}(x)} (\dot{\pi}_{(d_i)}^* \MC{L})_y$ with $(\dot{\pi}_{(d_i) \,!} \dot{\pi}_{(d_i)}^* \MC{L})_x$ (at least in the limit), but to see $(\pi_{(d_i)}^! \MC{L})_y = \MC{L}_{\pi_{(d_i)}(y)}$ (again at least in the limit). As this holds most notably for flat morphisms, we have such an equality when working with $\dot{\pi}_{(d_i)}$, but not for $\pi_{(d_i)}$. Thus we still need a similar comparison argument between $\dot{\pi}_{(d_i)}$ and $\pi_{(d_i)}$ as done in lemma \ref{lem:SheafFiberHomSpace}. 
Apart from this, proposition \ref{prop:SheafFiberRepresentation} and its corollary admit immediate translations to this alternative setting. \\
The main problem will only appear in the next section: When looking at the cohomology we wish to apply a K\"unneth formula to describe the cohomology of the product space in terms of the cohomology of the Rapoport-Zink space and of the Igusa variety. But instead of having a product of usual inverse images of sheaves, we will be forced to deal with a product of two exceptional inverse image sheaves. Unfortunately there seems to be no known K\"unneth formula for such a situation.
}

\subsection{Cohomology in the special fiber}\label{sec:CohomSpecialFiber}
We show the following spectral sequence (for a prime power $\ell^r$ with $\ell \neq p$)
\[E_2^{p, q} = \bigoplus_{t+s = q} Tor^p_{\MC{H}_r(\prod J_i)} \left(H^s_c\left(\prod_i \ov{\MB{M}}_{b_{\nu_i}}^{\preceq \mu_i}, \M{Z}/\ell^r\M{Z}\right), H^t_c\left(\ov{\op{Ig}}^{(\infty_i)}_U, \M{Z}/\ell^r\M{Z}\right) \right) \Rightarrow H_c^{p+q}\left(\ov{\MC{N}}^{(\nu_i)}_U, \M{Z}/\ell^r\M{Z}\right)\]
computing the cohomology of a Newton stratum in the special fiber at $(c_i)$ of the moduli space of global $G$-shtukas $\nabla_n^{\mmu}\MC{H}^1_U(C, G)$. A similar statement holds for more general coefficient sheaves.

\nota{}{
From now on, we abbreviate for any space $X$ its base-change to $\ACFq$ by $\ov{X}$. For example:
\begin{align*}
 \ov{\nabla^{\mmu}_n\MC{H}^1(C, G)} & \coloneqq \nabla^{\mmu}_n\MC{H}^1(C, G) \times_{C^n \setminus \Delta} ((C^n \setminus \Delta) \times \Sp \ACFq) \\
 \ov{\MC{N}}^{(\nu_i) \sharp}_U & \coloneqq \MC{N}^{(\nu_i) \sharp}_U \times_E \Sp \ACFq \\
 \ov{\op{Ig}}^{(d_i) \sharp}_U & \coloneqq \op{Ig}^{(d_i) \sharp}_U \times_E \Sp \ACFq \\
 \ov{\MB{M}}_{b_{\nu_i}}^{\preceq \mu_i \sharp} & \coloneqq \MB{M}_{b_{\nu_i}}^{\preceq \mu_i \sharp} \times_E \Sp \ACFq
\end{align*}
and similarly for non-perfect spaces or the corresponding adic spaces.
}
Assume for simplicity that $\MC{L}$ is now a sheaf of $\M{Z}/\ell^r\M{Z}$-modules over $\MC{N}^{(\nu_i)}_U$, which we consider as well as such a sheaf on $\MC{N}^{(\nu_i) \sharp}_U$. In all other matters the setup and notations of the previous section are retained. \\
Consider the morphism of sheaves $\dot{\pi}_!: \MC{F} \to \MC{L}$ base-changed to $\ov{\MC{N}}^{(\nu_i) \sharp}_U$. We have already seen that it is equivariant for the action of $\prod_i J_i$. Moreover as everything is already defined over $E$, it is equivariant for the action of $\Gamma = Gal(\ACFq/E)$. As the action of $\prod_i J_i$ is defined over $E'$, it actually commutes with the action of $\Gamma_{E'} = Gal(\ACFq/E') \subset \Gamma$. \\
Though we will not need it in the subsequent treatment, one can say even more about the interplay of the two actions: Consider the semi-direct product $\prod_i J_i \rtimes \Gamma$ given by the natural $\Gamma$-action on $\prod_i J_i$. Then $\MC{F}$ (and obviously $\MC{L}$) carries an action of $\prod_i J_i \rtimes \Gamma$ and $\dot{\pi}_!$ is equivariant with respect to this action\footnote{This contradict somewhat \cite[proposition 5.8]{MantoFoliation}, which states (among other things) that the $\Gamma$-action actually commutes with the $\prod_i J_i$-action. This difference is precisely the reason, why we can only make statements of the cohomology as $\Gamma_{E'}$-representations, but not as $\Gamma$-representations.}. \vspace{3mm} \\
We write $\MC{H}_r(\prod J_i)$ for the locally constant functions on $\prod_i J_i$ with values in $\M{Z}/\ell^r\M{Z}$ and let $\Lambda \cong \M{Z}/\ell^r\M{Z}$ be the trivial $\MC{H}_r(\prod J_i)$-representation. Then the functor of coinvariants on $\prod_i J_i$-modules has a left derived functor $\Lambda \otimes^L_{\MC{H}_r(\prod J_i)} (-)$.

\prop{\label{prop:CohomLiftToCovering}}{}{
There is an equality in the bounded derived category of $\M{Z}/\ell^r$-modules with $\Gamma_{E'}$-action
\[R\Gamma_c (\MC{N}^{(\nu_i) \sharp}_U, \MC{L}) = \Lambda \otimes^L_{\MC{H}_r(\prod J_i)} \left(\varinjlim_{(\theta_i, d_i)} R\Gamma_c\left(\prod_i \ov{\MB{M}}_{b_{\nu_i}}^{\circ \preceq \mu_i, \theta_i \sharp} \times_{\ACFq} \ov{\op{Ig}}^{(d_i) \sharp}_U, \dot{\pi}_{(d_i)}^* \MC{L}\right) \right).\]
}

\prooof
Corollary \ref{cor:SheafFiberRepresentation2} shows that we can apply \cite[corollary 5.4]{MantoFoliation} to get
\[R\Gamma_c(\ov{\MC{N}}^{(\nu_i) \sharp}_U, \MC{F}_{\prod J_i}) = \Lambda \otimes^L_{\MC{H}_r(\prod J_i)} R\Gamma_c(\ov{\MC{N}}^{(\nu_i) \sharp}_U, \MC{F})\]
where $\MC{F}_{\prod J_i}$ denotes the sheaf of $\prod_i J_i$-coinvariants of $\MC{F}$. By $\prod_i J_i$-equivariance of $\dot{\pi}_!$, this morphism factors via $\dot{\pi}_!: \MC{F}_{\prod J_i} \to \MC{L}$. In corollary \ref{cor:SheafFiberRepresentation2} we have seen that this map is an isomorphism over geometric points. Hence we may identify $\MC{F}_{\prod J_i}$ canonically with $\MC{L}$. Moreover 
\begin{align*}
 R\Gamma_c(\ov{\MC{N}}^{(\nu_i) \sharp}_U, \MC{F}) & = R\Gamma_c(\ov{\MC{N}}^{(\nu_i) \sharp}_U, \varinjlim_{(\theta_i, d_i)} \MC{F}^{(\theta_i, d_i)}) \\
 & = \varinjlim_{(\theta_i, d_i)} R\Gamma_c(\ov{\MC{N}}^{(\nu_i) \sharp}_U, \MC{F}^{(\theta_i, d_i)}) \\
 & = \varinjlim_{(\theta_i, d_i)} R\Gamma_c(\ov{\MC{N}}^{(\nu_i) \sharp}_U,  \dot{\pi}_{(d_i) \,!} \dot{\pi}_{(d_i)}^* \MC{L}) \\
 & = \varinjlim_{(\theta_i, d_i)} R\Gamma_c\left(\prod_i \ov{\MB{M}}_{b_{\nu_i}}^{\circ \preceq \mu_i, \theta_i \sharp} \times_{\ACFq} \ov{\op{Ig}}^{(d_i) \sharp}_U, \dot{\pi}_{(d_i)}^* \MC{L}\right).
\end{align*}
Putting all these identifications together, we get exactly the stated proposition. \exit

\thm{\label{thm:CohomKunnethDecomposition}}{}{
Let $\MC{L}$, $\MC{E}_1$ and $\MC{E}_2$ be \'etale sheaves of $\M{Z}/\ell^r\M{Z}$-modules over $\ov{\MC{N}}^{(\nu_i)}_U$, $\prod_i \ov{\MB{M}}_{b_{\nu_i}}^{\preceq \mu_i}$ and $\ov{\op{Ig}}^{(d_i)}_U$ together with a continuous action of $\Gamma_{E'}$. Consider them as well as sheaves on the perfection of the respective spaces. Denote by $pr_1, pr_2$ the two projections from $\prod_i \ov{\MB{M}}_{b_{\nu_i}}^{\preceq \mu_i} \times_{\ACFq} \ov{\op{Ig}}^{(d_i)}_U$ to the two factors. \\
Assume that $\MC{E}_1$ and $\MC{E}_2$ are equivariant for the $\prod_i J_i$-action on the underlying spaces and that there exists for all sufficiently large $(d_i, \theta_i)_i$ some $\Gamma_{E'}$-equivariant isomorphisms
\[\dot{\pi}_{(d_i)}^* \MC{L} \cong pr_1^* \MC{E}_1 \otimes pr_2^*\MC{E}_2\]
(over the perfection) which are compatible for different $(d_i, \theta_i)_i$. \\
Then there exists an isomorphism in the derived category of $\M{Z}/\ell^r\M{Z}$-modules with $\Gamma_{E'}$-action
\[R\Gamma_c\left(\ov{\MC{N}}^{(\nu_i)}_U, \MC{L}\right) = R\Gamma_c\left(\prod_i \ov{\MB{M}}_{b_{\nu_i}}^{\preceq \mu_i}, \MC{E}_1\right) \otimes^L_{\MC{H}_r(\prod J_i)} R\Gamma_c\left(\ov{\op{Ig}}^{(\infty_i)}_U, \MC{E}_2\right).\]
In other words, there exists a spectral sequence of $\Gamma_{E'}$-representations
\[E_2^{p, q} = \bigoplus_{t+s = q} Tor^p_{\MC{H}_r(\prod J_i)} \left(H^s_c\left(\prod_i \ov{\MB{M}}_{b_{\nu_i}}^{\preceq \mu_i}, \MC{E}_1\right), H^t_c\left(\ov{\op{Ig}}^{(\infty_i)}_U, \MC{E}_2\right) \right) \Rightarrow H_c^{p+q}\left(\ov{\MC{N}}^{(\nu_i)}_U, \MC{L}\right).\]
}

\prooof
By the K\"unneth formula for \'etale cohomology with compact support, we get
\begin{align*}
 R\Gamma_c(\ov{\MC{N}}^{(\nu_i) \sharp}_U, \MC{F}) & = \varinjlim_{(\theta_i, d_i)} R\Gamma_c\left(\prod_i \ov{\MB{M}}_{b_{\nu_i}}^{\circ \preceq \mu_i, \theta_i \sharp} \times_{\ACFq} \ov{\op{Ig}}^{(d_i) \sharp}_U, \dot{\pi}_{(d_i)}^* \MC{L}\right) \\
 & = \varinjlim_{(\theta_i, d_i)} R\Gamma_c\left(\prod_i \ov{\MB{M}}_{b_{\nu_i}}^{\circ \preceq \mu_i, \theta_i \sharp} \times_{\ACFq} \ov{\op{Ig}}^{(d_i) \sharp}_U, pr_1^* \MC{E}_1 \otimes pr_2^*\MC{E}_2\right) \\
 & = \varinjlim_{(\theta_i, d_i)} R\Gamma_c\left(\prod_i \ov{\MB{M}}_{b_{\nu_i}}^{\circ \preceq \mu_i, \theta_i \sharp}, \MC{E}_1\right) \otimes^L_{\M{Z}/\ell^r\M{Z}} R\Gamma_c\left(\ov{\op{Ig}}^{(d_i) \sharp}_U, \MC{E}_2\right) \\
 & = \left(\varinjlim_{(\theta_i)} R\Gamma_c\left(\prod_i \ov{\MB{M}}_{b_{\nu_i}}^{\circ \preceq \mu_i, \theta_i \sharp}, \MC{E}_1\right)\right) \otimes^L_{\M{Z}/\ell^r\M{Z}} \left(\varinjlim_{(d_i)} R\Gamma_c\left(\ov{\op{Ig}}^{(d_i) \sharp}_U, \MC{E}_2\right) \right).
\end{align*}
The next part is essentially a formal consequence of the previous proposition. Together with \cite[proposition 5.12]{MantoFoliation} this implies:
\begin{align*}
 R\Gamma_c (\MC{N}^{(\nu_i) \sharp}_U, \MC{L}) & = \Lambda \otimes^L_{\MC{H}_r(\prod J_i)} R\Gamma_c(\ov{\MC{N}}^{(\nu_i) \sharp}_U, \MC{F}) \\
 & = \Lambda \otimes^L_{\MC{H}_r(\prod J_i)} \left(\varinjlim_{(\theta_i)} R\Gamma_c\left(\prod_i \ov{\MB{M}}_{b_{\nu_i}}^{\circ \preceq \mu_i, \theta_i \sharp}, \MC{E}_1\right) \otimes^L_{\M{Z}/\ell^r\M{Z}} \varinjlim_{(d_i)} R\Gamma_c\left(\ov{\op{Ig}}^{(d_i) \sharp}_U, \MC{E}_2\right) \right) \\
 & = \left(\varinjlim_{(\theta_i)} R\Gamma_c\left(\prod_i \ov{\MB{M}}_{b_{\nu_i}}^{\circ \preceq \mu_i, \theta_i \sharp}, \MC{E}_1\right)\right) \otimes^L_{\MC{H}_r(\prod J_i)} \left(\varinjlim_{(d_i)} R\Gamma_c\left(\ov{\op{Ig}}^{(d_i) \sharp}_U, \MC{E}_2\right) \right).
\end{align*}
Taking cohomology commutes with passing to direct limits of spaces. Therefore one can identify $\varinjlim_{(\theta_i)} R\Gamma_c\left(\prod_i \ov{\MB{M}}_{b_{\nu_i}}^{\circ \preceq \mu_i, \theta_i \sharp}, \MC{E}_1\right)$ with the cohomology of $\prod_i \ov{\MB{M}}_{b_{\nu_i}}^{\preceq \mu_i \sharp}$. As the transition maps between the Igusa varieties are finite \'etale, taking cohomology commutes as well with taking inverse limits of Igusa varieties and we may identify the second factor with the cohomology on $\ov{\op{Ig}}^{(\infty_i) \sharp}_U$. \\
Finally by \cite[proposition A.4]{ZhuAffineGrass} $R\Gamma_c\left(\ov{\MC{N}}^{(\nu_i)}_U, \MC{L}\right) = R\Gamma_c\left(\ov{\MC{N}}^{(\nu_i) \sharp}_U, \MC{L}\right)$ and similarly for the other spaces. This allows us to descend from the perfections to the original spaces. \\
The only non-trivial steps to get compatibility with $\Gamma_{E'}$ were already dealt with in the previous proposition. \exit

\rem{}{
The theorem admits (at least in the formulation as an isomorphisms in the derived category) an immediate generalization by replacing the sheaves $\MC{L}$, $\MC{E}_1$ and $\MC{E}_2$ by arbitrary elements in the derived category of \'etale sheaves of $\M{Z}/\ell^r\M{Z}$-modules.
}

\cor{}{}{
There is a spectral sequence of $\Gamma_{E'}$-representations
\[E_2^{p, q} = \bigoplus_{t+s = q} Tor^p_{\MC{H}_r(\prod J_i)} \left(H^s_c(\prod_i \ov{\MB{M}}_{b_{\nu_i}}^{\preceq \mu_i}, \M{Z}/\ell^r\M{Z}), H^t_c(\ov{\op{Ig}}^{(\infty_i)}_U, \M{Z}/\ell^r\M{Z}) \right) \Rightarrow H_c^{p+q}(\ov{\MC{N}}^{(\nu_i)}_U, \M{Z}/\ell^r\M{Z})\]
}

\prooof
Apply the theorem for $\MC{F}$, $\MC{E}_1$ and $\MC{E}_2$ being the constant $\M{Z}/\ell^r\M{Z}$-sheaves on the respective spaces. \exit

\subsection{Further assumptions}
Unfortunately the formulas for the cohomology of the generic fiber will not hold in full generality, mainly because the comparison theorems for vanishing cycles are only known for proper schemes. Thus we make the following assumption, that will be in force throughout the rest of this paper.

\ass{\label{ass:CrucialProper}}{
All connected components of $\nabla^{\mmu}_n\MC{H}^1_{G(\M{A}^{int c_i})}(C, G) \times_{C^n \setminus \Delta} \Sp E[[\zeta_1, \ldots, \zeta_n]]$ (with trivial level structure) are proper over $\Sp E[[\zeta_1, \ldots, \zeta_n]]$.
}

\rem{}{
i) This assumption is satisfied e.g. for $G = GL_n$. To prove the valuative criterion, first note that any $GL_n$-torsor, i.e. vector bundle, over the generic fiber of a curve over a DVR extends to the whole curve. Indeed consider the double dual of the push-forward of the vector bundle on the generic fiber: By definition this sheaf coincides with the original vector bundle on the generic fiber and it is reflexive by construction, hence locally free as we are over a base of dimension $2$. Now the assertion follows from proposition \ref{Prop:BoundGlobHecke}. \\
ii) The equivalent condition was as well posed in the context of Shimura varieties of PEL-type, cf. \cite[section 8]{MantoFoliationPEL}. Note that in \cite{MantoFoliation} this is an immediate consequence from the setup and not explicitly mentioned. Only by a recent work of Lan and Stroh \cite{LanStrohCycles}, this condition was removed due to the existence of toroidal compactifications.
}
Moreover we will have to study all Newton strata at once, so it will be necessary to impose stronger assumptions on the base field $E$:

\ass{}{
$E'$ is a finite extension of $\Fq$ of some degree $s$ such that for every $\nu \in \MC{B}(G_{c_i})$ for which the associated Newton stratum is non-empty, there is a fundamental alcove $b_\nu$ defined over $E'$, which is decent for $s$, i.e. satisfies $(b_\nu\sigma)^s = z^{s\nu}$ (cf. \ref{def:JGroup}).
} 

Note that there are only finitely many such $\nu$, so the existence of such an $E'$ is clear. The decentness condition for $s$ has to be imposed due to the corresponding condition in \ref{prop:CohomLiftToCovering}. \vspace{3mm} \\
Finally let us stress the point, that whenever adic spaces occur in this section, they will be analytic (i.e. the analytic locus of a general adic space) as already considered in corollary \ref{cor:CoveringMorphismAnalytic}.

\subsection{Comparison of representations with torsion coefficients}\label{sec:TorsionCohomology}
We now start in earnest with dealing with representations appearing in the cohomology of the generic fiber $\nabla^{\mmu}_n\MC{H}^1(C, G)_\eta$. In this section we will focus on torsion coefficients $\M{Z}/\ell^r\M{Z}$. The crucial comparison between the generic and the special fiber is made using the vanishing cycles functor similarly to \cite[section 8.1]{MantoFoliation}, which then allows us to apply the results from section \ref{sec:CohomSpecialFiber}. \\
For the convenience of the reader, we recall first the construction of the $G(\M{A}^{c_i})$-action on the tower of $\nabla^{\mmu}_n\MC{H}^1_U(C, G)_\eta$ respectively on the cohomology $R\Gamma_c(\ov{\nabla^{\mmu}_n\MC{H}^1(C, G)}_\eta, \MC{F})$ (for any constant sheaf $\MC{F}$) from \cite[3.5]{Varsh} or \cite[2.20]{LafforgueVShtukas}:
\const{\label{const:ActionGenericCohomology}}{}{
First of all consider the situation where $g_x \in G(\M{A}_x) \subset G(\M{A}^{c_i})$ is concentrated in one point $x$ away from the $c_i$. Let $U \subset G(\M{A}^{int c_i})$ sufficiently small, which in this context means that $g_xU_xg_x^{-1} \subset G(\M{A}_x^{int})$, where $U_x$ denotes the factor of $U$ at the point $x$. 
Consider now any $S$-valued point in $\nabla^{\mmu}_n\MC{H}_U^1(C, G)_\eta$ given by a global $G$-shtuka $(\MS{G}, \varphi, \psi)$ over some scheme $S$. Wlog. assume that we may extend the (possibly trivial) level structure $\psi$ to a full trivialization $\widetilde{\psi}_x: \MF{L}_x(\MS{G}, \varphi) \to (L^+G_x, \sigma^*)$ locally at the point $x$ (as we may pass to some pro-\'etale cover, where such an extension exists and then descend at the very end back to our original scheme essentially in the same way as we descended $\pi_{\infty_i}$ to finite levels). Then define $g_x(\MS{G}, \varphi, \psi)$ to be the global $G$-shtuka obtained by changing the local $G$-shtuka at the point $x$ along the quasi-isogeny
\[\MF{L}_x(\MS{G}, \varphi) \xrightarrow{\widetilde{\psi}_x} (L^+G_x, \sigma^*) \xrightarrow{g_x} (L^+G_x, \sigma^*)\]
As $\widetilde{\psi}_x$ is well-defined up to an element in $U_x$, this quasi-isogeny is well-defined up to an element in $g_xU_xg_x^{-1}$. By our assumption that $g_xU_xg_x^{-1} \subset G(\M{A}_x^{int})$ it follows that $g_x(\MS{G}, \varphi, \psi)$ is well-defined and carries a $g_xU_xg_x^{-1}$-level structure. Thus we obtain a morphism
\[g_x: \nabla^{\mmu}_n\MC{H}^1_U(C, G)_\eta \to \nabla^{\mmu}_n\MC{H}^1_{g_xUg_x^{-1}}(C, G)_\eta\]
Note that this is the identity if $g_x \in U_x \subset G(\M{A}_x^{int})$. Thus if $g \in G(\M{A}^{c_i})$ is an arbitrary element and $U \subset G(\M{A}^{int c_i})$ satisfies $gUg^{-1} \subset G(\M{A}^{int c_i})$, then we obtain
\[g: \nabla^{\mmu}_n\MC{H}^1_U(C, G)_\eta \to \nabla^{\mmu}_n\MC{H}^1_{gUg^{-1}}(C, G)_\eta\]
by writing $g = \prod_{x \in C \setminus \{c_i\}_i} g_x$. Note that for almost all points $x$, the element $g_x$ will act trivially, so this a priori infinite product has only finitely many components that act non-trivially, hence our definition makes sense. \\
Now let $\MC{F}$ be any constant sheaf, for example $\M{Z}/\ell^r\M{Z}$ or $\M{Z}_\ell$. The action of $g$ induces an (iso)morphism of cohomology groups, where 
\[g^*: R\Gamma_c(\ov{\nabla^{\mmu}_n\MC{H}^1_{gUg^{-1}}(C, G)}_\eta, \MC{F}) \to R\Gamma_c(\ov{\nabla^{\mmu}_n\MC{H}^1_U(C, G)}_\eta, \MC{F})\]
These morphisms $g^*$ commute with the transition morphisms between the various sheaves and the various levels $U$ (of course for sufficiently small $U$). Hence they induce a morphism
\[g^*: \varinjlim_{U} R\Gamma_c(\ov{\nabla^{\mmu}_n\MC{H}^1_U(C, G)}_\eta, \MC{F}) \cong \varinjlim_{U} R\Gamma_c(\ov{\nabla^{\mmu}_n\MC{H}^1_{gUg^{-1}}(C, G)}_\eta, \MC{F}) \to \varinjlim_{U} R\Gamma_c(\ov{\nabla^{\mmu}_n\MC{H}^1_U(C, G)}_\eta, \MC{F})\]
Note that the action of $g$ on the level of schemes is already defined over $E'$. Thus the induced action on cohomology commutes with the canonical action of the Galois group $\Gamma_{E'}$. \\
The same constructions work not only over the generic fiber, but in fact over all of $\nabla^{\mmu}_n\MC{H}_U^1(C, G) \times_{C^n \setminus \Delta} \Sp E[[\zeta_1, \ldots, \zeta_n]]$ and in particular over the special fiber $\MB{X}^{\mmu}_U$. It extends as well to the perfections.
}

\warn{}{
The direct limit over the subgroups $U$ will in general not commute with an inverse limit over sheaves $\MC{F}$. So passing to the limit over all $U$ will only be useful in the next section, when dealing with $\M{Z}_\ell$-coefficients.
}

%
%
As already mentioned above, the vanishing cycles functor $\Psi_\eta$ is of central importance. For its definition in the category of schemes see \cite[XIII, 1.3.2.2]{GroSGAVII}. A similar construction works if one uses analytic adic spaces as generic fibers as explained in \cite[section 3.5]{HuberEtaleCohom}.
The right derived functors then define maps for each coefficient ring $\M{Z}/\ell^r\M{Z}$ of derived categories
\[R\Psi_\eta: D^+(X_\eta, \M{Z}/\ell^r\M{Z}) \to D^+(X_s, \M{Z}/\ell^r\M{Z})\]
if $X$ is a scheme with generic fiber $X_\eta$ and special fiber $X_s$, respectively
\[R\Psi_\eta^{an}: D^+(\MF{X}^{an}, \M{Z}/\ell^r\M{Z}) \to D^+(\MF{X}_s, \M{Z}/\ell^r\M{Z})\]
if $\MF{X}$ is (a perfection of) a locally noetherian formal scheme with associated analytic adic space $\MF{X}^{an}$ and special fiber $\MF{X}_s$.

\prop{}{}{
For all subgroups $U$ there is a canonical isomorphism
\[R\Gamma_c\left(\ov{\nabla^{\mmu}_n\MC{H}^1_U(C, G)}_\eta, \M{Z}/\ell^r\M{Z}\right) \cong R\Gamma_c\left(\ov{\MB{X}}^{\mmu}_U, R\Psi_\eta(\M{Z}/\ell^r\M{Z}_{\ov{\nabla^{\mmu}_n\MC{H}^1(C, G)}_\eta})\right)\]
in the derived category. Moreover it is compatible with the $\Gamma_{E'}$-action on both sides and for all $g \in G(\M{A}^{c_i})$
\[\begin{xy}
 \xymatrix {
  R\Gamma_c\left(\ov{\nabla^{\mmu}_n\MC{H}^1_{gUg^{-1}}(C, G)}_\eta, \M{Z}/\ell^r\M{Z}\right) \ar^-{\sim}[r] \ar^{g^*}[d] & R\Gamma_c\left(\ov{\MB{X}}^{\mmu}_{gUg^{-1}}, R\Psi_\eta(\M{Z}/\ell^r\M{Z}_{\ov{\nabla^{\mmu}_n\MC{H}^1(C, G)}_\eta})\right) \ar^-{g^*}[d] \\
  R\Gamma_c\left(\ov{\nabla^{\mmu}_n\MC{H}^1_U(C, G)}_\eta, \M{Z}/\ell^r\M{Z}\right) \ar^-{\sim}[r] & R\Gamma_c\left(\ov{\MB{X}}^{\mmu}_U, R\Psi_\eta(\M{Z}/\ell^r\M{Z}_{\ov{\nabla^{\mmu}_n\MC{H}^1(C, G)}_\eta})\right)
  }
\end{xy} \]
commutes.
}

\prooof
Note that our crucial properness assumption \ref{ass:CrucialProper} implies that for all subgroups $U \subset G(\M{A}^{int c_i})$, the space $\nabla^{\mmu}_n\MC{H}^1_U(C, G) \times_{C^n \setminus \Delta} \Sp E[[\zeta_1, \ldots, \zeta_n]]$ has only proper connected components. Thus \cite[XIII, 2.1.8]{GroSGAVII} there is for each $U$ and for each $r \geq 1$ an isomorphism
\[R\Gamma_c\left(\ov{\nabla^{\mmu}_n\MC{H}_U^1(C, G)}_\eta, \M{Z}/\ell^r\M{Z}\right) \cong R\Gamma_c\left(\ov{\MB{X}}^{\mmu}_U, R\Psi_\eta(\M{Z}/\ell^r\M{Z}_{\ov{\nabla^{\mmu}_n\MC{H}^1(C, G)}_\eta})\right)\]
Take now any $g \in G(\M{A}^{c_i})$. Then by \cite[XIII, 2.1.7]{GroSGAVII} applied to the smooth morphism $g$, shows that 
\[R\Psi_\eta(g^*): R\Psi_\eta (\M{Z}/\ell^r\M{Z}_{\ov{\nabla^{\mmu}_n\MC{H}_{gUg^{-1}}^1(C, G)}_\eta}) \to R\Psi_\eta (\M{Z}/\ell^r\M{Z}_{\ov{\nabla^{\mmu}_n\MC{H}_U^1(C, G)}_\eta})\]
(as the image under $R\Psi_\eta$ of $g^*$ on the generic fiber) coincides with the map coming from $g: \ov{\MB{X}}^{\mmu}_U \to \ov{\MB{X}}^{\mmu}_{gUg^{-1}}$. In particular both actions coincide after applying the functor $R\Gamma_c$.
The compatibility with the $\Gamma_{E'}$-action follows in the same way.
\exit

\prop{\label{prop:CohomDecompNewtonStrata}}{}{
The complexes 
\[R\Gamma_c\left(\ov{\nabla^{\mmu}_n\MC{H}^1_U(C, G)}_\eta, \M{Z}/\ell^r\M{Z}\right) \quad and \quad \bigoplus_{(\nu_i)} R\Gamma_c\left(\ov{\MC{N}}^{(\nu_i)}_U, R\Psi_\eta^{an}(\M{Z}/\ell^r\M{Z}_{\ov{\MF{N}}^{(\nu_i) \; an}_U})\right)\]
define the same virtual representation in the Grothendieck group of $\Gamma_{E'}$-representations of $\M{Z}/\ell^r\M{Z}$-modules, i.e.
\[\sum_i (-1)^i H^i\left(R\Gamma_c\left(\ov{\nabla^{\mmu}_n\MC{H}^1_U(C, G)}_\eta, \M{Z}/\ell^r\M{Z}\right)\right) \cong \sum_i (-1)^i H^i\left(\bigoplus_{(\nu_i)} R\Gamma_c\left(\ov{\MC{N}}^{(\nu_i)}_U, R\Psi_\eta^{an}(\M{Z}/\ell^r\M{Z}_{\ov{\MF{N}}^{(\nu_i) \; an}_U})\right)\right).\]
Moreover for any $g \in G(\M{A}^{c_i})$ the resulting morphisms in the derived category on both sides induce the same map on the alternating sum of cohomology groups (via the isomorphism above).
}

\prooof
By the previous proposition, it suffices to show all assertions for $R\Gamma_c\left(\ov{\nabla^{\mmu}_n\MC{H}^1_U(C, G)}_\eta, \M{Z}/\ell^r\M{Z}\right)$ instead of $R\Gamma_c\left(\ov{\MB{X}}^{\mmu}_U, R\Psi_\eta(\M{Z}/\ell^r\M{Z}_{\ov{\nabla^{\mmu}_n\MC{H}^1(C, G)}_\eta}) \right)$. As the Newton strata $\ov{\MC{N}}^{(\nu_i)}_U$ form a stratification of $\ov{\MB{X}}^{\mmu}_U$ into locally closed subsets, it follows
\begin{align*}
 & \sum_i (-1)^i H^i\left(R\Gamma_c\left(\ov{\MB{X}}^{\mmu}_U, R\Psi_\eta(\M{Z}/\ell^r\M{Z}_{\ov{\nabla^{\mmu}_n\MC{H}^1(C, G)}_\eta})\right)\right) \\
 & \hspace{15mm} = \sum_i (-1)^i H^i\left(\bigoplus_{(\nu_i)} R\Gamma_c\left(\ov{\MC{N}}^{(\nu_i)}_U, \left.R\Psi_\eta(\M{Z}/\ell^r\M{Z}_{\ov{\nabla^{\mmu}_n\MC{H}^1(C, G)}_\eta})\right|_{\ov{\MC{N}}^{(\nu_i)}_U}\right)\right)
\end{align*}
Now \cite[theorem 3.1]{BerkovichVanishingCyclesII}
\footnote{As the cited theorem is stated for Berkovich spaces and it is not immediately clear that its proof generalizes to our setting with adic spaces, we indicate another alternative way to see this isomorphism: By \cite[theorem 3.5.13]{HuberEtaleCohom} and after identifying $i^* \circ R^+j_*(j^*\M{Z}/\ell^r\M{Z}) = R\Psi_\eta\M{Z}/\ell^r\M{Z}$ for the vanishing cycles functor for schemes and $a^*(j^*\M{Z}/\ell^r\M{Z}) = \M{Z}/\ell^r\M{Z}$ as sheaves on the associated adic space, there is an isomorphism 
\[R\Psi_\eta(\M{Z}/\ell^r\M{Z}_{\ov{\nabla^{\mmu}_n\MC{H}^1_U(C, G)}_\eta}) \cong R\Psi_\eta^{an}(\M{Z}/\ell^r\M{Z}_{\ov{\MF{X}}^{\mmu \; an}_U})\]
So we may view the sheaf $R\Psi_\eta\M{Z}/\ell^r\M{Z}$ as coming from the formal scheme $\ov{\MF{X}}^{\mmu \; an}_U$. Thus it remains to see that the canonical morphism  
\[R\Psi_\eta^{an}(\M{Z}/\ell^r\M{Z}_{\ov{\MF{N}}^{(\nu_i) \; an}_U}) \to \left.\left(R\Psi_\eta^{an}(\M{Z}/\ell^r\M{Z}_{\ov{\MF{X}}^{\mmu \; an}_U})\right)\right|_{\ov{\MC{N}}^{(\nu_i)}_U}\]
is an isomorphism. This can be checked locally on stalks of points $x \in \ov{\MC{N}}^{(\nu_i)}_U$. But by \cite[theorem 3.5.8i)]{HuberEtaleCohom} this stalk does only depend on the formal neighborhood $\ov{\MF{N}}^{(\nu_i) \, \wedge x}_U$ (and the restriction of the sheaf on the associated adic space). So we can identify the stalks on both sides with $R\Psi_\eta^{an}(\M{Z}/\ell^r\M{Z}_{\ov{\MF{N}}^{(\nu_i) \, \wedge x \; an}_U})$.
}
implies that
\[\left.R\Psi_\eta(\M{Z}/\ell^r\M{Z}_{\ov{\nabla^{\mmu}_n\MC{H}^1(C, G)}_\eta})\right|_{\ov{\MC{N}}^{(\nu_i)}_U} \cong R\Psi_\eta^{an}(\M{Z}/\ell^r\M{Z}_{\ov{\MF{N}}^{(\nu_i) \; an}_U})\]

Putting all this together yields an isomorphism between the alternating sums of the cohomology groups of
\[R\Gamma_c\left(\ov{\MB{X}}^{\mmu}_U, R\Psi_\eta(\M{Z}/\ell^r\M{Z}_{\ov{\nabla^{\mmu}_n\MC{H}^1(C, G)}_\eta})\right) \quad \textnormal{and} \quad \bigoplus_{(\nu_i)} R\Gamma_c\left(\ov{\MC{N}}^{(\nu_i)}_U, R\Psi_\eta^{an}(\M{Z}/\ell^r\M{Z}_{\ov{\MF{N}}^{(\nu_i) \; an}_U})\right) \]
It is immediate that all isomorphisms respect the $\Gamma_{E'}$-action. Now consider any $g \in G(\M{A}^{c_i})$ (with $gUg^{-1} \subset G(\M{A}^{int c_i})$ as usual). As its action does not change the universal global $G$-shtuka on $\MB{X}^{\mmu}_U$ at the characteristic place (and would only change it via some quasi-isogeny anyway), the map $g: \MB{X}^{\mmu}_{gUg^{-1}} \to \MB{X}^{\mmu}_U$ respects the decomposition into Newton strata on both sides. Hence so does the corresponding map on cohomology. \exit
$\left. \right.$ \vspace{2mm} \\
So far everything we did concerning the cohomology was essentially a formal consequence of the existence of the Newton stratification in the special fiber. Now we use all our knowledge about the product decomposition, to relate these cohomology groups to the ones defined by Rapoport-Zink spaces and Igusa varieties. For this we first need to pass to perfections:

\prop{\label{prop:CohomCyclesPerfection}}{}{
Let $\MF{X}$ be a locally noetherian formal scheme with special fiber $X$ and associated adic space $\MF{X}^{an}$. Let $\MF{X}^\sharp$, $X^\sharp$ and $\MF{X}^{\sharp an}$ be their perfections. Then for any constructible sheaf $\MC{F}$ on $\MF{X}^{an}$, $R\Psi_{\eta}^{an}(\MC{F}_{\MF{X}^{\sharp an}})$ equals the pullback of $R\Psi_{\eta}^{an}(\MC{F}_{\MF{X}^{an}})$ to $X^\sharp$.
}

\prooof
Denote by $\pi_{\MF{X}}: \MF{X}^{an} \to X$ respectively $\pi_{\MF{X}^\sharp}: \MF{X}^{\sharp an} \to X^\sharp$ the specialization morphisms and by $\varepsilon$ the canonical morphism from the perfection to the original space. Then the vanishing cycles functor is nothing else than the pushforward along $\pi_{\MF{X}}$ respectively $\pi_{\MF{X}^\sharp}$. Applying again \cite[proposition A.4]{ZhuAffineGrass} we get:
\begin{align*}
 \varepsilon^* R\Psi_{\eta}^{an}(\MC{F}_{\MF{X}^{an}}) & = \varepsilon^* \pi_{\MF{X} *} \MC{F}_{\MF{X}^{an}}  = \varepsilon^* \pi_{\MF{X} *} \varepsilon_* \varepsilon^* \MC{F}_{\MF{X}^{an}} \\
 & = \varepsilon^* \varepsilon_* \pi_{\MF{X}^\sharp *}  \varepsilon^* \MC{F}_{\MF{X}^{an}} = \pi_{\MF{X}^\sharp *}  \varepsilon^* \MC{F}_{\MF{X}^{an}} = R\Psi_{\eta}^{an}(\MC{F}_{\MF{X}^{\sharp an}})
\end{align*}
\exit

So we may replace in proposition \ref{prop:CohomDecompNewtonStrata} the cohomology groups $R\Gamma_c\left(\ov{\MC{N}}^{(\nu_i)}_U, R\Psi_\eta^{an}(\M{Z}/\ell^r\M{Z}_{\ov{\MF{N}}^{(\nu_i) \; an}_U})\right)$ by their perfect analogues $R\Gamma_c\left(\ov{\MC{N}}^{(\nu_i) \sharp}_U, R\Psi_\eta^{an}(\M{Z}/\ell^r\M{Z}_{\ov{\MF{N}}^{(\nu_i) \sharp \; an}_U})\right)$.

\lem{}{}{
Let $\MF{X}$ and $\MF{Y}$ be two (perfections of) locally noetherian formal schemes with associated analytic adic spaces $\MF{X}^{an}$ and $\MF{Y}^{an}$. Let $\MC{E}$ be a constructible sheaf on $\MF{X}^{an}$ and $\MC{F}$ a constructible sheaf on $\MF{Y}^{an}$. Denote by $pr_1$ and $pr_2$ the two projections from $\MF{X} \times \MF{Y}$ to the two factors and by $pr_1^{an}$ and $pr_2^{an}$ their counterparts for adic spaces. Then
\[R\Psi_{\eta}^{an}\left(pr_1^{an *}\MC{E} \otimes pr_2^{an *}\MC{F} \right) \cong \left(pr_1^* R\Psi_{\eta}^{an}\MC{E}\right) \otimes \left(pr_2^* R\Psi_{\eta}^{an}\MC{F}\right)\]
}

\prooof
We use again that $\Psi_{\eta}^{an}$ is in the adic setting nothing else than the pushforward functor $\pi_{\MF{X} *}$ for the specialization morphism $\pi_{\MF{X}}: \MF{X}^{an} \to X$. So there exists an isomorphism
\[\Psi_{\eta}^{an}\left(pr_1^{an *}\MC{E} \otimes pr_2^{an *}\MC{F} \right) \cong \left(pr_1^* \Psi_{\eta}^{an}\MC{E}\right) \otimes \left(pr_2^* \Psi_{\eta}^{an}\MC{F}\right)\]
of sheaves on the special fiber of $\MF{X} \times \MF{Y}$. Take now the right derived functor on both sides. As pullbacks and tensor products of sheaves are right exact, one immediately obtains the formula above. \exit

\rem{}{
Under slightly stronger assumptions (and for Berkovich spaces), this statement can be found as lemma II.$4$ in a preprint version of \cite{HarrisTaylor}, which can be downloaded at \href{http://citeseerx.ist.psu.edu/viewdoc/download?doi=10.1.1.41.5795&rep=rep1&type=pdf}{citeseerx}. Note as well, that this lemma is much easier to prove for analytic adic spaces than for schemes or Berkovich spaces due to the simpler definition of $\Psi_{\eta}^{an}$.
}

\prop{\label{prop:ProductDecompVanishingCycles}}{}{
For all sufficiently large $(\theta_i)_i$ and $(d_i)_i$, there exists a canonical isomorphism of \'etale sheaves
\[\widehat{\dot{\pi}}_{(d_i)}^* R\Psi_{\eta}^{an}(\M{Z}/\ell^r\M{Z}_{\ov{\MF{N}}^{(\nu_i) \sharp \; an}_U}) \cong pr_1^*R\Psi_{\eta}^{an} (\M{Z}/\ell^r\M{Z}_{\prod \ov{\MC{M}}_{b_{\nu_i}}^{\circ \preceq \mu_i, \theta_i \sharp \; an}}) \otimes pr_2^*R\Psi_{\eta}^{an} (\M{Z}/\ell^r\M{Z}_{\ov{\MF{Ig}}^{(d_i) \sharp \; an}_U})\]
over $\prod_i \ov{\MB{M}}_{b_{\nu_i}}^{\circ \preceq \mu_i, \theta_i \sharp} \times \ov{\op{Ig}}^{(d_i) \sharp}_U$. Here $pr_1$ and $pr_2$ are the two canonical projections from $\prod_i \ov{\MB{M}}_{b_{\nu_i}}^{\circ \preceq \mu_i, \theta_i \sharp} \times \ov{\op{Ig}}^{(d_i) \sharp}_U$ to the two factors. \\
These isomorphisms are compatible with the transition morphisms for different $(\theta_i)_i$ and $(d_i)_i$. Moreover they are compatible for the $\Gamma_{E'}$-action on both sides and with the morphisms induced by elements in $G(\M{A}^{c_i})$.
}

\prooof
By proposition \ref{prop:FormalCoveringMorphEtaleEtale}, $\widehat{\dot{\pi}}_{(d_i)}$ is \'etale. Thus there is an isomorphism
\begin{align*}
 \widehat{\dot{\pi}}_{(d_i)}^* R\Psi_{\eta}^{an} (\M{Z}/\ell^r\M{Z}_{\ov{\MF{N}}^{(\nu_i) \sharp \; an}_U}) & \cong R\Psi_{\eta}^{an} (\M{Z}/\ell^r\M{Z}_{\prod \ov{\MC{M}}_{b_{\nu_i}}^{\circ \preceq \mu_i, \theta_i \sharp \; an} \times \ov{\MF{Ig}}^{(d_i) \sharp \; an}_U}) \\
 & = R\Psi_\eta^{an} \left(pr_1^{an *}\M{Z}/\ell^r\M{Z}_{\prod \ov{\MC{M}}_{b_{\nu_i}}^{\circ \preceq \mu_i, \theta_i \sharp \; an}} \otimes pr_2^{an *}\M{Z}/\ell^r\M{Z}_{\ov{\MF{Ig}}^{(d_i) \sharp \; an}_U}\right)
\end{align*}
essentially by definition of the vanishing sheaves. The isomorphism in the proposition now follows from the previous lemma for $\MC{E} = \MC{F} = \M{Z}/\ell^r\M{Z}$. \\
That this isomorphisms respects transition morphisms follows because the vanishing cycles functor $R\Psi_{\eta}^{an}$ commutes with these transition morphisms, as they are either open immersions (of Rapoport-Zink spaces) of finite \'etale covers (of Igusa varieties). 
Compatibility for the $\Gamma_{E'}$-action follows because $R\Psi_\eta^{an}$ respects actions of Galois groups and all morphisms used in the constructions are already defined over $E'$ (and even over $E$). Compatibility with elements in $G(\M{A}^{c_i})$ follows, because all $\widehat{\dot{\pi}}_{(d_i)}$ and all projections are equivariant for the $G(\M{A}^{c_i})$-action, because it just changes level structures away from the characteristic places.  \exit

\cor{\label{cor:CohomKunnethGeneric}}{}{
The complexes 
\[R\Gamma_c\left(\ov{\nabla^{\mmu}_n\MC{H}^1(C, G)}_\eta, \M{Z}_\ell\right)\]
and
\[\bigoplus_{(\nu_i)} \varinjlim_U \varprojlim_r R\Gamma_c \left(\prod_i \ov{\MB{M}}_{b_{\nu_i}}^{\preceq \mu_i},  R\Psi_\eta^{an}(\M{Z}/\ell^r\M{Z}_{\prod \ov{\MC{M}}_{b_{\nu_i}}^{\preceq \mu_i \; an}})\right) \otimes^L_{\MC{H}_r(\prod J_i)} \varinjlim_{(d_i)} R\Gamma_c \left(\ov{\op{Ig}}^{(d_i)}_U, R\Psi_\eta^{an}(\M{Z}/\ell^r\M{Z}_{\ov{\MF{Ig}}^{(d_i) \; an}_U})\right) \]
define the same virtual representation in the Grothendieck group of $G(\M{A}^{c_i}) \times \Gamma_{E'}$-representations of $\M{Z}_\ell$-modules.
}

\prooof
Applying theorem \ref{thm:CohomKunnethDecomposition} to $\MC{L} = R\Psi_\eta^{an}(\M{Z}/\ell^r\M{Z}_{\ov{\MF{N}}^{(\nu_i) \sharp \; an}_U})$ gives
\begin{align*}
& R\Gamma_c \left(\ov{\MC{N}}^{(\nu_i) \sharp}_U, R\Psi_\eta^{an}(\M{Z}/\ell^r\M{Z}_{\ov{\MF{N}}^{(\nu_i) \sharp \; an}_U})\right) \cong \\
& \hspace{0.5cm} \cong R\Gamma_c \left(\prod_i \ov{\MB{M}}_{b_{\nu_i}}^{\preceq \mu_i \sharp},  R\Psi_\eta^{an}(\M{Z}/\ell^r\M{Z}_{\prod \ov{\MC{M}}_{b_{\nu_i}}^{\preceq \mu_i \sharp \; an}})\right) \otimes^L_{\MC{H}_r(\prod J_i)} \varinjlim_{(d_i)} R\Gamma_c \left(\ov{\op{Ig}}^{(d_i) \sharp}_U, R\Psi_\eta^{an}(\M{Z}/\ell^r\M{Z}_{\ov{\MF{Ig}}^{(d_i) \sharp \; an}_U})\right). 
\end{align*}
But by proposition \ref{prop:CohomCyclesPerfection} and \cite[proposition A.4]{ZhuAffineGrass} we may identify 
\[R\Gamma_c \left(\ov{\MC{N}}^{(\nu_i)}_U, R\Psi_\eta^{an}(\M{Z}/\ell^r\M{Z}_{\ov{\MF{N}}^{(\nu_i) \; an}_U})\right) = R\Gamma_c \left(\ov{\MC{N}}^{(\nu_i) \sharp}_U, R\Psi_\eta^{an}(\M{Z}/\ell^r\M{Z}_{\ov{\MF{N}}^{(\nu_i) \sharp \; an}_U})\right)\]
and similarly for the Igusa varieties and Rapoport-Zink spaces. 
Thus proposition \ref{prop:CohomDecompNewtonStrata} yields an isomorphism between the virtual $\Gamma_{E'}$-representations associated to
\[R\Gamma_c\left(\ov{\nabla^{\mmu}_n\MC{H}^1_U(C, G)}_\eta, \M{Z}/\ell^r\M{Z}\right)\]
and
\[\bigoplus_{(\nu_i)} R\Gamma_c \left(\prod_i \ov{\MB{M}}_{b_{\nu_i}}^{\preceq \mu_i}, R\Psi_\eta^{an}(\M{Z}/\ell^r\M{Z}_{\prod \ov{\MC{M}}_{b_{\nu_i}}^{\preceq \mu_i \; an}})\right) \otimes^L_{\MC{H}_r(\prod J_i)} \varinjlim_{(d_i)} R\Gamma_c \left(\ov{\op{Ig}}^{(d_i)}_U, R\Psi_\eta^{an}(\M{Z}/\ell^r\M{Z}_{\ov{\MF{Ig}}^{(d_i) \; an}_U})\right) \]
that is compatible with the morphisms given by elements in $G(\M{A}^{c_i})$. As the transition morphisms for varying $r$ are compatible with all these actions, the same holds after passing to the inverse limit over all $\M{Z}/\ell^r\M{Z}$. Finally taking direct limits over all $U$ turns both sides into virtual $G(\M{A}^{c_i})$-representations, which are isomorphic by the compatibility assertions. \exit

\subsection{\texorpdfstring{A decomposition of $\ell$-adic representations}{A decomposition of l-adic representations}}\label{sec:LimitsOfMath}
The previous corollary \ref{cor:CohomKunnethGeneric} still has major disadvantages: First of all one passes to the limits only after taking the tensor product, while doing so for Rapoport-Zink spaces and Igusa varieties separately would be much more natural. This is rectified during this section. \\
Secondly one does not have a good control over the cohomology with values in vanishing cycles. So it would be helpful to either compute these sheaves explicitly or to identify their cohomology groups with $\ell$-adic cohomology of generic fibers of the respective schemes. However most natural conjectures in this direction either need general theorems in far greater generality than currently known or are even incorrect, despite appearances of such statements (of course in the corresponding world of mixed characteristic) in \cite{MantoFoliation} and \cite{MantoFoliationPEL}. 
%
\nota{}{
From now on we will drop the space as a subscript of the coefficient sheaf, e.g. writing $R\Gamma_c \left(\prod_i \ov{\MB{M}}_{b_{\nu_i}}^{\preceq \mu_i}, R\Psi_\eta^{an}(\M{Z}/\ell^r\M{Z})\right)$ instead of $R\Gamma_c \left(\prod_i \ov{\MB{M}}_{b_{\nu_i}}^{\preceq \mu_i}, R\Psi_\eta^{an}(\M{Z}/\ell^r\M{Z}_{\prod \ov{\MC{M}}_{b_{\nu_i}}^{\preceq \mu_i \; an}})\right)$. It should nevertheless be obvious at all times, what space is used. 
}

We first study the cohomology of Rapoport-Zink spaces in more detail. In particular we prove that it is a smooth representation and in fact compactly induced from a finite dimensional one. This however will take a bit of work. \\
The main idea is to compute this cohomology via a finite $\check{\op{C}}$ech complex, similarly to \cite[section 8.2.3-8.2.4]{MantoFoliation}. For this consider the open subset $U \coloneqq \prod_i \MB{M}_{b_{\nu_i}}^{\circ \preceq \mu_i, \theta_i}$. If $(\theta_i)_i$ is chosen large enough, then the translates of $U$ under the group $\prod_i J_i$ cover the whole Rapoport-Zink space $\prod_i \MB{M}_{b_{\nu_i}}^{\preceq \mu_i}$. This happens e.g. if $\pi_{(\infty_i)}$ stays surjective when restricted to $\prod_i \MB{M}_{b_{\nu_i}}^{\circ \preceq \mu_i, \theta_i}$. \\
Now observe that we constructed in claim $5$ of the proof of proposition \ref{prop:SheafJActionDefinition} an open subgroup $H \subset \prod_i J_i$, which acts trivially on $U$. 
For any positive integer $s$ define $\widetilde{(\prod J_i/H)}{}^s_{\neq}$ as the set of $s$-tuples $\delta = (\delta_1, \ldots, \delta_s)$ of distinct elements in $\prod J_i/H$ and $(\prod J_i/H)^s_{\neq} \subset \widetilde{(\prod J_i/H)}{}^s_{\neq}$ as the subset with $\delta_1 = 1$. Then abbreviate for any $\delta = (\delta_1, \ldots, \delta_s) \in \widetilde{(\prod J_i/H)}{}^s_{\neq}$ 
\[U_\delta = \delta_1 U \cap \ldots \cap \delta_s U.\]
We immediately obtain for any $\prod_i J_i$-equivariant $\M{Z}/\ell^r\M{Z}$-sheaf $\MC{F}$ a quasi-isomorphism in the derived category
\[R\Gamma_c\left(\prod \ov{\MB{M}}_{b_{\nu_i}}^{\preceq \mu_i}, \MC{F}\right) = \bigoplus_{\delta \in \widetilde{(\prod J_i/H)}{}^{\bullet}_{\neq}} R\Gamma_c\left(\ov{U}_\delta, \MC{F}\right) \coloneqq \bigoplus_{s} \bigoplus_{\delta \in \widetilde{(\prod J_i/H)}{}^{s}_{\neq}} (-1)^{s+1} R\Gamma_c\left(\ov{U}_\delta, \MC{F}\right)\]
(note that we may indeed sum only over $\widetilde{(\prod J_i/H)}{}^{\bullet}_{\neq}$ ignoring tuples with multiple equal entries, because we are using a cover by Zariski-open subschemes and not arbitrary \'etale covers).
On the left-hand side one has a natural action by $\prod_i J_i$ coming from the action on the spaces. On the right-hand side $\gamma \in \prod_i J_i$ induces isomorphisms isomorphisms $U_\delta \to U_{\gamma \delta}$, where $\gamma \delta = (\gamma \delta_1, \ldots, \gamma \delta_s)$ if $\delta = (\delta_1, \ldots, \delta_s) \in \widetilde{(\prod J_i/H)}{}^s_{\neq}$. The isomorphism above is then equivariant with respect to the $\Gamma_{E'} \times \prod_i J_i$-action on both sides. \\
Next consider $U_{\delta}$ for any $\delta = (\delta_1, \ldots, \delta_s) \in \widetilde{(\prod J_i/H)}{}^s_{\neq}$. Then the action by $\delta_1$ defines an isomorphism $U_{\delta} \to U_{\delta_1^{-1} \delta}$ and similarly on cohomology. Thus we obtain for each $\delta \in (\prod J_i/H)^s_{\neq}$ an isomorphism of $\Gamma_{E'} \times \prod_i J_i$-representations
\[\bigoplus_{\gamma \in \prod_i J_i/H} R\Gamma_c\left(\ov{U}_{\gamma \delta}, \MC{F}\right) \cong c\op{-}Ind_H^{\prod J_i}R\Gamma_c\left(\ov{U}_{\delta}, \MC{F}\right)\]
and 
\[R\Gamma_c\left(\prod \ov{\MB{M}}_{b_{\nu_i}}^{\preceq \mu_i}, \MC{F}\right) = \bigoplus_{\delta \in (\prod J_i/H)^{\bullet}_{\neq}} c\op{-}Ind_H^{\prod J_i}R\Gamma_c\left(\ov{U}_\delta, \MC{F}\right)\]
Moreover any element of $\gamma \in H$ defines an isomorphism between $U_{\delta}$ and $U_{\gamma \delta}$ with both $\delta, \gamma \delta \in (\prod J_i/H)^{\bullet}_{\neq}$. Hence setting $H_{\delta} = H \cap \delta_2H\delta_2^{-1} \cap \ldots \cap \delta_sH\delta_s^{-1}$ for $\delta = (1, \delta_2, \ldots, \delta_s) \in (\prod J_i/H)^s_{\neq}$, we may write
\begin{align*}
 \bigoplus_{\delta \in (\prod J_i/H)^{\bullet}_{\neq}} c\op{-}Ind_H^{\prod J_i}R\Gamma_c\left(\ov{U}_{\delta}, \MC{F}\right) & \cong \bigoplus_{\delta \in H\backslash (\prod J_i/H)^{\bullet}_{\neq}} c\op{-}Ind_{H_{\delta}}^H c\op{-}Ind_H^{\prod J_i}R\Gamma_c\left(\ov{U}_{\delta}, \MC{F}\right) \\
 & = \bigoplus_{\delta \in H\backslash (\prod J_i/H)^{\bullet}_{\neq}} c\op{-}Ind_{H_{\delta}}^{\prod J_i}R\Gamma_c\left(\ov{U}_{\delta}, \MC{F}\right)
\end{align*}

\prop{\label{prop:RZCohomDecomposition}}{}{
For every $\prod_i J_i$-equivariant constructible $\M{Z}/\ell^r\M{Z}$-sheaf $\MC{F}$ over $\prod_i \MB{M}_{b_{\nu_i}}^{\preceq \mu_i}$, there is a quasi-isomorphism in the derived category
\[R\Gamma_c\left(\prod \ov{\MB{M}}_{b_{\nu_i}}^{\preceq \mu_i}, \MC{F}\right) \cong \bigoplus_{\delta \in H\backslash (\prod J_i/H)^{\bullet}_{\neq}} c\op{-}Ind_{H_{\delta}}^{\prod J_i}R\Gamma_c\left(\ov{U}_{\delta}, \MC{F}\right)\]
of $\prod_i J_i$-representations. The direct sum on the right-hand side is actually finite and $R\Gamma_c\left(\prod \ov{\MB{M}}_{b_{\nu_i}}^{\preceq \mu_i}, \MC{F}\right)$ is smooth.
}

\prooof
Only the two last assertions are not yet shown. To see the finiteness of the direct sum, it suffices to prove that there are only finitely many $\delta$ with non-empty $U_{\delta}$ or even to see that $U$ intersects only finitely many of its translates under $\prod_i J_i$. 
Recall now that the universal quasi-isogeny $\beta^{univ}$ on $U$ is bounded by $(\theta_i)_i$. Hence, given some element $\gamma \in \prod_i J_i$, a necessary condition for $U \cap \gamma U$ to be non-empty is the existence of a point $x \in U$ such that both $\beta^{univ}|_x$ and $\gamma \circ \beta^{univ}|_x$ are bounded by $(\theta_i)_i$. But this implies that $\gamma$ is bounded by $(\theta_i \oplus -\theta_i)_i$. The set of all elements in $\prod J_i$ bounded by $\theta_i \oplus -\theta_i$ is compact, hence contains only finitely many $H$-cosets, which proves our claim. \\
Each $R\Gamma_c\left(\ov{U}_{\delta}, \MC{F}\right)$ is smooth, because it is a finitely generated module over a finite ring (which therefore has only a finite automorphism group as a $\M{Z}/\ell^r\M{Z}$-module). Taking compactly induced representations and direct sums preserve being smooth. Hence $R\Gamma_c\left(\prod \ov{\MB{M}}_{b_{\nu_i}}^{\preceq \mu_i}, \MC{F}\right)$ is indeed smooth. \exit

\cor{\label{prop:CohomRZSpaceSmooth}}{}{
$\varprojlim_r R\Gamma_c \left(\prod_i \ov{\MB{M}}_{b_{\nu_i}}^{\preceq \mu_i}, R\Psi_\eta^{an}(\M{Z}/\ell^r\M{Z})\right)$ exists in the derived category of smooth $\prod_i J_i$-representations in $\M{Z}/\ell^r\M{Z}$-modules.
}

\prooof
Applying the previous lemma to the inverse limit of sheaves $R\Psi_\eta^{an}(\M{Z}/\ell^r\M{Z})$, it remains to see that 
\[\varprojlim_r c\op{-}Ind_{H_{\delta}}^{\prod J_i}R\Gamma_c\left(\ov{U}_{\delta}, R\Psi_\eta^{an}(\M{Z}/\ell^r\M{Z})\right)\]
is smooth. Now inverse limits commute with compactly induced representations, which can be seen directly from the construction of $c\op{-}Ind_{H_{\delta}}^{\prod J_i}$. Hence we are reduced to see the smoothness of 
\[\varprojlim_r R\Gamma_c\left(\ov{U}_{\delta}, R\Psi_\eta^{an}(\M{Z}/\ell^r\M{Z})\right)\]
But this follows again from finiteness of cohomology with compact support on quasi-compact spaces. \exit

\rem{}{
Of course all representations above are still smooth, when considered as a $\Gamma_{E'} \times \prod_i J_i$-representation. This can be shown using the very same arguments.
}

\prop{\label{prop:CohomIgusaFiniteness}}{}{
$\varinjlim_{(d_i)} R\Gamma_c \left(\ov{\op{Ig}}^{(d_i)}_U, R\Psi_\eta^{an}(\M{Z}/\ell^r\M{Z})\right)$ can be represented by an infinite direct sum of complexes of admissible $\prod_i J_i$-representation in $\M{Z}_\ell$-modules. 
}

\prooof
We first deal with the problem, that Igusa varieties are not quasi-compact. For this write the moduli space of $G$-torsors $\MC{H}^1(C, G)$ as an increasing union $\varinjlim_j \MC{X}_j$ of open quasi-compact substacks $\MC{X}_j \subset \MC{H}^1(C, G)$. Then $\MC{X}_j \setminus \MC{X}_{j-1}$ is quasi-compact as well and we have seen in section \ref{subsec:GlobalModuli}, that $(\MC{X}_j \setminus \MC{X}_{j-1}) \times_{\MC{H}^1(C, G)} \MB{X}^{\mmu}_U$ exists as a quasi-compact DM-stack. It follows that the same quasi-compactness assertion holds for central leaves and as well for Igusa varieties $\op{Ig}^{(d_i)}_{U, \MC{X}_j} \coloneqq (\MC{X}_j \setminus \MC{X}_{j-1}) \times_{\MC{H}^1(C, G)} \op{Ig}^{(d_i)}_U$. Moreover the $\op{Ig}^{(d_i)}_{U, \MC{X}_j}$ define a partition of $\op{Ig}^{(d_i)}_U$ into locally closed subsets, compatible for varying $d_i$. Hence we may replace
\[\varinjlim_{(d_i)} R\Gamma_c \left(\ov{\op{Ig}}^{(d_i)}_U, R\Psi_\eta^{an}(\M{Z}/\ell^r\M{Z})\right) \qquad by \qquad \bigoplus_j \varinjlim_{(d_i)} R\Gamma_c \left(\ov{\op{Ig}}^{(d_i)}_{U, \MC{X}_j}, R\Psi_\eta^{an}(\M{Z}/\ell^r\M{Z})\right)\]
and we are left to show that the quasi-compact spaces $\op{Ig}^{(d_i)}_{U, \MC{X}_j}$ define admissible representations in their cohomology. For this note first that $R\Gamma_c \left(\ov{\op{Ig}}^{(d_i)}_{U, \MC{X}_j}, R\Psi_\eta^{an}(\M{Z}/\ell^r\M{Z})\right)$ exists in the derived category of finitely generated $\M{Z}/\ell^r\M{Z}$-modules. \\
Let now $H_{(d_i)} = \prod_i (J_i \cap I_{d_i}(b_{\nu_i}))$. Then the Hochschild-Serre spectral sequence for the transition morphisms of Igusa varieties degenerates because $H_{(d_i)}$ is a pro-$p$-group and we are considering coefficients with $\M{Z}/\ell^r\M{Z}$-coefficients. Hence we have an equality
\[R\Gamma_c \left(\ov{\op{Ig}}^{(d_i)}_{U, \MC{X}_j}, R\Psi_\eta^{an}(\M{Z}/\ell^r\M{Z})\right) = \left(\varinjlim_{(d_i)} R\Gamma_c \left(\ov{\op{Ig}}^{(d_i)}_{U, \MC{X}_j}, R\Psi_\eta^{an}(\M{Z}/\ell^r\M{Z})\right)\right)^{H_{(d_i)}}\]
where on the right-hand side the invariants under $H_{(d_i)}$ are taken. In particular it follows that these invariants exist in the derived category of finitely generated $\M{Z}/\ell^r\M{Z}$-modules. As the $H_{(d_i)}$ form a basis for the topology of $\prod_i J_i$, this implies the admissibility assertion. \exit

\lem{\label{lem:ExtTorComparison}}{}{
Let $J$ be a topological group and $H \subset J$ a compact open subgroup. Let $M$ be an element in the derived category of finite smooth $H$-representations and $N$ an element in the derived category of admissible $J$-representations in $\M{Z}_\ell$-modules. Then there exists a canonical isomorphism in the derived category of finite $\M{Z}_\ell$-modules
\[c\op{-}Ind_H^J M \otimes^L_{\MC{H}(J)} N \cong RHom_{J-smooth}(c\op{-}Ind_H^J (M^*), N).\]
Here $M^*$ denotes the dual of the finite representation $M$. If both $M$ and $N$ are bounded, then so is $c\op{-}Ind_H^J M \otimes^L_{\MC{H}(J)} N$.
}

\prooof
This is essentially a rewording of \cite[lemma 8.4]{MantoFoliation} in the derived setting. \exit

\prop{\label{prop:CommuteDirectInverseLimit}}{}{
For all sufficiently large $d_i'$ there exists a canonical isomorphism in the derived category of $G(\M{A}^{c_i}) \times \Gamma_{E'}$-representations
\begin{align*}
  & R\Gamma_c \left(\prod_i \ov{\MB{M}}_{b_{\nu_i}}^{\preceq \mu_i}, R\Psi_\eta^{an}(\M{Z}/\ell^r\M{Z})\right) \otimes^L_{\MC{H}(\prod J_i)} \varinjlim_{(d_i)} R\Gamma_c \left(\ov{\op{Ig}}^{(d_i)}_U, R\Psi_\eta^{an}(\M{Z}/\ell^r\M{Z})\right) \cong \\
  & \hspace{8mm} \cong \bigoplus_{\delta \in H\backslash (\prod J_i/H)^{\bullet}_{\neq}} R\Gamma_c\left(\ov{U}_{\delta}, R\Psi_\eta^{an}(\M{Z}/\ell^r\M{Z})\right) \otimes^L_{\MC{H}(H_{\delta})} R\Gamma_c \left(\ov{\op{Ig}}^{(d_i')}_U, R\Psi_\eta^{an}(\M{Z}/\ell^r\M{Z})\right),
\end{align*}
where $U_{\delta} \subset \prod_i \MB{M}_{b_{\nu_i}}^{\preceq \mu_i}$ are the open compact subschemes defined at the beginning of this section and $H_{\delta} \subset \prod J_i$ are compact subgroups. \\
Moreover $d_i'$ does not depend on the actual sheaves on both sides, but only on the geometry of the Rapoport-Zink space. The same formula holds when replacing both sides (or just one for that matter) with cohomology with $\M{Z}_{\ell}$-coefficients, i.e. after replacing the cohomology groups by their inverse limits over all $r$. 
}

\prooof
As the exact nature of the coefficient sheaves will play no role in this proof, let us abbreviate $\MC{F}_r \coloneqq R\Psi_\eta^{an}(\M{Z}/\ell^r\M{Z})$ for the sheaf over $\prod_i \ov{\MB{M}}_{b_{\nu_i}}^{\preceq \mu_i}$ and $\MC{F}'_r \coloneqq R\Psi_\eta^{an}(\M{Z}/\ell^r\M{Z})$ for the sheaf over $\ov{\op{Ig}}^{(d_i)}_U$.
By proposition \ref{prop:RZCohomDecomposition} and proposition \ref{prop:CohomIgusaFiniteness} we can write 
\begin{align*}
  & R\Gamma_c \left(\prod_i \ov{\MB{M}}_{b_{\nu_i}}^{\preceq \mu_i}, \MC{F}_r\right) \otimes^L_{\MC{H}(\prod J_i)} \varinjlim_{(d_i)} R\Gamma_c \left(\ov{\op{Ig}}^{(d_i)}_U, \MC{F}'_r\right) \cong \\
  & \hspace{2cm} \cong \left(\bigoplus_{\delta \in H\backslash (\prod J_i/H)^{\bullet}_{\neq}} c\op{-}Ind_{H_{\delta}}^{\prod J_i} R\Gamma_c\left(\ov{U}_{\delta}, \MC{F}_r\right) \right) \otimes^L_{\MC{H}(\prod J_i)} \varinjlim_{(d_i)} R\Gamma_c \left(\ov{\op{Ig}}^{(d_i)}_U, \MC{F}'_r\right) \\
  & \hspace{2cm} \cong \bigoplus_{\delta \in H\backslash (\prod J_i/H)^{\bullet}_{\neq}} \left(c\op{-}Ind_{H_{\delta}}^{\prod J_i} R\Gamma_c\left(\ov{U}_{\delta}, \MC{F}_r\right) \right) \otimes^L_{\MC{H}(\prod J_i)} \varinjlim_{(d_i)} R\Gamma_c \left(\ov{\op{Ig}}^{(d_i)}_U, \MC{F}'_r\right) 
\intertext{Now the smoothness respectively admissibility results imply together with lemma \ref{lem:ExtTorComparison}}
  & \hspace{2cm} \cong \bigoplus_{\delta \in H\backslash (\prod J_i/H)^{\bullet}_{\neq}} RHom_{\prod J_i-smooth}\left(c\op{-}Ind_{H_{\delta}}^{\prod J_i} R\Gamma_c\left(\ov{U}_{\delta}, \MC{F}_r\right)^*, \varinjlim_{(d_i)} R\Gamma_c \left(\ov{\op{Ig}}^{(d_i)}_U, \MC{F}'_r\right)\right)
\intertext{and by Frobenius reciprocity in the form of \cite[proposition 2.5]{BushHennLocalLanglands}, which can be applied after replacing $\varinjlim_{(d_i)} R\Gamma_c \left(\ov{\op{Ig}}^{(d_i)}_U, \MC{F}'_r\right)$ by its maximal smooth subrepresentation}
  & \hspace{2cm} \cong \bigoplus_{\delta \in H\backslash (\prod J_i/H)^{\bullet}_{\neq}} RHom_{H_{\delta}-smooth}\left( R\Gamma_c\left(\ov{U}_{\delta}, \MC{F}_r\right)^*, \varinjlim_{(d_i)} R\Gamma_c \left(\ov{\op{Ig}}^{(d_i)}_U, \MC{F}'_r\right)\right)
\intertext{We have already seen that there are open subgroups of $\prod J_i$ acting trivially on $\ov{U}_{\delta}$. Thus we may choose a sufficiently large tuple $(d_{i \delta})_i$ such that $H_{(d_{i \delta})} = \prod (J_i \cap I_{d_{i \delta}}(b_{\nu_i}))$ lies inside one of them. By finiteness of the direct sum over all $\delta$, we may as well choose one tuple $(d_i')_i$ for all $\delta$. As all elements in $R\Gamma_c\left(\ov{U}_{\delta}, \MC{F}_r\right)^*$ are by choice of the  $(d_i')_i$ invariant under $H_{(d_i')}$, one obtains}
  & \hspace{2cm} \cong \bigoplus_{\delta \in H\backslash (\prod J_i/H)^{\bullet}_{\neq}} RHom_{H_{\delta}-smooth}\left( R\Gamma_c\left(\ov{U}_{\delta}, \MC{F}_r\right)^*, \left(\varinjlim_{(d_i)} R\Gamma_c \left(\ov{\op{Ig}}^{(d_i)}_U, \MC{F}'_r\right)\right)^{H_{(d_i')}}\right) \\
  & \hspace{2cm} \cong \bigoplus_{\delta \in H\backslash (\prod J_i/H)^{\bullet}_{\neq}} RHom_{H_{\delta}-smooth}\left( R\Gamma_c\left(\ov{U}_{\delta}, \MC{F}_r\right)^*, R\Gamma_c \left(\ov{\op{Ig}}^{(d_i')}_U, \MC{F}'_r\right)\right)
\intertext{and finally by another use of lemma \ref{lem:ExtTorComparison}, but now for $J = H = H_{\delta}$}
  & \hspace{2cm} \cong \bigoplus_{\delta \in H\backslash (\prod J_i/H)^{\bullet}_{\neq}} R\Gamma_c\left(\ov{U}_{\delta}, \MC{F}_r\right) \otimes^L_{\MC{H}(H_{\delta})} R\Gamma_c \left(\ov{\op{Ig}}^{(d_i')}_U, \MC{F}'_r\right) 
\end{align*}
Independence of $d_i'$ of the sheaves is obvious and the very same argumentation works as well for cohomology with $\M{Z}_{\ell}$-coefficients. \exit

\lem{\label{lem:AbstractNonsenseTorInvLimit}}{}{
Let $H$ be a compact topological group, which is Hausdorff. Let $(M_r)_r$ and $(N_r)_r$ be projective systems of elements in the derived category of smooth $H$-representations on $\M{Z}_\ell$-modules of finite rank. Assume that both $\varprojlim_r M_r$ and $\varprojlim_r N_r$ satisfy the same finiteness assumptions. Then
\[(\varprojlim_r M_r) \otimes^L_{\MC{H}(H)} (\varprojlim_r N_r) \cong \varprojlim_r \left(M_r \otimes^L_{\MC{H}(H)} N_r \right)\]
}

\prooof
By finiteness and smoothness of the representations, there exists some open subgroup $H'$ acting trivially on $\varprojlim_r M_r$ and $\varprojlim_r N_r$. Thus 
\[(\varprojlim_r M_r) \otimes^L_{\MC{H}(H)} (\varprojlim_r N_r) \cong  (\varprojlim_r M_r) \otimes^L_{\MC{H}(H/H')} (\varprojlim_r N_r)\]
for the finite group $H/H'$. Similarly one obtains
\[\varprojlim_r \left(M_r \otimes^L_{\MC{H}(H)} N_r \right) \cong \varprojlim_r \left(M_r \otimes^L_{\MC{H}(H/H')} N_r \right)\]
Hence we may assume wlog that $H$ is a finite discrete group. \\
Consider first the case that the projective system $(N_r)_r$ is constantly equal to some $N$. Then observe that both $\varprojlim_r(M_r \otimes^L_{\MC{H}(H)} (-))$ (by \cite[\href{http://stacks.math.columbia.edu/tag/07KV}{Tag 07KV}]{stacks-project}(currently remark 15.64.17)) and $(\varprojlim_r M_r) \otimes^L_{\MC{H}(H)} (-)$ (by definition) are exact functors on the derived category. Moreover they obviously coincide on finite free $H$-modules. Hence they coincide on the derived category of the abelian category generated by finite free $H$-modules. But the finiteness assumption on $N$ together with finiteness of the group $H$ implies, that $N$ lies in this subcategory. Thus we obtain
\[(\varprojlim_r M_r) \otimes^L_{\MC{H}(H)} N \cong \varprojlim_r \left(M_r \otimes^L_{\MC{H}(H)} N \right)\]
Now applying this first for all $N = M_r$ and then for $N = \varprojlim_r N_r$ gives
\begin{align*}
 \varprojlim_r \left(M_r \otimes^L_{\MC{H}(H)} N_r \right) & = \varprojlim_r \varprojlim_{r'} \left(M_r \otimes^L_{\MC{H}(H)} N_{r'} \right) = \varprojlim_r \left(M_r \otimes^L_{\MC{H}(H)} (\varprojlim_{r'} N_{r'}) \right) \\
 & = (\varprojlim_r M_r) \otimes^L_{\MC{H}(H)} (\varprojlim_{r'} N_{r'})
\end{align*}
as desired. \exit

\rem{}{
Actually the following statement holds as well (and can easily be reduced to the lemma above): \\
Let $J$ be a topological group and $H \subset J$ an open compact subgroup. Let $(M_r)_r$ be a projective system of elements in the derived category of smooth $H$-representations on $\M{Z}_\ell$-modules of finite rank and $(N_r)_r$ a projective system of smooth $J$-representations on $\M{Z}_\ell$-modules of finite rank. Assume that both $\varprojlim_r M_r$ and $\varprojlim_r N_r$ satisfy the same finiteness assumptions. Then
\[(\varprojlim_r c\op{-}Ind_H^J M_r) \otimes^L_{\MC{H}(J)} (\varprojlim_r N_r) \cong \varprojlim_r \left(c\op{-}Ind_H^J M_r \otimes^L_{\MC{H}(J)} N_r \right)\]
}

\cor{\label{cor:MoveInverseLimits}}{}{
There is an isomorphism in the derived category of $G(\M{A}^{c_i}) \times \Gamma_{E'}$-representations 
\begin{align*}
  & \varprojlim_r \left(R\Gamma_c \left(\prod_i \ov{\MB{M}}_{b_{\nu_i}}^{\preceq \mu_i}, R\Psi_\eta^{an}(\M{Z}/\ell^r\M{Z})\right) \otimes^L_{\MC{H}(\prod J_i)} \varinjlim_{(d_i)} R\Gamma_c \left(\ov{\op{Ig}}^{(d_i)}_U, R\Psi_\eta^{an}(\M{Z}/\ell^r\M{Z})\right)\right) \cong \\
  & \hspace{12mm} \cong \left(\varprojlim_r R\Gamma_c \left(\prod_i \ov{\MB{M}}_{b_{\nu_i}}^{\preceq \mu_i}, R\Psi_\eta^{an}(\M{Z}/\ell^r\M{Z})\right)\right) \otimes^L_{\MC{H}(\prod J_i)} \left(\varinjlim_{(d_i)} \varprojlim_{r'} R\Gamma_c \left(\ov{\op{Ig}}^{(d_i)}_U, R\Psi_\eta^{an}(\M{Z}/\ell^{r'}\M{Z})\right)\right)
\end{align*}
}

\prooof
This follows from the previous two results: Abbreviate again $\MC{F}_r \coloneqq R\Psi_\eta^{an}(\M{Z}/\ell^r\M{Z})$ over $\prod_i \ov{\MB{M}}_{b_{\nu_i}}^{\preceq \mu_i}$ and $\MC{F}'_{r'} \coloneqq R\Psi_\eta^{an}(\M{Z}/\ell^{r'}\M{Z})$ over $\ov{\op{Ig}}^{(d_i)}_U$.
Then by proposition \ref{prop:CommuteDirectInverseLimit}
\begin{align*}
  & \varprojlim_r \left(R\Gamma_c \left(\prod_i \ov{\MB{M}}_{b_{\nu_i}}^{\preceq \mu_i}, \MC{F}_r\right) \otimes^L_{\MC{H}(\prod J_i)} \varinjlim_{(d_i)} R\Gamma_c \left(\ov{\op{Ig}}^{(d_i)}_U, \MC{F}'_r\right)\right) \cong \\
  & \hspace{2cm} \cong \varprojlim_r \bigoplus_{\delta \in H\backslash (\prod J_i/H)^{\bullet}_{\neq}} \left(R\Gamma_c\left(\ov{U}_{\delta}, \MC{F}_r\right) \otimes^L_{\MC{H}(H_{\delta})} R\Gamma_c \left(\ov{\op{Ig}}^{(d_i')}_U, \MC{F}'_r\right)\right) \\
  & \hspace{2cm} \cong \bigoplus_{\delta \in H\backslash (\prod J_i/H)^{\bullet}_{\neq}} \varprojlim_r \left(R\Gamma_c\left(\ov{U}_{\delta}, \MC{F}_r\right) \otimes^L_{\MC{H}(H_{\delta})} R\Gamma_c \left(\ov{\op{Ig}}^{(d_i')}_U, \MC{F}'_r\right)\right)
\intertext{Thus we may use lemma \ref{lem:AbstractNonsenseTorInvLimit} to obtain}
  & \hspace{2cm} \cong \bigoplus_{\delta \in H\backslash (\prod J_i/H)^{\bullet}_{\neq}} \left(\varprojlim_r R\Gamma_c\left(\ov{U}_{\delta}, \MC{F}_r\right)\right) \otimes^L_{\MC{H}(H_{\delta})} \left(\varprojlim_{r'} R\Gamma_c \left(\ov{\op{Ig}}^{(d_i')}_U, \MC{F}'_{r'}\right)\right)
\intertext{and once again proposition \ref{prop:CommuteDirectInverseLimit}, but now for cohomology with $\M{Z}_\ell$-coefficients}
  & \hspace{2cm} \cong \left(\varprojlim_r R\Gamma_c \left(\prod_i \ov{\MB{M}}_{b_{\nu_i}}^{\preceq \mu_i}, \MC{F}_r\right)\right) \otimes^L_{\MC{H}(\prod J_i)} \left(\varinjlim_{(d_i)} \varprojlim_{r'} R\Gamma_c \left(\ov{\op{Ig}}^{(d_i)}_U, \MC{F}'_{r'}\right)\right).
\end{align*}
\exit

\thm{\label{thm:FinalFormula}}{}{
Let $\nabla^{\mmu}_n\MC{H}^1(C, G)$ be a moduli space of global $G$-shtukas, such that all connected components of $\nabla^{\mmu}_n\MC{H}^1(C, G) \times_{C^n \setminus \Delta} \Sp E[[\zeta_1, \ldots, \zeta_n]]$ are proper over $\Sp E[[\zeta_1, \ldots, \zeta_n]]$. \\
Then there exists a canonical isomorphism between the virtual $G(\M{A}^{c_i}) \times \Gamma_{E'}$-representations
\[\sum_i (-1)^i H^i_c\left(\ov{\nabla^{\mmu}_n\MC{H}^1(C, G)}_\eta, \M{Q}_\ell\right)\]
and 
\[\sum_{(\nu_i)} \sum_{d, e, f} (-1)^{d+e+f} Tor_d^{\MC{H}(\prod J_i)} \left(H^e_c\left(\prod \ov{\MB{M}}_{b_{\nu_i}}^{\preceq \mu_i}, R\Psi_\eta^{an}\M{Q}_\ell \right), \varinjlim_U \varinjlim_{d_i} H^f_c \left(\ov{\op{Ig}}^{(d_i)}_U, R\Psi_\eta^{an}\M{Q}_\ell \right)\right) \]
Here $H^i_c\left(\prod \ov{\MB{M}}_{b_{\nu_i}}^{\preceq \mu_i}, R\Psi_\eta^{an}\M{Q}_\ell \right) \coloneqq \varprojlim_r H^i_c\left(\prod \ov{\MB{M}}_{b_{\nu_i}}^{\preceq \mu_i}, R\Psi_\eta^{an}(\M{Z}/\ell^r\M{Z})\right) \otimes_{\M{Z}_\ell} \M{Q}_\ell$ and similarly for Igusa varieties.
}

\prooof
We prove the statement already on the level of representations in $\M{Z}_\ell$-modules. Then essentially by definition $\sum_i (-1)^i H^i_c\left(\ov{\nabla^{\mmu}_n\MC{H}^1(C, G)}_\eta, \M{Z}_\ell\right)$ is the virtual representation associated to the complex 
\[R\Gamma_c\left(\ov{\nabla^{\mmu}_n\MC{H}^1(C, G)}_\eta, \M{Z}_\ell\right).\]
By corollary \ref{cor:CohomKunnethGeneric} we may replace it by the complex
\[\bigoplus_{(\nu_i)} \varinjlim_U \varprojlim_r \left( R\Gamma_c \left(\prod \ov{\MB{M}}_{b_{\nu_i}}^{\preceq \mu_i}, R\Psi_\eta^{an}(\M{Z}/\ell^r\M{Z})\right) \otimes^L_{\MC{H}_r(\prod J_i)} \varinjlim_{(d_i)} R\Gamma_c \left(\ov{\op{Ig}}^{(d_i)}_U, R\Psi_\eta^{an}(\M{Z}/\ell^r\M{Z}) \right)\right)\]
which by corollary \ref{cor:MoveInverseLimits} coincides with
\begin{align*}
 & \bigoplus_{(\nu_i)} \varinjlim_U \left( R\Gamma_c \left(\prod \ov{\MB{M}}_{b_{\nu_i}}^{\preceq \mu_i}, R\Psi_\eta^{an}\M{Z}_\ell\right) \otimes^L_{\MC{H}(\prod J_i)} \varinjlim_{d_i} R\Gamma_c \left(\ov{\op{Ig}}^{(d_i)}_U, R\Psi_\eta^{an}\M{Z}_\ell \right)\right) \\
 & \qquad = \bigoplus_{(\nu_i)} R\Gamma_c \left(\prod \ov{\MB{M}}_{b_{\nu_i}}^{\preceq \mu_i}, R\Psi_\eta^{an}\M{Z}_\ell\right) \otimes^L_{\MC{H}(\prod J_i)} \varinjlim_U \varinjlim_{d_i} R\Gamma_c \left(\ov{\op{Ig}}^{(d_i)}_U, R\Psi_\eta^{an}\M{Z}_\ell \right)
\end{align*}
But the virtual representation associated to this complex is nothing else than
\[\sum_{(\nu_i)} \sum_{d, e, f} (-1)^{d+e+f} Tor_d^{\MC{H}(\prod J_i)} \left(H^e_c\left(\prod \ov{\MB{M}}_{b_{\nu_i}}^{\preceq \mu_i}, R\Psi_\eta^{an}\M{Z}_\ell \right), \varinjlim_U \varinjlim_{d_i} H^f_c \left(\ov{\op{Ig}}^{(d_i)}_U, R\Psi_\eta^{an}\M{Z}_\ell \right)\right) \]
as desired. \exit

\rem{}{
One should expect, that one can deal with $G(\M{A}) \times \Gamma_{E'}$-representations, i.e. including the characteristic places, in a similar fashion: Extending Drinfeld's notion of full sets of sections, as introduced in \cite{DrinfeldElliptic}, to global respectively local $G$-shtukas should allow to define a tower of Rapoport-Zink spaces $\varinjlim_M \prod \MC{M}_{b_{\nu_i}, M}^{\preceq \mu_i}$ representing the moduli problem of full set of sections of level $M$ over the usual Rapoport-Zink space, cf. \cite[section 7.2]{MantoFoliation}. This tower should be compatible with the product decomposition, similarly to \cite[section 7.3]{MantoFoliation}, and therefore catch the action of $G(\M{A}_{c_i})$ on the generic fiber of the moduli space of global $G$-shtukas.
Then similar arguments as above will likely yield an equality of virtual $G(\M{A}) \times \Gamma_{E'}$-representations
\begin{align*}
 & \sum_i (-1)^i H^i_c\left(\ov{\nabla^{\mmu}_n\MC{H}^1(C, G)}_\eta, \M{Q}_\ell\right) = \\
 & \hspace{1cm} = \sum_{(\nu_i)} \sum_{d, e, f} (-1)^{d+e+f} Tor_d^{\MC{H}(\prod J_i)} \left(\varinjlim_M H^e_c\left(\prod \ov{\MB{M}}_{b_{\nu_i}, M}^{\preceq \mu_i}, R\Psi_\eta^{an}\M{Q}_\ell \right), \varinjlim_U \varinjlim_{d_i} H^f_c \left(\ov{\op{Ig}}^{(d_i)}_U, R\Psi_\eta^{an}\M{Q}_\ell \right)\right) 
\end{align*}
}
$\quad$ \vspace{7mm} \\
Finally we discuss some problems, which prevent us from obtaining formulas similar to the ones claimed in \cite{MantoFoliation}.

\prob{\label{prob:Cohomology1}}{
Ideally the cohomology of the moduli space of global $G$-shtukas could be expressed by the cohomology of the adic spaces associated to Rapoport-Zink spaces and the Igusa varieties. However the bad behavior of analytifications of stratifications comes back at us. Let us give a little bit more details here, though ignoring technical details like validity of K\"unneth formulas for adic spaces not locally of finite type...: \\
Recall from proposition \ref{prop:AdicNewtonStrataImage} that the formal completions of the Newton strata define adic subspaces $\MF{N}^{(\nu_i) \; an}_U$ of $\MF{X}^{\mmu \; an}_U$, which however do not form a stratification again. Moreover one can enlarge them to subsets $\MF{N}^{(\nu_i) \; an, strat}_U$, which form a partition, cf. lemma \ref{lem:AdicNewtonStrataPartition}. Thus the two complexes
\[R\Gamma_c\left(\ov{\MF{X}}^{\mmu \; an}_U, \M{Z}/\ell^r\M{Z}\right) \quad {\rm and} \quad \bigoplus_{(\nu_i)} R\Gamma_c\left(\ov{\MF{N}}^{(\nu_i) \; an, strat}_U, \M{Z}/\ell^r\M{Z}\right)\]
define isomorphic virtual representations in the Grothendieck group (and as well after passing to $\M{Z}_\ell$ or $\M{Q}_\ell$-coefficients).
On the other hand (similarly to our computations done in section \ref{sec:CohomSpecialFiber})
\[R\Gamma_c\left(\ov{\MF{N}}^{(\nu_i) \; an}_U, \M{Z}/\ell^r\M{Z}\right) \cong R\Gamma_c\left(\prod_i \ov{\MC{M}}_{b_{\nu_i}}^{\preceq \mu_i \; an}, \M{Z}/\ell^r\M{Z}\right) \otimes^L_{\MC{H}_r(\prod J_i)} \varinjlim_{(d_i)} R\Gamma_c\left(\ov{\MF{Ig}}^{(d_i) \; an}_U, \M{Z}/\ell^r\M{Z}\right).\]
So one is left to compare the cohomology with compact supports of the quasi-compact space $\MF{N}^{(\nu_i) \; an, strat}_U$ and its open non-quasi-compact subspace $\MF{N}^{(\nu_i) \; an}_U$. However already the set of global sections with compact support differ for these two spaces. So it would be rather surprising that their cohomology encode the same representation, when viewed in the Grothendieck group of $G(\M{A}^{c_i}) \times \Gamma_{E'}$-representations.
}

Unfortunately, it is not even clear how to remove even one of the sheaves of vanishing cycles in theorem \ref{thm:FinalFormula}:

\prob{\label{prob:Cohomology3}}{
Let us first discuss the case of Igusa varieties: In \cite[section 7.4]{MantoFoliation} it is claimed that the the sheaf of vanishing cycles $R\Psi_\eta^{an}(\M{Z}/\ell^r\M{Z}_{\ov{\MF{Ig}}^{(d_i) \; an}_U})$ on Igusa varieties is nothing else than the constant sheaf $\M{Z}/\ell^r\M{Z}$. The argument for this is based on the isomorphism
\[f^*R\Psi_\eta^{an}(\M{Z}/\ell^r\M{Z}_{\ov{\MF{Y}}^{an}}) \cong R\Psi_\eta^{an}(\M{Z}/\ell^r\M{Z}_{\ov{\MF{X}}^{an}})\]
for a smooth morphism $f: \MF{X} \to \MF{Y}$ of formal schemes. However this is applied to the structure morphism of the formal Igusa varieties, which are at most formally smooth, so that this isomorphism seems to be unknown in the situation here. \\
Note that the formal smoothness assumption is not true for the formal Igusa varieties used in this article, though it is in the situation of \cite{MantoFoliation} and can be obtained as well in our setting by changing the definitions slightly.
}

\prob{\label{prob:Cohomology4}}{
We now come to Rapoport-Zink spaces. One could hope for an isomorphism
\[R\Gamma_c\left(\prod \ov{\MB{M}}_{b_{\nu_i}}^{\preceq \mu_i}, R\Psi_\eta^{an}\M{Q}_\ell\right) \stackrel{??}{=} R\Gamma_c\left(\prod \ov{\MC{M}}_{b_{\nu_i}}^{\preceq \mu_i \; an}, \M{Q}_\ell\right)\]
using the notations of theorem \ref{thm:FinalFormula}.
However all known comparison theorems require the underlying formal scheme either to be locally of finite type over the base or to be a formal completion of a scheme locally of finite type along a closed subscheme. While the first situation certainly does not apply to Rapoport-Zink spaces, the second would yield (at least locally) a slightly different formula, cf. \cite[corollary 3.5]{BerkovichVanishingCyclesII}. \\
Note that in \cite[section 8.2.7-8.2.9]{MantoFoliation} (from where we will take the notations for this paragraph) it is tried to prove this equality by reduction to cohomology without supports. Apart from questionable use of Poincar\'e duality for non-compact analytic spaces, we encounter again the problem that analytification does not preserve stratifications (similarly to the discussion in problem \ref{prob:Cohomology1}): 
While the comparison theorem yields cohomology groups in the Berkovich space $U_\varepsilon^{cl \; rig}$ associated to the formal completion or the Rapoport-Zink space along $U_\varepsilon^{cl}$, the very last equality in the proof of \cite[theorem 8.7]{MantoFoliation} requires the use of $sp^{-1}U_\varepsilon^{cl}$ (where $sp: \MC{M}_M^{rig} \to \bar{\MC{M}}_M$ is the specialization morphism from Berkovich spaces). Now both spaces $U_\varepsilon^{cl \; rig}$ and $sp^{-1}U_\varepsilon^{cl}$ have the same underlying topological spaces, but are not the same Berkovich spaces, as their sets of admissible covers or equivalently their associated adic spaces differ.
}


\end{document}